# CALCULABILITÉ DE LA COHOMOLOGIE ÉTALE MODULO $\ell$

David A. Madore[†] et Fabrice Orgogozo[‡]

*À Jean-Louis Colliot-Thélène, qui nous a beaucoup appris*

Résumé. Soient $X$ un schéma algébrique sur un corps algébriquement clos et $\ell$ un nombre premier inversible sur $X$. D'après [SGA 4½, Th. finitude, 1.1], les groupes de cohomologie étale $\mathrm{H}^i(X, \mathbb{Z}/\ell\mathbb{Z})$ sont de dimension finie. Utilisant une variante $\ell$-adique des bons voisinages d'Artin et des résultats élémentaires sur la cohomologie des pro-$\ell$ groupes, on exprime la cohomologie de $X$ comme colimite bien contrôlée de celle de topos construits sur des $\mathsf{B}G$, où les $G$ sont des $\ell$-groupes finis calculables. On en déduit que les nombres de Betti modulo $\ell$ de $X$ sont algorithmiquement calculables (au sens de Church-Turing). La première partie du texte est consacrée à la démonstration de ce fait et de quelques compléments naturels. Elle s'appuie sur les outils de la seconde partie, dédiée à la géométrie algébrique effective. Enfin, dans la troisième partie, on présente un formalisme de calcul « universel » sur les éléments d'un corps.

Abstract. Let $X$ be an algebraic scheme over an algebraically closed field and $\ell$ a prime number invertible on $X$. According to [SGA 4½, Th. finitude, 1.1], the étale cohomology groups $\mathrm{H}^i(X, \mathbb{Z}/\ell\mathbb{Z})$ are finite-dimensional. Using an $\ell$-adic variant of Artin's good neighborhoods and elementary results on the cohomology of pro-$\ell$ groups, we express the cohomology of $X$ as a well controlled colimit of that of toposes constructed on $\mathsf{B}G$ where the $G$ are computable finite $\ell$-groups. From this, we deduce that the Betti numbers modulo $\ell$ of $X$ are algorithmically computable (in the sense of Church-Turing). The proof of this fact, along with certain related results, occupies the first part of this paper. This relies on the tools collected in the second part, which deals with computational algebraic geometry. Finally, in the third part, we present a "universal" formalism for computation on the elements of a field.

## Table des matières

[†]Institut Mines-Télécom, Télécom ParisTech, CNRS LTCI ; `david+math@madore.org`

[‡]CNRS, École Polytechnique ; `Fabrice.Orgogozo+math@normalesup.org`

## Introduction

L'objet principal de ce texte est de démontrer le théorème suivant, ainsi que la variante relative **0.9**.

**0.1. Théorème.** *Il existe un algorithme calculant la cohomologie étale* $\mathrm{H}^i(X, \mathbb{F}_\ell)$ *à coefficients dans* $\mathbb{F}_\ell$ *d'un schéma algébrique* $X$ *sur un corps algébriquement clos de caractéristique différente de* $\ell$*, ainsi que l'application* $\mathrm{H}^i(X, \mathbb{F}_\ell) \to \mathrm{H}^i(Y, \mathbb{F}_\ell)$ *déduite par fonctorialité d'un morphisme* $Y \to X$*.*

Bien entendu, il faut préciser l'énoncé et notamment ce qu'on entend par « calculer » : nous rappelons dans la partie II les faits essentiels dont nous aurons besoin sur les corps calculables et la calculabilité des opérations algébriques (et nous présentons dans la partie III un modèle de calcul « universel » qui permet d'obtenir des stratifications explicites comme expliqué en **0.7**) ; voir aussi **0.10** *infra*. Calculer le groupe $\mathrm{H}^i(X, \mathbb{F}_\ell)$ signifie notamment en calculer la dimension, en fonction de $i, \ell$ et des équations de $X$, et calculer l'application $\mathrm{H}^i(X, \mathbb{F}_\ell) \to \mathrm{H}^i(Y, \mathbb{F}_\ell)$ signifie en calculer la matrice dans une base déterminée par l'algorithme. (Voir **0.4** ci-dessous et ¶ **9.1.2**-**9.1.3** pour des précisions.)

L'énoncé précédent répond notamment à la question posée en [Poonen, Testa et Luijk 2012, hypothèse 7.4] et également considérée dans [Edixhoven et Couveignes 2011, chap. 1 et chap. 15, p. 401], où l'accent est mis sur la dépendance du temps d'exécution en $\ell$.

En caractéristique nulle ce résultat est déjà connu (voir par exemple [C. Simpson 2008, corollaire 2.5] ou [Poonen, Testa et Luijk 2012, § 7.2]). Pour une discussion du problème de la calculabilité des groupes d'homotopie ou d'homologie en topologie algébrique, voir par exemple [Sergeraert 1994].

**0.2.** Rappelons ([Deligne 1980, 5.2.2]) brièvement une définition *à la Čech* de ces groupes de cohomologie étale. Soit $X$ une variété algébrique sur un corps algébriquement clos dénombrable de caractéristique $p \neq \ell$ (par exemple $\overline{\mathbb{F}_p}$). Il existe un système projectif, indexé par les entiers naturels $\alpha$, de recouvrements étales $X_\alpha \twoheadrightarrow X$, cofinal au sens suivant : pour tout $U \twoheadrightarrow X$ étale, il existe une factorisation d'un $X_\alpha \twoheadrightarrow X$ à travers $U$. Notons $X_{\alpha\bullet}$ le cosquelette du morphisme $X_\alpha \to X$ c'est-à-dire le schéma simplicial $X_{\alpha n} \coloneqq X_\alpha \times_X \cdots \times_X X_\alpha$ ($n+1$ facteurs). La cohomologie de Čech $\check{\mathrm{H}}^i(X_{\alpha\bullet}, \mathbb{Z}/\ell\mathbb{Z})$ est, par définition, celle de l'ensemble simplicial $\pi_0(X_{\alpha\bullet})$. (Pour tout ensemble simplicial, on peut considérer le $\mathbb{Z}/\ell\mathbb{Z}$-module cosimplicial naturellement associé puis le complexe dont les dérivations sont les sommes alternées des faces (cf. [SGA 4, V, § 1.0 et § 2.3], [Milne 1980, III, § 2]).) Si $X$ est quasi-projective, il résulte de [Artin 1971, corollaire 4.2] que le groupe $\mathrm{H}^i(X, \mathbb{Z}/\ell\mathbb{Z})$ est isomorphe à la colimite des $\check{\mathrm{H}}^i(X_{\alpha\bullet}, \mathbb{Z}/\ell\mathbb{Z})$. Le problème auquel on est immédiatement confronté est que, donnés $k$, les $X_\alpha \to X$ et $\ell$, il n'est *a priori* pas évident de calculer deux entiers $\alpha \leq \beta$ tels que $\mathrm{H}^i(X, \mathbb{Z}/\ell\mathbb{Z}) = \mathrm{Im}\left(\check{\mathrm{H}}^i(X_{\alpha\bullet}, \mathbb{Z}/\ell\mathbb{Z}) \to \check{\mathrm{H}}^i(X_{\beta\bullet}, \mathbb{Z}/\ell\mathbb{Z})\right)$ (ces entiers existent car le système inductif $\check{\mathrm{H}}^i(X_{\alpha\bullet}, \mathbb{Z}/\ell\mathbb{Z})$ de groupes abéliens finis, ayant une colimite finie, est « essentiellement constant » au sens de la section 5 ; voir notamment **5.1**).

**0.3. Stratégie.** Inspirés par des travaux de Michael Artin ([SGA 4, XI]) et Gerd Faltings ([Faltings 1988, § 2]), nous nous ramenons au cas où l'on peut se restreindre dans la description précédente à des *revêtements* étales $X_\alpha \to X$ galoisiens de groupe un $\ell$-groupe. Nous montrons alors que le système inductif $\check{\mathrm{H}}^i(X_{\alpha\bullet}, \mathbb{Z}/\ell\mathbb{Z})$ — il s'agit maintenant de cohomologie de $\ell$-groupes — est « explicitement » essentiellement constant. La technique utilisée pour résoudre (ou plutôt ignorer) les divers problèmes d'extension que l'on rencontre est semblable à celle de [Schön 1991, chap. I-II] (également connue, indépendamment de [ibid.], de Ofer Gabber). La réduction à des espaces « $K(\pi, 1)$ » n'est possible que localement pour la topologie de la descente cohomologique universelle : si $X$ n'est pas lisse, on utilise un théorème de résolution des singularités de A. Johan de Jong. (En particulier, la *topologie des altérations* est suffisamment fine pour notre propos.)

Les ingrédients essentiels de la démonstration, présentés en § 1—8, sont résumés en § 10 qui récapitule l'algorithme de calcul des nombres de Betti.

**0.4.** Cette approche permet également de résoudre le problème posé à la fin de **0.2**. On montre en effet que l'on peut obtenir pour chaque $i$ des cocycles (hyper-Čech) pour une base de $\mathrm{H}^i(X, \mathbb{Z}/\ell\mathbb{Z})$ et, pour tout autre hyperrecouvrement $X_\bullet \to X$ pour la topologie des altérations, le moyen de développer l'image dans cette base d'un cocycle (hyper-Čech) par le morphisme $\check{\mathrm{H}}^i(X_\bullet, \mathbb{Z}/\ell\mathbb{Z}) \to \mathrm{H}^i(X, \mathbb{Z}/\ell\mathbb{Z})$ (voir **8.3**). Comme nous le fait remarquer Pierre Deligne, il n'est pas difficile d'en déduire un calcul de $\mathrm{R}\Gamma(X, \mathbb{Z}/\ell\mathbb{Z})$ (voir **8.4**).

**0.5. Remarque.** Dans l'énoncé du théorème **0.1**, le schéma $X$ est supposé décrit explicitement, c'est-à-dire par des équations le définissant (ou définissant un atlas d'ouverts affines : cf. § 16 pour les détails). Il est donc tentant de se demander quelles constructions naturelles, par exemple solutions de problèmes universels, peuvent être ainsi décrites explicitement par des équations. Cette question fait l'objet d'une vaste littérature ; signalons simplement ici les articles [Anderson 2002] (jacobienne ; voir aussi [Mumford 1975]) et [Haiman et Sturmfels 2004] (schéma de Hilbert). Bien entendu, nous utiliserons d'autres résultats de ce type au cours de la démonstration ; cf. par exemple ¶ **15.5**, où l'on « calcule » la normalisation.

Signalons maintenant quelques améliorations du théorème **0.1**. D'autres améliorations et compléments, ainsi que quelques questions, figurent dans la section 11.

**0.6.** Notons tout d'abord que la cas de la cohomologie à support compact, et plus généralement de la cohomologie relative (c'est-à-dire modulo un fermé), en résulte formellement ; cf. ¶ **11.2**. Nous avons énoncé **0.1** avec des coefficients dans $\mathbb{Z}/\ell\mathbb{Z}$ pour plus de simplicité, mais nous traiterons plus généralement le cas de la cohomologie à valeurs dans $\mathbb{Z}/\ell^n\mathbb{Z}$.

**0.7. Stratifications.** La première des affirmations suivantes, où $\ell$ est un nombre premier *fixé*, inversible sur les schémas considérés, découle immédiatement du modèle de calcul « universel » décrit dans la partie III (cf. **22.4**), et les deux suivantes d'une généralisation facile de celui-ci (cf. § 26 et notamment **26.6**) :

- si les équations de $X$ comportent des indéterminées, on peut calculer les équations des parties constructibles de l'espace (affine) de ces indéterminées correspondant à une stratification par la dimension : autrement dit, donnés $\ell$ et $n$ et un morphisme de type fini $\mathscr{X} \to \mathbb{A}^r_k$ (où $k$ est un corps algébriquement clos calculable, cf. **0.11**), on peut calculer une partition de $\mathbb{A}^r_k$ en un nombre fini de parties constructibles sur chacune desquelles $\mathrm{H}^i(X, \mathbb{Z}/\ell^n\mathbb{Z})$ (où $X$ désigne une fibre en un point géométrique de la partie) a un type d'isomorphisme[①] donné, qu'on peut calculer ; il en va de même, par exemple, du type de $\mathrm{H}^i(X, \mathbb{Z}/\ell^n\mathbb{Z}) \to \mathrm{H}^i(Y, \mathbb{Z}/\ell^n\mathbb{Z})$ (pour un morphisme $\mathscr{X} \to \mathscr{Y}$ sur $\mathbb{A}^r_k$) ;
- si $X$ est un schéma de type fini sur un ouvert de Spec $\mathbb{Z}$, on peut calculer, en fonction de $\ell$ et $n$ (et de la description de $X$), un $p_0$ tel que le type de $\mathrm{H}^i(X_p, \mathbb{Z}/\ell^n\mathbb{Z})$ ne dépende pas de $p$ lorsque $p \geq p_0$, où $X_p$ désigne une fibre géométrique au-dessus de $p$ ; il en va de même du type de $\mathrm{H}^i(X_p, \mathbb{Z}/\ell^n\mathbb{Z}) \to \mathrm{H}^i(Y_p, \mathbb{Z}/\ell^n\mathbb{Z})$ ;

[①] Par « type d'isomorphisme » d'un $\ell$-groupe abélien fini $V$, nous entendons ici bien sûr des $d_1 \leq \cdots \leq d_s$ tels que $V \simeq (\mathbb{Z}/\ell^{d_1}\mathbb{Z}) \times \cdots \times (\mathbb{Z}/\ell^{d_s}\mathbb{Z})$ ; et par type d'isomorphisme d'un morphisme entre deux tels groupes, la donnée d'une matrice à équivalence près.

- variante de l'affirmation précédente : donné un morphisme de type fini $\mathscr{X} \to \mathbb{A}^r_{\mathbb{Z}[\frac{1}{\ell}]}$, on peut calculer une partition de $\mathbb{A}^r_{\mathbb{Z}[\frac{1}{\ell}]}$ en un nombre fini de parties constructibles sur chacune desquelles $\mathrm{H}^i(X, \mathbb{Z}/\ell^n\mathbb{Z})$ a un type d'isomorphisme fixé, où $X$ désigne une fibre en un point géométrique de la partie.

Notons que l'*existence* d'une stratification (resp. d'un nombre premier $p_0$) comme ci-dessus est conséquence des résultats généraux de constructibilité et commutation aux changements de base. Les démonstrations de [Orgogozo 2013] (resp. [Katz et Laumon 1985, théorème 3.3.2]), de nature géométrique, devraient fournir une stratification (resp. un nombre premier $p_0$) explicite qui convienne pour *chaque* nombre premier $\ell$ inversible sur les schémas considérés.

Comme observé par Olivier Wittenberg ([Poonen, Testa et Luijk 2012, prop. 8.3]), le théorème **0.1** (étendu au cas des coefficients $\mathbb{Z}/\ell^n\mathbb{Z}$) a le corollaire suivant.

**0.8. Corollaire.** *Il existe un algorithme calculant la structure de la partie de torsion du groupe de cohomologie $\ell$-adique* $\mathrm{H}^i(X, \mathbb{Z}_\ell)$ *d'une variété* propre et lisse $X$ *sur un corps algébriquement clos de caractéristique différente de $\ell$.*

(Le rang des $\mathrm{H}^i(X, \mathbb{Z}_\ell)$ se calcule, dans le cas d'une variété sur un corps fini, en comptant les points de celles-ci et en utilisant une borne a priori, et dans le cas général en se ramenant au cas des corps finis : cf. [ibid.].) Notons par contre, qu'hormis dans le cas propre et lisse, permettant d'utiliser des arguments de poids, on ne sait malheureusement rien dire de la cohomologie *$\ell$-adique* (cf. **11.5**).

On déduira du théorème **0.1** le résultat suivant (voir **11.4**).

**0.9. Théorème.** *Soit $f : X \to S$ un morphisme de schémas algébriques sur un corps algébriquement clos $k$. Pour tout faisceau abélien $\mathscr{F}$ constructible sur $X$, de torsion inversible sur $k$, et tout entier $i \geq 0$, on peut explicitement calculer le faisceau $\mathrm{R}^i f_\star \mathscr{F}$, fonctoriellement en $\mathscr{F}$.*

Il en résulte formellement que l'on peut calculer $\mathrm{R}^i f_! \mathscr{F}$ : appliquer le théorème à une compactification de $f$ et au prolongement par zéro correspondant de $\mathscr{F}$.
(Idéalement, on aimerait plutôt déduire **0.1** du théorème **0.9**, démontré par dévissage.)

**0.10. Remarques sur la notion d'algorithme.** La « calculabilité » dans le titre de cet article, et le mot « algorithme calculant » dans l'énoncé du théorème **0.1** doivent se comprendre au sens de Church-Turing, c'est-à-dire le fait que les fonctions annoncées soient (générales) récursives, autrement dit calculables par un ordinateur idéalisé, par exemple une machine de Turing ou une machine à registres : cf. [Odifreddi 1989, définition I.1.7 et théorèmes I.4.3 et I.7.9]. Voir **0.11** sur la manière dont la machine doit manipuler les éléments du corps de base.
Soulignons que le fait de travailler avec des fonctions générales récursives nous permet d'effectuer des « recherches non bornées » (ce qu'on appelle aussi utiliser l'« opérateur $\mu$ de Kleene ») : si pour chaque $m$ il existe $n$ vérifiant une certaine propriété $\mathbf{P}(m, n)$ elle-même calculable, alors la fonction $\mu_n \mathbf{P}$ qui à $m$ associe le plus petit $n$ vérifiant $\mathbf{P}(m, n)$ est calculable (l'idée étant qu'on parcourt les $n$ jusqu'à en trouver un qui vérifie la propriété recherchée).
Nous utiliserons notamment librement ce résultat pour construire des objets géométriques : dès lors qu'un théorème garantit l'existence d'un objet géométrique (schéma, morphisme de schémas, ...) possédant une propriété *algorithmiquement testable*, on peut calculer algorithmiquement un tel objet, simplement en énumérant toutes les équations possibles pour les objets géométriques en question et en testant la propriété souhaitée jusqu'à en trouver un qui vérifie la condition voulue (cf. **12.8** pour plus de détails).
Pour une discussion sur la question de savoir dans quelle mesure on pourrait se passer de ce procédé, et si les fonctions dont on affirme la calculabilité seraient en fait *primitivement* récursives, voir ¶ **11.6** plus bas.

**0.11. Remarques générales sur les corps calculables et la géométrie algébrique effective.** Nous considérons dans cet article deux points de vue très liés, mais néanmoins différents, sur la calculabilité des opérations algébriques :

- Selon le premier point de vue, qui est le plus classique, on dispose de la notion de *corps calculable* (rappelée en **12.1**, et qui s'applique notamment à une clôture algébrique de $\mathbb{F}_p$ ou de $\mathbb{F}_p(T_1, \ldots, T_r)$), c'est-à-dire de corps dont les éléments sont codés par des entiers naturels et sur lequel on suppose qu'une machine de Turing puisse effectuer les opérations algébriques. Nous expliquons dans la partie II comment, donné un corps calculable $k$ (disposant d'un algorithme de factorisation et d'une $p$-base finie explicite, cf. **12.7**), on peut algorithmiquement manipuler les $k$-schémas algébriques de type fini et leurs morphismes (§ 16) et calculer ou tester la plupart des opérations et propriétés classiques de la géométrie algébrique. Si on adopte ce point de vue, on peut ignorer complètement la partie III de cet article.
- L'autre point de vue, que nous appellerons « modèle de calcul "universel" », dont l'exposition fait l'objet de la partie III, consiste à considérer un corps $k$ algébriquement clos *quelconque* et des machines qui peuvent manipuler (outre, évidemment, les entiers naturels) les éléments de $k$ sous forme de « boîtes noires » sur lesquelles elles peuvent effectuer les opérations algébriques et le test d'égalité. L'intérêt de ce point de vue est triple. D'abord, il se passe d'hypothèses sur le corps, et revient à travailler sur tous les corps algébriquement clos de caractéristique fixée (pour les corps non algébriquement clos, on renvoie à la discussion en § 25). Ensuite, il garantit que les fonctions calculables sont *uniquement* celles qui s'obtiennent par des opérations algébriques sur les corps (et pas, par exemple, celles qui pourraient dépendre du détail du codage d'un élément du corps par tel ou tel entier naturel : cf. **21.5**) : on observera donc au cours de la lecture de la partie II, que tous les algorithmes décrits sont de cette nature. Enfin, ce modèle de calcul donne automatiquement des stratifications « universelles », comme en **0.7**, sur lesquelles les fonctions calculables prennent des valeurs constantes — c'est d'ailleurs la manière dont on définit ce modèle de calcul (cf. **20.1**).

Il faut donc comprendre dans chaque énoncé du type « si $k$ est un corps algébriquement clos, alors telle opération est calculable » soit (selon le premier point de vue) « si $k$ est un corps algébriquement clos calculable au sens de **12.1**, alors telle opération est calculable », soit (selon le second) « si $k$ est un corps algébriquement clos, alors telle opération est calculable au sens des définitions **20.1** et **20.4** ». Si $k$ n'est pas supposé algébriquement clos, dans le point de vue « corps calculable », on aura encore besoin de supposer qu'on dispose d'un algorithme de factorisation et d'une $p$-base finie explicite (cf. **12.7**).

Encore un autre point de vue serait celui des *mathématiques constructives* à la Bishop (ou plus exactement de l'algèbre constructive), tel qu'exposé dans [Mines, Richman et Ruitenburg 1988] ou bien [Lombardi et Quitté 2011] : ce point de vue est probablement très proche de celui exposé dans notre partie III ; toutefois, n'étant pas assez familiers avec la logique intuitionniste utilisée par les mathématiques constructives pour dresser le parallèle (nous travaillons tout du long en logique classique), nous avons préféré nous contenter de quelques renvois vers des résultats analogues en mathématiques constructives là où ils pouvaient éclairer l'exposition.

**0.12. Leitfaden.** Des efforts ont été faits pour rendre la lecture de chacune des trois parties largement indépendantes des deux autres. L'ordre logique est III $\longrightarrow$ II $\longrightarrow$ I, mais comme on vient de l'expliquer on peut ignorer la partie III si on ne s'intéresse pas à des résultats tels que **0.7**, et par ailleurs le lecteur prêt à admettre la calculabilité des opérations classiques de la géométrie algébrique (ainsi que la représentation algorithmique des objets) pourra se contenter de lire la partie I. Inversement, la partie II peut servir de présentation autonome de certains résultats de géométrie algébrique effective, et la partie III décrit un formalisme indépendant de tout ce qui précède.

*Nous remercions chaleureusement Ofer Gabber auquel nous sommes redevables de plusieurs arguments importants et pour ses suggestions particulièrement utiles. À divers titres nous remercions également Pierre Deligne, Luc Illusie, Jean Lannes, Henri Lombardi, Ronald van Luijk, Jean-Pierre Serre, Lenny Taelman, et Olivier Wittenberg. Le second auteur sait également gré à Ahmed Abbes et Weizhe Zheng (郑维喆) de leurs invitations, respectivement à l'IHÉS et au centre Morningside (晨兴数学中心).*

# I. Cohomologie étale

## 1. $K(\pi, 1)$ PRO-$\ell$

Les résultats de cette section, essentiellement dus à Ofer Gabber[①], permettent d'établir une variante (**1.4.7**) de [FRIEDLANDER 1982, théorème 11.7] (dont l'énoncé est rappelé en **1.3.3**). Ceux qui préfèrent admettre les résultats de cette section pourront se contenter de lire les définitions **1.3.1** à **1.3.4** et **1.4.4**, ainsi que les propositions **1.4.7** et **1.4.12**.
Il est fort probable que les hypothèses nœthériennes faites dans cette section sur les schémas puissent être remplacées par des hypothèses de finitude plus faibles (cohérence, ouverture des composantes connexes).

### 1.1. Champs $\ell$-monodromiques : définitions.

**1.1.1.** Soit $\mathscr{C}$ un champ en groupoïdes sur un schéma localement nœthérien $S$, muni de la topologie étale ([GIRAUD 1971, II.1.2.1.3]). Rappelons que l'on note $\pi_0(\mathscr{C})$ le faisceau associé au préfaisceau

$$U/S \mapsto \{\text{classes d'isomorphisme d'objets de } \mathscr{C}(U)\}$$

— c'est aussi le faisceau des *sous-gerbes maximales* (=strictement pleines ; [ibid., III.2.1.3]) de $\mathscr{C}$ — et, pour chaque section locale $c_U \in \mathrm{Ob}\,\mathscr{C}(U)$, où $U/S$ est un ouvert étale, $\pi_1(\mathscr{C}, c_U)$ le faisceau en groupes $\underline{\mathsf{Aut}}(c_U)$ sur $U$. Le champ $\mathscr{C}$ est dit **constructible** si le faisceau $\pi_0(\mathscr{C})$ et les divers $\pi_1(\mathscr{C}, c)$ sont constructibles. Dans [SGA 1, XIII, § 0], un tel champ est dit 1-constructible. (Comparer avec [GIRAUD 1971, VII.2.2.1] et [ORGOGOZO 2003, § 2].)

**1.1.2.** Un **lien** sur $S$ est une section cartésienne (sur $S$) du champ associé au préchamp des faisceaux de groupes *à automorphisme intérieur près* ([GIRAUD 1971, IV.1.1]) ; un tel objet peut être représenté par un triplet constitué d'un recouvrement étale $S'$ de $S$, d'un faisceau en groupes $G'$ sur $S'$ et d'un isomorphisme extérieur $\varphi \in \mathsf{Isom\,ex}(p_1^\star G', p_2^\star G')$, où $p_1, p_2 : S' \times_S S' \rightrightarrows S'$ sont les deux projections, satisfaisant la condition de cocycle usuelle ([DELIGNE et al. 1982, II, appendice]).

**1.1.3.** À tout champ localement connexe-non vide (**gerbe**) $\mathscr{G}$ sur $S$, est associé un lien $\mathscr{L}$ ; si $c_{S'} \in \mathrm{Ob}\,\mathscr{G}(S')$ pour $S'$ couvrant $S$, on peut prendre $G' = \pi_1(\mathscr{G}, c_{S'})$ dans la description précédente : le lien « est » le faisceau des groupes d'automorphismes de sections locales, à conjugaison près. On dit que $\mathscr{G}$ est **lié** par le lien $\mathscr{L}$ ou encore que $\mathscr{G}$ est une $\mathscr{L}$-gerbe.

Si $\mathscr{G}$ est le champ $\mathsf{B}G$ des torseurs sous un $S$-faisceau en groupes $G$, son lien est celui naturellement associé à $G$, représenté par le triplet ($S' = S, G' = G, \varphi = \mathrm{Id}$). Le faisceau des automorphismes de ce lien est le faisceau des automorphismes *extérieurs* de $G$ et l'ensemble des classes d'isomorphisme de liens *localement* représentés par le $S$-schéma en groupes $G$ est naturellement isomorphe à $\mathrm{H}^1(S, \mathsf{Aut\,ex}(G))$.

**1.1.4.** Soit $\mathscr{L}$ un lien sur $S$. L'ensemble $\mathrm{H}^2(S, \mathscr{L})$ des classes d'équivalences de gerbes liées par $\mathscr{L}$ est naturellement muni d'une action libre et transitive du groupe de cohomologie $\mathrm{H}^2(S, Z(\mathscr{L}))$, où $Z(\mathscr{L})$ désigne le **centre** du lien $\mathscr{L}$. (Voir [ibid., loc. cit.] et [GIRAUD 1971, IV, 1.5.3, 3.1.1 et 3.3.3].)

**1.1.5. Définition.** Soit $\ell$ un nombre premier. Un objet sur $S$ du type suivant est dit **$\ell$-monodromique** s'il satisfait l'une des conditions suivantes.

- un faisceau d'*ensembles* $F$ : s'il est localement constant, constructible, et si pour tout point géométrique $s$ de $S$, le groupe de monodromie, image de $\pi_1(S, s)$ dans $\mathsf{Aut}(F_s)$, est un $\ell$-groupe ;
- un faisceau de *groupes* $G$ : s'il est $\ell$-monodromique en tant que faisceau d'ensembles et si ses fibres sont des $\ell$-groupes (finis) ;
- un lien $\mathscr{L}$ : s'il est *localement* représenté par un $\ell$-groupe fini (constant) $\mathsf{G}$ et si, sur chaque composante connexe, le $\mathsf{Aut\,ex}(\mathsf{G})$-torseur $\mathsf{Isom}_{\text{liens}}(\mathscr{L}, \text{lien}(\mathsf{G}))$ est $\ell$-monodromique, où $\text{lien}(\mathsf{G})$ désigne (abusivement) le lien du champ des $\mathsf{G}$-torseurs ;
- un champ $\mathscr{C}$ : si le faisceau d'ensembles $\pi_0(\mathscr{C})$ est $\ell$-monodromique et si, sur le revêtement correspondant de $S$, les sous-gerbes maximales sont à liens $\ell$-monodromiques.

[①] Les auteurs sont bien entendu seuls responsables des éventuelles imprécisions ou erreurs dans l'exposition.

Nous dirons également qu'un morphisme $Y/X$ est un $\ell$-**revêtement** s'il est fini étale galoisien d'ordre une puissance de $\ell$.

**1.1.6**. *Mise en garde.* Un faisceau $\ell$-monodromique provient par image inverse du topos $S_{\ell \text{ét}}$ considéré en 1.4.3 mais il n'en est pas de nécessairement de même d'un champ $\ell$-monodromique quelconque ; cela reflète le fait qu'une classe dans $\mathrm{H}^2(S, \mathbb{Z}/\ell\mathbb{Z})$ n'est pas nécessairement tuée par un revêtement galoisien d'ordre une puissance de $\ell$ (cf. 1.1.4).

**1.1.7. Remarque.** Les faisceaux $\ell$-monodromiques sont également utilisés en [Orgogozo 2003, 4.6], où l'on démontre la locale constance générique du type d'homotopie étale pro-$\ell$ des fibres d'un morphisme de type fini.

## 1.2. Champs $\ell$-monodromiques : sorites.

**1.2.1**. Un faisceau abélien extension (resp. un sous-quotient) de faisceaux abéliens $\ell$-monodromiques est également $\ell$-monodromique. La lissité (locale constance et constructibilité) résulte de [SGA 4, IX, 2.1(ii) et 2.6(ii)]. Que la monodromie d'une extension soit un $\ell$-groupe est conséquence immédiate du fait qu'un élément unipotent de $\mathrm{GL}_n(\mathbb{Z}/\ell^r\mathbb{Z})$ est de $\ell$-torsion, c'est-à-dire annulé par une *puissance* de $\ell$.

Réciproquement, tout faisceau abélien $\ell$-monodromique est extension successive du faisceau constant $\mathbb{Z}/\ell\mathbb{Z}$.

**1.2.2**. Un faisceau en groupes $G$ est $\ell$-monodromique si et seulement si son lien est $\ell$-monodromique. (Cette dernière propriété est tautologiquement équivalente au fait que la gerbe $\mathsf{B}G$ des $G$-torseurs soit $\ell$-monodromique.) Pour vérifier ce fait, on peut supposer $S$ connexe, et en choisir un point géométrique $s$. Les deux derniers termes de la suite exacte

$$1 \to G_s/Z(G_s) \to \mathsf{Aut}(G_s) \to \mathsf{Aut\,ex}(G_s) \to 1$$

reçoivent le groupe $\pi_1(S, s)$. Le noyau $G_s/Z(G_s)$ étant un $\ell$-groupe, l'image de $\pi_1(S, s)$ dans $\mathsf{Aut}(G_s)$ est un $\ell$-groupe si et seulement si l'image de $\pi_1(S, s)$ dans $\mathsf{Aut\,ex}(G_s)$ l'est.

**1.2.3**. *Gerbe quotient.* Par **sous-groupe normal** d'une gerbe $\mathscr{G}$ sur $S$ nous entendons la donnée pour chaque section locale $\sigma \in \mathrm{Ob}\,\mathscr{G}(U)$ d'un sous-faisceau en groupes distingués $N_\sigma \trianglelefteq \underline{\mathsf{Aut}}(\sigma) = \pi_1(\mathscr{G}, \sigma)$ et ceci de façon compatible aux restrictions et aux isomorphismes entre sections locales. La **gerbe quotient** de $\mathscr{G}$ par $N = (N_\sigma)_\sigma$ est, par définition, le champ associé à la catégorie fibrée $\mathscr{G}'$ ayant mêmes objets que $\mathscr{G}$ mais dont les homomorphismes sont les faisceaux quotients $\underline{\mathsf{Hom}}_{\mathscr{G}'}(\sigma, \tau) := \underline{\mathsf{Hom}}_{\mathscr{G}}(\sigma, \tau)/N_\sigma$. On la note $\mathscr{G}/N$. Le morphisme $\mathscr{G} \to \mathscr{G}/N$ est conservatif et couvrant, c'est-à-dire localement surjectif sur les objets et les flèches.

Cette construction s'applique en particulier aux centres $Z_\sigma = Z(\underline{\mathsf{Aut}}(\sigma))$ des faisceaux d'automorphismes locaux. Ici, $Z_\sigma$ ne « dépend pas » du choix de la section locale et se descend donc à $S$ : il existe un faisceau en groupes abéliens $Z$ et, pour chaque section locale $\sigma \in \mathrm{Ob}\,\mathscr{G}(U)$, un isomorphisme $Z(U) \xrightarrow{\sim} Z_\sigma$ commutant aux restrictions et aux isomorphismes entre sections locales. Lorsque $\mathscr{G}$ est une gerbe des $G$-torseurs, où $G$ est un faisceau en groupes, la gerbe quotient ainsi obtenue est équivalente à celle des $G'$-torseurs, où $G'$ est le quotient de $G$ par son centre.

Si $\sigma$ est une section globale de la gerbe quotient $\mathscr{G}/Z$, la gerbe $\mathscr{K}(\sigma)$ des relèvements de $s$ à $\mathscr{G}$ ([Giraud 1971, IV.2.5.1,2.5.4]) est naturellement une $Z$-gerbe (cf. [ibid., 2.5.6.(ii)]). Rappelons que $\mathscr{K}(\sigma)$ est triviale si et seulement si la section $\sigma$ se relève à $\mathscr{G}$.

**1.2.4. Lemme.** *Soit $\iota : \mathscr{H} \hookrightarrow \mathscr{G}$ un morphisme fidèle entre gerbes constructibles localement constantes sur un schéma localement nœthérien $S$. Si $\mathscr{G}$ est $\ell$-monodromique, il en est de même de $\mathscr{H}$.*

*Démonstration.* Supposons $S$ connexe, comme il est loisible de le faire, et choisissons en un point géométrique $s$. Pour chaque section locale $h$ de $\mathscr{H}$, le morphisme $\iota$ induit une injection $\underline{\mathsf{Aut}}_{\mathscr{H}}(h) \hookrightarrow \underline{\mathsf{Aut}}_{\mathscr{G}}(\iota h)$. Il en résulte que, localement, le morphisme $\mathrm{lien}(\iota) : \mathrm{lien}(\mathscr{H}) \to \mathrm{lien}(\mathscr{G})$ est représenté par une injection de groupes finis $\mathsf{H} \subseteq \mathsf{G}$ et que l'action du groupe fondamental $\pi_1(S, s)$ sur les liens de $\mathscr{H}$ et $\mathscr{G}$ se factorise à travers le morphisme horizontal ci-dessous :

$$\begin{array}{ccc} & \nearrow & \mathsf{Aut\,ex(H)} \\ & & \uparrow \\ \pi_1(S,s) & \longrightarrow & \mathsf{Aut(H \subseteq G)/(H/}Z\mathsf{(G) \cap H)} \\ & & \downarrow (\dagger) \\ & \searrow & \mathsf{Aut\,ex(G)} \end{array}$$

où $\mathsf{Aut(H \subseteq G)}$ désigne le groupe des automorphismes de $\mathsf{G}$ préservant globalement $\mathsf{H}$. Bien entendu, $\mathsf{G}$ étant un $\ell$-groupe, il en est de même de $\mathsf{H}$. Notons d'autre part que le noyau de la flèche (†) est un $\ell$-groupe : c'est un quotient du normalisateur de $\mathsf{H}$ dans $\mathsf{G}$. Il en résulte que l'image du morphisme horizontal est un $\ell$-groupe, car $\mathscr{G}$ est supposé $\ell$-monodromique. Partant, il en est de même de l'image du morphisme $\pi_1(S,s) \to \mathsf{Aut\,ex(H)}$. CQFD. □

**1.2.5. Lemme.** *Soient $\ell$ un nombre premier,* $\mathsf{G}$ *un $\ell$-groupe fini et $\mathscr{G}$ une gerbe $\ell$-monodromique sur un topos $S$, localement isomorphe à* $\mathsf{BG}$. *Si* $\mathrm{H}^2(S, \mathbb{Z}/\ell\mathbb{Z}) = 0$, *alors la gerbe $\mathscr{G}$ est triviale : elle a une section globale.*

En d'autres termes, $\mathscr{G}$ est une gerbe de torseurs sous un *faisceau* en groupes, localement isomorphe à $\mathsf{G}$.

*Démonstration.* Soit $\mathscr{G}'$ la gerbe quotient de $\mathscr{G}$ par son centre $Z$ (cf. **1.2.3**). Cette gerbe est localement isomorphe à $\mathsf{BG}'$, où $\mathsf{G}'$ est le quotient de $\mathsf{G}$ par son centre (non trivial, à moins que $\mathsf{G}$ ne le soit). Par récurrence sur l'ordre du groupe, elle a une section globale. La $Z$-gerbe $\mathscr{K}$ de ses relèvements à $\mathscr{G}$ (cf. *loc. cit.*) a une classe dans $\mathrm{H}^2(S, Z)$ (**1.1.4**), nécessairement nulle par hypothèse (utiliser **1.2.1**). Il en résulte que la gerbe $\mathscr{K}$ est triviale et, par conséquent, que la section globale de $\mathscr{G}'$ considérée se relève à la gerbe $\mathscr{G}$, qui est donc également triviale. □

**1.3. Courbes et polycourbes $\ell$-élémentaires.** Rappelons les définitions [SGA 4, XI.3.1 et XI.3.2].

**1.3.1.** On appelle **courbe élémentaire** sur $S$ un morphisme de schémas $f : X \to S$ qui peut être plongé dans un diagramme commutatif

$$\begin{array}{ccccc} X & \overset{j}{\hookrightarrow} & \overline{X} & \overset{i}{\hookleftarrow} & D \\ & {\scriptstyle f}\searrow & \downarrow{\scriptstyle \overline{f}} & \swarrow{\scriptstyle g} & \\ & & S & & \end{array}$$

satisfaisant aux conditions suivantes :

- **(i)** $j$ est une immersion ouverte et $X = \overline{X} - D$ ;
- **(ii)** $\overline{f}$ est une courbe (relative) projective lisse, à fibres géométriquement connexes ;
- **(iii)** $g$ est un revêtement étale à fibres non vides.

Donnés un morphisme projectif $\overline{f}$ et une immersion fermée $i$, les conditions (i)-(iii) sont algorithmiquement testables : pour la lissité de $\overline{f}$ et $g$, utiliser **17.2** ; pour la connexité géométrique (resp. non vacuité) des fibres de $\overline{f}$ (resp. de $g$), on peut supposer que $S$ est un corps en localisant en ses points maximaux (calculables par **15.2**), car le faisceau $\overline{f}_\star\underline{2}$ (resp. $g_\star\underline{2}$) est lisse, et utiliser **16.6**.

(Notons que si $S$ est *normal* intègre — comme il serait loisible de le supposer pour notre propos —, on pourrait également utiliser la connexité géométrique des fibres d'un morphisme propre $Y \to S$ dont la fibre générique est géométriquement connexe.)

**1.3.2.** On appelle **polycourbe élémentaire** sur $S$ un morphisme de schémas $X \to S$ admettant une factorisation en courbes élémentaires. (Lorsque cela ne semble pas prêter à confusion, on s'autorise à confondre $X$ avec le morphisme $X \to S$.) Notre terminologie, inspirée des courbes et polycourbes *hyperboliques* de [Mochizuki 1999], nous semble plus explicite que les « fibrations élémentaires » et « bons voisinages » de [SGA 4, XI]. (Voir aussi [Grothendieck 197?, IV.20 ou V.5 (« Sections hyperplanes et projections coniques ») § 13 (« Morphismes élémentaires et théorème de M. Artin »)], où un tel morphisme est dit « polyélémentaire ».)

**1.3.3.** Soient $X$ un schéma algébrique sur un corps parfait infini $k$ et $x$ un point fermé en lequel $X \to \mathrm{Spec}(k)$ est *lisse*. D'après [SGA 4, XI.3.3], il existe un voisinage ouvert de Zariski $U$ de $x$ qui est une polycourbe élémentaire sur $k$. (Si $k = \mathbb{C}$, l'espace topologique $U(\mathbb{C})$ est un $K(\pi, 1)$ et toute classe de cohomologie en degré $> 0$ est tuée par un revêtement fini étale $U' \to U$ ([SGA 4, XI.4.6]).) On a le raffinement suivant ([Friedlander 1982, théorème 11.7]) : il existe un voisinage *étale* de $x$ qui est un « $K(\pi, 1)$ pro-$\ell$ » (au sens expliqué en ¶ **1.4**), où $\pi$ est un pro-$\ell$-groupe extension itérée de pro-$\ell$-groupes libres de type fini. Dans cette section et la suivante, on donne une démonstration alternative de ce résultat, amélioré par la possibilité d'utiliser le théorème de résolution des singularités de A. J. de Jong.

**1.3.4.** Soit $\ell$ un nombre premier. Une courbe élémentaire $f : X \to S$ est dite **$\ell$-élémentaire** si le faisceau $\mathrm{R}^1 f_\star \mathbb{Z}/\ell\mathbb{Z}$ est $\ell$-monodromique (**1.1.5**). Lorsque $\ell$ est inversible sur les schémas considérés, ce qui est systématiquement le cas dans cet article, les faisceaux $\mathrm{R}^i f_\star \mathbb{Z}/\ell\mathbb{Z}$ sont automatiquement localement constants constructibles, de formation commutant aux changements de base (cf. [SGA 1, XIII.2.9] pour $i \leq 1$, et [SGA 4½, Th. finitude, appendice, 1.3.3] pour le $i \geq 0$ arbitraire) ; l'hypothèse précédente porte donc sur la *monodromie* de $\mathrm{R}^1 f_\star \mathbb{Z}/\ell\mathbb{Z}$ : si $S$ est connexe et $\overline{s}$ en est un point géométrique, on demande que l'image de $\pi_1(S, \overline{s})$ dans $\mathrm{GL}(\mathrm{H}^1(X_{\overline{s}}, \mathbb{Z}/\ell\mathbb{Z}))$ soit un $\ell$-groupe. Une *polycourbe* $f$ est dite $\ell$-élémentaire si elle admet une factorisation en courbes qui sont $\ell$-élémentaires ; nous déduirons de la proposition **1.3.7** ci-dessous que $\mathrm{R}^1 f_\star \mathbb{Z}/\ell\mathbb{Z}$ est alors $\ell$-monodromique.

Avant d'énoncer la proposition principale de cette section, commençons par deux lemmes.

**1.3.5. Lemme.** *Soient $X/S$ une courbe élémentaire (resp. $\ell$-élémentaire), $\ell$ un nombre premier inversible sur $S$ et $Y/X$ un $\ell$-revêtement connexe. Alors, le morphisme $Y \to S$ se factorise à travers un $\ell$-revêtement $S' \to S$ en une courbe élémentaire (resp. $\ell$-élémentaire) $Y \to S'$.*

*Démonstration.* Commençons par l'énoncé non respé. Soit $\overline{X}/S$ une compactification comme en **1.3.1**. Par modération et le lemme d'Abhyaṁkar relatif ([SGA 1, XIII.5.5]), il existe — *localement* sur $\overline{X}$ pour la topologie étale — un prolongement de $Y/X$ en un morphisme de Kummer généralisé relativement à $D \subseteq \overline{X}$ ([Grothendieck et Murre 1971, 1.3.9.c)]). Celui-ci est unique à isomorphisme unique près et sa source est lisse sur $S$. Si, en coordonnées locales, $D = V(t)$ et $Y \to X$ est $t = f(y)$, le diviseur $E$ d'équation $f/f'$ est étale sur $S$. (Si $f(y) = y^n$, avec $n$ inversible sur $S$, le diviseur $E$ est d'équation $y$, comme attendu.)[①] Par unicité (forte), on peut recoller ces morphismes construits étale-localement en un morphisme fini et plat $\overline{Y} \to \overline{X}$, prolongeant $Y \to X$, tel que $\overline{Y} \to S$ soit une courbe propre et lisse et $E = \overline{Y} - Y$ soit un diviseur fini étale sur $S$, à fibres non vides. (Cf. par exemple [Illusie, Laszlo et Orgogozo 2014, IX.2.1] pour une approche log-géométrique, lorsque $S$ est un schéma régulier, cas suffisant dans cet article.) Les fibres de $\overline{Y} \to S$ ne sont pas nécessairement géométriquement connexes mais sa factorisation de Stein répond à la question. (Rappelons que d'après [ÉGA $\mathrm{III}_2$, 7.8.10], la factorisation de Stein d'un morphisme propre et réduit est un revêtement *étale*.)

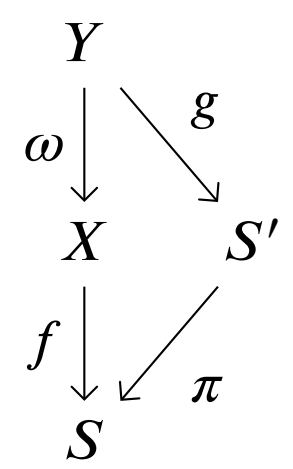

Montrons maintenant que si $f : X \to S$ est $\ell$-élémentaire, il en est de même de $g : Y \to S'$. L'égalité $\pi_\star \mathrm{R}^1 g_\star \mathbb{Z}/\ell\mathbb{Z} = \mathrm{R}^1(f\omega)_\star \mathbb{Z}/\ell\mathbb{Z}$, où $\pi$ (resp. $\omega$) est le $\ell$-revêtement $S' \to S$ (resp. $Y \to X$) nous ramène à montrer que le faisceau $\mathrm{R}^1(f\omega)_\star \mathbb{Z}/\ell\mathbb{Z}$ est $\ell$-monodromique car un faisceau $\mathscr{L}$ sur $S'$ est $\ell$-monodromique si $\pi_\star \mathscr{L}$ l'est. (Utiliser la surjectivité de la coünité $\pi^\star \pi_\star \mathscr{L} \to \mathscr{L}$.) Enfin, le faisceau $\omega_\star \mathbb{Z}/\ell\mathbb{Z}$ étant abélien $\ell$-monodromique, il est extension successive du faisceau constant $\mathbb{Z}/\ell\mathbb{Z}$ et, par conséquent, son image directe par $\mathrm{R}^1 f_\star$ est extension successive de sous-quotients du faisceau $\mathrm{R}^1 f_\star \mathbb{Z}/\ell\mathbb{Z}$, qui est $\ell$-monodromique par hypothèse. □

**1.3.6. Lemme.** *Soient $X \to S$ une courbe élémentaire, $\ell$ un nombre premier inversible sur $S$ et $\mathscr{C}$ un champ en groupoïdes sur $X$. Si $S$ est* strictement local *et si $\mathscr{C}$ est $\ell$-monodromique, il existe un $\ell$-revêtement étale $X' \to X$ tel que $\mathscr{C}' = (X' \to X)^\star \mathscr{C}$ soit isomorphe à un coproduit (fini) de champs de torseurs sous des $\ell$-groupes finis.*

Ceci nous permettra de ramener — *étale-localement sur $S$* — l'étude des champs $\ell$-monodromiques à celle des champs de torseurs sous un $\ell$-groupe fini.

① Notons que si $S$ est réduit, ce qui est suffisant pour nos applications, cela revient à prendre l'image inverse réduite de $D$.

*Démonstration.* D'après le lemme précédent (**1.3.5**) et l'hypothèse sur le champ $\mathscr{C}$, on peut supposer $\pi_0(\mathscr{C})$ constant. Quitte à remplacer $\mathscr{C}$ par une sous-gerbe maximale, on peut supposer que c'est une gerbe. Il suffit de montrer qu'elle a une section : elle sera alors isomorphe à une gerbe $\ell$-monodromique de torseurs et l'on pourra conclure par **1.2.2**. L'existence d'une section résulte de **1.2.5** et du fait que, par modération et commutation aux changements de base, les groupes $\mathrm{H}^2(X', \mathbb{Z}/\ell\mathbb{Z})$ sont nuls pour tout $\ell$-revêtement $X'$ de $X$. (Rappelons que $S$ est *strictement local.*) □

**1.3.7. Proposition.** *Soient $S$ un schéma localement nœthérien sur lequel un nombre premier $\ell$ est inversible et $f : X \to S$ une courbe $\ell$-élémentaire. Alors, pour tout champ $\ell$-monodromique $\mathscr{C}$ sur $X$, le champ image directe $f_\star\mathscr{C}$ est également $\ell$-monodromique. De plus, la formation de cette image directe commute aux changements de base $S' \to S$.*

*Démonstration.* La commutation aux changements de base et la lissité reviennent à montrer que si $S$ est strictement local de point fermé $s$ et $\overline{\eta}$ est un point générique géométrique, les morphismes $\mathscr{C}(X) \to \mathscr{C}(X_s)$ et $\mathscr{C}(X) \to \mathscr{C}(X_{\overline{\eta}})$ sont des équivalences. Par descente (finie étale) et les deux lemmes précédents (**1.3.5**, **1.3.6**), on se ramène au cas particulier où $\mathscr{C}$ est une gerbe de torseurs sous un $\ell$-groupe fini : cela résulte alors des théorèmes bien connus de spécialisation du groupe fondamental (modéré). La constructibilité de l'image directe est également déjà connue ([SGA 1, XIII, § 2]) et résulte d'ailleurs de la démonstration ci-dessous.

Commençons par montrer que pour tout champ $\mathscr{C}$ comme dans l'énoncé, le faisceau $\pi_0(f_\star\mathscr{C})$ est $\ell$-monodromique et traitant tout d'abord le cas d'un *faisceau* (c'est-à-dire d'un champ en catégories discrètes). Par passage à la limite, il suffit de montrer que si $S$ est un schéma (non nécessairement nœthérien) n'admettant pas de revêtement connexe d'ordre $\ell$ (c'est-à-dire : le topos $S_{\ell\text{ét}}$ défini en **1.4.3** est *local*) et $s$ est un point géométrique de $S$, alors le morphisme $\pi_1^{\text{pro-}\ell}(X \times_S S_{(s)}) \to \pi_1^{\text{pro-}\ell}(X)$ est *surjectif*. (Voir par exemple [Szamuely 2009, chap. 5] pour une théorie du groupe fondamental pour les schémas non nécessairement localement nœthériens.) Comme il s'agit de pro-$\ell$ groupes, cela est équivalent à montrer que la flèche induite par application du foncteur $\mathrm{H}^1(-, \mathbb{Z}/\ell\mathbb{Z})$ à $X \times_S S_{(s)} \to X$ est injective (cf. [Serre 1994, I, prop. 23]). Or, on a la suite exacte (Leray)

$$0 \to \mathrm{H}^1(S, f_\star\mathbb{Z}/\ell\mathbb{Z}) \to \mathrm{H}^1(X, \mathbb{Z}/\ell\mathbb{Z}) \to \mathrm{H}^0(S, \mathrm{R}^1 f_\star\mathbb{Z}/\ell\mathbb{Z}),$$

où : $f_\star\mathbb{Z}/\ell\mathbb{Z} = \mathbb{Z}/\ell\mathbb{Z}$ (par connexité des fibres géométriques), $\mathrm{H}^1(S, \mathbb{Z}/\ell\mathbb{Z}) = 0$ (par hypothèse sur $S$), et $\mathrm{R}^1 f_\star\mathbb{Z}/\ell\mathbb{Z}$ est constant sur $S$ (par hypothèse sur $f$) de sorte que $\mathrm{H}^0(S, \mathrm{R}^1 f_\star\mathbb{Z}/\ell\mathbb{Z}) = \mathrm{H}^1(X \times_S S_{(s)}, \mathbb{Z}/\ell\mathbb{Z})$. Supposons maintenant que $\mathscr{C}$ est un champ $\ell$-monodromique quelconque et montrons que $\pi_0(f_\star\mathscr{C})$ est $\ell$-monodromique. D'après ce qui précède, on peut supposer — quitte à remplacer $S$ par un $\ell$-revêtement — le faisceau $f_\star\pi_0(\mathscr{C})$ *constant* puis que le champ $\mathscr{C}$ est une *gerbe* ($\ell$-monodromique), notée $\mathscr{G}$ dorénavant. Soit $\mathscr{G}'$ son quotient par le centre $Z$ de $\mathscr{G}$ (**1.2.3**) et posons $\mathscr{L} \coloneqq \pi_0(f_\star\mathscr{G})$ et $\mathscr{L}' \coloneqq \pi_0(f_\star\mathscr{G}')$. Le faisceau en groupes abéliens $\mathrm{R}^1 f_\star Z$ agit sur le morphisme (localement constant) canonique $\mathscr{L} \to \mathscr{L}'$ et l'action est transitive sur les fibres : $(\mathrm{R}^1 f_\star Z)\backslash\mathscr{L} \to \mathscr{L}'$ est injective. (Ceci est l'analogue du fait bien connu (voir [Giraud 1971, III, 3.4.5] ou [Serre 1994, I, § 5.7]) que si $1 \to Z \to G \to G' \to 1$ est une extension centrale de groupes sur un topos $E$, le groupe abélien $\mathrm{H}^1(E, Z)$ agit sur $\mathrm{H}^1(E, G)$ par produit contracté, c'est-à-dire par $[T] \cdot [P] = [T \overset{Z}{\wedge} P]$, et l'application $\mathrm{H}^1(E, Z)\backslash\mathrm{H}^1(E, G) \to \mathrm{H}^1(E, G')$, de source le quotient à gauche de $\mathrm{H}^1(E, G)$ par $\mathrm{H}^1(E, Z)$, est *injective*.) Le morphisme $f$ étant $\ell$-monodromique, on peut supposer $\mathrm{R}^1 f_\star Z$ constant ; par récurrence, on peut supposer qu'il en est de même de $\mathscr{L}'$. La conclusion résulte alors formellement de l'observation suivante : si un $\ell$-groupe abélien fini $\mathsf{C}$ agit fidèlement et transitivement sur les fibres d'un morphisme $\mathscr{L} \to \mathscr{L}'$ de faisceaux lisses à but constant, le faisceau $\mathscr{L}$ est $\ell$-monodromique. En effet, l'action de la monodromie sur $\mathscr{L}$ se factorise à travers $\mathsf{C}$. Ceci achève la démonstration du fait que $\pi_0(f_\star\mathscr{C})$ est $\ell$-monodromique lorsque $f$ est une courbe $\ell$-élémentaire et $\mathscr{C}$ est $\ell$-monodromique.

Supposons dorénavant le faisceau $\mathscr{L} = \pi_0(f_\star\mathscr{C})$ trivial et considérons une sous-gerbe maximale $\mathscr{G}$ de $f_\star\mathscr{C}$. Il nous faut montrer que son lien est $\ell$-monodromique (**1.1.5**) ; on va utiliser le lemme **1.2.4**. On vérifie en effet que le morphisme $f^\star\mathscr{G} \to \mathscr{C}$ déduit par adjonction de l'inclusion $\mathscr{G} \hookrightarrow f_\star\mathscr{C}$ est fidèle : par commutation des images directes au passage aux fibres, on se ramène à observer que si $G$ est un faisceau localement constant sur $X$ et $x$ un point géométrique d'une fibre géométrique $X_s$, le morphisme $G(X_s) \to G_x$ est injectif. (On utilise la connexité des fibres géométriques du morphisme $f : X \to S$.) D'après *loc. cit.*, il en résulte que la gerbe $f^\star\mathscr{G}$ est $\ell$-monodromique. Pour démontrer que le lien de $\mathscr{G}$ est $\ell$-monodromique, on est donc ramené à vérifier qu'un faisceau $\mathscr{F}$ sur $S$ est $\ell$-monodromique si $f^\star\mathscr{F}$ l'est. La constructibilité résulte de la surjectivité de $f$, cf. p. ex. [SGA 4, IX.2.8], et la lissité n'est guère plus difficile. Que $\mathscr{F}$ soit $\ell$-monodromique résulte, dans le cas où les schémas sont connexes, de la surjectivité du morphisme $\pi_1(X) \to \pi_1(S)$ rappelée en 1.4.7 (premier paragraphe). □

## 1.4. Schémas $K(\pi, 1)$.

**1.4.1.** Commençons par rappeler quelques résultats de [Abbes et Gros 2011, § 9]. Pour chaque schéma $X$, on note $X_{\text{ét}}$ le topos étale, $X_{f\text{ét}}$ le topos *fini étale* (autrement dit, le topos associé à la sous-catégorie pleine des $X$-schémas étales finis munie de la topologie étale) et $\rho : X_{\text{ét}} \to X_{f\text{ét}}$ le morphisme évident ([ibid., § 9.2]). Lorsque $X$ est *cohérent*, le morphisme $\rho$ est un morphisme cohérent entre topos cohérents ([ibid., 9.11, 9.13]). Si l'on suppose de plus l'espace topologique quotient $\pi_0(X) = \operatorname{Spec} \operatorname{Idem}(X)$ discret (autrement dit : les composantes connexes de $X$ ouvertes), le morphisme d'adjonction $\mathrm{Id} \to \mathrm{R}\rho_\star\rho^\star$ induit un isomorphisme en degrés 0 et 1, comme on peut le vérifier par passage aux fibres, c'est-à-dire en supposant $X$ simplement connexe (ou encore : $X_{f\text{ét}}$ local). On dit ([ibid., 9.20]) que le schéma $X$ est un $K(\pi, 1)$ s'il satisfait la condition d'acyclicité suivante : pour tout entier $n \geq 1$ inversible sur $X$ et tout faisceau $\mathscr{L}_f$ de $\mathbb{Z}/n\mathbb{Z}$-modules sur $X_{f\text{ét}}$, le morphisme d'adjonction $\mathscr{L}_f \to \mathrm{R}\rho_\star\rho^\star\mathscr{L}_f$ est un isomorphisme.

**1.4.2.** *Traduction dans le cas connexe.* Si $X$ est connexe, le choix d'un point géométrique $x$ de $X$ identifie le topos $X_{f\text{ét}}$ au topos $\mathsf{B}\pi_1(X, x)$ des $\pi_1(X, x)$-ensembles continus (cf. [ibid., § 9.7]). En particulier, la cohomologie de $X_{f\text{ét}}$ n'est autre que la cohomologie (« galoisienne ») du groupe profini $\pi_1(X, x)$ ([ibid., § 9.7.6]). Si $X$ est un $K(\pi, 1)$, alors pour tout chaque entier $n$ inversible sur $X$ et chaque faisceau étale $\mathscr{L}$ de $\mathbb{Z}/n\mathbb{Z}$-modules localement constant constructible sur $X$, la flèche canonique $\mathrm{R}\Gamma(\pi_1(X, x), \mathscr{L}_x) \to \mathrm{R}\Gamma(X, \mathscr{L})$ est un isomorphisme. La réciproque est vraie lorsque $X$ est cohérent et ceci est encore équivalent au fait qu'un revêtement universel $\widetilde{X}$ de $X$ n'a pas de cohomologie à valeur dans $\mathbb{Z}/n\mathbb{Z}$ en degré $\geq 1$ (ou, de façon équivalente, $> 1$).

**1.4.3.** *Variante pro-$\ell$.* Soit $\ell$ un nombre premier inversible sur un schéma cohérent $X$ à composantes connexes ouvertes et considérons la sous-catégorie pleine des $X$-schémas finis étales dont la monodromie (en tant que faisceau d'ensembles représenté) est, sur chaque composante connexe, un $\ell$-groupe, autrement dit, dont le groupe de Galois sur $X$ d'une clôture galoisienne est un $\ell$-groupe. Munie de la topologie étale, cette catégorie donne lieu à un site, dont on note $X_{\ell\text{ét}}$ le topos associé, naturellement équipé d'un morphisme $\rho_\ell : X_{\text{ét}} \to X_{\ell\text{ét}}$. Par construction, la cohomologie de $X_{\ell\text{ét}}$ est, dans le cas où $X$ est connexe, la cohomologie du groupe fondamental pro-$\ell$ de $X$ (pointé en un point géométrique).

**1.4.4.** Un schéma $X$ comme ci-dessus est un $K(\pi, 1)$ **pro-$\ell$** si pour tout faisceau abélien de $\ell$-torsion ([SGA 4, IX.1.1]) $\mathscr{L}_\ell$ sur $X_{\ell\text{ét}}$, l'unité d'adjonction $\mathscr{L}_\ell \to \mathrm{R}\rho_{\ell\star}\rho_\ell^\star\mathscr{L}_\ell$ est un isomorphisme. En conséquence, si $x$ est un point géométrique de $X$, supposé connexe, et $\mathscr{L}$ est un faisceau abélien $\ell$-monodromique (donc constructible) sur $X$, on a $\mathrm{R}\Gamma(\pi_1^{\text{pro-}\ell}(X, x), \mathscr{L}_x) \xrightarrow{\sim} \mathrm{R}\Gamma(X, \mathscr{L})$ et, réciproquement, ceci caractérise les $K(\pi, 1)$ pro-$\ell$.

**1.4.5.** *Mise en garde.* Notons qu'avec nos définitions, un schéma $K(\pi, 1)$ n'est pas nécessairement un $K(\pi, 1)$ pro-$\ell$ : le noyau du morphisme de pro-$\ell$-complétion (« sous-groupe $\ell$-résiduel ») d'un groupe profini n'est en général pas pro-$\ell'$ et, *a fortiori*, pas nécessairement acyclique pour les coefficients de $\ell$-torsion.

**1.4.6.** *Exemple.* Toute courbe affine lisse $C$ sur un corps algébriquement clos est un $K(\pi, 1)$ pro-$\ell$ pour chaque nombre premier $\ell$. En effet, ni la courbe (affine), ni le pro-$\ell$ groupe fondamental (libre, cf. p. ex. [Wingberg 1984, théorème 1.1] lorsque $\ell$ est inversible sur $C$) n'ont de cohomologie en degré $>$ 1. (On utilise également le fait général que l'unité d'adjonction $\mathrm{Id} \to \mathrm{R}\rho_{\ell\star}\rho_\ell^\star$ est un isomorphisme en degré $\leq 1$ pour les faisceaux de $\ell$-torsion.)

**1.4.7. Proposition.** *Soient $\ell$ un nombre premier inversible sur un corps algébriquement clos $k$ et $X$ une polycourbe $\ell$-élémentaire sur* $\mathrm{Spec}(k)$*. Alors, le schéma $X$ est un $K(\pi,1)$ pro-$\ell$ (cf. **1.4.4**), et le pro-$\ell$ complété du groupe fondamental de chaque composante connexe de $X$ est extension itérée de pro-$\ell$ groupes libres de type fini.*

*Démonstration.* Soit $f : X \to Y$ une polycourbe $\ell$-élémentaire, où $Y$ est une courbe affine connexe lisse sur corps algébriquement clos de caractéristique $\neq \ell$. On note $\overline{\eta}$ un point générique géométrique de $Y$ et on souhaite tout d'abord montrer que la suite de pro-$\ell$ groupes

$$(\dagger) \qquad 1 \to \pi_1^{\text{pro-}\ell}(X_{\overline{\eta}}) \to \pi_1^{\text{pro-}\ell}(X) \to \pi_1^{\text{pro-}\ell}(Y) \to 1$$

est *exacte.* (On omet ici la notation de points bases.) L'exactitude à droite résulte de la surjection $\pi_1(X) \twoheadrightarrow \pi_1(Y)$ ([SGA 1, XIII.4.1]). Soit $\mathsf{G}$ un $\ell$-groupe fini. Rappelons (cf. par exemple [SGA 4, XII, § 1]) que l'on a une suite exacte d'ensembles pointés :

$$1 \to \mathrm{H}^1(Y, f_\star \mathsf{G}) \to \mathrm{H}^1(X, \mathsf{G}) \to \mathrm{H}^0(Y, \mathrm{R}^1 f_\star \mathsf{G}),$$

analogue non abélien de la suite exacte de bas degré usuelle (Leray). Ici, on a l'égalité $f_\star \mathsf{G} = \mathsf{G}$ et le faisceau $\mathrm{R}^1 f_\star \mathsf{G}$ est localement constant de fibre en $\overline{\eta}$ isomorphe à $\mathrm{H}^1(X_{\overline{\eta}}, \mathsf{G})$ de sorte que le noyau de la flèche (d'ensembles pointés) $\mathrm{H}^1(X, \mathsf{G}) \to \mathrm{H}^1(X_{\overline{\eta}}, \mathsf{G})$ est l'image de $\mathrm{H}^1(Y, \mathsf{G})$ : (†) est exacte au centre (voir p. ex. [SGA 1, V, § 6] pour la traduction en terme de groupes fondamentaux).

Notons $\mathscr{C}$ le champ sur $Y$ image directe par $f$ du champ des $\mathsf{G}$-torseurs sur $X$. Rappelons ([Giraud 1971, V.3.1.6]) l'interprétation champêtre de la suite exacte d'ensembles pointés précédente. Le terme de droite s'identifie à l'ensemble des sous-gerbes maximales de $\mathscr{C}$ : à une classe $c \in \mathrm{H}^0(Y, \mathrm{R}^1 f_\star \mathsf{G})$ est associée la gerbe $D(c)$ sur $Y$ dont la fibre en $V \to Y$ est la catégorie des $\mathsf{G}$-torseurs sur $U = X \times_Y V$ dont la classe dans $\mathrm{H}^1(U, \mathsf{G})$ s'envoie sur la restriction $c_{|V} \in \mathrm{H}^0(V, \mathrm{R}^1 f_\star \mathsf{G})$. La gerbe $D(c)$ est triviale si et seulement si $c$ est l'image d'une classe $[T] \in \mathrm{H}^1(X, \mathsf{G})$ et, dans ce cas, $D(c)$ est équivalente au champ des $f_\star \mathsf{Aut}(T)$-torseurs sur $Y$. D'après **1.3.7**, il existe un $\ell$-revêtement étale de $Y' \to Y$ tel que le faisceau $\pi_0(\mathscr{C}')$ des sous-gerbes maximales de $\mathscr{C}' = (Y' \to Y)^\star \mathscr{C}$ soit fini constant et que chacune des sous-gerbes maximales soient *localement* liées par le lien d'un $\ell$-groupe constant. (On utilise la commutation au changement de base de la formation de l'image directe $f_\star \mathsf{BG}$.) D'après **1.2.5**, ces sous-gerbes maximales, sur la *courbe affine* $Y'$, sont triviales. Il en résulte que le foncteur de restriction $\mathscr{C}'(Y') = \mathscr{C}(Y') \to \mathscr{C}(\overline{\eta})$ est essentiellement surjectif : tout $\mathsf{G}$-torseur sur $X_{\overline{\eta}}$ s'obtient par restriction à partir d'un $\mathsf{G}$-torseur sur $X' = X \times_Y Y'$. Ceci suffit pour achever la démonstration de l'exactitude de (†) (cf. [SGA 1, V.6.8 et XIII.4.3]). On utilise le fait qu'une clôture galoisienne *sur* $X$ d'un $\mathsf{G}$-torseur sur $X'$ est un revêtement d'ordre une puissance de $\ell$ : si $\mathsf{H}'' \trianglelefteq \mathsf{H}' \trianglelefteq \mathsf{H}$ sont des groupes finis avec $[\mathsf{H} : \mathsf{H}']$ et $[\mathsf{H}' : \mathsf{H}'']$ des puissances d'un nombre premier $\ell$, le sous-groupe $\bigcap_{h \in \mathsf{H}/\mathsf{H}''} h\mathsf{H}'' h^{-1}$ est d'indice une puissance de $\ell$ dans $\mathsf{H}$.

Pour montrer que $X$ est un $K(\pi,1)$ pro-$\ell$, commençons par constater que si $\mathscr{F}$ est un faisceau abélien $\ell$-monodromique, et $f$ une polycourbe $\ell$-élémentaire, les faisceaux $\mathrm{R}^j f_\star \mathscr{F}$ sont $\ell$-monodromiques. La stabilité par extension des faisceaux abéliens $\ell$-monodromiques (**1.2.1**) nous permet en effet de nous ramener, par récurrence, au cas où $f$ est une courbe $\ell$-élémentaire, auquel cas cela résulte de la proposition **1.3.7**. La conclusion résulte alors de la suite exacte d'homotopie précédente, de l'égalité $\mathrm{R}\Gamma(X, \mathscr{F}) = \mathrm{R}\Gamma(Y, \mathrm{R}f_\star \mathscr{F})$, de l'isomorphisme $\mathrm{R}\Gamma(\pi_1^{\text{pro-}\ell}(Y, \overline{\eta}), (\mathrm{R}f_\star \mathscr{F})_{\overline{\eta}}) \overset{\sim}{\to} \mathrm{R}\Gamma(Y, \mathrm{R}f_\star \mathscr{F})$ (cf. **1.4.6**) et enfin de l'isomorphisme $(\mathrm{R}f_\star \mathscr{F})_{\overline{\eta}} = \mathrm{R}\Gamma(X_{\overline{\eta}}, \mathscr{F}) \overset{\sim}{\leftarrow} \mathrm{R}\Gamma(\pi_1^{\text{pro-}\ell}(X_{\overline{\eta}}), \mathscr{F})$, obtenu par récurrence sur la dimension relative. □

**1.4.8.** Avant d'énoncer le résultat d'abondance des $K(\pi,1)$ pro-$\ell$ ci-dessous, rappelons que P. Deligne a défini dans [Deligne 1974, § 5.3.5] la topologie de la descente cohomologique universelle ([SGA 4, V^bis, § 3.3]) ; pour une introduction pédagogique à cette notion, le lecteur pourra se référer à [Conrad 2003]. Nous appellerons **topologie des altérations** la topologie définie par la prétopologie engendrée par les recouvrements étales et les altérations (= morphisme propre et surjectif induisant un morphisme fini au-dessus d'un ouvert partout dense et envoyant tout point maximal sur un point maximal.) La topologie de la descente cohomologique universelle est plus fine que la topologie des altérations.

**1.4.9. Proposition.** *Soit $\ell$ un nombre premier inversible sur un corps algébriquement clos $k$.* Localement pour la topologie des altérations, *tout $k$-schéma algébrique est une polycourbe $\ell$-élémentaire (et en* *particulier, d'après* **1.4.7**, *localement un* $K(\pi, 1)$ *pro-*$\ell$*). De plus, si le schéma est supposé lisse, c'est même vrai localement pour la topologie* étale.

*Démonstration.* Soit $X$ un schéma comme dans l'énoncé, que l'on peut supposer intègre. D'après [A. J. de Jong 1996, 4.1], il est localement (par une altération) *lisse* sur $k$. D'après [SGA 4, XI.3.3] (rappelé en **1.3.3**), on peut donc supposer le $k$-schéma lisse connexe $X$ est une polycourbe élémentaire. Factorisons $X \to \mathrm{Spec}(k)$ à travers une courbe élémentaire $f : X \to Y$. Il existe un revêtement étale $h : Y' \to Y$ trivialisant $\mathrm{R}^1 f_\star \mathbb{Z}/\ell\mathbb{Z}$, c'est-à-dire tel que le faisceau $h^\star \mathrm{R}^1 f_\star \mathbb{Z}/\ell\mathbb{Z} = \mathrm{R}^1 f'_\star \mathbb{Z}/\ell\mathbb{Z}$ soit trivial (i.e., constant) donc, en particulier, $\ell$-monodromique. Un tel revêtement est calculable : quitte à se placer au-dessus du point générique du schéma normal $Y$, ceci est expliqué en § 2. (Le cas considéré ici est celui, plus facile, d'une courbe ; cf. **2.4** et **2.7**.) Étant obtenu par changement de base de $f$, le morphisme $f'$ est une courbe élémentaire, et même $\ell$-élémentaire par construction. Par récurrence sur la dimension du schéma considéré, on peut supposer qu'il existe, localement pour la topologie étale au voisinage de chaque point, un morphisme $Y'' \to Y'$ tel que le morphisme $g''$ du diagramme commutatif ci-dessous soit une polycourbe $\ell$-élémentaire. Le schéma $X''$ obtenu par changement de base convient.

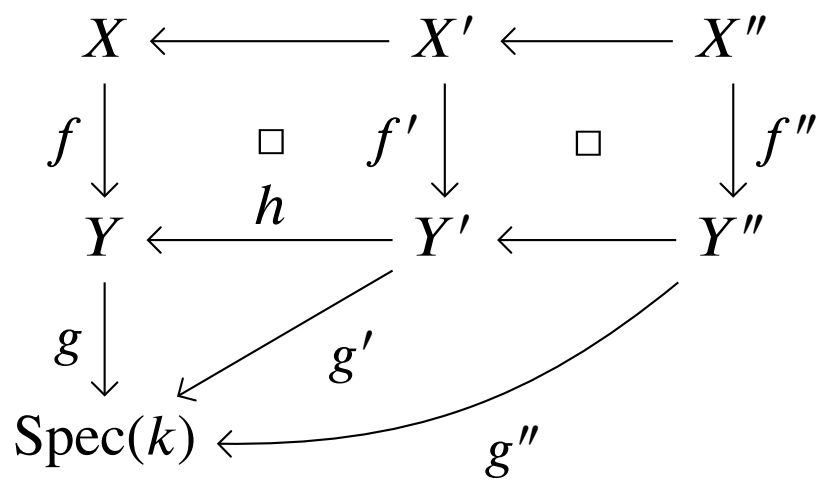

□

**1.4.10.** Notons que l'utilisation du théorème de A. J. de Jong fait perdre l'éventuelle primitive récursivité de notre algorithme (puisqu'on utilise ce théorème en énumérant tous les morphismes sur $X$ jusqu'à trouver une altération qui en résolve les singularités, sans aucune borne sur le temps d'exécution autre que le fait que cette énumération terminera, cf. **11.6** et **12.8**). Il est probable que, dans le cas lisse, on puisse préserver l'éventuelle primitive récursivité de notre algorithme par une analyse précise de [SGA 4, XI, § 2-3].

**1.4.11.** Si $k$ est algébrique séparable sur un corps ${}_0k$ et $X/k$ obtenu par extension des scalaires d'un ${}_0k$-schéma algébrique ${}_0X$, tout recouvrement $\{Y_\alpha \to X\}$ par des polycourbes $\ell$-élémentaires est dominé par un recouvrement du même type *défini sur* ${}_0k$[①]. En effet, si $Y_\alpha \to X$ se descend en un morphisme de ${}_1k$-schémas algébriques ${}_1Y_\alpha \to {}_1X$, où ${}_1k/{}_0k$ est étale, il suffit de considérer le ${}_0k$-morphisme composé ${}_1Y_\alpha \to {}_0X$. Le schéma ${}_1Y_\alpha \otimes_{{}_0k} k$ est isomorphe à $Y_\alpha \otimes_k ({}_1k \otimes_{{}_0k} k)$ ; c'est un coproduit de polycourbes $\ell$-élémentaires, fini surjectif au-dessus de $Y_\alpha$.

Nous utiliserons la proposition précédente sous la forme suivante (cf. [SGA 4, $\mathrm{V}^{\mathrm{bis}}$.5.1] ou bien [Deligne 1974, 5.3.3.1 et § 6.2]).

**1.4.12. Corollaire.** *Soient* $\ell$ *un nombre premier inversible sur un corps algébriquement clos* $k$ *et* $X$ *un* $k$*-schéma algébrique. Alors, il existe un* $X$*-schéma simplicial* $X_\bullet$ *tel que chaque* $X_i$ *soit un coproduit fini de polycourbes* $\ell$*-élémentaires et tel que la flèche d'adjonction* $\mathrm{R}\Gamma(X_{\text{ét}}, \mathbb{Z}/\ell\mathbb{Z}) \to \mathrm{R}\Gamma(\mathrm{Tot}\, X_{\bullet\text{ét}}, \mathbb{Z}/\ell\mathbb{Z})$ *soit un isomorphisme. De plus :*

- *un morphisme* $Y \to X$ *de* $k$*-schémas algébriques peut être coiffé par un morphisme simplicial* $Y_\bullet \to X_\bullet$ *du type précédent ;*
- *si* $k$ *est algébrique séparable sur* ${}_0k$ *et* $X/k$ *obtenu par extension des scalaires d'un* ${}_0k$*-schéma algébrique* ${}_0X$*, on peut supposer pour tout entier* $r$ *que la flèche* $X_{\bullet \leq r} \to X$ *provient par extension des scalaires d'une flèche* ${}_0X_{\bullet \leq r} \to {}_0X$.

[①] Les indices sont mis à gauche pour éviter la confusion avec les indices simpliciaux utilisés ci-dessous.

Le dernier point résulte de l'observation **1.4.11** et de la construction des hyperrecouvrements rappelée ci-dessous, effectuée sur ${}_0k$. Pour la définition du *topos total* Tot $X_{\bullet \text{ét}}$, voir les références citées en § 4. (Voir également **4.1.1** pour une variante « $n$-approchée ».) Pour la définition d'un *hyperrecouvrement*, voir par exemple [SGA 4, V.7.3.1.1] ou bien [Deligne 1974, 5.3.5] (ou, de nouveau, [Conrad 2003] pour une approche pédagogique).

**1.4.13.** Par la suite, nous utiliserons implicitement le fait que les objets et les flèches ci-dessus sont calculables en tout étage : la construction décrite en [SGA 4, $\mathrm{V}^{\mathrm{bis}}$.5.1.3] (voir aussi [Deligne 1974, prop. 6.2.4]) ne fait intervenir que des limites et coproduits finis et explicites de schémas, que l'on « améliore » en utilisant **1.4.9**. Traitons le problème de la construction d'une flèche $Y_\bullet \to X_\bullet$ plus en détail ; le cas de $X_\bullet$ en est un cas particulier. Par produit fibré, il suffit de montrer le lemme suivant.

**Lemme.** *Tout hyperrecouvrement $X_\bullet$ de $X$ peut être dominé par un hyperrecouvrement $X'_\bullet$ à composantes des coproduits de polycourbes $\ell$-élémentaires, calculables en tout étage.*

Ici, « hyperrecouvrement » = « hyperrecouvrement pour la topologie des altérations ».

*Démonstration.* Fixons un entier $r \geq 0$ et supposons donné un morphisme simplicial ($r$-tronqué) $X'_{\bullet\leq r} \to \mathrm{sq}_r X_\bullet$, où $X'_{\bullet\leq r}$ est scindé à composantes des coproduits de polycourbes $\ell$-élémentaires. D'après **1.4.9**, on peut « couvrir » le produit fibré de $X_{r+1}$ et $\mathrm{cosq}_r(X'_{\bullet\leq r})_{r+1}$ au-dessus de $\mathrm{cosq}_r(\mathrm{sq}_r X_\bullet)_{r+1}$ par un coproduit de polycourbes $\ell$-élémentaires, noté $N$. Pour le calcul du cosquelette, on utilise la formule :

$$\mathrm{cosq}_r(Z)_s = \lim_{\substack{k\leq r \\ [k]\to[s]}} Z_k$$

(la limite étant comprise dans la catégorie des $X$-schémas, et calculable d'après **16.3**). D'après [SGA 4, $\mathrm{V}^{\mathrm{bis}}$.5.1.3], il existe un schéma simplicial $r+1$-tronqué scindé (calculable) $X'_{\bullet\leq r+1}$ prolongeant $X'_{\bullet\leq r}$, s'envoyant sur $\mathrm{sq}_{r+1}(X_\bullet)$, tel que $X'_{r+1}$ soit un coproduit de $N$ et de composantes scindées de $X'_{\bullet\leq r}$ ; c'est donc un coproduit de polycourbes $\ell$-élémentaires. □

## 2. Calculabilité du $\mathrm{H}^1$

L'objectif de cette section est de démontrer le théorème suivant.

**2.1. Théorème.** *Soient $X$ un schéma* normal*, de type fini sur un corps algébriquement clos $k$, et $G$ un groupe fini constant d'ordre une puissance d'un nombre premier $\ell$ inversible sur $X$. Alors, on peut calculer $\mathrm{H}^1(X, G)$, c'est-à-dire produire une liste de représentants des classes d'isomorphie de $G$-torseurs sur $X$.*

**2.2.** Il suffit de trouver une extension finie galoisienne du corps des fractions $K$ de $X$ qui domine tous ces $G$-torseurs, ou encore une extension les trivialisant tous. En effet, si $L/K$ est une telle extension, de groupe de Galois $\pi$, l'ensemble $\mathrm{H}^1(X, G)$ est naturellement un sous-ensemble de l'ensemble fini $\mathrm{H}^1(L/K, G) = \mathsf{Hom}(\pi, G)/G$ des classes d'isomorphisme de $G$-torseurs sur $K$ trivialisés par $L/K$. Si $\varphi \in \mathsf{Hom}(\pi, G)$, on peut construire explicitement le $G$-torseur $A_\varphi = \mathsf{Hom}_{\pi\text{-}\mathsf{Ens}}(G, L)$ sur $K$ correspondant par la théorie de Galois-Grothendieck, et tout $G$-torseur sur $X$ est obtenu par normalisation de $X$ dans un tel $A_\varphi$ (voir **15.5** pour la calculabilité de la normalisation). Parmi ces $X$-schémas, en nombre fini, il faut vérifier quels sont ceux qui sont étales et galoisiens de groupe $G$ sur $X$ (cf. **17.2** et **17.3**).

**2.3. Réduction au cas abélien.** Notons $Z$ le centre (non trivial) du $\ell$-groupe fini $G$. La suite exacte

$$1 \to Z \to G \to G/Z \to 1$$

induit ([Giraud 1971, V.2.3]) une suite exacte d'*ensembles pointés*

$$\mathrm{H}^1(X, Z) \to \mathrm{H}^1(X, G) \to \mathrm{H}^1(X, G/Z).$$

Supposons, par récurrence, que l'on sache trouver une extension $L$ du corps des fonctions de $X$ trivialisant les $G/Z$-torseurs. Soit $X_L$ le normalisé de $X$ dans $L$. L'image inverse sur $X_L$ de chaque $G$-torseur sur $X$ provient d'un $Z$-torseur sur $X_L$. (On utilise le fait que la restriction au point générique induit une injection sur les $\mathrm{H}^1$ car les schémas considérés sont normaux.) Ceci nous ramène au cas particulier où $G$ est un $\ell$-groupe abélien et, finalement, au cas où $G = \mathbb{Z}/\ell\mathbb{Z}$, ce qu'on supposera maintenant.

### 2.4. Cas d'une courbe.

**2.4.1.** Si $X$ est une *courbe* lisse sur le corps algébriquement clos $k$, le cardinal de $\mathrm{H}^1(X, \mathbb{Z}/\ell\mathbb{Z})$ est connu. On peut donc effectivement produire tous les $\mathbb{Z}/\ell\mathbb{Z}$-torseurs sur $X$ en un temps fini.

**2.4.2. Remarque.** Signalons comment l'on pourrait se ramener au cas d'une courbe projective (lisse). Si $\overline{X}$ est la complétion projective de $X$ et $\overline{X} - X = \{c_1, \dots, c_r\}$ sont les points à l'infini, il existe (Riemann-Roch) une fonction $f \in K^\times$ s'annulant exactement en ces points. Notons $L = K(\sqrt[\ell]{f})$. D'après le lemme d'Abhyaṁkar, le tiré en arrière d'un $\mathbb{Z}/\ell\mathbb{Z}$-torseur sur $X$ au normalisé $X_L$ de $X$ dans $L$ s'étend à la complétion projective $\overline{X}_L$. Pour trouver $f$, on note qu'il existe un entier $n$ explicite tel que, pour chaque $i \in \{1, \dots, r\}$, il existe une fonction $f_i \in \mathscr{L}((n+1)c_i) - \mathscr{L}(nc_i)$, que l'on peut calculer algorithmiquement, cf. [Heß 2002] ; la fonction $f = \sum_i f_i$ convient. (Une variante de cet argument est également possible lorsque $X$ est un $k$-schéma algébrique normal, quitte à l'altérer pour en faire le complémentaire d'un diviseur à croisements normaux dans un $k$-schéma projectif lisse (de Jong).)

Signalons pour terminer que si $X$ est une courbe projective lisse, on sait ([Serre 1975, chap. VI, nº12]) que tout revêtement connexe de groupe $\mathbb{Z}/\ell\mathbb{Z}$ est induit par une isogénie de la jacobienne de $X$ (que l'on peut construire explicitement ; cf. p. ex. [Anderson 2002]). Ceci fournit une approche différente de 2.4.1 pour la détermination effective de $\mathrm{H}^1(X, \mathbb{Z}/\ell\mathbb{Z})$, qui peut probablement être rendue primitivement récursive.

**2.5. Fibration en courbes.** ① Supposons dorénavant le schéma $X$ de dimension $d \geq 2$ et démontrons 2.1 par récurrence sur l'entier $d$. On peut supposer le schéma $X$ intègre et, quitte à le modifier — c'est-à-dire le remplacer par un $X$-schéma propre et birationnel —, on peut également supposer qu'il existe un $k$-schéma de type fini intègre $S$ de dimension $d-1$ et un morphisme $X \to S$ faisant de $X$ une courbe relative sur $S$ à fibre générique lisse et géométriquement connexe. (Cf. p. ex. [A. J. de Jong 1996, 4.11-12].) Il résulte de ce qui précède — appliqué à la courbe $X_{\overline{\eta}}$ où $\overline{\eta}$ est un point générique géométrique de $S$ — qu'il existe un $S$-schéma étale $S'$ intègre de point générique $\eta'$ et un revêtement $Y' \to X_{S'}$ tel que si $T$ est un $G$-torseur sur $X$ et $T'$ le torseur obtenu par le changement de base $Y' \to X$, la fibre générique géométrique $T'_{\overline{\eta'}}$ est le $G$-torseur trivial. Quitte à changer les notations, on peut supposer $S = S'$ et $X = Y'$.

(Du point de vue algorithmique, signalons que l'on peut déterminer $\eta$ par 15.2 et que les corps $\kappa(\eta)$ et $\kappa(\overline{\eta})$ conservent les bonnes propriétés de calculabilité que l'on peut imposer à $k$ par 12.5.)

**2.6.** Quitte à remplacer $S$ par un ouvert étale, on peut supposer que le morphisme $X \to S$ est une *courbe élémentaire*. Vérifions que, sous cette hypothèse, tout $G$-torseur $T$ sur $X$ à fibre générique géométrique triviale provient de $S$ : ceci nous permettra de conclure, par récurrence, car $\dim(S) = d - 1 < d$. Supposons, comme il a été fait ci-dessus, que $G$ est le groupe abélien $\mathbb{Z}/\ell\mathbb{Z}$, et notons $f$ le morphisme $X \to S$. Considérons la suite exacte

$$0 \to \mathrm{H}^1(S, f_\star \mathbb{Z}/\ell\mathbb{Z}) \to \mathrm{H}^1(X, \mathbb{Z}/\ell\mathbb{Z}) \to \mathrm{H}^0(S, \mathrm{R}^1 f_\star \mathbb{Z}/\ell\mathbb{Z}).$$

Comme rappelé en 1.3.4, le morphisme d'adjonction $\mathbb{Z}/\ell\mathbb{Z} \to f_\star \mathbb{Z}/\ell\mathbb{Z}$ est un isomorphisme, le faisceau $\mathrm{R}^1 f_\star \mathbb{Z}/\ell\mathbb{Z}$ est lisse sur $S$ et la fibre générique géométrique de $\mathrm{R}^1 f_\star \mathbb{Z}/\ell\mathbb{Z}$ est isomorphe à $\mathrm{H}^1(X_{\overline{\eta}}, \mathbb{Z}/\ell\mathbb{Z})$. Il en résulte que la flèche $\mathrm{H}^0(S, \mathrm{R}^1 f_\star \mathbb{Z}/\ell\mathbb{Z}) \to \mathrm{H}^1(X_{\overline{\eta}}, \mathbb{Z}/\ell\mathbb{Z})$ est une injection et, par conséquent, que tout $\mathbb{Z}/\ell\mathbb{Z}$-torseur sur $X$ trivialisé par $X_{\overline{\eta}}$ provient de $S$. (Voir aussi 1.3.7 et 1.4.7.) Ceci achève la démonstration du théorème 2.1.

**2.7. Remarque.** Supposons $X$ obtenu par extension algébrique des scalaires à partir d'un schéma $X_0$ sur $k_0$ et $T^{(1)}, \dots, T^{(N)}$ des représentants des classes d'isomorphie de $G$-torseurs sur $X$. Il existe une sous-extension étale $k_1/k_0$ de $k$, calculable (16.11), telle que les $G$-torseurs précédents soient définis sur $k_1$, c'est-à-dire que chaque $T^{(i)}$ soit $X$-isomorphe à $T_1^{(i)} \times_{X_1} X$, où $X_1 = X_0 \times_{k_0} k_1$ et $T_1^{(i)}$ est un $G$-torseur sur $X_1$. Dans ce cas, l'action naturelle de $\mathsf{Aut}(k/k_0)$ sur $\mathrm{H}^1(X, G)$ se factorise à travers une

① Cette méthode nous a été suggérée par Ofer Gabber.

action explicite du quotient fini $\mathrm{Gal}(k_1/k_0)$ : si (plus généralement) $g \in \mathsf{Aut}_{k_0}(X_1)$, l'action (à droite) de $g$ sur $[T_1] \in \mathrm{H}^1(X_1, G)$ est donnée par $[T_1] \cdot g = [T_1 \times_{X_1,g} X_1]$.

## 3. Série $\ell$-centrale descendante et groupe fondamental

### 3.1. Une filtration.

**3.1.1.** Soient $\ell$ un nombre premier et $G$ un pro-$\ell$ groupe. Pour tout sous-groupe $H$ de $G$, notons $F(H)$ l'adhérence du sous-groupe $H^\ell \cdot (H, G)$ de $G$, où $(X, Y)$ désigne le sous-groupe engendré par les commutateurs $x^{-1}y^{-1}xy$ pour $x$ dans $X$ et $y$ dans $Y$ : autrement dit, il s'agit du sous-groupe fermé de $G$ engendré par les puissances $\ell$-ièmes des éléments de $H$ et des commutateurs de ces éléments avec ceux de $G$. Rappelons ([Dixon et al. 1999, 1.15] ou [Neukirch, A. Schmidt et Wingberg 2000, définition 3.8.1]) que la **série $\ell$-centrale descendante** de $G$ est la filtration définie de la façon suivante : $G^{[1]} = G$, et $G^{[n+1]} = F(G^{[n]})$. Ces sous-groupes sont caractéristiques dans $G$ et on note $G^{(n)}$ le quotient $G/G^{[n]}$. En particulier, $G^{[2]}$ est le sous-groupe de Frattini $\Phi(G)$ de $G$ ([Serre 1978-79, § 3.6] ou [Rotman 1995, théorème 5.48(ii)]) et le $\mathbb{F}_\ell$-espace vectoriel $G^{(2)}$ est naturellement le dual de $\mathsf{Hom}_{\mathrm{cont}}(G, \mathbb{Z}/\ell\mathbb{Z})$.

Pour des raisons typographiques, nous noterons parfois $F^nG$ pour $G^{[n]}$, notamment en § 6.

**3.1.2. Remarques.** On vérifie sans peine que $G^{[n+1]}$ est le plus petit sous-groupe fermé normal de $G$ contenu dans $G^{[n]}$ tel que $G^{[n]}/G^{[n+1]}$ soit $\ell$-élémentaire abélien et contenu dans le centre de $G/G^{[n+1]}$.

Pour notre propos, nous avons une certaine liberté dans le choix de la filtration : on aurait aussi bien pu considérer, par exemple, la filtration de Frattini itérée $\Phi^nG$, où $\Phi(H) = \overline{H^\ell(H, H)} \subseteq F(H)$. À ce propos, signalons que la filtration de Frattini itérée $\Phi^nG$ et la filtration $\ell$-centrale descendante $F^nG$ définie ci-dessus sont équivalentes au sens où il existe une fonction $\tau(d, n)$ (explicitement calculable, et même, primitivement récursive) telle que $F^{\tau(d,n)}G \subseteq \Phi^nG \subseteq F^nG$ si $G$ a $d$ générateurs (comme pro-$\ell$-groupe). (Esquisse de démonstration : l'inclusion $\Phi^nG \subseteq F^nG$ est évidente. L'inclusion $F^{\tau(d,n)}G \subseteq \Phi^nG$ se montre en combinant le (i) du lemme 6.1 ci-dessous avec un l'analogue du (ii) pour la filtration $\Phi^nG$. Cet analogue du (ii) revient à majorer l'ordre de $L/\Phi^nL$ où $L$ est le pro-$\ell$-groupe libre à $d$ générateurs : ceci peut se faire par récurrence sur $n$ en utilisant le fait que $\#(L/\Phi L) = \ell^d$ et le théorème de l'indice de Schreier [Lubotzky et Segal 2003, prop. 16.4.5] pour calculer le nombre de générateurs de $\Phi^nL$.)

Signalons également que si $G$ est topologiquement de type fini, le passage à l'adhérence est superflu, les sous-groupes considérés étant déjà fermés, tant pour $F$ que pour $\Phi$ ([Dixon et al. 1999, 1.20]).

### 3.2. Revêtements et topos associés.

**3.2.1.** Soient $X$ un schéma *connexe*, séparé de type fini sur un corps algébriquement clos $k$, $x$ un point géométrique, $\ell$ un nombre premier inversible sur $k$ et $\pi_X$ le groupe fondamental pro-$\ell$ de $(X, x)$. Pour chaque entier $n \geq 1$, notons $X^{[n]} \to X$ un revêtement étale de $X$ correspondant au sous-groupe d'indice fini $\pi_X^{[n]}$ de $\pi_X$ ; nous dirons que c'est un **revêtement $\ell$-étale $n$-approché universel** de $X$. On note $X_{\ell\text{ét}}^{(n)}$, voire simplement $X^{(n)}$, le *topos* des faisceaux sur $X_{\ell\text{ét}}$ trivialisés par le revêtement étale $X^{[n]} \to X$. Le morphisme naturel $X_{\ell\text{ét}} \to X_{\ell\text{ét}}^{(n)}$ s'identifie au morphisme de topos $\mathsf{B}\pi_X \to \mathsf{B}\pi_X^{(n)}$ déduit de la surjection $\pi_X \twoheadrightarrow \pi_X^{(n)}$ : l'image inverse est le foncteur des $\pi_X^{(n)}$-ensembles vers les $\pi_X$-ensembles continus (obtenu par composition) et l'image directe est le foncteur « invariants sous $\pi_X^{[n]}$ » ; cf. p. ex. [SGA 4, IV.4.5.1]. (On laisse le soin au lecteur de définir $X_{\ell\text{ét}}^{(n)}$ sous des hypothèses plus générales ; nous n'en aurons pas usage.) Par construction, pour chaque $n \geq 2$, le topos $X_{\ell\text{ét}}^{(n)}$ est ponctuel si et seulement si $\pi_X$ est trivial, c'est-à-dire si $X$ ne possède pas de revêtement étale connexe d'ordre $\ell$.

**3.2.2.** *Mutatis mutandis*, les définitions précédentes s'étendent au cas d'un schéma (algébrique sur un corps algébriquement clos) *non nécessairement connexe*. Pour définir $X^{[n]}$, on écrit $X$ comme coproduit de ses composantes connexes (ouvertes) et l'on procède de manière évidente : si $X = \coprod_{c\in\pi_0(X)} X_c$, on pose $X^{[n]} \coloneqq \coprod_{c\in\pi_0(X)} X_c^{[n]}$. La définition du topos $X^{(n)}$ est inchangée et celle du groupe $\pi_X^{(n)}$ devrait être remplacée par la définition d'un *groupoïde* $\Pi_X^{(n)}$.

Notons que $X^{(1)}_{\ell\text{ét}}$ est le topos discret des faisceaux sur l'ensemble (fini) $\pi_0(X)$.

**3.2.3.** Donnons une description alternative de $X^{(n)}$. Commençons par observer qu'il est équivalent au topos des faisceaux sur le site **$\ell$-étale $n$-approché** de $X$ dont les objets sont les $U \to X$ finis étales qui sont, composante connexe par composante connexe, quotients d'un revêtement $\ell$-étale $n$-approché universel de $X$ ; une famille est couvrante si son image recouvre $X$. Fixons un revêtement $\ell$-étale $n$-approché universel $X^{[n]}$ de $X$, et notons $\Pi_{X^{[n]}/X}$ le groupoïde totalement disconnexe (au sens de [Gabriel et Zisman 1967, 6.1.4]) dont les objets sont les $c \in \pi_0(X)$, de groupe d'automorphismes $\mathsf{Aut}_{X_c}(X^{[n]}_c)$. Le foncteur $\mathscr{F} \mapsto \mathscr{F}(X^{[n]})$ induit une équivalence entre le topos $X^{(n)}$ et le topos des préfaisceaux sur $\Pi_{X^{[n]}/X}$. Les foncteurs $\mathscr{F} \mapsto \mathscr{F}(X^{[n]}_c)$ sont des *points* du topos $X^{(n)}$. (Noter qu'un schéma $X^{[n]}_c$ n'est pas nécessairement local au sens de la topologie $\ell$-étale $n$-approchée considérée.)

**3.2.4.** Nous dirons qu'un faisceau d'ensembles constructible $\mathscr{F}$ sur $X^{(n)}$ (resp. un morphisme $\mathscr{F} \to \mathscr{G}$) est **calculable** (en fonction des données) s'il en est ainsi du $X$-schéma $\ell$-étale $n$-approché $X^{\mathscr{F}}$ (resp. du morphisme $X^{\mathscr{F}} \to X^{\mathscr{G}}$) qui le représente (ce schéma ou morphisme étant lui-même décrit au sens de la section 16) ; si le faisceau $\mathscr{F}$ est abélien, on demande de plus que les $X$-morphismes $+ : X^{\mathscr{F}} \times_X X^{\mathscr{F}} \to X^{\mathscr{F}}$ (addition), $[-1] : X^{\mathscr{F}} \to X^{\mathscr{F}}$ (opposé) et $e : X \to X^{\mathscr{F}}$ (neutre) soient calculables. L'ensemble des sections sur un ouvert $\ell$-étale $n$-approché donné $U$ d'un faisceau calculable $\mathscr{F}$, ainsi que les flèches de fonctorialité, sont calculables : il suffit de calculer l'ensemble des sections du morphisme fini étale de schémas $X^{\mathscr{F}} \to X$ au-dessus de $U$, et les flèches induites ; cf. **17.3**. Réciproquement, l'action de la monodromie sur $\mathscr{F}(X^{[n]})$ permet de reconstruire $X^{\mathscr{F}}$.

**3.2.5.** Il résulte immédiatement de ce qui précède que l'on peut calculer le noyau d'un morphisme calculable $\mathscr{F} \to \mathscr{G}$ de faisceaux constructibles abéliens. Il en est de même du conoyau : il est représenté par le plus grand quotient $X'$ de $X^{\mathscr{G}}$ tel que le morphisme de schémas en groupes $X^{\mathscr{F}} \to X'$ se factorise à travers la section nulle $X \to X'$.

**3.2.6.** Nous dirons qu'un faisceau d'ensembles constructible $\mathscr{F}$ sur $X^{(n)}$ est **induit** s'il est isomorphe à un faisceau image directe étale $\pi_{\text{ét}\star} E$, où $\pi$ est le morphisme (fini étale) $X^{[n]} \to X$ et $E$ est un faisceau constant constructible sur chaque composante connexe. Tout faisceau d'ensemble constructible $\mathscr{F}$ sur $X^{(n)}$ est sous-objet d'un faisceau induit : l'unité $\mathscr{F} \to \pi_{\text{ét}\star}\pi^{\star}_{\text{ét}}\mathscr{F}$ est un monomorphisme. Si $\mathscr{F}$ est *calculable*, il en est de même de l'injection précédente.

Rappelons que l'image inverse d'un faisceau représentable est représentée par le produit fibré évident et que si $\mathscr{G}$ est un faisceau d'ensembles constructible sur un schéma $X'$, représenté par un $X'$-schéma $Y'$, le faisceau étale image directe $f_{\text{ét}\star}\mathscr{G}$ par un morphisme *fini étale* $f : X' \to X$ est représenté par même schéma $Y'$, vu sur $X$. On utilise ici l'égalité $f_{\text{ét}\star} = f_{\text{ét}!}$ et la description de ce dernier foncteur faite par exemple en [SGA 4, IV.11.3.1].

### 3.3. Calculabilité de $X^{[n]}$.

**3.3.1.** Soient $k$ un corps algébriquement clos et $X$ un $k$-schéma algébrique *normal*, non nécessairement connexe. Comme précédemment, on a fixé un nombre premier $\ell$ inversible sur $k$. On se propose de montrer que l'on peut calculer un revêtement $\ell$-étale $n$-approché universel $X^{[n]}$ de $X$. Une décomposition de $X$ en composantes connexes étant calculable (cf. **16.6**), on suppose dorénavant le schéma $X$ connexe. Observons maintenant qu'une fois calculé $X^{[n]} \to X$, on peut en déduire $\pi^{(n)}_X$ qui est isomorphe au groupe d'automorphismes $\mathsf{Aut}_X(X^{[n]})$, lui-même en bijection avec $\pi_0(X^{[n]} \times_X X^{[n]})$ par son action sur la composante connexe diagonale (cf. **17.3** ou **16.6**). Réciproquement, si $\pi^{(n)}_X$ — ou même simplement son cardinal — est connu, on peut construire $X^{[n]}$, au pire par recherche non bornée (**12.8**) d'un revêtement étale ayant le bon groupe de Galois. Ci-dessous, nous allons donc nous contenter de calculer le groupe $\pi^{(n)}_X$ dans le cas où $X$ est normal connexe (en fait, il n'est pas difficile de se convaincre qu'on construit bien $X^{[n]}$).

**3.3.2. Remarque.** Par « calcul » d'un groupe fini, on entend la détermination du cardinal et d'une table de multiplication du groupe. (En particulier, on peut en déterminer des générateurs — par exemple, le groupe tout entier — et, pour tout ensemble fini de générateurs, une présentation finie associée.)

**3.3.3.** *$n = 2$.* Comme signalé ci-dessus, le groupe abélien $\pi_X^{(2)}$ est isomorphe au dual du $\mathbb{F}_\ell$-espace vectoriel de dimension finie $\mathrm{H}^1(X, \mathbb{F}_\ell)$. D'après le théorème 2.1, on peut calculer le rang de cet espace vectoriel. En particulier, si $\pi_X$ est un pro-$\ell$-groupe libre de type fini, on peut calculer son rang (nombre minimal de pro-générateurs).

**3.3.4.** *Récurrence sur $n$.* On suppose que l'on sait calculer le groupe $\pi_X^{(n)}$ et par conséquent un revêtement $\ell$-étale $n$-approché universel $X^{[n]} \to X$. Considérons maintenant un revêtement $\ell$-étale 2-approché universel $X^{[n][2]}$ de $X^{[n]}$. C'est également un revêtement galoisien de $X$ car le sous-groupe de Frattini $\Phi(\pi_X^{[n]})$ est caractéristique dans $\pi_X$ donc distingué. Notons $G$ le groupe de Galois de $X^{[n][2]}$ sur $X$, que l'on peut calculer : c'est le groupe d'automorphismes d'un $X$-schéma explicite. Abstraitement, il est isomorphe à $\pi_X/\Phi(\pi_X^{[n]})$ ; le groupe $\pi_X^{(n+1)} = \pi_X/\pi_X^{[n+1]}$ — que l'on cherche à calculer — en est donc un quotient. Plus précisément, par fonctorialité ([Dixon et al. 1999, 1.16 (i)]), $\pi_X^{(n+1)} = G^{(n+1)}$. Ceci permet de conclure.

**3.4. Fonctorialité.** De même que l'on ne peut associer à un espace topologique (connexe localement simplement connexe) un revêtement universel de façon fonctorielle, on ne peut espérer choisir $X^{[n]}$ fonctoriellement en $X$. Tout comme $X \rightsquigarrow X_{\text{ét}}$ (cf. p. ex. [SGA 4, VII.1.4]), la construction $X \rightsquigarrow X_{\ell\text{ét}}^{(n)}$ définit seulement un *pseudo*-foncteur ([SGA 1, VI, § 8] ou [Borceux 1994, 7.5.1]) covariant de la catégorie des schémas vers la 2-catégorie des topos.

**3.4.1.** Soit $f : (X, x) \to (Y, y)$ un morphisme de schémas pointés, entre $k$-schémas algébriques, où $k$ est un corps algébriquement clos sur lequel un nombre premier $\ell$ fixé est inversible. (On ne suppose pas que le morphisme $f$ est un $k$-morphisme.) Soient $n \geq 1$ un entier et $X^{[n]}$ (resp. $Y^{[n]}$) un revêtement $\ell$-étale $n$-approché universel de $X$ (resp. $Y$). Il résulte de la théorie des revêtements et de la fonctorialité de $G \mapsto G^{(n)}$ qu'il existe un morphisme $f^{[n]} : X^{[n]} \to Y^{[n]}$ au-dessus de $f : X \to Y$, c'est-à-dire rendant commutatif le diagramme suivant :

$$\begin{array}{ccc} X^{[n]} & \overset{f^{[n]}}{\dashrightarrow} & Y^{[n]} \\ \downarrow & & \downarrow \\ X & \underset{f}{\longrightarrow} & Y. \end{array}$$

Si $X$ et $Y$ sont connexes, ce morphisme $f^{[n]}$ est unique à translation près par $\mathsf{Aut}_Y(Y^{[n]})$. Si $f$ est un $k$-morphisme explicite entre schémas normaux, on peut le calculer.

Rigidification : choisissons des points géométriques $x^{[n]} \to X^{[n]}$ et $y^{[n]} \to Y^{[n]}$ au-dessus de $x$ et $y$ respectivement. Alors, il existe un *unique* morphisme pointé $f_{xy}^{[n]} : (X^{[n]\circ}, x^{[n]}) \to (Y^{[n]\circ}, y^{[n]})$ au-dessus de $f$, où le terme de gauche d'une paire $(Z^\circ, z)$ désigne la composante connexe de $Z$ contenant $z$.

**3.4.2.** *Version simpliciale.* Soit $X_\bullet$ un $k$-schéma algébrique simplicial tronqué — c'est-à-dire un foncteur d'une catégorie $\mathbf{\Delta}_{\leq r}$ ① ($r \in \mathbb{N}$) vers les $k$-schémas algébriques —, et $P_\bullet \to X_\bullet$ un morphisme simplicial (tronqué), tel que chaque $P_i$ soit coproduit fini de points géométriques. Il existe pour chaque $i$ un relèvement de $P_i \to X_i$ à un revêtement $\ell$-étale $n$-approché universel $X_i^{[n]}$ de $X_i$. Notons $X_{P_i}^{[n]\circ}$ le coproduit $\coprod_{p_i \in P_i} X_{p_i}^{[n]\circ}$, où $X_{p_i}^{[n]\circ}$ désigne la composante connexe de $X_i^{[n]}$ contenant l'image de $p_i$. D'après le paragraphe précédent, les $X_{P_i}^{[n]\circ}$ s'organisent *de façon unique* en un schéma simplicial $X_{P_\bullet}^{[n]\circ}$ de telle sorte que $P_\bullet \to X_{P_\bullet}^{[n]\circ}$ soit une factorisation simpliciale de $P_\bullet \to X_\bullet$.

En d'autres termes, quitte à introduire des « multiplicités », on peut relever simplicialement les points géométriques d'un schéma à un revêtement $\ell$-étale $n$-approché universel.

① C'est la catégorie des ensembles $\{0, \dots, s\} \subseteq \mathbb{N}$, où $s \leq r$, munis des applications croissantes ; elle est notée $(\Delta)_r$ dans [Deligne 1974, 5.1.1].

## 4. Cohomologie $\ell$-étale $n$-approchée d'un schéma simplicial

On vérifie ici le fait — énoncé en 4.3.1 et intuitivement évident compte tenu de ce qui précède (3.3) —, que l'on sait calculer la cohomologie des topos obtenus à partir d'un schéma simplicial tronqué par application du pseudo-foncteur $S \rightsquigarrow S^{(n)}$ (cf. ¶ 3.4). Il est logiquement possible, et peut-être préférable, de ne lire cette section qu'après la section 8. Pour un raccourci, cf. 4.2.3.

**4.1. Généralités.** Nous renvoyons le lecteur à [Illusie 1971-1972, VI, § 5.1 et § 6.2] ou [SGA 4, VI, § 7] pour les détails, et [Deligne 1974, § 5.1] pour un résumé (qui est le point de départ de la théorie) dans le cas des espaces topologiques.

**4.1.1.** *Cohomologie d'un topos simplicial.* Soient $X_\bullet$ un $k$-schéma algébrique simplicial et $n$ un entier $\geq 1$. Notons $X_\bullet^{(n)}$ le topos simplicial (c'est-à-dire fibré sur $\mathbf{\Delta}$) qui s'en déduit par application du pseudo-foncteur $S \rightsquigarrow S^{(n)}$ et $\operatorname{Tot} X_\bullet^{(n)}$ le topos total associé, noté $\operatorname{Top} X_\bullet^{(n)}$ dans [Illusie 1971-1972, VI, § 5.1]. On peut voir un objet $\mathscr{F}_\bullet$ de $\operatorname{Tot} X_\bullet^{(n)}$ comme la donnée pour chaque ouvert $U$ d'un $X_i$ (ouvert : objet du site $\ell$-étale $n$-approché défini en 3.2.3) d'un ensemble $\mathscr{F}_\bullet(U) = \mathscr{F}_i(U)$, fonctoriellement en un sens que nous ne répétons pas ([ibid., VI, § 5.2] ou [Deligne 1974, 5.1.7] ; essentiellement, la fonctorialité est « cosimpliciale en $i$ et faisceautique en $U$ »).

Rappelons par contre un procédé de calcul de la cohomologie de $\operatorname{Tot} X_\bullet^{(n)}$ à valeurs dans un faisceau abélien $\mathscr{F}_\bullet$. Le cas du $\mathrm{H}^0$ (sections globales) est particulièrement simple : c'est $\lim_{i\in\mathbf{\Delta}} \Gamma(X_i^{(n)}, \mathscr{F}_i) = \operatorname{Ker}\big(\Gamma(X_0^{(n)}, \mathscr{F}_0) \rightrightarrows \Gamma(X_1^{(n)}, \mathscr{F}_1)\big)$. Soit $u_\bullet : P_\bullet \to X_\bullet^{(n)}$ un morphisme simplicial tel que pour chaque étage $i$, le topos $P_i$ soit discret (c'est-à-dire coproduit de topos ponctuels) et le morphisme $u_i : P_i \to X_i^{(n)}$ d'image inverse conservative. Notons $\mathscr{F}_i^\bullet = (\mathscr{F}_i^j)$ la résolution flasque (« de Godement ») du faisceau $\mathscr{F}_i$ associée. Le système des $\Gamma(X_i^{(n)}, \mathscr{F}_i^j)$ est cosimplicial en $i$ et différentiel gradué en $j$ ; il fournit un complexe double, la différentielle en $i$ étant la somme alternée usuelle. La cohomologie cherchée est celle du complexe simple associé :

$$\mathrm{R}\Gamma(\operatorname{Tot} X_\bullet^{(n)}, \mathscr{F}_\bullet) \simeq \operatorname{Tot} \Gamma(X_i^{(n)}, \mathscr{F}_i^j).$$

Le terme $\Gamma(X_i^{(n)}, \mathscr{F}_i^j)$ de droite n'est autre que l'ensemble des sections globales de $\mathscr{F}_i^j$, vu comme faisceau étale sur le schéma $X_i$.

On en déduit notamment que $\tau_{<r}\mathrm{R}\Gamma(\operatorname{Tot} X_\bullet^{(n)}, \mathscr{F}_\bullet)$ ne dépend que du schéma simplicial tronqué $X_{\bullet\leq r}$ et de la restriction de $\mathscr{F}_\bullet$ correspondante. Bien que cela ne soit pas absolument nécessaire pour les résultats de cet article, nous précisons en 4.1.2 ci-dessous cette observation.

**4.1.2.** *Variante tronquée (cf. [Gabber 2001]).* Pour tout $r \in \mathbb{N} \cup \{+\infty\}$ et tout faisceau abélien $F_{\bullet\leq r}$ sur $\mathbf{\Delta}_{\leq r}^{\mathrm{op}}$ (c'est-à-dire : groupe abélien cosimplicial tronqué $\mathbf{\Delta}_{\leq r} \to \mathsf{Ab}$), les sections globales dérivées $\mathrm{R}\Gamma(\mathbf{\Delta}_{\leq r}^{\mathrm{op}}, F_{\bullet\leq r})$ sont calculées par le complexe normalisé $N F_{\bullet\leq r} \in \operatorname{Ob} C^{[0,r]}(\mathsf{Ab})$ de la correspondance de Dold-Kan. Pour une définition de ce complexe dans un contexte non tronqué, cf. par exemple [Dold et Puppe 1961, § 3]. Les sections globales dérivées d'un *complexe* (de groupes abéliens cosimpliciaux tronqués) se calculent en prenant le *complexe simple* obtenu par ce procédé. Si $X_{\bullet\leq r}$ est un $k$-schéma algébrique tronqué et $\mathscr{F}_{\bullet\leq r}$ est un faisceau sur $\operatorname{Tot} X_{\bullet\leq r}^{(n)}$, on en déduit en poussant par $X_{\bullet\leq r}^{(n)} \to \mathbf{\Delta}_{\leq r}$ que pour toute résolution $\mathscr{F}_{\bullet\leq r}^\bullet$ à $\mathscr{F}_{i\leq r}^j$ acycliques sur $X_i^{(n)}$, on a

$$\mathrm{R}\Gamma(\operatorname{Tot} X_{\bullet\leq r}^{(n)}, \mathscr{F}_{\bullet\leq r}) = \operatorname{Tot} N_i \Gamma(X_{i\leq r}^{(n)}, \mathscr{F}_i^j)$$

où, pour chaque $j$, on note $N_i\Gamma(X_{i\leq r}^{(n)}, \mathscr{F}_i^j)$ le complexe normalisé déduit du groupe cosimplicial tronqué $i \mapsto \Gamma(X_{i\leq r}^{(n)}, \mathscr{F}_i^j)$. Ceci est compatible avec la description non « normalisée » du paragraphe précédent par le théorème d'Eilenberg-Mac Lane, [ibid., 3.22].

Si $X_\bullet$ (resp. $\mathscr{F}_\bullet$) est une extension de $X_{\bullet\leq r}$ (resp. de $\mathscr{F}_{\bullet\leq r}$) en un $k$-schéma simplicial non tronqué (resp. en un faisceau sur $\operatorname{Tot} X_\bullet^{(n)}$), on a donc un triangle distingué

$$(\text{complexe dans } D^{>r}(\mathbb{Z})) \to \mathrm{R}\Gamma(\operatorname{Tot} X_\bullet^{(n)}, \mathscr{F}_\bullet) \to \mathrm{R}\Gamma(\operatorname{Tot} X_{\bullet\leq r}^{(n)}, \mathscr{F}_{\bullet\leq r}) \overset{+1}{\to}$$

de sorte qu'en particulier la flèche

$$\mathrm{H}^d(\operatorname{Tot} X_\bullet^{(n)}, \mathscr{F}_\bullet) \to \mathrm{H}^d(\operatorname{Tot} X_{\bullet\leq r}^{(n)}, \mathscr{F}_{\bullet\leq r})$$

est un isomorphisme pour $d < r$ et injective pour $d = r$.

4.2. **Résolution de Godement explicite.** Montrons maintenant que les considérations précédentes permettent de calculer les $\mathrm{H}^i(\mathrm{Tot}\, X^{(n)}_{\bullet\leq r}, \mathscr{F}_{\bullet\leq r})$ lorsque les objets sont donnés explicitement.

**4.2.1.** On reprend les notations de **4.1.1** et l'on suppose chacun des étages $X_i$ de $X_{\bullet\leq r}$ *normaux*, afin de pouvoir appliquer les résultats de **3.3**. Observons tout d'abord que, si $P$ est un point géométrique d'un des $X_i$ (obtenu par exemple par le procédé décrit en **16.10**), l'ensemble de ses images par toutes les composées de morphismes de bord et de dégénérescence est fini, de cardinal borné par celui de l'ensemble de toutes les applications croissantes $\{0, \dots, m\} \to \{0, \dots, n\}$ avec $m, n \leq r$. On peut donc choisir un ensemble *fini* $\{P_j\}$ de points géométriques formant un morphisme simplicial $P_{\bullet\leq r} \to X_{\bullet\leq r}$ et tel que chaque composante connexe de chacun des $X_i$ contienne au moins un des $P_j$. Comme expliqué en **3.4.2**, on peut relever $P_{\bullet\leq r}$ en un morphisme de schémas simpliciaux $P_{\bullet\leq r} \to X^{[n]\circ}_{P_{\bullet\leq r}}$, qui nous permet de calculer les images directes et inverses par le morphisme de topos simpliciaux $u_{\bullet\leq r} : P_{\bullet\leq r} \to X^{(n)}_{\bullet\leq r}$, où l'on voit maintenant (abusivement) les $P_i$ comme des topos discrets. Le calcul de l'image inverse $u_i^\star \mathscr{F}$ d'un faisceau $\mathscr{F}$ sur $X_i^{(n)}$ est évident : sa fibre en $p_i$ est $\mathscr{F}(X^{[n]\circ}_{p_i})$. L'image directe $u_{i_\star} E$, pour la topologie $\ell$-étale $n$-approchée, est le faisceau induit (**3.2.6**) image directe usuelle (=étale) par le morphisme fini étale $X^{[n]\circ}_{P_i} \to X_i$ du faisceau constant sur chaque composante connexe correspondant à $E$. Donné $E$, ce faisceau est calculable, étant représenté par un coproduit explicite $X^{[n]\circ}_{P_i,E}$ de copies de $X^{[n]\circ}_{p_i}$. De plus, on peut calculer les $X$-morphismes (addition, neutre, opposé) faisant de $X^{[n]\circ}_{P_i,E}$ un $X$-schéma en groupes fini.

**4.2.2.** Pour tout faisceau $\mathscr{F}_{\bullet\leq r} \in \mathrm{Ob}\,\mathrm{Tot}\, X^{(n)}_{\bullet\leq r}$, le morphisme d'adjonction $\mathscr{F}_{\bullet\leq r} \to \mathscr{G}_{\bullet\leq r} := u_{\bullet\leq r_\star} u^\star_{\bullet\leq r} \mathscr{F}_{\bullet\leq r}$ est le début de la résolution flasque de Godement considérée ci-dessus. Supposons $\mathscr{F}_{\bullet\leq r}$ calculable, c'est-à-dire représentable par un schéma simplicial tronqué en groupes abéliens finis $X^{\mathscr{F}}_{\bullet\leq r}$ au-dessus de $X_{\bullet\leq r}$ qui est calculable (cf. **3.2.4**). Comme expliqué ci-dessus, le faisceau $\mathscr{G}_{\bullet\leq r}$ — induit étage par étage — est également calculable, ainsi que la flèche $\mathscr{F}_{\bullet\leq r} \to \mathscr{G}_{\bullet\leq r}$. D'après **3.2.5**, le conoyau de cette flèche est calculable : on a montré que, donné $\mathscr{F}_{\bullet\leq r}$ — et les données auxiliaires, non canoniques, $P_{\bullet\leq r}$, etc. —, on peut calculer une résolution « de Godement » $\mathscr{C}^{(n)}_P(\mathscr{F}_{\bullet\leq r})$ jusqu'à des degrés arbitrairement grands.

**4.2.3.** *Variante par recherche non bornée.* Si l'on s'autorise à être moins explicite (et perdre l'éventuelle primitive récursivité), on peut procéder plus simplement, c'est-à-dire sans avoir recourt aux schémas simpliciaux $P_{\bullet\leq r}$ et $X^{[n]\circ}_{P_{\bullet\leq r}}$. Fixons comme ci-dessus un corps $k$, un nombre premier $\ell$, un schéma simplicial tronqué $X_{\bullet\leq r}$, un anneau $\Lambda$, et deux entiers $n$ et $d$. Il existe une résolution tronquée par un complexe d'induits $\mathscr{F}_i^j$ (cosimplicial en $i \leq r$, différentiel gradué en $j \leq d$) du faisceau constant de valeur $\Lambda$ sur $X_{\bullet\leq r}$. Cette résolution tronquée est, par hypothèse, acyclique étage par étage et calcule la cohomologie tronquée : $\tau_{<d}\mathrm{R}\Gamma(\mathrm{Tot}\, X^{(n)}_{\bullet\leq r}, \Lambda) = \tau_{j<d}\, \mathrm{Tot}\, N_i \Gamma(X_i, \mathscr{F}_i^j)$. Pour calculer la cohomologie de $\mathrm{Tot}\, X^{(n)}_{\bullet\leq r}$ à valeurs dans $\Lambda$ en degrés strictement inférieurs à $d$, il suffit donc de parcourir les morphismes $\Lambda \to \mathscr{F}^\bullet_\bullet$, où $\mathscr{F}_i^j$ est un complexe tronqué d'induits, et, si c'est une résolution tronquée (c'est-à-dire acyclique en degrés $\leq d$) — fait que l'on sait vérifier (cf. **3.2.5**) —, de calculer le terme de droite correspondant. On peut procéder de même pour calculer les flèches de fonctorialité en $n$ et en $X_{\bullet\leq r}$.

**4.2.4. Remarque.** Ce qui précède peut être vu comme une variante du fait bien connu que l'on sait calculer en chaque degré la cohomologie d'un groupe fini agissant sur un $\Lambda$-module explicite (voir par exemple la définition donnée en ([Serre 1994, I, § 2.2]). Sur le problème de la détermination de l'*algèbre* de cohomologie $\mathrm{H}^\star(G, \Lambda)$ d'un groupe fini $G$, voir par exemple [Carlson 2001]. (D'après le théorème de Venkov-Evens ([Evens 1991, 7.4.6], c'est une $\Lambda$-algèbre de type fini.)

**4.2.5.** La construction précédente d'un résolution de Godement est fonctorielle en $n$ en un sens évident : une fois choisi $P_{\bullet\leq r} \to X^{[n]\circ}_{P_{\bullet\leq r}}$ comme ci-dessus, on a pour tout $m \leq n$ un choix naturel de $P_{\bullet\leq r} \to X^{[m]\circ}_{P_{\bullet\leq r}}$, qui permet de calculer — jusqu'à des degrés arbitrairement grands — la flèche $\rho_m^\star \mathscr{C}^{(m)}_P(\mathscr{F}_{\bullet\leq r}) \to \rho_n^\star \mathscr{C}^{(n)}_P(\mathscr{F}_{\bullet\leq r})$, où $\rho_?$ est le morphisme de topos $\mathrm{Tot}\, X_{\bullet\leq r} \to \mathrm{Tot}\, X^{(?)}_{\bullet\leq r}$.

**4.2.6**. Soit maintenant $f_{\bullet\leq r} : X_{\bullet\leq r} \to Y_{\bullet\leq r}$ un morphisme de $k$-schémas normaux simpliciaux tronqués. Expliquons brièvement comment calculer la flèche induite sur la cohomologie des topos totaux $\ell$-étales $n$-approchés. On commence par produire une flèche simpliciale de points conservatifs $P^X_{\bullet\leq r} \to P^Y_{\bullet\leq r}$ au-dessus de $X_{\bullet\leq r} \to Y_{\bullet\leq r}$, que l'on relève arbitrairement aux schémas $X_i^{[n]}$ et $Y_i^{[n]}$, puis on construit l'unique morphisme simplicial $X^{[n]}_{P^X_{\bullet\leq r}} \to Y^{[n]}_{P^Y_{\bullet\leq r}}$ au-dessus de $X_{\bullet\leq r} \to Y_{\bullet\leq r}$ les respectant. On peut alors calculer un morphisme de résolution de Godement $\mathscr{C}^{(n)}_{P^Y}(\mathscr{F}_{\bullet\leq r}) \to \mathscr{C}^{(n)}_{P^X}(\mathscr{F}_{\bullet\leq r})$ au-dessus du morphisme de topos $\mathrm{Top}\, X^{(n)}_{\bullet\leq r} \to \mathrm{Top}\, Y^{(n)}_{\bullet\leq r}$.

**4.3**. Pour référence ultérieure, résumons les rappels et observations précédentes sous la forme suivante.

**4.3.1**. **Proposition.** *Soient $k$ un corps algébriquement clos, $\ell$ un nombre premier inversible sur $k$, $r \geq 0$ un entier et $X_{\bullet\leq r}$ un $k$-schéma algébrique normal simplicial tronqué. Pour tout $\ell$-groupe abélien fini $\Lambda$ et tout triplet d'entiers $d \geq 0, n \geq m \geq 1$, on peut calculer le complexe $\tau_{\leq d}\mathrm{R}\Gamma(\mathrm{Tot}\, X^{(n)}_{\bullet\leq r}, \Lambda)$ et la flèche $\tau_{\leq d}\mathrm{R}\Gamma(\mathrm{Tot}\, X^{(m)}_{\bullet\leq r}, \Lambda) \to \tau_{\leq d}\mathrm{R}\Gamma(\mathrm{Tot}\, X^{(n)}_{\bullet\leq r}, \Lambda)$. De plus, donné un morphisme $X_{\bullet\leq r} \to Y_{\bullet\leq r}$ de $k$-schémas algébriques normaux simpliciaux tronqués, on peut calculer la flèche $\tau_{\leq d}\mathrm{R}\Gamma(\mathrm{Tot}\, Y^{(n)}_{\bullet\leq r}, \Lambda) \to \tau_{\leq d}\mathrm{R}\Gamma(\mathrm{Tot}\, X^{(n)}_{\bullet\leq r}, \Lambda)$.*

Par « calcul » d'un complexe ou d'une flèche, on entend par là que l'on peut trouver des représentants explicites au sens de 13.3.

**4.3.2**. **Remarque.** Considérant le cas particulier $m = 1$, on voit que l'on peut calculer les flèches

$$\tau_{\leq d}\check{\mathrm{R}}\Gamma(X_{\bullet\leq r}, \Lambda) \to \tau_{\leq d}\mathrm{R}\Gamma(\mathrm{Tot}\, X^{(n)}_{\bullet\leq r}, \Lambda),$$

où le terme de gauche est le tronqué du complexe associé au $\Lambda$-module cosimplicial tronqué des $\Gamma(X_i, \Lambda)$. Cela résulte du fait que $X^{(1)}_{\bullet\leq r} = \pi_0(X_{\bullet\leq r})$ (voir 3.2.2).

## 5. Systèmes essentiellement constants

**5.1**. Soit $\mathscr{A}$ une catégorie abélienne satisfaisant les conditions AB3 et AB5 ([Grothendieck 1957, I]) — qui garantissent l'existence et l'exactitude des colimites filtrantes — et $A_\bullet = (A_i)$ un système inductif de $\mathscr{A}$ indexé par $\mathbb{N}$. Notons $A_\infty$ la colimite de $A_\bullet$ et, pour chaque $j \leq k$ dans $\mathbb{N} \cup \{\infty\}$, posons $A(j,k) := \mathrm{Im}(A_j \to A_k)$. Si $j^- \leq j \leq k \leq k^+$, on a naturellement des flèches $A(j^-, k) \hookrightarrow A(j,k) \twoheadrightarrow A(j,k^+)$. Si $A_\infty$ est *nœthérien*, il existe un $j$ tel que $A(j,\infty) = A_\infty$ ; si de plus $A_j$ est nœthérien, il existe un $k \geq j$ tel que $A(j,k) \xrightarrow{\sim} A(j,\infty) = A_\infty$.

**5.2**. **Définition.** Soient $\mathscr{A}$ une catégorie abélienne, $i_0$ un entier, et $c : \mathbb{N}_{\geq i_0} \to \mathbb{N} \times \mathbb{N}$ une fonction. On dit qu'un système inductif $A_\bullet = (A_i)$ de $\mathscr{A}$, indexé par $\mathbb{N}_{\geq i_0}$, est $c$-**essentiellement constant** si la colimite $A_\infty = \mathrm{colim}_i\, A_i$ est représentable dans $\mathscr{A}$ et si pour chaque $i \in \mathbb{N}_{\geq i_0}$, tel que $c(i) = (j,k)$, l'inégalité $i \leq j \leq k$ est satisfaite et la flèche canonique $A(j,k) \to A_\infty$ est un *isomorphisme.*

**5.3**. Variante : soient $N \geq i_0$ deux entiers et $\varphi : \mathbb{N}_{\geq i_0} \to \mathbb{N}$ une fonction telle que $\varphi(j) \geq j$ pour tout $j \in \mathbb{N}_{\geq i_0}$. On dit qu'un système inductif $A_{\bullet \geq i_0}$ est $(N,\varphi)$-**essentiellement constant** si :

- (**i**) pour chaque $j$, le système inductif $A(j,k)_{k\geq\varphi(j)}$, à flèches de transition *a priori* épimorphiques, est constant ;
- (**ii**) le système inductif $A(j,\varphi(j))_{j\geq N}$, à flèches de transition *a priori* monomorphiques, est constant.

Explicitement : (i) $\mathrm{Ker}(A_j \to A_{\varphi(j)}) \xrightarrow{\sim} \mathrm{Ker}(A_j \to A_k)$ pour $k \geq \varphi(j)$ et (ii) $\mathrm{Im}(A_N \to A_{\varphi(j)}) \xrightarrow{\sim} \mathrm{Im}(A_j \to A_{\varphi(j)})$ pour $j \geq N$.

Nous laissons le soin au lecteur de vérifier que, donné $c$, on peut calculer une paire $(N,\varphi)$ telle que tout système inductif $c$-essentiellement constant soit $(N,\varphi)$-essentiellement et que, réciproquement, donnée $(N,\varphi)$, on peut calculer un $c$.

**5.4**. Il résulte des observations précédentes que si $A_\bullet$ est un système inductif constitué d'objets nœthériens et à colimite nœthérienne, il existe une telle fonction. Nous dirons, de façon un peu vague, qu'un système inductif est **explicitement essentiellement constant** s'il est $(N,\varphi)$-essentiellement constant pour un entier $N$ et une fonction $\varphi$ *calculables* en fonction des données. Cette notion apparaît, avec un but semblable, dans [Schön 1991], puis [Rubio et Sergeraert 2002, § 2.1], où un tel système inductif est appelé « module de Schön ». Voir également [Grothendieck 1956, p. 3][①].

**5.5**. Notons qu'un tel système inductif est essentiellement constant au sens usuel : il appartient à l'image essentielle du plongement de $\mathscr{A}$ dans la catégorie abélienne (cf. p. ex. [Kashiwara et Schapira 2006, 8.6.5 (i)]) des ind-objets $\mathrm{Ind}(\mathscr{A})$. La proposition clef suivante est le pendant « explicite », du fait que la catégorie $\mathscr{A}$ est naturellement une sous-catégorie *épaisse* de $\mathrm{Ind}(\mathscr{A})$ (cf. p. ex. [ibid., 8.6.11]).

**5.6. Proposition.** *Soit* $0 \to A'_\bullet \to A_\bullet \to A''_\bullet \to 0$ *une suite exacte de systèmes inductifs. Si deux des trois termes sont explicitement essentiellement constants, il en est de même du troisième.*

Cette proposition est élémentaire et bien connue (cf. [Schön 1991, lemme 5, p. 4] ou [Rubio et Sergeraert 2002, théorème 2.3]). Pour la commodité du lecteur, nous en donnons une démonstration, dans le cas d'une catégorie de modules pour simplifier l'exposition.

*Démonstration.* Supposons d'abord que $(A_n)$ et $(A''_n)$ soient respectivement $(N,\varphi)$- et $(N'',\varphi'')$-essentiellement constants. La première condition de 5.3 est vérifiée de $(A'_n)$ pour la fonction $\varphi$ : en effet, un élément de $A'_n$ qui s'annule dans $A'_m$ pour $m \geq \varphi(n)$ s'annule en particulier dans $A_m$ (d'après la même condition sur $(A_n)$) donc s'annule dans $A_{\varphi(n)}$ donc dans $A'_{\varphi(n)}$. Soit maintenant $N' = \varphi''(N)$ (qui est supérieur ou égal à $N$) : si $x$ appartient à $A'_n$, son image dans $A'_{\varphi(n)}$ vue dans $A_{\varphi(n)}$ est l'image d'un élément $y$ de $A_N$ (d'après la deuxième condition sur le système $(A_n)$) : l'image de ce $y$ dans $A''_N$ s'annule dans $A''_{\varphi(n)}$, c'est-à-dire appartient à $\mathrm{Ker}(A''_N \to A''_{\varphi(n)})$, et la première condition sur $(A''_n)$ entraîne que cette image s'annule dans $A''_{N'}$, donc l'image de $y$ dans $A_{N'}$ provient d'un élément de $A'_{N'}$, qui par construction a la même image dans $A'_{\varphi(n)}$ que l'élément $x$ qu'on s'était fixé. On a donc montré que le système $(A'_n)$ était $(N',\varphi)$-essentiellement constant (pour $N' = \varphi''(N)$).

Supposons maintenant que $(A_n)$ et $(A'_n)$ soient respectivement $(N,\varphi)$- et $(N',\varphi')$-essentiellement constants. Soit $z \in A''_n$ avec $n \geq N'$, et soit $m \geq n$ tel que l'image de $z$ s'annule dans $A''_m$ : alors, si $y$ est un relèvement quelconque de $z$ à $A_n$, l'image de $y$ dans $A''_m$ s'annule, donc l'image de $y$ dans $A_m$ provient d'un $x \in A'_m$ ; puisque $(A'_n)$ est essentiellement constant, il existe $x_0 \in A'_n$ tel que $x_0$ et $x$ aient même image dans $A'_{\varphi'(m)}$ ; alors $y' := y - x_0$ (vu comme élément de $A_n$) a une image nulle dans $A_{\varphi'(m)}$ : donc l'image de $y'$ dans $A_{\varphi(n)}$ est déjà nulle (d'après la première condition sur $(A_n)$), mais ceci implique que l'image de $z$ dans $A''_{\varphi(n)}$ est nulle. Ceci montre la première condition pour $(A''_n)$, pour la fonction $\varphi''$ égale à $\max(\varphi, N')$. S'agissant de la seconde condition, si $z \in A''_n$ et si $y$ en est un relèvement quelconque à $A_n$, il existe un $\tilde{y}$ dans $A_N$ tel que $y$ et $\tilde{y}$ aient la même image dans $A_{\varphi''(n)}$, et alors l'image $\tilde{z}$ de $\tilde{y}$ dans $A''_N$ a la même image que $z$ dans $A''_{\varphi''(n)}$. On a donc montré que le système $(A''_n)$ était $(N,\varphi'')$-essentiellement constant (pour $\varphi'' = \max(\varphi, N')$).

Enfin, supposons que $(A'_n)$ et $(A''_n)$ soient respectivement $(N',\varphi')$- et $(N'',\varphi'')$-essentiellement constants. Soit $y \in A_n$ qui s'annule dans $A_m$ pour $m \geq n$ : alors en particulier son image dans $A''_m$ s'annule, donc elle s'annule déjà dans $A''_{\varphi''(n)}$ (d'après la première propriété sur $(A''_n)$) ; donc l'image de $y$ dans $A_{\varphi(n)}$ provient d'un élément $x$ de $A'_{\varphi''(n)}$ ; si $m \geq \varphi''(n)$, l'image de $x$ dans $A'_m$ s'annule et s'annule donc déjà (d'après la première propriété sur $(A'_n)$) dans $A'_{\varphi'(\varphi''(n))}$. Ceci montre la première condition sur $(A_n)$ pour la fonction $\varphi : n \mapsto \varphi'(\varphi''(n))$. Enfin, soit $y \in A_n$ : son image dans $A''_n$ a la même image dans $A''_{\varphi''(n)}$ qu'un certain élément $\tilde{z} \in A''_{N''}$, donc si $\tilde{y}$ est un relèvement quelconque de $\tilde{z}$ à $A_{N''}$, les éléments $y$ et $\tilde{y}$ (de $A_n$ et $A_{N''}$ respectivement) ont même image dans $A_{\varphi''(n)}$, donc la différence entre ces images provient d'un élément $x \in A'_{\varphi''(n)}$ ; ce dernier a la même image dans $A'_{\varphi(n)} = A'_{\varphi'(\varphi''(n))}$ qu'un certain élément $\tilde{x} \in A'_{N'}$ : si on appelle $N$ le maximum de $N'$ et $N''$ alors la somme des images de $\tilde{x}$ et $\tilde{y}$ dans $A_{\varphi(n)}$ est la même que celle de $y$. On a donc montré que le système $(A_n)$ était $(N,\varphi)$-essentiellement constant pour $\varphi = \varphi' \circ \varphi''$ et $N = \max(N', N'')$. ☐

---

[①] Nous remercions Luc Illusie de nous avoir communiqué cette référence.

**5.7.** Soient $\mathscr{A}$ une catégorie abélienne, et $(\mathrm{E}_{r,\lambda}^{\star,\star})_{\lambda\in\mathbb{N}}$ un système inductif (indexé par $\lambda$) de suites spectrales ($r \geq r_0$) d'objets de $\mathscr{A}$, supposées dans le premier quadrant, dont on note, pour chaque indice $\lambda$, l'aboutissement $\mathrm{E}_\lambda^\star$. (Suivant p. ex. [ÉGA III$_1$, $\mathbf{0}_{\mathrm{III}}$, § 11.1], on considère que cet objet filtré de $\mathscr{A}$ fait partie de la donnée.)

**5.8. Corollaire.** *Soit $m$ un entier tel que les systèmes inductifs $(\mathrm{E}_{r_0,\lambda}^{p,q})_\lambda$ soient explicitement essentiellement constants pour chaque paire d'entiers $p, q$ d'entiers tels que $p + q \leq 2m + 1$. Alors, pour chaque $0 \leq d \leq m$, le système inductif $(\mathrm{E}_\lambda^d)_\lambda$ est explicitement essentiellement constant.*

*Démonstration.* Pour chaque indice $\lambda$, le calcul de $\mathrm{E}_{r,\lambda}^{p,q}$ ne fait intervenir que des flèches entre sous-quotients de $\mathrm{E}_{r_0,\lambda}^{p',q'}$ avec $p' + q' \leq p + q + (r - r_0)$. Comme d'autre part $\mathrm{E}_{\infty,\lambda}^{p,q} = \mathrm{E}_{r,\lambda}^{p,q}$ si $r > p + q + 1$, il résulte de la proposition précédente (5.6) que les systèmes inductifs $(\mathrm{E}_\infty^{p,q})_\lambda$ pour $p + q < m$ sont explicitement essentiellement constants. Enfin, comme pour chaque $0 \leq d \leq m$ l'aboutissement $\mathrm{E}_\lambda^d$ est une extension itérée de ces $\mathrm{E}_{\infty,\lambda}^{p,q}$ la conclusion résulte d'une nouvelle application de *loc. cit.* □

**5.9. Proposition.** *Soient $A_\bullet \to B_\bullet$ un morphisme de systèmes inductifs et $\tau : \mathbb{N} \to \mathbb{N}$ une fonction strictement croissante. Supposons qu'il existe un diagramme commutatif*

$$\begin{array}{ccc} A_\bullet & \longrightarrow & B_\bullet \\ \downarrow & \overset{h_\bullet}{\swarrow} & \downarrow \\ A_{\tau(\bullet)} & \longrightarrow & B_{\tau(\bullet)}. \end{array}$$

*Alors, si $B_\bullet$ est $(N, \varphi)$-essentiellement constant, le système inductif $A_\bullet$ est $(\tau N, \tau\varphi)$-essentiellement constant. En particulier, lorsque $\tau$ est calculable, le système inductif $A_\bullet$ est* explicitement *essentiellement constant si $B_\bullet$ l'est.*

*Démonstration.* Soit $j$ un entier. Remplaçons, dans le carré commutatif de l'énoncé, le système inductif $A_\bullet$ (resp. $B_\bullet$, etc.) par le système $A'_\bullet \coloneqq A(j, \bullet \geq \varphi(j))$ — à morphismes de transition épimorphiques — (resp. $B'_\bullet \coloneqq B(j, \bullet \geq \varphi(j))$ — constant par hypothèse —, etc.). La commutativité du diagramme montre alors que $h'_\bullet$ est un isomorphisme ; le système inductif $(A'_k)_{k\geq\varphi(j)}$ est donc constant pour $k \geq \tau\varphi(j)$. Posons $\psi = \tau\varphi$. Le même argument, appliqué à $A''_\bullet \coloneqq A(j, \psi(j))_j$, etc. montre que ce système est constant pour $j \geq \tau(N)$. (On utilise le fait que $B$ est $(N, \psi)$-essentiellement constant car $\psi \geq \varphi$.) □

## 6. Approximation d'un pro-$\ell$-groupe par ses quotients finis

Pour $\pi$ un $\ell$-groupe fini, on rappelle qu'on a défini en § 3 la filtration $\ell$-centrale descendante par $F^1\pi = \pi$ et $F^{n+1}\pi = (F^n\pi)^\ell \cdot (\pi, F^n\pi)$ (groupe topologiquement engendré).

**6.1. Lemme.** *Il existe deux fonctions calculables $\varphi_\ell$ et $\psi_\ell$ telles que :*

- ***(i)*** *si $\pi$ est un $\ell$-groupe fini d'ordre $\leq n$ alors $F^{\varphi_\ell(n)}\pi = 1$, et*
- ***(ii)*** *si $\pi$ est un $\ell$-groupe fini à $d$ générateurs tel que $F^n\pi = 1$, alors $\#\pi \leq \psi_\ell(d, n)$.*

Il résulte de la démonstration que $\ell^{(d+1)^n}$ convient pour $\psi_\ell(d, n)$, et que 1 plus la valuation $\ell$-adique de $n$ convient pour $\varphi_\ell(n)$.

*Démonstration.* Pour ce qui est de $\varphi$ : pour chaque $\ell$-groupe fini il existe un $r$ tel que $F^r\pi = 1$ (cf. [Neukirch, A. Schmidt et Wingberg 2000, proposition 3.8.2]) ; or comme $F^{i+1}\pi$ est défini en fonction de $\pi$ et de $F^i\pi$, deux termes consécutifs de la suite $F^i\pi$ ne peuvent pas être égaux sauf à ce que cette suite stationne, et on vient de dire que ceci ne se produit que pour $F^i\pi = 1$ : il en résulte que la valuation $\ell$-adique de l'ordre de $F^i\pi$ doit décroître strictement jusqu'à atteindre 0, donc la valuation $\ell$-adique de $n$ (plus 1, puisque la filtration $F^i\pi$ est numérotée à partir de 1) convient pour $\varphi_\ell(n)$.

Pour ce qui est de $\psi$ : d'après [Lubotzky et Segal 2003, théorème 3.5.1], si $L$ est le pro-$\ell$-groupe libre sur $d \geq 2$ générateurs et $N$ un sous-groupe distingué ouvert de $L$ d'indice $\ell^s > 1$, si on note $N' = N^\ell \cdot (N, L)$, on a $(N : N') \leq \ell^{(d-1)s+1}$, de sorte que $(L : N') \leq \ell^{ds+1} \leq \ell^{(d+1)s}$ ; en appliquant ceci à $N = F^rL$ et par récurrence sur $r$ on en déduit $(L : F^rL) \leq \ell^{(d+1)^r}$ (on vérifie immédiatement que cette inégalité fonctionne encore pour $d = 1$ et $r = 1, 2$). Par conséquent, si $\pi$ est un $\ell$-groupe fini à $d$ générateurs tel que $F^n\pi = 1$, en considérant $L \twoheadrightarrow \pi$ la surjection donnée par ces $d$ générateurs, on a une surjection $L/F^nL \twoheadrightarrow \pi$, donc $\#\pi \le \ell^{(d+1)^n}$. □

**6.2.** Il résulte de ce lemme que pour chaque $\ell$ et chaque $d$, on peut calculer (au sens de 3.3.2) le $\ell$-groupe fini $L^{(n)}$ quotient du pro-$\ell$ groupe libre $L$ à $d$ générateurs : parmi les groupes $\pi$ comme en (ii) ci-dessus, c'est celui ayant le plus gros cardinal. (Notons qu'ici, il est *a priori* trivial de déterminer une présentation finie de $L^{(n)}$ : c'est le quotient du pro-$\ell$-groupe libre $L$ par $L^{[n]}$. Par Frattini, le nombre minimal de générateurs de $L^{(n)}$ est $d$ ; par Golod-Šafarevič ([Serre 1994, I, § 4.4]), le nombre de relations entre ces générateurs est $> \frac{1}{4}d^2$.) En conséquence, on peut — pour chaque entier $n$ — déterminer explicitement le système projectif tronqué $L^{(\bullet\le n)}$ et un système compatible de $d$ générateurs.

**6.3. Proposition** (Lemme d'Artin-Rees-Frattini effectif). *Il existe une fonction $\tau_\ell$ calculable telle que, si $1 \to \pi' \to \pi \to \pi'' \to 1$ est une suite exacte courte de pro-$\ell$-groupes, où $\pi', \pi''$ ont respectivement $d', d''$ générateurs, on a $F^{\tau_\ell(d',d'',n)}\pi \cap \pi' \subseteq F^n\pi' \subseteq F^n\pi \cap \pi'$ pour tout $n$.*

*Démonstration.* Il est évident que $F^n\pi' \subseteq F^n\pi \cap \pi'$. On souhaite montrer que, réciproquement, $F^n\pi' \supseteq F^{\tau(n)}\pi \cap \pi'$ pour une certaine fonction $\tau$ explicitement calculable (dépendant du nombre $d', d''$ de générateurs de $\pi', \pi''$, mais pas d'autres données).

Expliquons pourquoi on peut supposer que $\pi''$ est libre (en tant que pro-$\ell$-groupe) : il existe en tout cas un morphisme surjectif $L \twoheadrightarrow \pi''$ où $L$ est le pro-$\ell$-groupe libre sur $d''$ générateurs ; et quitte à relever à $\pi$ les images par ce morphisme de chacun des générateurs, on peut le factoriser comme la composée d'un morphisme $s\colon L \to \pi$ et de la surjection donnée $\pi \twoheadrightarrow \pi''$. Soit $\hat\pi = \pi \times_{\pi''} L$ l'ensemble des éléments de $\pi \times L$ dont les deux composantes ont même image dans $\pi''$ (la première projection est donc un morphisme surjectif $\hat\pi \twoheadrightarrow \pi$ qui se restreint à l'identité sur $\pi'$) : ce $\hat\pi$, qui s'inscrit dans une suite exacte $1 \to \pi' \to \hat\pi \to L \to 1$, se décrit aussi comme le produit semidirect $\hat\pi = \pi' \rtimes_* L$ par l'action de $L$ sur $\pi'$ donnée par $z * x = s(z)\, x\, s(z)^{-1}$. Si on a montré la conclusion voulue pour la suite exacte $1 \to \pi' \to \hat\pi \to L \to 1$, la même vaut encore pour $1 \to \pi' \to \pi \to \pi'' \to 1$ (puisque l'image de $F^N\hat\pi$ dans $\pi$ est contenue dans, et même égale à, $F^N\pi$).

On peut donc bien supposer que $\pi''$ est libre, et qu'il existe une section $s\colon \pi'' \to \pi$, qui fait de $\pi$ le produit semidirect $\pi = \pi' \rtimes_* \pi''$ où $*$ désigne l'action de $\pi''$ sur $\pi'$ définie par $z * x = s(z)\, x\, s(z)^{-1}$, et $\pi = \pi' \rtimes_* \pi''$.

Fixons $n$. On veut montrer qu'il existe $N$ tel que $F^n\pi' \supseteq F^N\pi \cap \pi'$, et expliquer pourquoi $N$ se calcule sous la forme $\tau(d', d'', n)$ en fonction de $d', d''$ et $n$.

L'action de $\pi''$ sur $\pi'$ stabilise $F^n\pi'$, donc définit une action sur $\pi'/F^n\pi'$, et on a $(\pi/F^n\pi') = (\pi'/F^n\pi') \rtimes_* \pi''$ pour cette action quotient.

Comme $\pi'/F^n\pi'$ est fini, $\mathsf{Aut}(\pi'/F^n\pi')$ est lui-même fini, et comme $\pi'' \to \mathsf{Aut}(\pi'/F^n\pi')$ (donné par $*$) est continu, et que les $F^m\pi''$ forment un système fondamental de voisinages de l'unité dans $\pi''$, il existe $m$ tel que $F^m\pi''$ agisse trivialement sur $\pi'/F^n\pi'$ (cf. [Neukirch, A. Schmidt et Wingberg 2000, proposition 3.8.2]). On peut être plus précis : on a $\#(\pi'/F^n\pi') \le \psi(d', n)$ avec les notations du lemme, donc $\#\mathsf{Aut}(\pi'/F^n\pi') \le \psi(d', n)!$, donc $\varphi(\psi(d', n)!\,)$ convient (en considérant l'image de $\pi''$ dans $\mathsf{Aut}(\pi'/F^n\pi')$) — ce qui nous importe est qu'un $m$ qui convient puisse être calculé en fonction de $d'$ et $n$.

L'action de $\pi''$ sur $\pi'/F^n\pi'$ passe donc au quotient par $F^m\pi''$, c'est-à-dire définit une action de $\pi''/F^m\pi''$ sur $\pi'/F^n\pi'$, et on a $(\pi/F^n\pi')/s(F^m\pi'') = (\pi'/F^n\pi') \rtimes_* (\pi''/F^m\pi'')$ pour cette action quotientée.

Notons $\bar\pi$ ce $\ell$-groupe fini $\pi/((F^n\pi')\cdot s(F^m\pi'')) = (\pi'/F^n\pi') \rtimes_* (\pi''/F^m\pi'')$. Son ordre est majoré par $\psi(d', n) \times \psi(d'', m)$ (et rappelons que $m = \varphi(\psi(d', n)!\,)$ convient).

Il existe alors $N \ge n, m$ tel que $F^N\bar\pi = 1$ : précisément, $\varphi(\psi(d', n) \times \psi(d'', m))$ convient pour $N$. On a alors $F^N\pi \subseteq (F^n\pi')\cdot s(F^m\pi'')$, donc $F^N\pi \cap \pi' \subseteq F^n\pi'$, ce qu'on voulait démontrer. □

**6.4. Corollaire.** *On reprend les hypothèses et les notations de la proposition. Soient $\tilde\pi'^{(n)} = \pi'/\big(\pi' \cap \pi^{[n]}\big)$ le noyau de la surjection naturelle $\pi^{(n)} \twoheadrightarrow \pi''^{(n)}$ et $\Lambda$ un groupe abélien. Pour tout entier $j$, si le système inductif $\mathrm{H}^j(\pi'^{(n)}, \Lambda)$ est explicitement essentiellement constant, il en est de même de $\mathrm{H}^j(\tilde\pi'^{(n)}, \Lambda)$. (La fonction explicitant ce fait intervenir uniquement $d', d''$ et celle explicitant le fait que $\mathrm{H}^j(\pi'^{(n)}, \Lambda)$ est essentiellement constant.)*

*Démonstration.* Cela résulte de la proposition précédente, réécrite sous la forme d'un diagramme commutatif (pour chaque $n$)

$$\begin{array}{ccc} \tilde{\pi}'^{(n)} & \twoheadleftarrow & \pi'^{(n)} \\ \big\uparrow\!\!\uparrow & \nearrow\!\!\!\nearrow\  & \big\uparrow\!\!\uparrow \\ \tilde{\pi}'^{(\tau n)} & \twoheadleftarrow & \pi'^{(\tau n)} \end{array}$$

et de 5.9. □

**6.5. Proposition.** *Soient $L$ un pro-$\ell$-groupe libre à $d$ générateurs topologiques, $n_0 \geq 1$ un entier, $\Lambda$ un $\ell$-groupe abélien fini et $V$ un $\Lambda$-module de type fini muni d'une action explicite de $L^{(n_0)}$. Pour tout entier $i$, le système inductif $\mathrm{H}^i(L^{(n)}, V)$, $n \geq n_0$, est explicitement essentiellement constant.*

Par « action explicite », on entend la donnée d'une présentation explicite (13.3) de $V$ et de $d$ éléments de $\mathsf{Aut}(V)$ (13.4) satisfaisant des relations explicites décrivant $L^{(n_0)}$ (cf. 6.2). Le module $V$ est naturellement muni, pour chaque $n \geq n_0$, de l'action de $L^{(n)}$ déduite de la surjection $L^{(n)} \twoheadrightarrow L^{(n_0)}$.

*Démonstration.* Distinguons trois cas :

$i = 0$. Le système inductif $\mathrm{H}^0(L^{(n)}, V)$ étant *constant*, il est $(n_0, \mathrm{Id})$-essentiellement constant (5.3).

$i = 1$. Rappelons (cf. p. ex. [Serre 1994, I, § 2.6 b)]) que les flèches $\mathrm{H}^1(L^{(n)}, V) \to \mathrm{H}^1(L^{(n+1)}, V)$ sont *injectives*, de sorte que la propriété 5.3 (i) est satisfaite pour $\varphi = \mathrm{Id}$. (Notons que $V$ est fixe par $L^{[n]}/L^{[n+1]}$.) Il reste à trouver $N \geq n_0$ tel que la flèche (injective) $\mathrm{H}^1(L^{(N)}, V) \to \mathrm{H}^1(L, V)$ soit un isomorphisme ou, de façon équivalente, tel que l'on ait l'égalité $\#\mathrm{H}^1(L^{(N)}, V) = \#\mathrm{H}^1(L, V)$. La conclusion résulte du fait que ces cardinaux sont calculables. Pour le terme de gauche c'est clair : on sait calculer le système projectif $L^{(n)}$ ; pour le terme de droite, rappelons ([Ogg 1962, p. 188]) que l'on a une suite exacte

$$0 \to \mathrm{H}^0(L, V) \to V \to V^d \to \mathrm{H}^1(L, V) \to 0$$

si bien que l'on a l'égalité (formule « d'Euler-Poincaré », due à Ogg et Šafarevič) $\#\mathrm{H}^1(L, V) = (\#V)^{d-1} \times \#\mathrm{H}^0(L, V)$.

$i \geq 2$. La colimite $\mathrm{H}^i(L, V)$ étant nulle (cf. [Serre 1994, I, § 3.4]), il suffit de trouver $\varphi : \mathbb{N} \to \mathbb{N}$ telle que $\mathrm{H}^i(L^{(n)}, V) \to \mathrm{H}^i(L^{(\varphi(n))}, V)$ soit nulle pour chaque $n \geq n_0$ et de poser, par exemple, $N = n_0$. Une telle fonction $\varphi$ existe et est calculable car les objets et les flèches le sont.

□

**6.6. Remarque.** Bien que cela ne soit pas nécessaire — sauf pour ne pas perdre la primitive récursivité — signalons que l'on peut être plus précis. Avec les notations de l'énoncé, on a pour chaque $N \geq n \geq n_0$ un morphisme de la suite exacte

$$0 \to \mathrm{H}^1(L^{(n)}, V) \to \mathrm{H}^1(L, V) \to \mathrm{H}^1(L^{[n]}, V)$$

vers la suite exacte

$$0 \to \mathrm{H}^1(L^{(N)}, V) \to \mathrm{H}^1(L, V) \to \mathrm{H}^1(L^{[N]}, V),$$

où $L^{[n]}$ (resp. $L^{[N]}$) agit trivialement sur $V$ et les flèches sont les flèches de fonctorialité évidentes. Par chasse au diagramme, l'injection $\mathrm{H}^1(L^{(N)}, V) \hookrightarrow \mathrm{H}^1(L, V)$ est un isomorphisme si la flèche « verticale » $\mathrm{H}^1(L^{[n]}, V) \to \mathrm{H}^1(L^{[N]}, V)$ est nulle. Il suffit pour cela que l'on ait l'inclusion $L^{[N]} \subseteq L^{[n][2]}$. À $n$ fixé, un tel $N$ peut être obtenu à partir de la suite exacte

$$1 \to L^{[n]} \to L \to L^{(n)} \to 1$$

par application de la proposition 6.3 et de l'estimation du rang du groupe libre $L^{[n]}$ par la formule de l'indice de Schreier ([Serre 1977, I, § 3.4] ou [Rotman 1995, théorème 11.45]). Considérons maintenant le cas $i \geq 2$. La suite spectrale de Hochschild-Serre associée à la suite exacte précédente dégénère en $\mathrm{E}_3$, car $\mathrm{H}^j(L^{[n]}, V)$ est nul pour $j > 1$. Comme l'aboutissement est nul — pour la même raison — en degré cohomologique $> 1$, la flèche $d_2 : \mathrm{H}^{i-2}\big(L^{(n)}, \mathrm{H}^1(L^{[n]}, V)\big) \to \mathrm{H}^i(L^{(n)}, V)$ est surjective pour chaque $i \geq 2$. Ceci est bien entendu valable pour chaque $N \geq n$. Il en résulte que pour tuer la flèche $\mathrm{H}^i(L^{(n)}, V) \to \mathrm{H}^i(L^{(N)}, V)$, il suffit de tuer $\mathrm{H}^1(L^{[n]}, V) \to \mathrm{H}^1(L^{[N]}, V)$. C'est ce que l'on a fait ci-dessus.

**6.7. Remarque.** Il serait intéressant de calculer la plus petite fonction $\varphi$ telle que les flèches $\mathrm{H}^i(L^{(n)}, \mathbb{Z}/\ell\mathbb{Z}) \to \mathrm{H}^i(L^{(\varphi(d,i,n))}, \mathbb{Z}/\ell\mathbb{Z})$ soient nulles. (Lorsque $d = 1$ — cas d'un pro-$\ell$-groupe abélien libre — la fonction $n \mapsto n+1$ convient.) Nous ignorons la réponse à cette question, mais nous indiquons un argument, duquel nous sommes redevables à Jean-Pierre Serre, qui montre que si $i = 2$ et si $L^{[n]}$ désigne maintenant la filtration de Frattini itérée $\Phi^n L$ plutôt que la filtration $F^n L$ considérée ci-dessus (cf. 3.1.2 à ce sujet), et bien sûr $L^{(n)} = L/\Phi^n L$, alors la flèche fonction $n \mapsto n+1$ convient, autrement dit que la flèche $\mathrm{H}^i(L^{(n)}, \mathbb{Z}/\ell\mathbb{Z}) \to \mathrm{H}^i(L^{(n+1)}, \mathbb{Z}/\ell\mathbb{Z})$ est nulle :

En comparant les suites spectrales de Hochschild-Serre associées aux suites exactes $1 \to L^{[n]} \to L \to L^{(n)} \to 1$ et $1 \to L^{[n]}/L^{[n+1]} \to L^{(n+1)} \to L^{(n)} \to 1$, le morphisme évident de la seconde suite spectrale vers la première donne :

$$\begin{array}{ccccc}
\mathrm{H}^0(L^{(n)}, \mathrm{H}^1(L^{[n]}/L^{[n+1]}, \mathbb{Z}/\ell\mathbb{Z})) & \longrightarrow & \mathrm{H}^2(L^{(n)}, \mathbb{Z}/\ell\mathbb{Z}) & \longrightarrow & \mathrm{H}^2(L^{(n+1)}, \mathbb{Z}/\ell\mathbb{Z}) \\
\downarrow & & \Vert & & \downarrow \\
\mathrm{H}^0(L^{(n)}, \mathrm{H}^1(L^{[n]}, \mathbb{Z}/\ell\mathbb{Z})) & \longrightarrow & \mathrm{H}^2(L^{(n)}, \mathbb{Z}/\ell\mathbb{Z}) & \longrightarrow & \mathrm{H}^2(L, \mathbb{Z}/\ell\mathbb{Z}) = 0.
\end{array}$$

Or la flèche canonique de $\mathrm{H}^1(L^{[n]}/L^{[n+1]}, \mathbb{Z}/\ell\mathbb{Z}) = \mathrm{Hom}(L^{[n]}/L^{[n+1]}, \mathbb{Z}/\ell\mathbb{Z})$ vers $\mathrm{H}^1(L^{[n]}, \mathbb{Z}/\ell\mathbb{Z}) = \mathrm{Hom}(L^{[n]}, \mathbb{Z}/\ell\mathbb{Z})$ est un isomorphisme puisque nous avons pris la filtration où $L^{[n+1]}$ est le Frattini de $L^{[n]}$. Comme tout élément de $\mathrm{H}^2(L^{(n)}, \mathbb{Z}/\ell\mathbb{Z})$ se relève à $\mathrm{H}^0(L^{(n)}, \mathrm{H}^1(L^{[n]}, \mathbb{Z}/\ell\mathbb{Z}))$, on en conclut que son image dans $\mathrm{H}^2(L^{(n+1)}, \mathbb{Z}/\ell\mathbb{Z})$ est nulle.

## 7. Calcul de la cohomologie d'une polycourbe $\ell$-élémentaire

**7.1.** Soit $X$ une polycourbe $\ell$-élémentaire sur $\mathrm{Spec}(k)$, où $k$ est un corps algébriquement clos, que l'on peut supposer factorisée en courbes $\ell$-élémentaires ($X = X_m \to X_{m-1} \to \cdots \to X_1 \to \mathrm{Spec}\, k$ où $\dim X_i = i$). D'après 1.4.7, c'est un $K(\pi, 1)$ pro-$\ell$, où $\pi$ est le pro-$\ell$ complété du groupe fondamental de $X$, qui est extension itérée de pro-$\ell$ groupes libres de type fini. L'objectif de cette section est de montrer que pour chaque $\ell$-groupe abélien fini $\Lambda$ et chaque $d \geq 0$ on peut déterminer $(N_d, \varphi_d)$ tels que le système inductif $\big(\mathrm{H}^d(\pi^{(n)}, \Lambda)\big)_{n \geq 0}$ — dont $\mathrm{H}^d(X, \Lambda)$ est la colimite — est $(N_d, \varphi_d)$-essentiellement constant au sens de 5.3.

**7.2. Dévissage.** On raisonne par récurrence sur la dimension $m$ de $X$. Par hypothèse (cf. 1.4.7, démonstration), le groupe $\pi$ s'insère dans une suite exacte $1 \to \pi' \to \pi \to \pi'' \to 1$, où $\pi''$ est un pro-$\ell$ groupe libre (non abélien) et $\pi'$ est une extension itérée de tels groupes. Cette suite exacte est d'*origine géométrique*, c'est-à-dire déduite de morphismes calculables de schémas (normaux connexes) comme ci-dessous, par application du foncteur « groupe fondamental pro-$\ell$ ».

$$\begin{array}{rccl}
\pi & X & \longleftarrow & X_{\overline{\eta}} \quad \pi' \\
 & \downarrow & & \downarrow \\
[\text{pro-}\ell\text{-libre}]\ \pi'' & Y & \longleftarrow & \overline{\eta} \\
 & \downarrow & & \\
 & k & &
\end{array}$$

(où $Y = X_1$ est une courbe $\ell$-élémentaire et $\overline{\eta}$ un point générique géométrique de celle-ci ; soulignons que $X_{\overline{\eta}} \to \overline{\eta}$ est encore une polycourbe $\ell$-élémentaire, cette fois de dimension $m-1$). Notons que l'on peut calculer le nombre de pro-générateurs (3.3.3) de $\pi''$ et $\pi'$, qui apparaissent dans le lemme d'Artin-Rees-Frattini effectif 6.3. (Pour $\pi'$, on peut procéder par récurrence ou bien utiliser la calculabilité du $\mathrm{H}^1$.)

Fixons $j$. D'après 3.3 et 3.4.1, on peut calculer pour chaque $n \geq 1$ la suite exacte $1 \to \tilde{\pi}'^{(n)} \to \pi^{(n)} \to \pi''^{(n)} \to 1$ (de groupes finis) considérée en 6.4 et, en particulier, calculer $\tilde{\pi}'^{(n)} = \pi'/(\pi' \cap \pi^{[n]})$. L'hypothèse de récurrence permet d'affirmer que le système inductif $\mathrm{H}^j(\pi'^{(n)}, \Lambda)$ est explicitement essentiellement constant. D'après *loc. cit.*, il en est de même de $V_n := \mathrm{H}^j(\tilde{\pi}'^{(n)}, \Lambda)$. (Comme rappelé en 4.2.4, on sait calculer chacun de ces différents groupes de cohomologie.)

**7.3. Cas d'une courbe.** Soit $V = \operatorname{colim}_n V_n$ ; c'est un $\Lambda$-module de type fini. Il résulte du caractère explicitement essentiellement constant de la colimite que l'on peut calculer $V$ ainsi que l'action induite d'un quotient explicite $\pi''^{(n_0)}$ de $\pi''$. Fixons $i$. D'après 6.5, on peut calculer un couple $(M, \psi)$ tel que le système $\mathrm{H}^i(\pi''^{(n)}, V)$, $n \geq n_0$, soit $(M, \psi)$-essentiellement constant. On veut montrer que, quitte à changer $M$ et $\psi$, il en est de même du système inductif $\mathrm{H}^i(\pi''^{(n)}, V_n)$. Par hypothèse, il existe un entier $N \geq n_0$ et une fonction strictement croissante $\varphi : \mathbb{N} \to \mathbb{N}$ tels que $(V_n)_n$ soit $(N, \varphi)$-essentiellement constant ; en particulier, le morphisme $(V_n)_{n \geq N} \to (V_{\varphi(n)})_{n \geq N}$ se factorise à travers le morphisme $(V_n)_n \to (V)_n$, où $(V)_n$ est le système inductif constant de valeur $V$. Passant à la cohomologie, on en déduit un diagramme commutatif

$$\begin{array}{ccc} \mathrm{H}^i(\pi''^{(\bullet)}, V_\bullet) & \longrightarrow & \mathrm{H}^i(\pi''^{(\bullet)}, V) \\ \downarrow & \overset{h_\bullet}{\swarrow} & \downarrow \\ \mathrm{H}^i(\pi''^{(\varphi\bullet)}, V_{\varphi\bullet}) & \longrightarrow & \mathrm{H}^i(\pi''^{(\varphi\bullet)}, V). \end{array}$$

D'après 5.9, le système inductif $\mathrm{H}^i(\pi''^{(\bullet)}, V_\bullet)$ est $(\varphi M, \varphi\psi)$-essentiellement constant.

**7.4. Suite spectrale de Hochschild-Serre.** Revenons maintenant au calcul de la cohomologie du schéma $X$. On a

$$\mathrm{R}\Gamma(X, \Lambda) = \mathrm{R}\Gamma(\pi, \Lambda) = \mathrm{R}\Gamma\big(\pi'', \mathrm{R}\Gamma(\pi', \Lambda)\big),$$

que l'on approxime par

$$\mathrm{R}\Gamma(\pi^{(n)}, \Lambda) = \mathrm{R}\Gamma\big(\pi''^{(n)}, \mathrm{R}\Gamma(\tilde{\pi}'^{(n)}, \Lambda)\big)$$

D'après [Serre 1994, I, § 2.6], on a pour chaque entier $\lambda \geq 1$ une suite spectrale

$$\mathrm{E}_{2,\lambda}^{i,j} = \mathrm{H}^i\big(\pi''^{(\lambda)}, \mathrm{H}^j(\tilde{\pi}'^{(\lambda)}, \Lambda)\big) \Rightarrow \mathrm{H}^{i+j}(\pi^{(\lambda)}, \Lambda).$$

Il résulte de 5.8 et de ce qui précède que pour chaque entier $d \geq 0$ on peut calculer $(N_d, \varphi_d)$ tels que le système inductif $\mathrm{H}^d(\pi^{(\bullet)}, \Lambda)$ soit $(N_d, \varphi_d)$-essentiellement constant. En particulier, on peut trouver deux entiers $\alpha \leq \beta$ tels que

$$\mathrm{H}^d(\pi, \Lambda) = \operatorname{Im}\big(\mathrm{H}^d(\pi^{(\alpha)}, \Lambda) \to \mathrm{H}^d(\pi^{(\beta)}, \Lambda)\big).$$

Ces objets sont donc algorithmiquement calculables.

**7.5. Synthèse.** Résumons la situation de cette section sous la forme du diagramme suivant avec des notations légèrement différentes pour expliciter la récurrence :

$$\begin{array}{ccccc}
\pi = \pi^{\langle 0\rangle} & \pi^{\langle 1\rangle} & \cdots & \pi^{\langle m-2\rangle} & \pi^{\langle m-1\rangle} \\[2em]
X_m & X^{\langle 1\rangle}_{m-1} & \cdots & X^{\langle m-2\rangle}_2 & X^{\langle m-1\rangle}_1 \\
\downarrow & \downarrow & & \downarrow & \downarrow{\scriptstyle \pi''^{\langle m-1\rangle} = \pi^{\langle m-1\rangle}} \\
X_{m-1} & X^{\langle 1\rangle}_{m-2} & \cdots & X^{\langle m-2\rangle}_1 & \operatorname{Spec} k^{\langle m-1\rangle} \\
\downarrow & \downarrow & & \downarrow{\scriptstyle \pi''^{\langle m-2\rangle}} & \\
\vdots & \vdots & & \operatorname{Spec} k^{\langle m-2\rangle} & \\
\downarrow & \downarrow{\scriptstyle \pi''^{\langle 1\rangle}} & & & \\
X_1 & \operatorname{Spec} k^{\langle 1\rangle} & & & \\
\downarrow{\scriptstyle \pi''^{\langle 0\rangle}} & & & & \\
\operatorname{Spec} k & & & &
\end{array}$$

Si $X = X_m \to X_{m-1} \to \cdots \to X_1 \to \operatorname{Spec} k$ est la polycourbe $\ell$-élémentaire de départ, on appelle $X^{\langle i\rangle}_{m-i} \to X^{\langle i\rangle}_{m-i-1} \to \cdots \to X^{\langle i\rangle}_1 \to \operatorname{Spec} k^{\langle i\rangle}$ sa fibre au-dessus d'un point générique géométrique de $X_i$, et $\pi^{\langle i\rangle}$ le groupe fondamental pro-$\ell$ de cette fibre $X^{\langle i\rangle}_{m-i}$, ainsi que $\pi''^{\langle i\rangle}$ celui de $X^{\langle i\rangle}_1$. Les $\pi''^{\langle i\rangle}$ sont des groupes pro-$\ell$-libres dont on peut calculer le nombre de générateurs ; les $\pi^{\langle i\rangle}$ s'inscrivent dans des suites exactes $1 \to \pi^{\langle i+1\rangle} \to \pi^{\langle i\rangle} \to \pi''^{\langle i\rangle} \to 1$, permettant de calculer leur nombre de générateurs, et des fonctions explicitant le fait que les $\mathrm{H}^d(\pi^{\langle i\rangle}, \Lambda)$ sont essentiellement constants.

## 8. Descente

**8.1**. Soient $k$ un corps algébriquement clos, $X$ un $k$-schéma algébrique (supposé décrit comme en **16.2**), et $\ell$ un nombre premier inversible sur $k$. D'après **1.4.12**, il existe un $X$-schéma simplicial $X_\bullet$ calculant la cohomologie étale de $X$ à coefficients dans le $\ell$-groupe abélien fini $\Lambda$, et dont les constituants sont des coproduits de polycourbes $\ell$-élémentaires. Par descente cohomologique et **1.4.4**, les flèches ci-dessous sont des isomorphismes :

$$\mathrm{R\Gamma}(X_{\text{ét}}, \Lambda) \xrightarrow{\sim} \mathrm{R\Gamma}(\mathrm{Tot}\, X_{\bullet\text{ét}}, \Lambda) \xleftarrow{\sim} \mathrm{R\Gamma}(\mathrm{Tot}\, X_{\bullet\ell\text{ét}}, \Lambda),$$

où $\mathrm{Tot}\, X_{\bullet\text{ét}}$ (resp. $\mathrm{Tot}\, X_{\bullet\ell\text{ét}}$) désigne le topos total associé au système simplicial des topos $X_{i\,\text{ét}}$ (resp. $X_{i\,\ell\text{ét}}$), $i \geq 0$. Le second isomorphisme résulte du fait que les images directes entre topos simpliciaux se calculent étage par étage, si bien que l'adjonction est un isomorphisme si elle l'est sur chaque étage (cf. [Deligne 1974, 5.2.5] ou [Illusie 1971-1972, VI.5.8.1 (iii)]).

**8.2**. Fixons un entier $d \geq 0$ puis un entier $r > d$. D'après les observations précédentes et **1.4.13**, il existe un $X$-schéma simplicial $X_\bullet$, à tronqué (=squelette) $X_{\bullet\leq r}$ calculable, tel que $\mathrm{H}^d(X_{\text{ét}}, \Lambda) = \mathrm{H}^d(\mathrm{Tot}\, X_{\bullet\ell\text{ét}}, \Lambda)$ et tel que les $X_i$ soient des coproduits finis de $k$-polycourbes $\ell$-élémentaires.

*Pour alléger les notations, on omet dorénavant les indices « $_{\ell\text{ét}}$ » et « $_{\text{ét}}$ ».*

Pour chaque entier $\lambda \geq 2$, on a une suite spectrale ([ibid., VI.6.2.3.2])

$$\mathrm{E}_{1,\lambda}^{i,j} = \mathrm{H}^j(X_i^{(\lambda)}, \Lambda) \Rightarrow \mathrm{H}^{i+j}(\mathrm{Tot}\, X_\bullet^{(\lambda)}, \Lambda),$$

où les topos $X_i^{(\lambda)}$ sont comme définis en **3.2**. Il résulte donc de l'exactitude des colimites filtrantes ([ÉGA $\mathrm{III}_1$, **0**.11.1.8]) et des isomorphismes $\mathrm{colim}_\lambda\, \mathrm{H}^j(X_i^{(\lambda)}, \Lambda) \xrightarrow{\sim} \mathrm{H}^j(X_i, \Lambda)$, que la cohomologie de $\mathrm{Tot}\, X_\bullet$ est colimite de la cohomologie des $\mathrm{Tot}\, X_\bullet^{(\lambda)}$. D'autre part, d'après les résultats de § 7 (cas d'une polycourbe $\ell$-élémentaire), on peut calculer $(N_1, \varphi_1)$ tels que les systèmes inductifs $(\mathrm{E}_{1,\lambda}^{i,j})_\lambda$ soient $(N_1, \varphi_1)$-essentiellement constants pour chaque $i, j \geq 0$ tels que $i + j \leq 2d + 1$. Il en résulte (**5.8**) que l'on peut calculer $(N_\infty, \varphi_\infty)$ tels que le système inductif $\mathrm{H}^d(\mathrm{Tot}\, X_\bullet^{(\lambda)}, \Lambda)$ soit $(N_\infty, \varphi_\infty)$-essentiellement constants ; en particulier, on peut calculer deux entiers $\mu \leq \nu$ tels que l'on ait

$$\mathrm{H}^d(X, \Lambda) = \mathrm{Im}\left(\mathrm{H}^d(\mathrm{Tot}\, X_\bullet^{(\mu)}, \Lambda) \to \mathrm{H}^d(\mathrm{Tot}\, X_\bullet^{(\nu)}, \Lambda)\right).$$

Comme expliqué en **4.1.2**, on a $\mathrm{H}^d(\mathrm{Tot}\, X_\bullet^{(\mu)}, \Lambda) \xrightarrow{\sim} \mathrm{H}^d(\mathrm{Tot}\, X_{\bullet\leq r}^{(\mu)}, \Lambda)$, et de même pour $\nu$. Que l'on puisse trouver une présentation explicite (**13.3**) de $\mathrm{H}^d(X, \Lambda)$ résulte alors de la proposition **4.3.1**.

### 8.3. Hyper-Čech.

**8.3.1**. Soit $X'_\bullet \to X$ un hyperrecouvrement pour la topologie des altérations. Vérifions que l'on peut calculer les morphismes $\check{\mathrm{H}}^d(X'_\bullet, \Lambda) \to \mathrm{H}^d(X, \Lambda)$, comme annoncé en **0.4**.

D'après le lemme **1.4.13**, on peut calculer (en tout étage) un hyperrecouvrement $X_\bullet \to X$ comme ci-dessus, se factorisant à travers un morphisme $X_\bullet \to X'_\bullet$. Fixons $d \geq 0$. La flèche $\check{\mathrm{H}}^d(X'_\bullet, \Lambda) \to \mathrm{H}^d(X, \Lambda)$ étant la composée des flèches $\check{\mathrm{H}}^d(X'_\bullet, \Lambda) \to \check{\mathrm{H}}^d(X_\bullet, \Lambda)$ et $\mathrm{H}^d(X_\bullet, \Lambda) \to \mathrm{H}^d(X, \Lambda)$ — la première étant trivialement calculable pour des schémas simpliciaux donnés (par calculabilité fonctorielle du $\pi_0$) —, on est ramené au cas particulier où $X'_\bullet = X_\bullet$. La conclusion résulte alors d'une part du fait que, comme observé en **4.3.2**, on a $\check{\mathrm{H}}^d(X_\bullet, \Lambda) = \mathrm{H}^d(\mathrm{Tot}\, X_{\bullet\leq r}^{(1)}, \Lambda)$ pour $r > d$ et, d'autre part, de la calculabilité des flèches $\mathrm{H}^d(\mathrm{Tot}\, X_{\bullet\leq r}^{(\mu)}, \Lambda) \to \mathrm{H}^d(\mathrm{Tot}\, X_{\bullet\leq r}^{(\nu)}, \Lambda)$ pour $\mu \leq \nu \leq \infty$.

**8.3.2**. Il résulte de ce qui précède que, donnés deux hyperrecouvrements $X'_\bullet \to X_\bullet$ de $X$ pour la topologie des altérations, on sait vérifier si la flèche $\mathrm{H}^d(X_\bullet, \Lambda) \to \mathrm{H}^d(X, \Lambda)$ identifie la cohomologie de $X$ au quotient de $\mathrm{H}^d(X_\bullet, \Lambda)$ par le noyau (calculable) de $\mathrm{H}^d(X_\bullet, \Lambda) \to \mathrm{H}^d(X'_\bullet, \Lambda)$. (D'autre part, on sait qu'il existe deux tels hyperrecouvrements.) En particulier, si les $X_{\alpha\bullet}$ sont comme en **0.2**, le système inductif $\check{\mathrm{H}}^d(X_{\alpha\bullet}, \Lambda)$ est explicitement essentiellement constant (mais le « explicitement » utilise une recherche non bornée).

**8.3.3. Remarque.** Notons que si $V \to U$ est un morphisme de $k$-schémas algébriques se factorisant à travers un revêtement $\ell$-étale $n$-approché universel $U^{[n]}$ de $U$, le morphisme de topos $V^{(n)} \to U^{(n)}$ se factorise à travers $V^{(n)} \to V^{(1)}$, dont le but est naturellement équivalent au topos discret des faisceaux sur $\pi_0(V)$. D'autre part, pour chaque $X_\bullet$ comme en **8.2** et chaque entier $n \geq 1$, on devrait sans aucun doute pouvoir fabriquer en utilisant les techniques usuelles de construction d'hyperrecouvrements (cf. **1.4.13** et **4.2.1**) — donc, en particulier, sans nouvelle recherche non bornée — un hyperrecouvrement $\tilde{X}_\bullet$ de $X$ au-dessus de $X_\bullet$ tel que les $\tilde{X}_i \to X_i$ se factorisent par un revêtement $n$-approché universel de $X_i$. Que $\mathrm{R}\Gamma(X_\bullet^{(n)}, \Lambda) \to \mathrm{R}\Gamma(X, \Lambda)$ se factorise à travers $\check{\mathrm{R}}\Gamma(\tilde{X}_\bullet, \Lambda) \to \mathrm{R}\Gamma(X, \Lambda)$ entraîne que l'on peut obtenir (sans nouvelle recherche non bornée) des cocycles hyper-Čech pour une base des $\mathrm{H}^i(X, \Lambda)$.

### 8.4. Calcul de $\mathrm{R}\Gamma(X, \Lambda)$.

**8.4.1.** Soient $X$ et $\Lambda$ comme ci-dessus et $X_{\alpha\bullet}$ un système projectif (indexé par les entiers)[①] d'hyperrecouvrements de $X$ tel que pour chaque entier $i$, on ait l'égalité $\mathrm{colim}_\alpha\, \check{\mathrm{H}}^i(X_{\alpha\bullet}, \Lambda) \xrightarrow{\sim} \mathrm{H}^i(X, \Lambda)$, ou encore un quasi-isomorphisme $\mathrm{hocolim}_\alpha\, \check{\mathrm{R}}\Gamma(X_{\alpha\bullet}, \Lambda) \xrightarrow{\sim} \mathrm{R}\Gamma(X, \Lambda)$, où $\check{\mathrm{R}}\Gamma(X_{\alpha\bullet}, \Lambda)$ est le complexe de Čech déduit du $\Lambda$-module cosimplicial $\Gamma(X_{\alpha\bullet}, \Lambda)$. (Rappelons que dans une catégorie abélienne satisfaisant la condition AB5 de Grothendieck, la cohomologie d'une colimite homotopique est la colimite des groupes de cohomologie.) Le complexe $\mathrm{R}\Gamma(X, \Lambda)$ appartenant à $\mathsf{D}^b_c(\Lambda)$, il résulte du lemme classique [SGA 4½, Rapport, 4.7] qu'il existe un complexe $\mathscr{K}$ de $\Lambda$-modules de type fini, concentré en degrés $[0, 2\dim(X)]$ et, pour $\alpha$ suffisamment grand, un morphisme de (vrais) complexes $\mathscr{K} \to \check{\mathrm{R}}\Gamma(X_{\alpha\bullet}, \Lambda)$ tel que la flèche composée $\mathscr{K} \to \mathrm{R}\Gamma(X, \Lambda)$ soit un quasi-isomorphisme. Pour calculer un tel $\mathscr{K}$, il suffit de parcourir les morphismes $\mathscr{K} \to \check{\mathrm{R}}\Gamma(X_{\alpha\bullet}, \Lambda)$ et de s'arrêter lorsqu'on en a trouvé un induisant le quasi-isomorphisme recherché (en degré $0 \leq i \leq 2\dim(X)$). C'est possible car on sait calculer les flèches $\check{\mathrm{H}}^i(X_{\alpha\bullet}, \Lambda) \to \mathrm{H}^i(X, \Lambda)$.

**8.4.2.** La remarque précédente devrait même permettre de calculer $\mathrm{R}\Gamma(X, \Lambda)$ sans plus de recherches non bornées que celles faites jusqu'à **8.2**.

## 9. Fonctorialité

### 9.1. Fonctorialité sur $\mathrm{Spec}(k)$.

**9.1.1.** Soient $k$ un corps algébriquement clos, $f : Y \to X$ un morphisme de $k$-schémas algébriques (supposé décrit comme en **16.2**), et $\Lambda$ un $\ell$-groupe abélien fini, avec $\ell$ inversible sur $k$. D'après **1.4.12**, il existe un morphisme simplicial $Y_\bullet \to X_\bullet$ au-dessus de $f$, calculable jusqu'à des étages arbitrairement élevés (et dépendant de $\ell$ mais pas de $\Lambda$), donnant lieu à un diagramme commutatif

$$\begin{array}{ccccccc}
\mathrm{R}\Gamma(\mathrm{Tot}\, X_{\bullet \ell\text{ét}}, \Lambda) & \xrightarrow{\sim} & \mathrm{R}\Gamma(\mathrm{Tot}\, X_{\bullet\text{ét}}, \Lambda) & \longrightarrow & \mathrm{R}\Gamma(\mathrm{Tot}\, Y_{\bullet\text{ét}}, \Lambda) & \xleftarrow{\sim} & \mathrm{R}\Gamma(\mathrm{Tot}\, Y_{\bullet \ell\text{ét}}, \Lambda) \\
 & & \uparrow \wr & & \uparrow \wr & & \\
 & & \mathrm{R}\Gamma(X, \Lambda) & \longrightarrow & \mathrm{R}\Gamma(Y, \Lambda). & &
\end{array}$$

D'après ce qui précède (§ 8), il existe deux entiers explicites $\mu \leq \nu$ tels que pour chaque $d \leq 2\max\{\dim(X), \dim(Y)\} < r$, on ait $\mathrm{H}^d(Z, \Lambda) = \mathrm{Im}\left(\mathrm{H}^d(\mathrm{Tot}\, Z^{(\mu)}_{\bullet \leq r}, \Lambda) \to \mathrm{H}^d(\mathrm{Tot}\, Z^{(\nu)}_{\bullet \leq r}, \Lambda)\right)$ pour $Z = X$ ou $Y$. Le morphisme $\mathrm{H}^d(\mathrm{Tot}\, X^{(\nu)}_{\bullet \leq r}, \Lambda) \to \mathrm{H}^d(\mathrm{Tot}\, Y^{(\nu)}_{\bullet \leq r}, \Lambda)$ étant calculable (cf. **4.2.6**), on peut trouver une présentation explicite du morphisme $\mathrm{H}^d(f, \Lambda) : \mathrm{H}^d(X, \Lambda) \to \mathrm{H}^d(Y, \Lambda)$ (au sens de **13.3**). Si $\Lambda$ est le corps $\mathbb{Z}/\ell\mathbb{Z}$, cela revient bien entendu à calculer le rang de $\mathrm{H}^d(f, \mathbb{Z}/\ell\mathbb{Z})$.

Notons que l'on pourrait bien entendu utiliser la présentation tirée de **8.3**, pour obtenir le même résultat (voir aussi **9.2**).

---

[①] Lorsque $k$ est dénombrable, ce qui suffit pour notre propos, l'existence d'un tel système projectif est élémentaire (voir [Deligne 1980, 5.2.2], cité en **0.2**). Pour $k$ quelconque, on peut améliorer le résultat classique selon lequel la catégorie des hyperrecouvrements *à homotopie près* est cofiltrante (cf. p. ex. [Artin et Mazur 1969, 8.13]) en la « rigidifiant » ([Friedlander 1982, § 4] ; comparer avec **3.4.2**). D'après O. Gabber (communication personnelle) on a des résultats semblables pour des sites généraux, sans hypothèse de finitude (ni, notamment, d'existence de suffisamment de points).

**9.1.2**. *Amélioration.* Vérifions maintenant que pour toute collection finie $f_1, \dots, f_r$ de $k$-morphismes explicites $Y \to X$, on peut calculer des présentations explicites des $\mathrm{H}^d(f_i, \Lambda)$ relativement à de *mêmes* présentations explicites de $\mathrm{H}^d(X, \Lambda)$ et $\mathrm{H}^d(Y, \Lambda)$. Il résulte en effet de **1.4.12** et **1.4.13** qu'il existe pour chaque $\alpha \in \{1, \dots, r\}$ un morphisme $Y_{\alpha\bullet} \to X_\bullet$ comme en *loc. cit.* au-dessus de $f_\alpha$. En considérant le produit fibré des $Y_{\alpha\bullet}$ au-dessus de $Y$ et en réappliquant la construction de *loc. cit.*, on en déduit qu'il existe un diagramme commutatif

$$\begin{array}{ccc} Y_\bullet & \overset{f_{1\bullet}}{\underset{f_{r\bullet}}{\vdots}} & X_\bullet \\ \downarrow & & \downarrow \\ Y & \overset{f_1}{\underset{f_r}{\vdots}} & X \end{array}$$

où les flèches verticales *ne dépendent pas de l'indice* $\alpha \in \{1, \dots, r\}$. La conclusion en résulte aussitôt.

(En particulier, donnés deux morphismes $f, g : Y \to X$, on peut décider si $\mathrm{H}^d(g, \Lambda) = \mathrm{H}^d(f, \Lambda)$.)

Notons que si $Y = X$, on peut supposer que les $\Lambda$-modules explicites $\mathrm{H}^d(X, \Lambda)$ et $\mathrm{H}^d(Y, \Lambda)$ sont *égaux*. Pour s'en convaincre, il suffit par exemple de rajouter l'identité $Y \to X$ aux morphismes $f_1, \dots, f_r$ et de composer avec l'inverse de l'isomorphisme $\mathrm{H}^d(X, \Lambda) \to \mathrm{H}^d(Y, \Lambda)$ qui s'en déduit.

On peut reformuler la fonctorialité établie sous la forme suivante.

**9.1.3**. Soit $\mathscr{G}$ un graphe fini orienté avec arêtes multiples possibles. Supposons donné un étiquetage de $\mathscr{G}$ par la catégorie des $k$-schémas algébriques, c'est-à-dire un étiquetage de ses sommets par des $k$-schémas algébriques et un étiquetage des arêtes par des $k$-morphismes (entre les schémas correspondants). On peut calculer un étiquetage du graphe opposé $\mathscr{G}^{\mathrm{op}}$ par la catégorie des $\Lambda$-modules finis, déduit du précédent par application du foncteur $\mathrm{H}^d(-, \Lambda)$.

On peut déduire cet énoncé du précédent en considérant le coproduit $X = \mathrm{Spec}(k) \amalg \coprod_s X_s$, où $X_s$ parcourt les étiquettes des sommets $s$ de $\mathscr{G}$, et les endomorphismes $X \to X$ envoyant chaque $X_{s'}$, sauf un $X_s$, sur $\mathrm{Spec}(k)$ et déterminé par l'étiquette d'une arête sur ce dernier $X_s$.

## 9.2. Action galoisienne.

**9.2.1**. Soient ${}_0k$ un corps (calculable et disposant d'un algorithme de factorisation et d'une $p$-base finie explicite : cf. **12.7**) et ${}_0X$ un ${}_0k$-schéma algébrique explicite. Fixons une clôture algébrique $k$ de ${}_0k$ et notons $X$ le $k$-schéma obtenu par extension des scalaires. Nous allons montrer que l'on peut calculer une extension finie galoisienne ${}_1k/{}_0k$ telle que l'action du groupe de Galois de ${}_0k$ sur $\mathrm{H}^*(X, \Lambda)$ se factorise à travers $\Gamma = \mathrm{Gal}({}_1k/{}_0k)$ et calculer la représentation du groupe fini correspondante. Toute extension étale de la clôture parfaite de ${}_0k$ dans $k$ se descendant explicitement à ${}_0k$, on peut supposer le corps ${}_0k$ parfait (**12.5**, (v)).

**9.2.2**. Fixons $d$. Comme on l'a vu en **8.3**, il existe un hyperrecouvrement tronqué pour la topologie des altérations $X_{\bullet\leq r} \to X$ tel que $\check{\mathrm{H}}^d(X_{\bullet\leq r}, \Lambda) \to \mathrm{H}^d(X, \Lambda)$ soit *surjective* ; par cofinalité (observée en **1.4.11**), il existe un tel hyperrecouvrement *défini sur* ${}_0k$. La calculabilité de l'action du groupe de Galois sur $\mathrm{H}^d(X, \Lambda)$ se déduit donc de celle de l'action sur $\check{\mathrm{H}}^d(X_{\bullet\leq r}, \Lambda)$ et de la calculabilité de la flèche. Plus précisément, si $c \in \mathrm{H}^d(X, \Lambda)$ est l'image d'une $d$-chaîne $z \in \mathrm{H}^0(X_d, \Lambda)$, la classe de cohomologie $\gamma \cdot c$, où $\gamma \in \mathrm{Gal}(k/{}_0k)$, est l'image de la $d$-chaîne $\gamma \cdot z$ déduite de l'action du groupe de Galois sur $\pi_0(X_d) = \pi_0({}_0X_d \otimes_{{}_0k} k)$.

Notons que l'hyperrecouvrement tronqué $X_{\bullet\leq r} \to X$ est défini sur une sous-extension galoisienne finie ${}_1k/{}_0k$, l'action précédente se factorise à travers le quotient fini $\Gamma = \mathrm{Gal}({}_1k/{}_0k)$. Ce dernier est calculable car ${}_0k$ est un corps calculable avec un algorithme de factorisation ([Fried et Jarden 2008, 19.3.2]).

## 10. Structure de l'algorithme et exemple simple

**10.1. Structure générale.** Récapitulons brièvement comment les différents éléments qui ont été présentés s'emboîtent pour fournir, en principe, un algorithme permettant de calculer $\mathrm{H}^d(X, \Lambda)$ pour $X$ un schéma algébrique sur un corps algébriquement clos de caractéristique différente de $\ell$ et $\Lambda$ un $\ell$-groupe abélien fini.

Dans un premier temps, on calcule, jusqu'à un certain niveau $r$, un hyperrecouvrement $X_\bullet$ de $X$ tel qu'explicité en **1.4.12**. Plus exactement, on calcule $X_0 \to X$ qui recouvre $X$ et dont les composantes sont des polycourbes $\ell$-élémentaires (et en particulier, des $K(\pi, 1)$ pro-$\ell$) : ceci se fait au moyen de la proposition **1.4.9** (compte tenu des remarques qui suivent au sujet de la constructivité) ; puis de même $X_1 \to X_0 \times_X X_0$ et ainsi de suite comme expliqué en **1.4.13**. Cette construction des $X_i$ doit être menée pour $i \leq r$ avec $r$ qui dépend uniquement du degré en lequel on veut calculer la cohomologie (comme expliqué en **4.1.2**, en fait $r = d + 1$ suffit).

D'après les résultats de la section 3 — s'appuyant sur la calculabilité du nombre de $\mathbb{Z}/\ell\mathbb{Z}$-torseurs établie en **2.1** —, on sait calculer, pour chaque niveau d'approximation fini $\lambda \geq 1$, et pour tout $k$-schéma algébrique normal $Y$, un revêtement $\ell$-étale $\lambda$-approché universel $Y^{[\lambda]}$, ainsi que le groupe de Galois correspondant $\pi_Y^{[\lambda]}$ si $Y$ est connexe.

D'après § 4, on sait calculer, pour chaque niveau d'approximation fini $\lambda \geq 1$, le groupe $\mathrm{H}^d(\mathrm{Tot}\, X_\bullet^{(\lambda)}, \Lambda)$ (isomorphe à $\mathrm{H}^d(\mathrm{Tot}\, X_{\bullet \leq r}^{(\lambda)}, \Lambda)$ puisque $r$ a été choisi assez grand), et même la flèche $\mathrm{H}^d(\mathrm{Tot}\, X_\bullet^{(\mu)}, \Lambda) \to \mathrm{H}^d(\mathrm{Tot}\, X_\bullet^{(\nu)}, \Lambda)$ pour deux entiers $\mu \leq \nu$.

Il s'agit donc de calculer de tels entiers pour que l'image de cette flèche soit le groupe $\mathrm{H}^d(X, \Lambda)$ recherché. Comme expliqué en § 8, ceci résulte de **5.8** appliqué à la suite spectrale $\mathrm{E}_{1,\lambda}^{i,j} = \mathrm{H}^j(X_i^{(\lambda)}, \Lambda) \Rightarrow \mathrm{H}^{i+j}(\mathrm{Tot}\, X_\bullet^{(\lambda)}, \Lambda)$, une fois connus des fonctions explicitant, pour chaque $i$ et $j$, le fait que $\mathrm{H}^j(X_i^{(\lambda)}, \Lambda)$ est essentiellement constant. De telles bornes sont obtenues en § 7.

**10.2. Esquisse d'exemple.** Pour illustrer la manière dont l'algorithme s'exécuterait, nous esquissons le calcul de $\mathrm{H}^i(\mathbb{P}_k^1, \mathbb{Z}/\ell\mathbb{Z})$ pour les petits $i$ en en suivant les différentes étapes. (Le cas encore plus simple d'une courbe *affine* lisse $X$ consiste essentiellement à calculer $X^{[\lambda]}$ par **3.3** pour des petites valeurs de $\lambda$ et à appliquer la proposition **6.5** : de toute façon, on est ramené à § 2.)

Notons $U = \mathbb{P}_k^1 \setminus \{\infty\}$ et $U' = \mathbb{P}_k^1 \setminus \{0\}$ deux ouverts de Zariski qui recouvrent $\mathbb{P}_k^1$ et dont on note $V \coloneqq U \times_{\mathbb{P}_k^1} U'$ l'intersection. Soit $X_0 \coloneqq U \amalg U' \to \mathbb{P}_k^1$, vu comme un $\mathbb{P}_k^1$-schéma simplicial 0-tronqué ; son cosquelette est donné par $X_p = X_0 \times_{\mathbb{P}_k^1} \cdots \times_{\mathbb{P}_k^1} X_0$ (avec $p + 1$ facteurs), qui est le coproduit $U \amalg V \amalg \cdots \amalg V \amalg U'$ de $U$, $U'$ et $2^{p+1} - 2$ copies de $V$ qu'on imaginera étiquetés par les $2^{p+1}$ mots binaires $w$ de longueur $p + 1$, et pour $0 \leq i \leq p + 1$, le morphisme $X_{\delta_{p,i}} : X_{p+1} \to X_p$ envoie par le morphisme évident la composante $X_{p+1,w}$ étiquetée $w$ sur celle $X_{p,w'}$ étiquetée par le mot $w'$ égal à $w$ privé de son $i$-ième bit. Il s'agit manifestement d'un hyperrecouvrement. Comme $U, U', V$ sont des (poly)courbes $\ell$-élémentaires sur $\mathrm{Spec}\, k$, il est possible que la recherche de courbes élémentaires effectuée par l'algorithme retourne cet hyperrecouvrement.

Examinons maintenant comment se déroulerait le calcul de la cohomologie de $\mathrm{Tot}\, X_\bullet^{(\lambda+1)}$ (en fonction d'un entier $\lambda + 1 \geq 1$) à valeurs dans le faisceau constant $\mathbb{Z}/\ell\mathbb{Z}$. Comme $V = \mathbb{G}_\mathrm{m}$, le topos $V^{(\lambda+1)}$ est le topos des $\mathbb{Z}/\ell^\lambda\mathbb{Z}$-ensembles, tandis que $U^{(\lambda+1)}$ et $U'^{(\lambda+1)}$ sont, bien sûr, celui des ensembles. (Et pour un $\lambda$ donné, l'algorithme est capable d'effectuer ce calcul en suivant **3.3**.) Un faisceau abélien de $\mathrm{Tot}\, X_\bullet^{(\lambda+1)}$ est donc (cf. **4.1.1**) la donnée pour chaque mot binaire $w$ de longueur $p+1$ (pour $p \geq 0$) d'un groupe abélien $A_{p,w}$, muni d'une action de $\mathbb{Z}/\ell^\lambda\mathbb{Z}$ sauf si $w$ est l'un des mots $00 \cdots 0$ ou $11 \cdots 1$, ainsi que de morphismes $A_{p,w} \to A_{p+1,w'}$ pour chaque mot $w'$ obtenu en insérant un bit dans $w$, vérifiant les compatibilités évidentes. Si comme topos discret utilisé en § 4 on prend $P_\bullet = X_\bullet^{(1)}$, alors un faisceau abélien de $\mathrm{Tot}\, P_\bullet$ correspond à de telles données sans l'action de $\mathbb{Z}/\ell^\lambda\mathbb{Z}$ : à un tel objet est associé un complexe de différentielle (†) $\bigoplus_{w \in \{0,1\}^{p+1}} A_{p,w} \to \bigoplus_{w \in \{0,1\}^{p+2}} A_{p+1,w}$ somme alternée des morphismes $A_{p,w} \to A_{p+1,w'}$.

Les foncteurs $u^\star$ associant à un $\mathbb{Z}/\ell^\lambda\mathbb{Z}$-ensemble son ensemble sous-jacent, et $u_\star$ son adjoint à droite $X \mapsto \mathsf{Hom}_{\mathsf{Ens}}(\mathbb{Z}/\ell^\lambda\mathbb{Z}, X)$, définissent pour tout $\mathbb{Z}/\ell^\lambda\mathbb{Z}$-module une résolution de Godement, analogue (et quasi-isomorphe) à l'une des résolutions habituelles définissant la cohomologie des groupes (par exemple [Serre 1994, I, § 2.2] ou [Neukirch, A. Schmidt et Wingberg 2000, I, § 2]). Ces foncteurs forment un morphisme $\mathrm{Tot}\, P_\bullet \to \mathrm{Tot}\, X_\bullet^{(\lambda+1)}$. À partir du faisceau constant $A_{p,w} = \mathbb{Z}/\ell\mathbb{Z}$, l'algorithme (4.2) va donc calculer, en bas degrés, le complexe double dont les colonnes sont sommes de copies de la résolution qu'on vient de dire (et pour $w = 00\cdots 0, 11\cdots 1$, de la résolution triviale), et dont les différentielles horizontales sont données par (†). Le calcul du $\mathrm{H}^i(\mathrm{Tot}\, X_\bullet^{(\lambda+1)}, \mathbb{Z}/\ell\mathbb{Z})$ est alors donné par la cohomologie du complexe simple associé à ce complexe double.

On peut prédire quel sera le résultat de ce calcul en utilisant la suite spectrale de descente décrite en 8.2 (et qui est la « première » suite spectrale associée au complexe double décrit ci-dessus) : il est facile de se convaincre que cette suite spectrale dégénère en $\mathrm{E}_2$, les seuls termes non nuls étant $\mathrm{E}_2^{0,0} = \mathbb{Z}/\ell\mathbb{Z}$ et $\mathrm{E}_2^{1,q} = \mathrm{H}^q(\mathbb{Z}/\ell^\lambda\mathbb{Z}, \mathbb{Z}/\ell\mathbb{Z})$ si $q \geq 2$. Ainsi, $\mathrm{H}^n(\mathrm{Tot}\, X_\bullet^{(\lambda+1)}, \mathbb{Z}/\ell\mathbb{Z})$ vaut $\mathbb{Z}/\ell\mathbb{Z}$ si $n = 0$, $0$ si $n = 1$, et $\mathrm{H}^{n-1}(\mathbb{Z}/\ell^\lambda\mathbb{Z}, \mathbb{Z}/\ell\mathbb{Z})$ si $n \geq 2$ (c'est-à-dire $\mathbb{Z}/\ell\mathbb{Z}$ dès que $\lambda \geq 1$). Pour $\lambda \leq \mu$, les flèches $\mathrm{H}^n(\mathrm{Tot}\, X_\bullet^{(\lambda+1)}, \mathbb{Z}/\ell\mathbb{Z}) \to \mathrm{H}^n(\mathrm{Tot}\, X_\bullet^{(\mu+1)}, \mathbb{Z}/\ell\mathbb{Z})$ correspondent bien aux morphismes fonctoriels (d'inflation) $\mathrm{H}^{n-1}(\mathbb{Z}/\ell^\lambda\mathbb{Z}, \mathbb{Z}/\ell\mathbb{Z}) \to \mathrm{H}^{n-1}(\mathbb{Z}/\ell^\mu\mathbb{Z}, \mathbb{Z}/\ell\mathbb{Z})$, et grâce à 6.5–6.7 on sait que cette flèche sera un isomorphisme pour $n \leq 2$ et nulle pour $n \geq 3$ dès que $1 \leq \lambda < \mu$. Muni de cette borne, l'algorithme calcule $\mathrm{H}^n(\mathbb{P}^1_k, \mathbb{Z}/\ell\mathbb{Z})$ comme l'image de $\mathrm{H}^n(\mathrm{Tot}\, X_\bullet^{(2)}, \mathbb{Z}/\ell\mathbb{Z}) \to \mathrm{H}^n(\mathrm{Tot}\, X_\bullet^{(3)}, \mathbb{Z}/\ell\mathbb{Z})$, ce qui donne bien le résultat attendu.

Soulignons sur cet exemple le fait suivant : si on remplace les coefficients $\mathbb{Z}/\ell\mathbb{Z}$ par le groupe $\boldsymbol{\mu}_\ell$ des racines $\ell$-ièmes de l'unité (qui lui est non canoniquement isomorphe) et le groupe $\mathbb{Z}/\ell^\lambda\mathbb{Z}$ par le groupe $\boldsymbol{\mu}_{\ell^\lambda}$ (de nouveau non canoniquement isomorphe) des automorphismes du revêtement étale $\mathbb{G}_\mathrm{m} \to \mathbb{G}_\mathrm{m}$ donné par $z \mapsto z^{\ell^\lambda}$, le calcul ne fait plus intervenir aucun choix arbitraire, et on voit donc $\mathrm{H}^2(\mathbb{P}^1_k, \boldsymbol{\mu}_\ell)$ comme $\mathrm{H}^1(\boldsymbol{\mu}_{\ell^\lambda}, \boldsymbol{\mu}_\ell)$ (l'action de $\boldsymbol{\mu}_{\ell^\lambda}$ étant triviale), lui-même isomorphe à $\mathsf{Hom}(\boldsymbol{\mu}_{\ell^\lambda}, \boldsymbol{\mu}_\ell) = \mathbb{Z}/\ell\mathbb{Z}$ sans faire de choix arbitraire.

## 11. Compléments

On considère ici quelques résultats qui sont des prolongements naturels de notre théorème principal et on énonce en 11.5 quelques questions, dont certaines sont probablement hors de portée. Cette section se termine (11.6) par quelques précisions métamathématiques sur la nature des algorithmes que l'on espère pouvoir obtenir.

Ci-dessous $k$ est un corps algébriquement clos et $\Lambda$ un anneau commutatif fini de cardinal inversible sur $k$.

**11.1. Cohomologie d'un schéma simplicial.** Soit $X_\bullet$ un $k$-schéma algébrique simplicial. D'après [Deligne 1974, § 6.4]et les résultats de § 1.4, il existe un hyperrecouvrement $X_{\bullet\bullet} \to X_\bullet$ par des polycourbes $\ell$-élémentaires tel que $\mathrm{H}^*(X_\bullet, \Lambda) = \mathrm{H}^*(X_{\bullet\bullet}, \Lambda)$. On peut donc procéder comme dans le cas non simplicial pour calculer le terme de gauche en tout degré donné à l'avance. Noter qu'il n'est pas absolument nécessaire d'utiliser une variante bisimpliciale de arguments précédents : d'après le théorème de Cartier-Eilenberg-Zilber ([ibid., 6.4.2.2] ou [Illusie 1971-1972, I, § 1.2]), on a $\mathrm{R}\Gamma(X_{\bullet\bullet}, \Lambda) = \mathrm{R}\Gamma(\delta X_{\bullet\bullet}, \Lambda)$, où $\delta X_{\bullet\bullet}$ est le schéma simplicial *diagonal* déduit de $X_{\bullet\bullet}$. Le même argument est valable si l'on veut calculer la cohomologie de $X_\bullet$ à valeurs dans un faisceau à composantes localement constantes (étage par étage).

Tout espace algébrique (au sens d'Artin) étant localement pour la topologie étale un schéma, on peut probablement utiliser le résultat précédent pour en calculer la cohomologie étale.

**11.2. Cohomologie relative.** On présente ici deux constructions de la cohomologie relative d'un morphisme $Y \to X$ : la première (simpliciale) est valable en toute généralité, la seconde (11.2.5) s'applique uniquement à la cohomologie relativement à un sous-schéma fermé.

Rappelons que la cohomologie relative s'inscrit dans un triangle distingué

$$\mathrm{R}\Gamma(X/Y, \Lambda) \to \mathrm{R}\Gamma(X, \Lambda) \to \mathrm{R}\Gamma(Y, \Lambda) \overset{+1}{\to} .$$

Lorsque $\Lambda$ est un corps, la calculabilité de la dimension du groupe de cohomologie relative $\mathrm{H}^i(X/Y, \Lambda)$ résulte donc immédiatement du théorème 0.1 : c'est une extension de $\mathrm{Ker}\left(\mathrm{H}^i(X, \Lambda) \to \mathrm{H}^i(Y, \Lambda)\right)$ par $\mathrm{Im}\left(\mathrm{H}^{i-1}(X, \Lambda) \to \mathrm{H}^{i-1}(Y, \Lambda)\right)$.

**11.2.1**. *Cône d'un morphisme simplicial.* Soit $f_\bullet : Y_\bullet \to X_\bullet$ un morphisme de topos simpliciaux (dont on notera $\mathrm{Tot}\, f_\bullet : \mathrm{Tot}\, Y_\bullet \to \mathrm{Tot}\, X_\bullet$ le morphisme entre les topos totaux, cf. **4.1.1**) ; on pourra penser au cas d'un morphisme déduit, par passage aux topos étales, d'un morphisme de $k$-schémas algébriques simpliciaux. Notons $Cf_\bullet$ son **cône** : $Cf_n = X_n \amalg Y_{<n}$, où $Y_{<n} = \coprod_{-1 \leq i < n} Y_i$ et, conformément à l'usage, $Y_{-1}$ est le topos final $e$. Explicitons le morphisme de fonctorialité $Cf_\varphi : Cf_m \to Cf_n$ déduit d'une application croissante $\varphi : [n] = \{0, \ldots, n\} \to [m] = \{0, \ldots, m\}$. Sur le facteur $X_m$, c'est le morphisme composé $X_m \to X_n \hookrightarrow Cf_n$, la première flèche étant $X_\varphi$. Soit maintenant $-1 \leq \mu < m$ et considérons l'application $\varphi_\mu : [\nu] \to [\mu]$ déduite de $\varphi$ par changement de base $[\mu] \hookrightarrow [m]$. Sur le facteur $Y_\mu$, le morphisme $Cf_\varphi$ est le morphisme composé $Y_\mu \to Y_\nu \hookrightarrow Cf_n$, la première flèche étant $Y_{\varphi_\mu}$ si $\nu < n$, et $Y_\mu \to Y_\nu \to X_\nu \hookrightarrow Cf_n$, où $Y_\nu \to X_\nu$ est $f_\nu$. Par exemple, le morphisme face $d_i : Cf_{n+1} \to Cf_n$, déduit de l'unique injection croissante $[n] \to [n+1]$ d'image ne contenant pas $i$, est induit par les identités $Y_j \to Y_j$ pour $j < i$ (composée avec $Y_n \to X_n$ si $j = n$) et les faces $d_i : Y_l \to Y_{l-1}$ pour $i \leq l \leq n$ et $d_i : X_{n+1} \to X_n$.

Soit $\varphi_\bullet : \mathscr{F}_\bullet \to \mathscr{G}_\bullet$ un $f_\bullet$-morphisme entre faisceaux abéliens $\mathscr{F}_\bullet$ sur $\mathrm{Tot}\, X_\bullet$ et $\mathscr{G}_\bullet$ sur $\mathrm{Tot}\, Y_\bullet$, c'est-à-dire un morphisme $f_\bullet^\star \mathscr{F}_\bullet \to \mathscr{G}_\bullet$, objet du topos flèche $\mathrm{Fl}(\mathrm{Tot}\, f_\bullet)$. Les coproduits et les flèches du paragraphe précédent, calculés dans la catégorie des paires $(X, \mathscr{F})$ — cf. [Deligne 1974, 6.3.1.b)], en remplaçant « *espace topologique* » par « *topos* »[①] —, permettent de définir un faisceau *abélien* simplicial sur $\mathrm{Tot}\, Cf_\bullet$, cône de $\varphi_\bullet$, que nous notons $C\varphi_\bullet$. (Notons que l'objet final de cette catégorie de paires est le faisceau nul sur $e$.)

**11.2.2**. La cohomologie relative de $\mathrm{Tot}\, X_\bullet$ modulo $\mathrm{Tot}\, Y_\bullet$ est un cas particulier de celle définie, pour tout morphisme $S \to T$ de topos, en [Illusie 1971-1972, III, § 4] (voir aussi [Deligne 1980, § 4.3.4]). On la note $\mathrm{R}\Gamma(\mathrm{Tot}\, X_\bullet / \mathrm{Tot}\, Y_\bullet, \varphi_\bullet)$ ; ce n'est en général *pas* la cohomologie de $\mathrm{Fl}(\mathrm{Tot}\, f_\bullet)$ à valeurs dans $\varphi_\bullet$. (Cette dernière étant isomorphe à $\mathrm{R}\Gamma(\mathrm{Tot}\, X_\bullet, \mathscr{F}_\bullet)$ d'après [Illusie 1971-1972, III.4.2].) Implicite en [Deligne 1974, § 6.3] est la formule suivante :

$$(\star) \qquad \mathrm{R}\Gamma(\mathrm{Tot}\, X_\bullet / \mathrm{Tot}\, Y_\bullet, \varphi_\bullet) = \mathrm{R}\Gamma(\mathrm{Tot}\, Cf_\bullet, C\varphi_\bullet).$$

En degré cohomologique nul, cette formule résulte immédiatement de la description explicite de $d_0$ et $d_1$, et du fait que le terme de gauche est $\mathrm{Ker}(\mathrm{H}^0(\mathrm{Tot}\, X_\bullet, \mathscr{F}_\bullet) \to \mathrm{H}^0(\mathrm{Tot}\, Y_\bullet, \mathscr{G}_\bullet))$. (On pourra aussi comparer à [Anderson 1987, § 1], où le cas des *ensembles* simpliciaux et des coefficients constants est considéré.)

**11.2.3**. Justifions brièvement $(\star)$. Notons $\mathscr{R}$ le topos $\mathrm{Fl}(\mathrm{Tot}\, f_\bullet)$ et $\mathscr{C}$ le topos $\mathrm{Tot}\, Cf_\bullet$. Il n'y a pas de morphisme naturel $\mathscr{C} \to \mathscr{R}$ — comme on le voit par exemple en considérant l'unique endomorphisme du topos simplicial vide — mais un morphisme $a : \mathscr{C}^\times \to \mathscr{R}$, où $\mathscr{C}^\times$ est le *sous-topos ouvert* de $\mathscr{C}$, défini par la condition : $i^\star = \varnothing$, où $i : \mathrm{Fl}(e_\bullet) \to \mathscr{C}$ est le morphisme évident. Le morphisme image inverse $a^\star$ est la variante ensembliste du cône d'un morphisme de faisceau abélien considérée en **11.2.1** : sur chaque $Y_{-1} = e$, on considère le faisceau initial (=vide). Notant $j : \mathscr{C}^\times \hookrightarrow \mathscr{C}$, on vérifie sans peine l'égalité $C\varphi_\bullet = j_! a^\star \varphi_\bullet$ pour tout faisceau abélien sur $\mathscr{R}$. La formule $(\star)$ résulte alors, par dérivation et exactitude de $j_! a^\star$, de la formule en degré cohomologique nul.

**11.2.4**. Soient $f : Y \to X$ un morphisme de $k$-schémas algébriques et $f_\bullet : Y_\bullet \to X_\bullet$ le morphisme induit entre les schémas simpliciaux constants. Il résulte de ce qui précède et de **11.1** l'on peut calculer le groupe de cohomologie étale $\mathrm{H}^i(X/Y, \Lambda)$ : c'est $\mathrm{H}^i(Cf_\bullet, C\Lambda_\bullet)$. Ceci s'applique en particulier au calcul de la cohomologie à support compact : si $X$ est *propre* sur $k$, la cohomologie de $X/F$ — où $F$ est un fermé de $X$ — est la cohomologie à support compact de $U = X - F$ :

$$\mathrm{H}^i_c(U, \Lambda) = \mathrm{H}^i(X/F, \Lambda).$$

**11.2.5**. Esquissons une autre approche du calcul de la cohomologie de $X$ relative à un fermé $F$. Notons $Z = X \amalg_F X$ la somme amalgamée (pincement) de deux copies de $X$ le long de $F$ ; si $X = \mathrm{Spec}(A)$ et $F = \mathrm{Spec}(A/I)$, le schéma $Z$ est le spectre du sous-anneau $B = \{(a_1, a_2) \in A^2 : a_1 \equiv a_2 \bmod I\}$ de $A^2$. Voir [Ferrand 2003, § 5 et § 7, notamment le théorème 7.1] pour une discussion de l'existence de coproduits sous des hypothèses bien plus générales et **16.4** pour une approche effective.

[①] Voir aussi [SGA 4, $\mathrm{V}^{\mathrm{bis}}$, 4.3.0].

Notons $i_1$ et $i_2$ les deux immersions fermées canoniques de $X$ dans $Z$, et $i$ l'immersion fermée de $F$ dans $Z$. Considérons les morphismes d'adjonctions $\rho_\alpha : i_{\alpha\star}\Lambda \to i_\star\Lambda$ (pour $\alpha \in \{1,2\}$) et $\delta$ le morphisme « différence » $(i_\star\Lambda)^2 \to i_\star\Lambda$, $(a_1, a_2) \mapsto a_1 - a_2$. Le complexe concentré en degrés $[0,1]$

$$[i_{1\star}\Lambda \oplus i_{2\star}\Lambda \overset{\delta\circ(\rho_1,\rho_2)}{\to} i_\star\Lambda]$$

calcule la cohomologie $\mathrm{R}\Gamma(Z,\Lambda)$ de $Z$ (« Mayer-Vietoris »). Or, il est *isomorphe* au complexe

$$[i_{1\star}\Lambda \oplus i_{2\star}\Lambda \overset{0\oplus\rho_2}{\to} i_\star\Lambda]$$

dont la cohomologie est la somme directe $\mathrm{R}\Gamma(X,\Lambda)\bigoplus\mathrm{R}\Gamma(X/F,\Lambda)$. La cohomologie de $X/F$ se déduit donc fonctoriellement de la cohomologie (usuelle) de $X$ et de celle de $Z = X \amalg_F X$.

(Notons que la rétraction $Z = X \amalg_F X \to X = X \amalg_X X$ des deux inclusions $X \rightrightarrows Z$ induit la décomposition en somme directe ci-dessus.)

**11.3. Structure d'algèbre graduée.**

**11.3.1.** Soient $X$ et $Y$ deux $k$-schémas de type fini, $\Lambda$ comme au début de la section 11 et $i, j$ deux entiers. Comme expliqué en **8.3**, on peut notamment construire deux hyperrecouvrements (étales) $X_\bullet \to X$ et $Y_\bullet \to Y$ tels que les flèches $\check{\mathrm{H}}^i(X_\bullet,\Lambda) \to \mathrm{H}^i(X,\Lambda)$ et $\check{\mathrm{H}}^j(Y_\bullet,\Lambda) \to \mathrm{H}^j(Y,\Lambda)$ soient *surjectives*. Utilisant la formule de Künneth triviale $\pi_0(X_\bullet \times_k Y_\bullet) = \pi_0(X_\bullet)\times\pi_0(Y_\bullet)$ et l'existence d'un homotopisme d'Eilenberg-Zilber *explicite* ([Mac Lane 1963, VIII, théorèmes 8.1 et 8.8])

$$\mathsf{Hom}(\pi_0(X_\bullet),\Lambda) \otimes \mathsf{Hom}(\pi_0(Y_\bullet),\Lambda) \to \mathsf{Hom}(\pi_0(X_\bullet \times_k Y_\bullet),\Lambda),$$

on en déduit un morphisme

$$\check{\mathrm{H}}^i(X_\bullet,\Lambda) \otimes \check{\mathrm{H}}^j(Y_\bullet,\Lambda) \to \check{\mathrm{H}}^{i+j}(X_\bullet \times_k Y_\bullet,\Lambda)$$

relevant le morphisme de Künneth $\mathrm{H}^i(X,\Lambda) \otimes \mathrm{H}^j(Y,\Lambda) \to \mathrm{H}^{i+j}(X \times_k Y,\Lambda)$ ([Milne 1980, V.1.19] ; comparer avec [Deligne 1974, 8.1.25]). Ce dernier est donc calculable compte tenu de ce qui précède. (On utilise le fait que la flèche naturelle $\check{\mathrm{H}}^{i+j}(X_\bullet \times_k Y_\bullet,\Lambda) \to \mathrm{H}^{i+j}(X \times_k Y,\Lambda)$ l'est.)

**11.3.2.** Lorsque $X = Y$, on en déduit la structure de $\Lambda$-algèbre graduée commutative sur $\mathrm{H}^*(X,\Lambda)$ par composition avec la flèche de fonctorialité $\mathrm{H}^*(X \times_k X,\Lambda) \to \mathrm{H}^*(X,\Lambda)$ induite par la diagonale $X \to X \times_k X$ : on peut donc calculer le produit

$$\mathrm{H}^*(X,\Lambda) \otimes \mathrm{H}^*(X,\Lambda) \to \mathrm{H}^*(X,\Lambda)$$

$$c_1 \otimes c_2 \mapsto c_1 \smile c_2.$$

**11.3.3.** Le résultat précédent a une application immédiate aux cycles algébriques. Supposons dorénavant $X$ propre, lisse, connexe de dimension $d_X$, et fixons un isomorphisme $t_X : \mathrm{H}^{2d_X}(X,\Lambda) \xrightarrow{\sim} \Lambda$. Par la dualité de Poincaré, toute forme linéaire $\varphi : \mathrm{H}^i(X,\Lambda) \to \Lambda$ est de la forme $b \mapsto t_X(a_\varphi \smile b)$ pour une unique classe $a_\varphi \in \mathrm{H}^{2d_X-i}(X,\Lambda)$, que l'on peut calculer. En particulier, tout morphisme $Z \to X$ de source une variété propre lisse connexe de dimension $d_Z$ induit une classe $c_Z \in \mathrm{H}^{2(d_X-d_Z)}(X,\Lambda)$, correspondant à la forme linéaire composée $\mathrm{H}^{2d_Z}(X,\Lambda) \to \mathrm{H}^{2d_Z}(Z,\Lambda) \xrightarrow{\sim} \Lambda$, où le dernier morphisme est un morphisme trace pour $Z$.

La classe de cycle $c_Z$ est ainsi définie à multiplication par un élément de $\Lambda^\times$ près (lié au choix arbitraire de traces pour $X$ et $Z$), ambiguïté que l'on devrait pouvoir lever. Ceci nous permet cependant de répondre à la question : la classe d'un cycle algébrique lisse de $X$ est-elle triviale ?

**11.4. Images directes.** On montre que pour tout morphisme $f : X \to S$ entre $k$-schémas algébriques, tout faisceau constructible $\mathscr{F}$ de $\Lambda$-modules sur $X$ et tout entier $i$, on peut calculer le faisceau $\mathrm{R}^i f_\star\mathscr{F}$, fonctoriellement en $\mathscr{F}$.

**11.4.1.** Il conviendrait de donner une description explicite de la catégorie des faisceaux constructibles sur un schéma (explicite) $X$ et de vérifier que quelques opérations usuelles (noyau, conoyau, etc.) sont bien calculables. Signalons simplement que différentes approches sont possibles : stratifications et flèches de recollement (cf. [SGA 4, IV.9.3]), description des générateurs ou cogénérateurs usuels ([SGA 4, IX 2.9 (ii) et 2.14 (ii)]), espaces algébriques ([Artin 1973, chap. VII], ou [Milne 1980, chap. V, § 1]). Cette dernière est probablement la plus économique.

**11.4.2**. *Effaçabilité et dévissages.* Soient $S$ un schéma nœthérien et $f : X \to S$ un morphisme de type fini. D'après une variante de [SGA 4½, Arcata, IV.3.5], les foncteurs $\mathrm{R}^i f_\star$ pour $i > 0$ sont *effaçables* dans la catégorie des faisceaux constructibles sur $X$ : pour chaque faisceau (abélien) constructible $\mathscr{F}$ sur $X$, il existe un plongement de $\mathscr{F}$ dans un faisceau *constructible* $\widetilde{\mathscr{F}}$ tel que le morphisme $\mathrm{R}^i f_\star \mathscr{F} \to \mathrm{R}^i f_\star \widetilde{\mathscr{F}}$ soit nul. De plus, il est formel de vérifier que tout morphisme $\mathscr{F}_1 \to \mathscr{F}_2$ de faisceaux constructibles s'insère dans un diagramme

$$\begin{array}{ccc} \mathscr{F}_1 & \longrightarrow & \mathscr{F}_2 \\ \downarrow & & \downarrow \\ \widetilde{\mathscr{F}_1} & \longrightarrow & \widetilde{\mathscr{F}_2} \end{array}$$

où $\widetilde{\mathscr{F}_1}$ et $\widetilde{\mathscr{F}_2}$ effacent respectivement $\mathrm{R}^i f_\star \mathscr{F}_1$ et $\mathrm{R}^i f_\star \mathscr{F}_2$.

Fixons un entier $n > 1$ et supposons que l'on sache calculer, fonctoriellement en $\mathscr{F}$, les $\mathrm{R}^i f_\star \mathscr{F}$ pour chaque entier $i < n$, et chaque faisceau constructible $\mathscr{F}$ de $\Lambda$-modules sur $X$. Considérons un faisceau constructible $\mathscr{F}$ sur $X$ et $\mathscr{F} \hookrightarrow \widetilde{\mathscr{F}}$ un monomorphisme effaçant $\mathrm{R}^n f_\star \mathscr{F}$. Notant $\mathscr{G}$ le faisceau quotient $\widetilde{\mathscr{F}}/\mathscr{F}$, on a la suite exacte $\mathrm{R}^{n-1} f_\star \widetilde{\mathscr{F}} \to \mathrm{R}^{n-1} f_\star \mathscr{G} \to \mathrm{R}^n f_\star \mathscr{F} \to 0$. La calculabilité (fonctorielle) de $\mathrm{R}^n f_\star \mathscr{F}$ se ramène donc à la détermination explicite d'un monomorphisme $\mathscr{F} \hookrightarrow \widetilde{\mathscr{F}}$ comme ci-dessus. La possibilité de plonger tout $\Lambda$-faisceau constructible dans l'image directe par un morphisme fini d'un faisceau constant sur chaque composante connexe nous ramène au problème suivant : trouver un morphisme fini surjectif $\pi : X' \to X$ tel que la flèche de fonctorialité (déduite de l'unité de l'adjonction $\pi^\star \dashv \pi_\star$)

$$\mathrm{R}^n f_\star \Lambda \to \mathrm{R}^n f'_\star \Lambda$$

soit nulle, où $f' = f \circ \pi'$.

**11.4.3**. *Cas d'un morphisme propre.* Quitte à énumérer les morphismes $\pi$ (pour faire une recherche non bornée : cf. 12.8), on se ramène au problème de *tester* la nullité d'une flèche comme ci-dessus. Or, on peut calculer une stratification explicite $S = \bigcup_\iota S_\iota$ (réunion disjointe) telle que les faisceaux $\mathrm{R}^n f_\star \Lambda$ et $\mathrm{R}^n f'_\star \Lambda$ ci-dessus soient lisses sur les $S_\iota$ ; c'est un corollaire immédiat des démonstrations « géométriques » de la constructibilité des images directes. (Voir [Orgogozo 2013] pour des variantes sur ce thème.) Tester si $\mathrm{R}^n f_\star \Lambda \to \mathrm{R}^n f'_\star \Lambda$ est nulle revient donc à tester si sa fibre l'est en tout point géométrique $\overline{\eta}$ localisé en un point maximal $\eta$ des $S_\iota$. Lorsque $f$ est *propre* (ou simplement cohomologiquement propre pour les faisceaux de torsion), cette fibre n'est autre que la flèche de fonctorialité

$$\mathrm{H}^n(X_{\overline{\eta}}, \Lambda) \to \mathrm{H}^n(X'_{\overline{\eta}}, \Lambda).$$

D'après le théorème 0.1, on peut décider si une telle flèche est nulle ou non ; la conclusion résulte alors du fait que les points maximaux $\eta$ sont en nombre fini.

Ceci démontre le théorème 0.9 dans le cas particulier où $f$ est *propre*, ou bien lorsque $S = \mathrm{Spec}(k)$.

**11.4.4**. Notons que l'on peut améliorer légèrement le résultat d'effaçabilité précédent : donné $f : X \to S$ propre (et $\Lambda$), on peut calculer $X' \to X$ fini surjectif

$$\begin{array}{ccc} X & \overset{\pi}{\dashleftarrow} & X' \\ {\scriptstyle f}\downarrow & \swarrow {\scriptstyle f'} & \\ S & & \end{array}$$

tel que le morphisme $\mathrm{R}^{\geq 1} f_\star \Lambda \to \mathrm{R}^{\geq 1} f'_\star \Lambda$ soit nul, où l'on note $\mathrm{R}^{\geq 1} = \tau_{\geq 1}\mathrm{R}$ pour simplifier. Il suffit d'itérer suffisamment la construction (cf. par exemple [Bhatt 2011, § 2]).

Déduisons de cette observation que, si $f$ n'est plus nécessairement propre (mais toujours de type fini), l'existence d'un morphisme *fini* surjectif $\pi : X' \to X$ tel que les $\mathrm{R}^i f_\star \Lambda \to \mathrm{R}^i f'_\star \Lambda$ soient nuls pour $i > 0$ se déduit de l'existence d'une *altération a* effaçant la cohomologie. Supposons en effet que l'on ait un diagramme commutatif

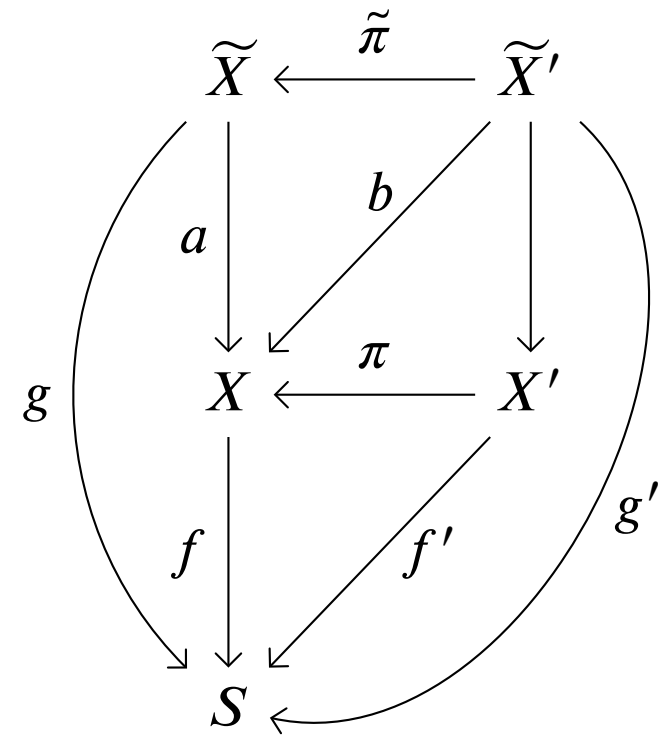

où :

- $a$ est une altération effaçant la cohomologie de $f$ : la flèche de fonctorialité $\mathrm{R}^i f_\star \Lambda \to \mathrm{R}^i g_\star \Lambda$ est nulle pour chaque $i > 0$ ;
- $\tilde{\pi}$ un morphisme fini surjectif effaçant la cohomologie du morphisme *propre* $a$ : la flèche de fonctorialité $\mathrm{R}^{\geq 1} a_\star \Lambda \to \mathrm{R}^{\geq 1} b_\star \Lambda$ est nulle ;
- $\widetilde{X}' \to X' \to X$ est la *factorisation de Stein* de $b$ ; en particulier, $\pi$ est *fini* surjectif.

(D'après ce qui précède, donné $a$, on sait *calculer* un tel diagramme.) Il résulte de la seconde hypothèse que l'on a une factorisation diagonale du carré commutatif ci-dessous

$$\begin{array}{ccc} \Lambda & \longrightarrow & \mathrm{R}^0 b_\star \Lambda \ = \pi_\star \Lambda \\ \downarrow & \nearrow & \downarrow \\ \mathrm{R} a_\star \Lambda & \longrightarrow & \mathrm{R} b_\star \Lambda \end{array}$$

où l'égalité du coin supérieur droit est conséquence de la troisième hypothèse. Fixons $n > 0$ et appliquons le foncteur $\mathrm{R}^n f_\star$. Le carré précédent devient

$$\begin{array}{ccc} \mathrm{R}^n f_\star \Lambda & \longrightarrow & \mathrm{R}^n f'_\star \Lambda \\ {\scriptstyle 0}\downarrow & \nearrow & \downarrow \\ \mathrm{R}^n g_\star \Lambda & \longrightarrow & \mathrm{R}^n g'_\star \Lambda. \end{array}$$

La nullité de la flèche verticale de gauche correspond à la première hypothèse. La flèche horizontale supérieure est donc nulle. CQFD.

**11.4.5.** Les observations précédentes et le théorème de résolution des singularités [A. J. DE JONG 1996, 4.1] ramènent le calcul d'un morphisme fini surjectif effaçant la cohomologie d'un morphisme (non nécessairement propre) $f : X \to S$ au cas particulier où : $X$ est le complémentaire d'un diviseur à croisements normaux stricts $D$ dans un schéma projectif lisse $\overline{X}$ sur $k$, et où l'on s'autorise à effacer par une *altération* $a$. Le cas propre étant déjà connu, il suffit de montrer que l'on peut calculer un diagramme commutatif

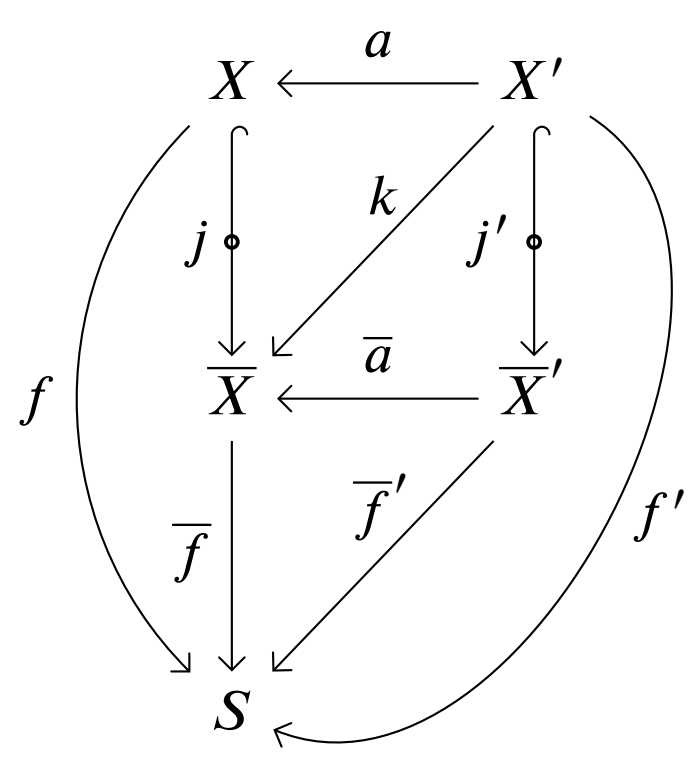

où $\overline{a}$ est une altération et le carré commutatif ci-dessous

$$\begin{array}{ccc} \mathrm{R}^n\overline{f}_\star\Lambda & \longrightarrow & \mathrm{R}^n\overline{f}'_\star\Lambda \\ \downarrow & \nearrow & \downarrow \\ \mathrm{R}^n f_\star\Lambda & \longrightarrow & \mathrm{R}^n f'_\star\Lambda \end{array}$$

se factorise diagonalement comme indiqué. En effet, le cas propre appliqué au morphisme $\overline{f}'$ permet de construire une altération (et même un morphisme fini surjectif) $\overline{X}'' \to \overline{X}'$ effaçant $\mathrm{R}^n\overline{f}'_\star\Lambda$ ; la conclusion est alors immédiate.

Comme au paragraphe précédent (11.4.4), il suffit de montrer que l'on peut calculer une altération $\overline{a}$ effaçant $\mathrm{R}^{\geq 1}j_\star\Lambda$, c'est-à-dire telle que le morphisme $\mathrm{R}^{\geq 1}j_\star\Lambda \to \mathrm{R}^{\geq 1}k_\star\Lambda$ soit nul. En effet, on a sous cette hypothèse un diagramme commutatif et une factorisation

$$\begin{array}{ccccc} \Lambda & \longrightarrow & k_\star\Lambda = \overline{a}_\star\Lambda & \longrightarrow & \mathrm{R}\overline{a}_\star\Lambda \\ \downarrow & \nearrow & & & \downarrow \\ \mathrm{R}j_\star\Lambda & & \longrightarrow & & \mathrm{R}k_\star\Lambda \end{array}$$

induisant après application du foncteur $\mathrm{R}^n\overline{f}_\star$ le diagramme commutatif

$$\begin{array}{ccccc} \mathrm{R}^n\overline{f}_\star\Lambda & \longrightarrow & \mathrm{R}^n\overline{f}_\star(\overline{a}_\star\Lambda) & \longrightarrow & \mathrm{R}^n\overline{f}_\star(\mathrm{R}\overline{a}_\star\Lambda) = \mathrm{R}^n\overline{f}'_\star\Lambda \\ \downarrow & \nearrow & & & \downarrow \\ \mathrm{R}^n f_\star\Lambda & & \longrightarrow & & \mathrm{R}^n f'_\star\Lambda \end{array}$$

désiré. Pour effacer $\mathrm{R}^{\geq 1}j_\star\Lambda$, il suffit — par itération, cf. 11.4.4, premier paragraphe — de savoir effacer chaque $\mathrm{R}^i j_\star\Lambda$ pour $i > 0$ (en nombre fini). Pour tout entier $N > 0$, il existe un morphisme fini surjectif $\overline{X}' \to \overline{X}$ tel que, Zariski-localement sur $\overline{X}$, le tiré en arrière du diviseur $D = \overline{X} - X$ soit une puissance $N$-ième dans $\overline{X}'$. En effet, il est possible de trouver, Zariski-localement, une extension finie du corps des fractions de $\overline{X}$ telle que la clôture intégrale de $\overline{X}$ convienne (cf. [Deligne 1980, § 1.7.9]) ; il suffit alors de considérer la clôture intégrale de $\overline{X}$ dans une extension composée. Vérifions que $\overline{X}' \to \overline{X}$ efface les $\mathrm{R}^i j_\star\Lambda$. Soit $d$ un point géométrique de $D$, $U = X \times_{\overline{X}} \overline{X}_{(d)}$ le complémentaire de $D_{(d)}$ dans le schéma régulier strictement local $\overline{X}_{(d)}$ et $V$ l'ouvert correspondant dans $\overline{X}'_{(d)} = \overline{X}_{(d)} \times_{\overline{X}} \overline{X}'$. La fibre en $d$ de $\mathrm{R}^i j_\star\Lambda \to \mathrm{R}^i k_\star\Lambda$ s'identifie à l'application de fonctorialité $\mathrm{H}^i(U,\Lambda) \to \mathrm{H}^i(V,\Lambda)$. Cette dernière est nulle par construction et pureté lorsque $N \cdot \Lambda = \{0\}$. Notons que l'on peut bien itérer cette construction car, quitte à altérer (ce qui est licite comme on l'a vu), on peut supposer que $\overline{X}'$ est régulier et que l'ouvert $X'$ image inverse de $X$ est le complémentaire d'un diviseur à croisements normaux stricts. Ceci achève la démonstration du théorème 0.9.

### 11.5. Questions.

**11.5.1.** Calculer $\mathrm{H}^i(X,\mathscr{K})$, voire $\mathrm{R}\Gamma(X,\mathscr{K})$, pour $\mathscr{K}$ un complexe borné constructible de $\Lambda$-modules. Variante relative : étendre le théorème 0.9 au calcul de $\mathrm{R}f_\star\mathscr{K}$.

**11.5.2.** (Théorème de changement de base propre effectif.) Donnés $X \to S$ un morphisme *propre* et $s \in S$, construire un voisinage étale $U$ de $s$ tel que toute élément de $\mathrm{H}^*(X_s,\Lambda)$ se relève à $X_U$.

**11.5.3.** Calculer les cycles proches $\psi^r_f(\Lambda)$ d'un morphisme $f : X \to S = \mathrm{Spec}(k[t]_{(t)})$.

**11.5.4.** Calculer les nombres de Betti $\ell$-adiques dans le cas non nécessairement propre et lisse[①].

[①] Ce problème nous semble actuellement hors de portée.

**11.6. Primitive récursivité, ou existence de bornes algorithmiques.** Comme nous l'avons signalé en **0.10** (voir aussi **12.8** pour une explication plus détaillée), la notion de « calculabilité » que nous avons utilisée est la notion classique de calculabilité au sens de Church-Turing, qui permet notamment d'effectuer des « recherches non bornées », c'est-à-dire énumérer des objets (toujours ramenables aux entiers naturels) jusqu'à en trouver un, si on sait qu'il existe, vérifiant une propriété algorithmiquement testable. L'utilisation de ce procédé, sans aucune borne *a priori* sur la longueur des recherches en question, fait perdre tout contrôle sur la complexité de nos algorithmes.

Il nous semble cependant plausible que de telles bornes puissent être trouvées. Plus exactement, nous pensons que les fonctions calculées algorithmiquement dans le présent article sont au moins *primitivement récursives*, c'est-à-dire calculables par un algorithme dont toutes les boucles peuvent être bornées *a priori* au sens où on doit avoir calculé un majorant sur le nombre d'exécutions de toute boucle avant d'entrer dans celle-ci : cf. [Odifreddi 1989, définition I.1.6 et proposition I.5.8]. Ceci interdit l'utilisation des « recherches non bornées » et correspond à la façon la plus naturelle de les interdire[1]. (Pour rendre plus parlante la notion de fonction primitivement récursive, on peut décrire un langage de programmation qui ne permet pas d'appels récursifs de fonctions et dans lequel toutes les boucles sont des boucles bornées par la valeur d'une variable à l'entrée de la boucle : tel est le langage « BlooP » décrit dans [Hofstadter 1999, chapitre XIII], qu'on pourra consulter pour une description agréable à lire de la différence entre fonctions primitivement récursives et générales récursives, ces dernières y étant définies par le langage « FlooP ».) Par ailleurs, les fonctions primitivement récursives sont une « classe de complexité » car une fonction primitivement récursive est une fonction qui peut être calculée algorithmiquement avec une complexité (en espace ou en temps) elle-même donnée par une fonction primitivement récursive : cf. [Odifreddi 1999, VIII.8.8].

Les fonctions que nous calculons sont peut-être même « élémentaires » (ou « élémentairement récursives ») au sens de Kalmár, c'est-à-dire de temps d'exécution borné par une tour d'exponentielles, cf. [ibid., définition VIII.7.1] ; de nouveau, il revient au même de dire qu'elle est calculable algorithmiquement avec une complexité elle-même élémentaire (cf. [ibid., théorème VIII.7.6]).

Cette supposition est motivée, entre autres, par un slogan proposé par certains logiciens, cf. [Friedman 1999, conjecture 1] ou [Avigad 2003, « grand conjecture » page 258], selon lequel tout théorème mathématique « ordinaire »[2], qui peut s'énoncer dans le langage de l'arithmétique, est en fait prouvable dans des systèmes formels faibles de l'arithmétique ; or, dès que ces systèmes prouvent qu'un algorithme termine, ils prouvent en fait qu'il appartient à une classe de complexité bien comprise. On pense notamment au système PRA, « Primitive Recursive Arithmetic » (défini par exemple en [S. G. Simpson 2009, IX, § 3]) et qui démontre la terminaison précisément des fonctions primitivement récursives, (cf. [Hájek et Pudlák 1998, corollaire IV.3.7]), et au système EA, « Elementary [function] Arithmetic », défini en [Avigad 2003, § 2], et qui démontre la terminaison des fonctions Kalmár-élémentaires, (cf. [ibid., théorème 2.2 et remarque qui suit]). Si notre théorème (affirmant qu'un certain algorithme termine en calculant la cohomologie étale) est démontrable dans ces systèmes faibles, c'est que l'algorithme est primitivement récursif voire Kalmár-élémentaire.

Pour s'en convaincre, il faudrait chercher à reprendre tous les résultats d'existence utilisés ici pour faire apparaître une borne explicite sur les objets construits (de manière à placer les algorithmes dans une des hiérarchies décrites dans [Odifreddi 1999, chapitre VIII]). Nous n'avons pas eu le courage de mener cet exercice, mais nous avons au moins cherché à limiter autant les appels aux « recherches non bornées » (notamment en **6.6**, ou encore **15.5** pour le calcul de la normalisation). Il reste que nous n'avons pas réussi à l'éviter dans la construction d'un hyperrecouvrement d'un schéma $X$ par des polycourbes $\ell$-élémentaires (proposition **1.4.9** et suite) : s'il est probable que, dans le cas où $X$ est lisse on puisse sans trop de mal construire explicitement les polycourbes en question en suivant la démonstration de [SGA 4, XI, § 2-3], le cas général nécessiterait aussi de revoir les résultats de [A. J. de Jong 1996] sous un angle algorithmique.

[1] Il n'est malheureusement pas possible de dire que toute fonction intuitivement calculable sans recherche non bornée est primitivement récursive, car il existe différentes sortes de récursion garantissant la terminaison qui ne peuvent pas s'exprimer sous forme primitivement récursive : la fonction d'Ackermann en est un exemple ; on pourra consulter [Odifreddi 1999, VIII.9, notamment les définitions VIII.9.1 et VIII.9.3] pour des notions plus générales. De toute manière, l'argument général du problème de l'arrêt (i.e., un principe diagonal) ne permet pas qu'on puisse formaliser la notion intuitive de « fonction calculable sans recherche non bornée » (donc en particulier, terminant toujours).

[2] Il va de soi que l'affirmation est fausse sans le qualificatif « ordinaire », et que celui-ci ne peut pas être défini rigoureusement. Voir [Smoryński 1985] (ainsi que les autres articles de cet auteur dans le même recueil [Harrington et al. 1985]) pour un aperçu généraliste de ces questions, ou bien l'introduction de [Friedman 201?] pour une présentation plus systématique.

## II. Algèbre commutative et géométrie algébrique effectives

L'objet de cette partie est de vérifier la calculabilité des propriétés et opérations algébriques et géométriques utilisées dans la première partie.

### 12. Corps et extensions de corps

**12.1. Définition.** On appelle **corps calculable** la donnée d'une partie calculable (=récursive) $\mathfrak{K}$ de $\mathbb{N}$, d'une relation d'équivalence calculable $\equiv$ sur $\mathfrak{K}$, d'éléments $0_{\mathfrak{K}}$ et $1_{\mathfrak{K}}$ de $\mathfrak{K}$, et de fonctions calculables $(+)_{\mathfrak{K}} : \mathfrak{K}\times\mathfrak{K} \to \mathfrak{K}$ et $(-)_{\mathfrak{K}} : \mathfrak{K} \to \mathfrak{K}$ et $(\times)_{\mathfrak{K}} : \mathfrak{K}\times\mathfrak{K} \to \mathfrak{K}$ et $({}^{-1})_{\mathfrak{K}} : \{z \in \mathfrak{K} : z \not\equiv 0_{\mathfrak{K}}\} \to \mathfrak{K}$, telles que ces opérations passent au quotient par $\equiv$ et définissent sur $\mathfrak{K}/\equiv$ une structure de corps. On notera généralement $K = (\mathfrak{K}/\equiv)$ et on dira abusivement que $K$ « est » un corps calculable pour sous-entendre qu'on s'est donné une structure de corps calculable dont $K$ est le quotient ; s'il faut désambiguïser, on pourra dire que $\mathfrak{K}$ est l'ensemble d'**étiquettes** qui servent à décrire les éléments de $K$.

Une **extension calculable** $\mathfrak{K} \to \mathfrak{L}$ de corps calculables est la donnée d'une fonction calculable $f : \mathfrak{K} \to \mathfrak{L}$ telle que $x \equiv_{\mathfrak{K}} y$ implique $f(x) \equiv_{\mathfrak{L}} f(y)$ et que l'application $K \to L$ définie par passage au quotient soit une extension (un morphisme) de corps. (On dira aussi que $K$ est un sous-corps calculable de $L$.) Si de plus $K$ est une partie récursive de $L$, c'est-à-dire s'il existe une fonction calculable qui, donné $y \in \mathfrak{L}$, décide s'il existe $x$ tel que $f(x) \equiv_{\mathfrak{L}} y$ (et, cf. **12.8**, on peut alors supposer[①] qu'elle calcule ce $x$), alors on dit que $K$ est un sous-corps calculable **reconnaissable** de $L$ (ou que $L$ est une extension calculable reconnaissable de $K$).

On dit qu'un corps calculable $K$ admet un **algorithme de factorisation** lorsqu'il existe un algorithme qui, donné un polynôme à coefficients dans $K$, calcule sa factorisation en polynômes irréductibles.

**12.2. Remarques. (0)** Nous renvoyons notamment à [Fröhlich et Shepherdson 1956] et [Stoltenberg-Hansen et Tucker 1999] pour des généralités sur les corps calculables. Il existe différentes variantes autour de la définition, essentiellement sans importance dans le cadre dans lequel nous nous plaçons (par exemple, quitte à n'utiliser que le plus petit élément — pour l'ordre de $\mathfrak{K}$ en tant que partie de $\mathbb{N}$ — de chaque classe d'équivalence, on peut omettre la relation d'équivalence et demander directement que $\mathfrak{K}$ soit une partie récursive de $\mathbb{N}$ munie d'opérations qui en font un corps) ; celle proposée ci-dessus (équivalente à celle de [ibid., 2.1.5]), nous semble la plus naturelle et celle qui se transpose le plus agréablement, par exemple, au cas où on remplacerait les fonctions récursives par des fonctions seulement primitivement récursives (cf. [Jacobsson et Stoltenberg-Hansen 1985, § 1]).

Il est notamment utile de rappeler les faits suivants.

**(1)** Si $K$ est un corps calculable, alors $K(T)$ (où $T$ est une indéterminée), ainsi que $K[X]/(f)$ (où $f \in K[X]$ est un polynôme irréductible) sont des extensions calculables et reconnaissables de $K$. (C'est-à-dire qu'il y a une façon standard de faire de $K(T)$ ou de $K[X]/(f)$ des corps calculables et de l'extension une extension calculable reconnaissable, et c'est de cette structure qu'on parlera toujours ; par exemple, un élément de $K[X]/(f)$, si $d = \deg f$, est décrit comme un $d$-uplet $(c_0, \dots, c_{d-1})$ d'éléments de $K$, ou plus précisément de l'ensemble $\mathfrak{K}$ d'étiquettes des éléments de $K$, représentant la classe modulo $f$ du polynôme $\sum c_i X^i$, l'addition se faisant terme à terme et la multiplication se faisant en terminant par le reste de la division euclidienne par $f$, laquelle est évidemment calculable.)

**(2)** Si $K$ est un corps calculable, alors « la » clôture algébrique de $K$ est une extension calculable de $K$ ([Rabin 1960, théorème 7] ; cf. [Stoltenberg-Hansen et Tucker 1999, corollaire 3.1.11]), mais non reconnaissable en général (cf. le point suivant). **(3)** L'existence d'un algorithme de factorisation pour un corps calculable $K$ équivaut à l'existence d'un algorithme qui décide si un polynôme

[①] Ce serait la bonne définition à prendre si on voulait remplacer la notion de calculabilité par celle de fonction primitivement récursive.

admet une racine, ou encore à l'existence d'un algorithme qui reconnaît si un élément de la clôture algébrique de $K$ (calculable comme on vient de le dire en (2)) appartient à $K$ ([Rabin 1960, théorème 8] : cf. aussi [Stoltenberg-Hansen et Tucker 1999, proposition 3.2.2] et [Miller 2010, théorème 2.5(2)]). Bien entendu, **(4)** tout corps calculable algébriquement clos admet un algorithme de factorisation.

De plus, **(5)** si $K$ est un corps calculable admettant un algorithme de factorisation, alors c'est aussi le cas de $K(T)$ (où $T$ est une indéterminée) et de $K[X]/(f)$ si $f \in K[X]$ est un polynôme irréductible *séparable* (cf. [Stoltenberg-Hansen et Tucker 1999, théorèmes 3.2.3 et 3.2.4] ; et **(6)** la nécessité de l'hypothèse « séparable » pour le point précédent est montrée dans [Fröhlich et Shepherdson 1956, théorème 7.12]).

**(7)** Si $K$ est un corps *parfait* calculable admettant un algorithme de factorisation (notamment si $K$ est un corps calculable algébriquement clos), et si $K(x_1, \dots, x_n)$ est une extension de type fini de $K$, puisqu'on peut extraire de $x_1, \dots, x_n$ une base de transcendance séparante ([Matsumura 1989, remarque précédent le théorème 26.3]), le point (5) montre que le corps calculable $K(x_1, \dots, x_n)$ admet lui aussi un algorithme de factorisation (cf. [Stoltenberg-Hansen et Tucker 1999, 3.2.6]).

**12.3**. On rappelle (voir notamment [Bourbaki 2007, V, § 13] ou [ÉGA IV$_1$, **0**, § 21] ou encore [Fried et Jarden 2008, § 2.7]) qu'une *p-base* (resp. une famille *p-libre*) finie d'un corps $K$ de caractéristique $p > 0$ (sous-entendu : sur $K^p$) est une famille $b_1, \dots, b_r \in K$ tels que les produits $b_1^{i_1} \cdots b_r^{i_r}$ pour $0 \leq i_u < p$ forment une base (resp. une famille libre) du $K^p$-espace vectoriel $K$ ; il existe une $p$-base finie de $K$ si et seulement si $K$ est de degré fini sur $K^p$, auquel cas ce degré vaut $p^r$ où $r$ est le cardinal de la $p$-base : on appelle $r$ le *p-rang* ou *exposant d'imperfection* de $K$ (sous-entendu : sur $K^p$).

**12.4**. **Proposition.** *Soit $K$ un corps calculable de caractéristique $p > 0$ et de p-rang fini. Il revient au même de se donner :*

(i) *le p-rang $r$ de $K$ et un algorithme décidant si des éléments $a_1, \dots, a_s$ de $K$ sont linéairement indépendants sur $K^p$ ;*
(ii) *le p-rang $r$ de $K$ et un algorithme décidant si des éléments $a_1, \dots, a_s$ de $K$ sont p-libres;*
(iii) *une p-base $b_1, \dots, b_r$ de $K$ ;*
(iv) *une p-base $b_1, \dots, b_r$ de $K$ et un algorithme exprimant un élément $x$ de $K$ sous la forme $\sum_{\underline{i}} \xi_{\underline{i}}^p b_1^{i_1} \cdots b_r^{i_r}$ où $\xi_{\underline{i}} \in K$ pour $\underline{i} = (i_1, \dots, i_r)$ vérifiant $0 \leq i_u < p$ pour tout $u$.*

« Il revient au même de se donner » signifie qu'on peut exprimer n'importe laquelle de ses données en fonction de n'importe quelle autre de façon algorithmique et uniforme — c'est-à-dire par un algorithme indépendant de $K$ et des autres données.

(On pourra comparer ces équivalences avec [Richman 1981b, théorème 1] qui en est l'analogue dans le cadre de l'algèbre constructive.)

On dira qu'on a sur un corps calculable une ***p*-base finie explicite** en référence à n'importe laquelle de ces données.

*Démonstration.* Il est évident que connaître (i) permet de connaître (ii) (tester la $p$-liberté revient, par définition de ce terme, à tester l'indépendance linéaire de certaines puissances). Connaissant (ii), on peut connaître (iii) en énumérant les éléments de $K$ et en ajoutant ceux qui sont $p$-libres avec les précédents jusqu'à atteindre le $p$-rang de $K$. Connaissant (iii) on obtient (iv) en énumérant toutes les écritures possibles de $x$ sur la $p$-base jusqu'à en obtenir une qui convient. Enfin, connaissant (iv), on a un isomorphisme explicite de $K^p$-espaces vectoriels entre $K$ et $(K^p)^{\oplus p^r}$ (somme de $p^r$ copies de $K^p$), ce qui permet donc facilement de tester l'indépendance linéaire d'une famille (il s'agit simplement de calculer des déterminants). ☐

Comme on l'a souligné en 12.2(5–6), si $K$ est un corps calculable admettant un algorithme de factorisation, ces propriétés valent pour $L := K[X]/(f)$ (avec $f \in K[X]$ irréductible) lorsque $f$ est séparable, mais pas nécessairement dans le cas général. Si on fait l'hypothèse qu'on dispose sur $K$ d'une $p$-base finie explicite, cette difficulté n'existe plus :

**12.4.1**. **Lemme.** *Soit $K$ un corps calculable pour lequel on dispose d'un algorithme de factorisation et d'une p-base finie explicite. Soit $a \in K$ n'appartenant pas à $K^p$. Alors sur le corps $L := K(\sqrt[p]{a})$*

*(extension calculable reconnaissable de $K$ comme rappelé en 12.2(1)), on dispose d'un algorithme de factorisation et d'une $p$-base finie explicite.*

*Démonstration.* Puisque $a \notin K^p$ (autrement dit, le singleton $a$ est $p$-libre), d'après 12.4(ii), on peut explicitement construire une $p$-base de $K$ contenant l'élément $a$ (comme dans la démonstration du fait que (ii) permet de trouver (iii), en partant de $a$), disons $a, b_2, \ldots, b_r$. Alors $a^{1/p}, b_2, \ldots, b_r$ constitue une $p$-base explicite de $L$.

Pour montrer que $L$ dispose d'un algorithme de factorisation, d'après [Fröhlich et Shepherdson 1956, 7.3] ou son amélioration citée dans [Stoltenberg-Hansen et Tucker 1999, 3.2.5], il suffit de montrer qu'on peut décider si un élément de $L$ est une puissance $p$-ième : or d'après 12.4(iv) on sait l'écrire sur la $p$-base $a^{1/p}, b_2, \ldots, b_r$, et il suffit de vérifier que seul le coefficient devant 1 est non nul dans cette écriture. □

Le résultat suivant a pour objet de convaincre que tous les corps que nous serons amenés à considérer sont calculables avec un algorithme de factorisation :

**12.5. Proposition.** *Soit $K$ un corps calculable pour lequel on dispose d'un algorithme de factorisation et (si $K$ est de caractéristique $p > 0$) d'une $p$-base finie explicite. Soit $L$ l'extension de $K$ définie par l'une des opérations suivantes :*

(i) *l'ajout d'un transcendant : $L = K(T)$ où $T$ est une indéterminée,*
(ii) *l'ajout d'un élément algébrique : $L = K[X]/(f)$ où $f \in K[X]$ est irréductible (non supposé séparable), donné,*
(iii) *le passage à la*[①] *clôture algébrique $L = K^{\mathrm{alg}}$ de $K$,*
(iv) *le passage à la clôture séparable $L = K^{\mathrm{sép}}$ de $K$,*
(v) *(dans le cas où $K$ est de caractéristique $p > 0$) le passage à la clôture parfaite $L = K^{1/p^\infty}$ de $K$,*

*alors $L$ est une extension calculable reconnaissable de $K$, et on dispose d'un algorithme de factorisation et d'une $p$-base finie explicite pour $L$.*

Plus précisément, on va esquisser des algorithmes explicites qui, donnés des algorithmes qui calculent les opérations sur $K$ et la factorisation des polynômes de $K[X]$ et une $p$-base finie explicite de $K$, et donnés le cas où on se place, et le polynôme $f$ dans le cas (ii), présentent $L$ comme une extension calculable reconnaissable de $K$, permettent de factoriser les polynômes de $L[X]$, et fournissent une $p$-base de $L$.

(On pourra comparer avec [Mines et Richman 1982, théorème 3.9], analogue de (ii) ci-dessus mais dans le cadre de l'algèbre constructive.)

*Démonstration.* Traitons chacun des cas séparément.

(i) Le corps $L = K(T)$ est une extension calculable reconnaissable de $K$ comme on l'a rappelé en 12.2(1), et dispose d'un algorithme de factorisation d'après 12.2(5). En ajoutant $T$ à la $p$-base de $K$ on obtient une $p$-base de $L$ (cf. [Fried et Jarden 2008, lemme 2.7.2 et sa démonstration]).

Pour le cas (ii), on peut distinguer le cas où $f$ est séparable et celui où il est purement inséparable : en effet, il est algorithmique d'écrire un polynôme irréductible $f \in K[X]$ sous la forme $h(X^{p^e})$ avec $h$ irréductible et séparable, ce qui ramène l'extension $K[X]/(f)$ aux deux extensions $E \coloneqq K[X]/(h)$ avec $h$ séparable puis $L = E[X]/(X^{p^e} - a)$ avec $a$ la classe de $X$ modulo $h$. Par ailleurs, pour les extensions purement inséparables, on peut encore se ramener au cas où le polynôme est de la forme $X^p - a$ (quitte à écrire une racine $(p^e)$-ième comme extractions successives de racines $p$-ièmes).

Dans le cas (ii) avec $f$ séparable, le corps $L$ est une extension calculable reconnaissable de $K$ comme on l'a rappelé en 12.2(1), et dispose d'un algorithme de factorisation d'après 12.2(5). Une $p$-base de $K$ est encore une $p$-base de $L$ (cf. [ibid., lemme 2.7.3]).

Dans le cas (ii) avec $f$ de la forme $X^p - a$ a été traité en 12.4.1.

(iii) La clôture algébrique $L$ de $K$ est une extension calculable de $K$ comme on l'a rappelé en 12.2(2), et dispose d'un algorithme de factorisation (12.2(4)). La $p$-base vide convient pour $L$. De plus, comme $K$ était supposé disposer d'un algorithme de factorisation, on peut reconnaître $K$ dans $L$ d'après 12.2(3).

[①] On conviendra que, s'agissant d'un corps calculable, « la » clôture algébrique désigne celle qui est construite explicitement par l'algorithme de Rabin : cf. 12.2(2) ; et de même « la » clôture séparable désigne la clôture séparable dans cette clôture algébrique.

(iv) La clôture séparable $L$ de $K$ se voit comme un sous-corps de la clôture algébrique : pour reconnaître si un élément de cette dernière appartient à $L$, il suffit de calculer son polynôme minimal sur $K$ (quitte à énumérer tous les polynômes de $K[X]$ jusqu'à en trouver un qui annule l'élément considéré[①], cf. **12.8**, puis le factoriser) et vérifier s'il est séparable. On dispose d'un algorithme de factorisation puisque, d'après **12.2**(3), il suffit pour cela de savoir identifier un élément de $L$ dans la clôture algébrique commune de $K$ et $L$, et on vient d'expliquer que c'est possible. Enfin, une $p$-base de $K$ est encore une $p$-base de $L$ (cf. [Fried et Jarden 2008, lemme 2.7.3]).

(v) Si $K$ est un corps calculable (sans autre hypothèse pour l'instant), alors on peut construire $L = K^{1/p^\infty}$ extension calculable de $K$ explicitement selon sa définition : on définit $\mathfrak{L}$ comme l'ensemble des couples $(e, a)$ où $e \in \mathbb{N}$ et $a$ est un élément de $K$ : ce couple représente alors la racine $(p^e)$-ième de $a$ dans $L$, et on peut définir $(e, a) \equiv_{\mathfrak{L}} (e', a')$ (disons pour $e' \geq e$) lorsque $a^{p^{e'-e}} = a'$ dans $K$, ce qui est bien une relation calculable. Pour ajouter ou multiplier $(e, a)$ et $(e', a')$ (disons pour $e' \geq e$), on remplace $(e, a)$ par $(e', a^{p^{e'-e}})$ et on effectue l'opération entre $a^{p^{e'-e}}$ et $a'$, qui est calculable.

Avec l'hypothèse supplémentaire que $K$ dispose d'un algorithme de factorisation, on peut tester si un élément de $K$ a sa racine $p$-ième dans $K$ (et le cas échéant la calculer) : on peut donc considérer uniquement les couples $(e, a) \in \mathfrak{L}$ « canoniques », définis comme ceux pour lesquels $a \notin K^p$ si $e > 0$, convertir un couple $(e, a) \in \mathfrak{L}$ quelconque en un couple « canonique », et on voit alors clairement que $L$ est une extension calculable reconnaissable de $K$ (c'est d'ailleurs essentiellement ce qui est fait dans [Steel 2005, § 2.1]).

Pour montrer que $L$ dispose d'un algorithme de factorisation, il suffit clairement de montrer qu'on peut factoriser dans $L$ les polynômes $f$ de $K[X]$ (quitte à appliquer l'isomorphisme entre $K^{1/p^e}$ et $K$ pour un $e$ assez grand) : comme $K$ a un algorithme de factorisation, on peut évidemment supposer $f$ irréductible dans $K[X]$, et l'écrire sous la forme $h(X^{p^e})$ avec $h \in K[X]$ irréductible séparable : ceci se réécrit $(h_1(X))^{p^e}$ avec $h_1$ le polynôme de $L[X]$ dont les coefficients sont les racines $(p^e)$-ièmes de ceux de $h$ ; en utilisant de nouveau l'isomorphisme entre $K^{1/p^e}$ et $K$ il est clair que $h_1$ est irréductible dans $L[X]$, et on a la factorisation voulue. Enfin, la $p$-base vide convient pour $L$. □

**12.6. Remarque.** L'énoncé ci-dessus considère des extensions algébriques du type $L = K[X]/(f)$ avec $f$ irréductible (i.e., $(f)$ maximal dans $K[X]$). Il n'y aura pas de difficulté, dès qu'on saura manipuler les idéaux d'une algèbre de polynômes à plusieurs variables (13), à y exprimer des extensions du type $L = K[Z_1, \ldots, Z_d]/\mathfrak{m}$ avec $\mathfrak{m}$ un idéal maximal de $K[Z_1, \ldots, Z_d]$ comme une tour d'extensions monogènes. Il suffit en effet d'utiliser un algorithme d'élimination, cf. [Eisenbud 1995, § 15.10.4], pour calculer les intersections $\mathfrak{m} \cap K[Z_1, \ldots, Z_i]$, qui définissent autant d'extensions de corps $K_i \subseteq K_{i+1}$ avec $K_i = K[Z_1, \ldots, Z_i]/(\mathfrak{m} \cap K[Z_1, \ldots, Z_i])$, algébriques engendrées par un seul élément dont on connaît le polynôme minimal, et $K = K_0 \subseteq K_1 \subseteq \cdots \subseteq K_d = L$.

Plus généralement, si $\mathfrak{p}$ est un idéal premier de $K[Z_1, \ldots, Z_d]$, l'extension $L := \mathrm{Frac}(K[Z_1, \ldots, Z_d]/\mathfrak{p})$ de $K$ se ramène également explicitement à des extensions comme envisagées ci-dessus : en effet, il suffit de considérer un ensemble maximal $E$ de variables $Z_i$ telles que $\mathfrak{p} \cap K[E] = \{0\}$ (de nouveau calculable par élimination), si bien que $K(E) = \mathrm{Frac}(K[E])$ est une extension transcendante pure de $K$ et que $L$ en est une extension algébrique du type considéré au paragraphe précédent (on peut aussi invoquer le lemme de normalisation de Noether [ibid., théorème 13.3] ou [Serre 1965, III(D)2, théorème 2], mais ce n'est pas nécessaire ici car on a recherche simplement la finitude générique, c'est-à-dire la finitude de l'extension de corps).

**12.7. Convention.** On fera souvent l'abus de langage consistant à écrire qu'une certaine opération algébrique ou géométrique est calculable sans autre précision sur le corps $k$ dans lequel vivent les données : sauf mention du contraire, il faut en fait comprendre : pour tout corps *$k$ calculable, disposant d'un algorithme de factorisation* (cf. **12.2**) *et une $p$-base finie explicite* (**12.4**), l'opération en question est calculable. Lorsque $k$ est algébriquement clos, bien entendu, seule l'hypothèse « calculable » est utile (cf. **12.2**(4)) ; par ailleurss, d'après le résultat qu'on vient de montrer, les corps que nous sommes amenés à manipuler vérifient toujours les hypothèses qui viennent d'être dites.

[①] Il va de soi que sur une description réellement explicite de la clôture algébrique on n'aurait pas besoin de faire quelque chose d'aussi absurde !

En fait, si on a choisi d'utiliser le point de vue de la partie III (cf. **0.11**), il faudra même comprendre (ce qui est plus fort, cf. **21.4** et **21.5**) que l'opération décrite est calculable au sens du modèle de calcul « universel » définie en **20.1** et **20.4**.

(Dans le langage de l'algèbre constructive, l'analogue de nos hypothèses « calculable, disposant d'un algorithme de factorisation et une $p$-base finie explicite » serait les corps « discrets pleinement factoriels » : voir [Mines, Richman et Ruitenburg 1988, théorème VII.3.3].)

**12.8. Remarque.** Le fait de considérer des fonctions récursives (générales, par opposition par exemple aux fonctions primitivement récursives) signifie que si $\mathbf{P}(m, n)$ est une propriété calculable (des entiers naturels) et que pour chaque $m$ il existe $n$ vérifiant $\mathbf{P}(m, n)$, alors la fonction $\mu_n\mathbf{P}$ qui à $m$ associe le plus petit $n$ vérifiant $\mathbf{P}(m, n)$ est calculable (il s'agit de l'« opérateur $\mu$ de Kleene », cf. par exemple [Odifreddi 1989, définitions I.1.5 et I.1.7]). Concrètement, cela signifie que les fonctions calculables (=récursives, donc) peuvent effectuer des « recherches non bornées » dans les entiers naturels, c'est-à-dire parcourir tous les $n$ jusqu'à en trouver un qui vérifie la propriété $\mathbf{P}(m, n)$ demandée, même si on ne dispose d'aucune borne de complexité sur le temps qu'une telle recherche pourra prendre (on demande simplement qu'elle termine pour tout $m$ si on veut que la fonction soit totale).

Comme les corps calculables (**12.1**) sont étiquetés par les entiers naturels, cette notion de « recherche non bornée » s'applique aussi bien à eux. C'est la raison pour laquelle, dans la définition d'un corps calculable, on pourrait par exemple se passer d'imposer que la fonction $({}^{-1})_{\mathfrak{K}} : \{z \in \mathfrak{K} : z \not\equiv 0_{\mathfrak{K}}\} \to \mathfrak{K}$ soit calculable (elle l'est automatiquement puisque pour calculer $x^{-1}$ on peut parcourir tous les $y \in \mathfrak{K}$ jusqu'à en trouver un qui vérifie $x \times_{\mathfrak{K}} y \equiv_{\mathfrak{K}} 1_{\mathfrak{K}}$). Notons d'ores et déjà que si on préfère à la notion de corps calculable le modèle de calcul « universel » présenté dans la partie III, bien qu'on ne puisse plus énumérer les « éléments de $\mathbb{A}^1$ », la possibilité d'effectuer des recherches subsiste grâce aux résultats de la section 23.

Une fois définis les schémas et autres objets géométriques en 16, on pourra effectuer de même ce type de « recherches non bornées » sur l'ensemble des schémas, morphismes de schémas, ou tout autre objet géométrique du même type : il s'agira d'énumérer toutes les données susceptibles de décrire, par exemple, un schéma, et pour chacune d'elles de tester si elle vérifie la propriété $\mathbf{P}$ considérée. Ainsi, il découle du théorème de A. J. de Jong sur la résolution des singularités par altérations ([A. J. de Jong 1996, 4.1]), et du fait que cette propriété soit algorithmiquement testable, qu'on peut en fait calculer une telle résolution ! Il va de soi que cette façon de procéder fait perdre toute utilisabilité pratique à nos algorithmes — et comme observé en **11.6** on *devrait* pouvoir s'en passer.

*Sous-remarque :* On pourrait même être un peu plus général dans les recherches par test non bornées : il n'est pas nécessaire que la propriété $\mathbf{P}(m, n)$ soit décidable (=calculable, =récursive), il suffit qu'elle soit semi-décidable (=récursivement énumérable), autrement dit qu'il existe une machine de Turing qui termine en répondant « vrai » lorsque $\mathbf{P}(m, n)$ est vrai, sans imposer qu'elle termine quand $\mathbf{P}(m, n)$ n'est pas vrai — on peut alors calculer *un* $n$ vérifiant $\mathbf{P}(m, n)$ (pas nécessairement le plus petit) en lançant en parallèle la vérification de tous les $\mathbf{P}(m, n)$ jusqu'à ce que l'une d'entre elles termine. (Ceci n'est néanmoins pas compatible avec l'invocation de 23 dans le cas du modèle de calcul universel, s'il s'agit d'énumérer des éléments du corps, mais ce type de recherche est surtout utile pour parcourir des entiers $n$ ; nous n'utiliserons de toute façon pas cette opération.)

## 13. Modules de type fini sur une $k$-algèbre de type fini

Soit $S = k[Z_1, \ldots, Z_d]$ où $k$ est un corps calculable (nous n'utiliserons pas ici l'ensemble des hypothèses **12.7** puisque nous n'aurons jamais affaire à des factorisations de polynômes : la simple calculabilité de $k$ suffit). Le but de cette section est de montrer qu'on peut manipuler algorithmiquement les $S$-modules de type fini, et, plus généralement, les modules de type fini sur une algèbre quotient de $S$.

**13.1. Algorithmes fondamentaux.** Un morphisme $\varphi : S^m \to S^n$ de $S$-modules est représenté par une matrice $n \times m$ d'éléments de $S$. On sait (par exemple en utilisant des bases de Gröbner) répondre algorithmiquement aux deux questions fondamentales suivantes :

- *(Test d'appartenance à un sous-module.)* Donné $\varphi : S^m \to S^n$ et un élément $y \in S^n$, décider si $y$ appartient à l'image $M$ de $\varphi$ et, le cas échéant, en calculer un antécédent. (Ou, de façon équivalente : donnés $x_1, \dots, x_m \in S^n$ et $y \in S^n$, décider si $y$ peut être écrit comme combinaison linéaire de $x_1, \dots, x_m$ et, le cas échéant, en trouver les coefficients. Pour le faire, on peut calculer une base de Gröbner du sous-module $M$ engendré par $x_1, \dots, x_m$, cf. [EISENBUD 1995, algorithmes 15.7 et 15.9], [DECKER et LOSSEN 2006, problème 2.16] ou [BECKER et WEISPFENNING 1993, lemme 6.7 et discussion qui suit], cette base de Gröbner étant elle-même écrite comme combinaison des $x_i$, et ceci permet de décider si un élément $y$ appartient à $M$ et, le cas échéant, l'écrire comme combinaison des élément de la base, dont des $x_1, \dots, x_m$.)

- *(Calcul de syzygies.)* Donné $\varphi : S^m \to S^n$, calculer un ensemble de générateurs du noyau de $\varphi$ (ou, si on préfère, donnés des éléments $x_1, \dots, x_m$ de $S^n$, calculer un système de générateurs des syzygies entre eux, cf. [EISENBUD 1995, algorithme 15.10] ou [BECKER et WEISPFENNING 1993, théorème 6.4]).

(De façon alternative, on pourra se référer à [MINES, RICHMAN et RUITENBURG 1988, VIII.1.5 combiné à III.2.6] ou bien [LOMBARDI et QUITTÉ 2011, VII.1.10] pour une présentation de ces faits en algèbre constructive et sans utiliser de bases de Gröbner.)

Soulignons que les algorithmes en question n'utilisent que des opérations rationnelles et des tests d'égalité dans le corps $k$.

**13.2. Noyaux, conoyaux et images sur les anneaux de polynômes.** Une **présentation explicite** d'un $S$-module $M$ de type fini est la donnée d'un morphisme $\varphi : G \to F$, où $F, G$ sont deux $S$-modules libres de type fini, tel que $M = \operatorname{Coker} \varphi$ (concrètement, $\varphi$ est fourni sous la forme d'une matrice $n \times m$ d'éléments de $S$ où $m$ et $n$ sont les rangs de $G$ et $F$ respectivement).

Une **présentation explicite** d'un morphisme $\alpha : P \to Q$ de $S$-modules, où $P = \operatorname{Coker}(G_P \xrightarrow{\varphi_P} F_P)$ et $Q = \operatorname{Coker}(G_Q \xrightarrow{\varphi_Q} F_Q)$ sont deux $S$-modules explicitement présentés comme ci-dessus, est la donnée d'un morphisme $\alpha_F : F_P \to F_Q$ tel que $\alpha$ se déduise de $\alpha_F$ par passage au quotient (autrement dit, tel qu'il existe $\alpha_G : G_P \to G_Q$ vérifiant $\alpha_F \varphi_P = \varphi_Q \alpha_G$).

Dans le contexte ci-dessus, donnés $P = \operatorname{Coker}(G_P \xrightarrow{\varphi_P} F_P)$ et $Q = \operatorname{Coker}(G_Q \xrightarrow{\varphi_Q} F_Q)$, on sait tester si un $\alpha_F : F_P \to F_Q$ est la présentation explicite d'un morphisme $\alpha : P \to Q$ ; et de plus, on sait calculer algorithmiquement l'image, le conoyau, et le noyau, d'un morphisme explicitement présenté : les algorithmes à cet effet sont bien connus, mais nous les rappelons brièvement ci-dessous pour la commodité du lecteur.

- Pour tester si $\alpha_F : F_P \to F_Q$ passe au quotient et définit (une présentation explicite d')un morphisme $\alpha : P \to Q$ (où $P = \operatorname{Coker}(G_P \xrightarrow{\varphi_P} F_P)$ et $Q = \operatorname{Coker}(G_Q \xrightarrow{\varphi_Q} F_Q)$), il s'agit, grâce à l'algorithme de test d'appartenance évoqué en **13.1**, de tester si les éléments images de la base de $G_P$ par $\alpha_F \varphi_P$ sont dans l'image de $\varphi_Q$ : cf. [DECKER et LOSSEN 2006, problème 4.1 et début de § 4.2]).
- Donné $\alpha : L \to F$ un morphisme entre $S$-modules libres de type fini, on sait calculer une présentation explicite de $\operatorname{Im} \alpha$ (sous la forme $\operatorname{Im} \alpha = \operatorname{Coker}(H \to L)$). En effet, ceci revient exactement à calculer des générateurs des syzygies entre les images par $\alpha$ des éléments de la base de $L$ (c'est-à-dire, des colonnes de la matrice décrivant $\alpha$), algorithme déjà évoqué en **13.1** comme calcul de syzygies.
- Si $\alpha : L \to Q$ est un morphisme de $S$-modules, où $L$ est toujours libre mais cette fois $Q$ est défini par une présentation explicite $Q = \operatorname{Coker}(G_Q \xrightarrow{\varphi_Q} F_Q)$, le morphisme $\alpha$ étant explicitement présenté par la donnée de $\alpha_F : L \to F_Q$, alors on peut encore calculer une présentation explicite de $\operatorname{Im} \alpha$. En effet, il est isomorphe à $(\operatorname{Im} \alpha_F + \operatorname{Im} \varphi_Q)/\operatorname{Im} \varphi_Q$, or $\operatorname{Im} \alpha_F + \operatorname{Im} \varphi_Q$ est l'image du morphisme $L \oplus G_Q \xrightarrow{(\alpha_F, \varphi_Q)} F_Q$ de modules libres, cas traité par le paragraphe précédent, et si $H \to L \oplus G_Q$ est la présentation de cette image, alors $H \to L$ définit la présentation de $\operatorname{Im} \alpha$ recherchée (cf. [ibid., problème 4.2]).
- Donné un morphisme $\alpha : P \to Q$ explicitement présenté entre modules explicitement présentés $P = \operatorname{Coker}(G_P \xrightarrow{\varphi_P} F_P)$ et $Q = \operatorname{Coker}(G_Q \xrightarrow{\varphi_Q} F_Q)$, il est facile de calculer une présentation explicite de $\operatorname{Coker} \alpha$, à savoir $\operatorname{Coker} \alpha = \operatorname{Coker}(F_P \oplus G_Q \xrightarrow{(\alpha_F, \varphi_Q)} F_Q)$ (cf. [ibid., § 4.2.1]). Mais on peut également calculer une présentation de $\operatorname{Im}\alpha$ sous la forme $\operatorname{Coker}(H \to F_P)$ avec $H$ libre (en effet, il s'agit de $\operatorname{Im}(F_P \to Q)$, cas qu'on a traité au paragraphe précédent), et aussi de $\operatorname{Ker}\alpha$ (il s'agit de l'image de $H \to P$, de nouveau le cas qu'on a traité). (Cf. [ibid., problème 4.3] ou [Eisenbud 1995, proposition 15.32].)

On notera bien évidemment qu'on obtient non seulement une présentation explicite de l'image ou du noyau d'un morphisme explicitement présenté $\alpha : P \to Q$, mais aussi une présentation explicite de l'inclusion canonique $\operatorname{Im}\alpha \to Q$ ou $\operatorname{Ker}\alpha \to P$. Ceci permettra aisément de se convaincre, par exemple, qu'on peut calculer des sommes ou intersections de sous-modules.

**13.3. Modules sur les algèbres de type fini quelconques.** Si $R = S/I$ est une algèbre de type fini sur un corps $k$ où $S = k[Z_1, \dots, Z_d]$ et $I$ l'idéal engendré par $h_1, \dots, h_r \in S$, un $R$-module n'est autre qu'un $S$-module annulé par $I$. On peut appeler **$R$-module explicitement présenté** le conoyau d'un morphisme $\varphi : R^m \to R^n$ (décrit par une matrice $n \times m$ d'éléments de $R$) ou, de façon équivalente, comme le conoyau d'un morphisme $\tilde{\varphi} : S^m \to S^n$ tel que $h_j e_i$ soit dans l'image de $\tilde{\varphi}$ pour tout $1 \leq j \leq r$ et tout $1 \leq i \leq n$ (en notant $e_1, \dots, e_n$ la base canonique de $S^n$) : ce critère est algorithmiquement testable, et on passe de façon évidente d'une description à l'autre (dans un sens en reprenant la matrice de $\tilde{\varphi}$ comme matrice de $\varphi$, et dans l'autre en relevant de façon quelconque la matrice de $\varphi$ et en ajoutant des colonnes $h_j e_i$). Les morphismes de $R$-modules sont simplement des morphismes de $S$-modules qui s'avèrent être des $R$-modules, et ce qui précède montre qu'on sait calculer l'image, le conoyau et le noyau d'un morphisme explicitement présenté de $R$-modules.

Notons par ailleurs que pour manipuler des *sous-modules* d'un module explicitement présenté, on peut représenter ceux-ci comme l'image d'un morphisme (ou, ce qui revient au même, par un ensemble d'éléments engendrant le sous-module) ou comme le noyau d'un morphisme : on a vu qu'on peut passer d'une représentation à l'autre. Dès lors, il est clair qu'on peut calculer des sommes ou intersections de sous-modules, de tester l'inclusion ou l'égalité entre deux sous-modules (tout se ramène facilement à tester la nullité d'un sous-module, ce qui est facile si on le décrit comme engendré par certains éléments).

(En particulier, on sait tester l'inclusion et l'égalité d'idéaux.)

**13.4. Produits tensoriels et Hom de modules ; transporteurs et annulateurs.** On continue de noter $R = S/I$ une $k$-algèbre de type fini où $S = k[Z_1, \dots, Z_d]$.

Si $P = \operatorname{Coker}(G_P \xrightarrow{\varphi_P} F_P)$ et $Q = \operatorname{Coker}(G_Q \xrightarrow{\varphi_Q} F_Q)$ sont deux $R$-modules explicitement présentés, on peut calculer une présentation explicite de $P \otimes_R Q$, à savoir $P \otimes Q = \operatorname{Coker}((G_P \otimes F_Q) \oplus (F_P \otimes G_Q) \to F_P \otimes F_Q)$ où les produits tensoriels de modules libres sont triviaux à calculer et la flèche est $(\varphi_P \otimes \operatorname{Id}) \oplus (\operatorname{Id} \otimes \varphi_Q)$.

De même, grâce au fait qu'on sait calculer les noyaux, on peut calculer une présentation explicite du module $\operatorname{Hom}_R(P, Q) = \operatorname{Ker}(\operatorname{Hom}_R(F_P, Q) \to \operatorname{Hom}_R(G_P, Q))$ où $\operatorname{Hom}_R(F_P, Q)$ admet la présentation explicite évidente $\operatorname{Coker}(\operatorname{Hom}_R(F_P, G_Q) \to \operatorname{Hom}_R(F_P, F_Q))$ et de même pour $\operatorname{Hom}_R(G_P, Q)$, et où la flèche entre eux est donnée par le morphisme $\operatorname{Hom}_R(F_P, F_Q) \to \operatorname{Hom}_R(G_P, F_Q)$ de composition à gauche par $\varphi_P$.

Si $N, N'$ sont des sous-modules d'un module $M$ explicitement présenté, on peut calculer l'idéal « transporteur » $(N : N') \coloneqq \{f \in R : fN' \subseteq N\}$ (c'est-à-dire, en calculer des générateurs) : en effet, on peut le voir comme le noyau du morphisme $R \to \operatorname{Hom}_R(N', M/N)$. De même, si $J$ est un idéal de $R$ et $N$ un sous-module de $M$, on peut calculer le sous-module $(N : J) \coloneqq \{x \in M : Jx \subseteq N\}$ de $M$ (car on peut le voir comme noyau de $M \to \operatorname{Hom}_R(J, M/N)$). (Comparer [ibid., exercice 15.41].)

En particulier, on sait calculer l'idéal annulateur $(0 : M)$ d'un module $M$ explicitement présenté.

**13.5. Tor et Ext.** On continue de noter $R = S/I$ une $k$-algèbre de type fini où $S = k[Z_1, \dots, Z_d]$.

Si $M$ est un $R$-module défini par une présentation explicite $M = \operatorname{Coker}(G \to F)$, on peut aisément en calculer la troncation à un ordre quelconque d'une résolution libre : en effet, il suffit de poser $F_0 = F$ et $F_1 = G$ et récursivement construire le noyau de $F_{i+1} \to F_i$ comme l'image d'un morphisme $F_{i+2} \to F_{i+1}$ de modules libres (par l'algorithme de syzygies de 13.1 si on est sur l'algèbre de polynômes $S$, ou par les techniques générales présentées ci-dessus).

Si $F_i$ est une résolution libre d'un $R$-module $P$ explicitement présenté, et si $Q$ est un $R$-module explicitement présenté, en calculant le complexe $F_i \otimes_R Q$, puis l'homologie de celui-ci, on peut calculer $\mathrm{Tor}_i^R(P,Q)$ pour $i$ arbitraire (mais fixé). De même, en calculant la cohomologie du complexe $\mathsf{Hom}_R(F_i, Q)$ on calcule $\mathrm{Ext}^i_R(P,Q)$ pour $i$ arbitraire (mais fixé).

**13.6. Présentation d'une algèbre finie comme module de type fini.** Soit $A$ un anneau quelconque. Une $A$-algèbre de présentation finie en tant que $A$-module est, en particulier, de présentation finie en tant que $A$-algèbre : la proposition facile suivante explicite ce fait (on pourra comparer avec [T. de Jong 1998, § 3], ou [Lombardi et Quitté 2011, exercice IV.15] dans le cas libre) :

**13.6.1. Proposition.** *Soit $A$ un anneau et $B$ une $A$-algèbre engendrée comme $A$-module par les éléments $1, x_1, \ldots, x_n$ dont le module des syzygies (c'est-à-dire le sous-$A$-module $N$ de $A^{n+1}$ formé des $(c_0, \ldots, c_n) \in A^{n+1}$ tels que $c_0 + c_1x_1 + \cdots + c_nx_n = 0 \in B$) soit de type fini engendré par $y_1, \ldots, y_m$.*

*Soit $\Phi \colon A[T_1, \ldots, T_n] \to B$ le morphisme de $A$-algèbres envoyant $T_i$ sur $x_i$ et $J$ son noyau (c'est-à-dire l'idéal des relations algébriques entre $x_1, \ldots, x_n$). On considère chaque $y_i$ comme un élément de $J$ en identifiant $(c_0, \ldots, c_n) \in A^{n+1}$ au polynôme $c_0 + c_1T_1 + \cdots + c_nT_n$ de degré 1. De plus, pour chaque $1 \leq i, j \leq n$, écrivons $x_ix_j = b_0^{(i,j)} + b_1^{(i,j)}x_1 + \cdots + b_n^{(i,j)}x_n$ pour certains $b_u^{(i,j)} \in A$, et soit $q_{i,j} = b_0^{(i,j)} + b_1^{(i,j)}T_1 + \cdots + b_n^{(i,j)}T_n - T_iT_j$, polynôme de degré 2 appartenant à $J$ (relation quadratique).*

*Alors $J$ est engendré, en tant qu'idéal de $A[T_1, \ldots, T_n]$, par les $y_i$ et par les $q_{i,j}$.*

*Démonstration.* Soit $J'$ l'idéal de $A[T_1, \ldots, T_n]$ engendré par les $y_i$ et par les $q_{i,j}$. Fixons un ordre total quelconque sur les monômes en $T_1, \ldots, T_n$ qui raffine l'ordre partiel donné par le degré total (par exemple, l'ordre lexicographique gradué) : il s'agit donc d'un bon ordre sur les monômes. Soit $z$ un élément de $J$ n'appartenant pas à $J'$ et dont le monôme initial (c'est-à-dire, le plus grand monôme intervenant dans $z$ avec un coefficient non nul) soit le plus petit possible. Si $z$ est de degré $\geq 2$, ce monôme initial est divisible par un $T_iT_j$, disons $z = aT_iT_je + u$ où $a \in A$, où $e$ est un monôme et où $u$ ne fait intervenir que des monômes plus petits que $T_iT_je$ : alors $z' := z + aeq_{i,j}$ est congru à $z$ modulo $J'$, donc appartient à $J$ mais non à $J'$, et son monôme initial est strictement plus petit que celui de $z$, une contradiction. Si $z$ est de degré $\leq 1$, alors $z$ appartient à $N$ donc $z$ est engendré par $y_1, \ldots, y_m$, de nouveau une contradiction. ☐

**13.6.2.** En particulier, si $R = S/I$ est une algèbre de présentation finie sur un corps $k$ et si $B$ est un $R$-module explicitement présenté, et si on dispose sur $B$ d'une multiplication décrite, par exemple, sous forme de tous les produits $x_ix_j$ pour $x_i, x_j$ parcourant le système de générateurs donné par la présentation de $B$ comme $R$-module, alors on peut calculer une présentation de $B$ comme $R$-algèbre, et donc comme $k$-algèbre.

Ceci s'applique notamment pour montrer que si $J$ est un idéal de $R$, on peut calculer une présentation de $\mathsf{End}_R(J)$ en tant que $R$-algèbre ou en tant que $k$-algèbre.

## 14. Algèbres de type fini sur un corps : description algorithmique

Comme dans la section précédente, $k$ est ici un corps calculable (nous n'utiliserons pas ici l'ensemble des hypothèses 12.7 puisque nous n'aurons jamais affaire à des factorisations de polynômes : la simple calculabilité de $k$ suffit).

**14.1. Algèbres de type fini sur $k$.** Une algèbre de type fini sur un corps $k$ sera représentée comme un quotient $A = k[X_1, \ldots, X_m]/I$ d'une algèbre de polynômes, c'est-à-dire par la donnée d'un ensemble (fini !) de générateurs de $I$. Un élément de $R$ sera représenté par un polynôme dans $k[X_1, \ldots, X_m]$ qui le relève : on peut ainsi calculer algorithmiquement les sommes et produits dans $A$, et le fait de pouvoir tester l'appartenance à $I$ permet de tester la nullité, donc l'égalité, d'éléments de $A$. Remarquons aussi qu'on sait tester l'inversibilité (un $f \in R$ est inversible si et seulement si, pour n'importe quel $\tilde{f} \in k[X_1, \ldots, X_m]$ le relevant, l'idéal $I + (\tilde{f})$ obtenu en adjoignant $\tilde{f}$ aux générateurs décrivant $I$, est égal à l'idéal unité de $k[X_1, \ldots, X_m]$, chose qu'on sait tester).

**14.2. Morphismes entre algèbres de type fini.** Si $A = k[X_1, \dots, X_m]/I$ et $B = k[Y_1, \dots, Y_n]/J$ sont deux algèbres de type fini sur $k$ comme ci-dessus, on représentera un morphisme de $k$-algèbres $A \to B$ comme la donnée de $m$ éléments $h_1, \dots, h_m$ de $Y$ tels que $f_i(h_1, \dots, h_m) = 0$ pour tout $i$ si $f_1, \dots, f_r$ sont les générateurs choisis pour représenter $I$ (cette condition est évidemment testable algorithmiquement).

Dans ces conditions, on peut aussi présenter $B$ comme $A$-algèbre de la manière suivante : $B = A[Y_1, \dots, Y_n]/\tilde{J}$, où $\tilde{J}$ est l'idéal $J + (x_i - h_i)$ décrit en adjoignant aux générateurs décrivant $J$ les relations supplémentaires identifiant l'image $x_i$ de $X_i$ dans $A$ avec l'élément $h_i$ de $B$ (ou plus exactement, n'importe quel polynôme dans $k[Y_1, \dots, Y_n]$ le relevant). Et réciproquement, donnée une présentation $B = A[Y_1, \dots, Y_n]/J$ où $J$ est décrit par des générateurs, on peut écrire $B = k[X_1, \dots, X_m, Y_1, \dots, Y_n]/\hat{J}$ où $\hat{J}$ est l'idéal obtenu en relevant les générateurs décrivant $J$ et en y adjoignant ceux de $I$ ; et le morphisme $A \to B$ est alors évident. On pourra donc librement choisir entre la représentation d'un morphisme entre $k$-algèbres de type fini et celle d'une algèbre sur une autre algèbre.

Dans les conditions ci-dessus, on peut tester algorithmiquement si le morphisme $A \to B$ est *surjectif*. En effet, sur la description où $A = k[X_1, \dots, X_m]/I$ et $B = k[X_1, \dots, X_m, Y_1, \dots, Y_n]/\hat{J}$, il s'agit de tester, pour chaque $Y_i$, si $Y_i$ est congru modulo $\hat{J}$ à un élément de $k[X_1, \dots, X_m]$, or ceci peut se faire testant, au moyen d'une base de Gröbner, si chaque $Y_i$ appartient à l'idéal initial (cf. [Eisenbud 1995, § 15.2]) de $\hat{J}$ pour un ordre monomial pour lequel $Y_i$ est supérieur à tout monôme en les $X_1, \dots, X_m$. (En effet, si pour chaque $i \in \{1, \dots, n\}$ il existe $u_i \in k[X_1, \dots, X_m]$ tel que $Y_i - u_i \in \hat{J}$ alors pour un tel ordre monomial le terme initial de $Y_i - u_i$, à savoir $Y_i$ d'après l'hypothèse faite sur l'ordre monomial, appartient à l'idéal initial de $\hat{J}$ ; et réciproquement, si chaque $Y_i$ appartient à l'idéal initial de $\hat{J}$, disons que $Y_1 < \cdots < Y_n$, alors chaque $Y_i$ est congru modulo $\hat{J}$ a un polynôme en $X_1, \dots, X_m, Y_1, \dots, Y_{i-1}$, donc en les $X_1, \dots, X_m$.)

On peut aussi tester algorithmiquement si le morphisme $A \to B$ est *injectif* ou même calculer son noyau (un idéal de $A$, qu'on peut représenter comme l'image d'un idéal de $k[X_1, \dots, X_m]$). En effet, si $A = k[X_1, \dots, X_m]/I$ et $B = k[X_1, \dots, X_m, Y_1, \dots, Y_n]/\hat{J}$ où $I \subseteq \hat{J}$, le noyau du morphisme $A \to B$ est (l'image modulo $I$ de) l'intersection $\hat{J} \cap k[X_1, \dots, X_m]$, laquelle se calcule par un algorithme d'élimination (cf. [ibid., § 15.10.4]). Ceci permet naturellement de calculer aussi une présentation de l'image $\operatorname{Im} \varphi \simeq A/\operatorname{Ker} \varphi$ d'un morphisme $\varphi\colon A \to B$ (il n'y a pas redondance avec le paragraphe précédent, car cette présentation de l'image ne permet pas trivialement de savoir si elle est $B$ tout entier sauf, justement, à utiliser ce qui précède).

**14.3. Lemme.** *Soit $\varphi\colon A \to B$ un morphisme de $k$-algèbres (où, ici, $k$ est un anneau quelconque). Pour $f \in A$, on note comme d'habitude $A[\frac{1}{f}]$ la $k$-algèbre $A[T]/(Tf - 1)$ localisée de $A$ en inversant $f$ (et $B[\frac{1}{f}] = B \otimes_A A[\frac{1}{f}]$). Alors l'ensemble des $f \in A$ tels que le morphisme $A[\frac{1}{f}] \to B[\frac{1}{f}]$ déduit de $\varphi$ soit injectif (resp. soit surjectif, resp. soit un isomorphisme) est un idéal de $A$.*

*Démonstration.* Soient $N = \operatorname{Ker} \varphi$ et $Q = \operatorname{Coker} \varphi$, vus comme $A$-modules. La suite exacte $0 \to N \to A \to B \to Q \to 0$, tensorisée par le $A$-module plat $A[\frac{1}{f}]$, donne $0 \to N[\frac{1}{f}] \to A[\frac{1}{f}] \to B[\frac{1}{f}] \to Q[\frac{1}{f}] \to 0$ : on voit donc que l'ensemble des $f$ tels que $A[\frac{1}{f}] \to B[\frac{1}{f}]$ soit injective (resp. surjective, resp. bijective) est l'ensemble des $f$ tels que $N[\frac{1}{f}] = 0$ (resp. $Q[\frac{1}{f}] = 0$, resp. $N[\frac{1}{f}] = 0$ et $Q[\frac{1}{f}] = 0$). Or si $M$ est un $A$-module, dire de $f \in A$ que $M[\frac{1}{f}] = 0$ signifie que chaque élément de $M$ est annulé par une puissance de $f$ (pouvant dépendre de l'élément), c'est-à-dire que $f$ est dans l'intersection des radicaux des annulateurs de tous les éléments $z \in M$ — sous cette forme, il est clair que l'ensemble des $f$ en question est bien un idéal de $A$. ☐

**14.4. Proposition.** *Dans les conditions du lemme ci-dessus (mais en reprenant pour $k$ un corps calculable), si $A$ et $B$ sont des $k$-algèbres de type fini décrites par une présentation, alors on peut algorithmiquement calculer les idéaux indiqués par le lemme qui précède, à condition de savoir calculer le radical d'un idéal de $A$ (ce qui sera possible d'après* **15.2** *au prix des hypothèses* **12.7** *sur le corps $k$).*

*Démonstration.* On note comme précédemment $N = \operatorname{Ker} \varphi$ et $Q = \operatorname{Coker} \varphi$, vus comme, respectivement, un idéal de $A$ et un $A$-module (qui n'est pas, en général, de type fini). Par ailleurs, on suppose $A = k[X_1, \dots, X_m]/I$ et $B = k[X_1, \dots, X_m, Y_1, \dots, Y_n]/\hat{J}$ où $I \subseteq \hat{J}$.

On a déjà expliqué qu'on pouvait calculer l'idéal $\hat{J} \cap k[X_1, \dots, X_m]$ de $k[X_1, \dots, X_m]$ par un algorithme d'élimination : les générateurs ainsi obtenus, lus modulo $I$, engendrent $N$, et fournissent donc une description de celui-ci comme idéal de $A$, ou en particulier, comme $A$-module de type fini. L'idéal des $f \in A$ tels que $A[\frac{1}{f}] \to B[\frac{1}{f}]$ soit injectif, soit $N[\frac{1}{f}] = 0$, est le radical de l'annulateur de $N$, qu'on peut calculer d'après 13.4 et l'hypothèse faite sur la calculabilité du radical.

Pour calculer l'idéal des $f \in A$ tels que $A[\frac{1}{f}] \to B[\frac{1}{f}]$ soit surjectif, on peut supposer que $A \xrightarrow{\varphi} B$ est injectif (quitte à quotienter par $N$, qu'on sait calculer, pour remplacer $A$ par l'image de $\varphi$), autrement dit que $I = \hat{J} \cap k[X_1, \dots, X_m]$. Nous ferons donc cette hypothèse.

Pour chaque variable $Y_i$, dont on note $y_i$ l'image dans $B$, on peut calculer l'idéal formé des $f \in A$ tels qu'il existe $r \geq 0$ pour lequel $f^r y_i \in A$ : en effet, il s'agit du radical de l'idéal $(A : y_i)$ des $f \in A$ tels que $f y_i \in A$, et ce dernier est calculable en travaillant dans le sous-$A$-module de type fini de $B$ engendré par 1 et $y_i$ (dont les relations sont connues : c'est l'intersection de $\hat{J}$ avec $k[X_1, \dots, X_m, Y_i]$). On peut donc aussi calculer l'intersection de ces $n$ idéaux, c'est-à-dire l'ensemble des $f \in A$ tels que pour chaque $i$ il existe $r \geq 0$ vérifiant $f^r y_i \in A$. Mais il est clair que cet idéal est aussi l'idéal des $f \in A$ tels que pour chaque $h \in B$ il existe $r \geq 0$ vérifiant $f^r h \in A$ : c'est bien l'idéal des $f \in A$ tels que $A[\frac{1}{f}] \to B[\frac{1}{f}]$ soit surjectif. □

Vu au niveau du morphisme de schémas $\operatorname{Spec} B \to \operatorname{Spec} A$, l'idéal des $f$ tels que $A[\frac{1}{f}] \to B[\frac{1}{f}]$ soit injectif définit le plus grand ouvert au-dessus duquel $\operatorname{Spec} B \to \operatorname{Spec} A$ est schématiquement dominant, tandis que l'idéal des $f$ tels que $A[\frac{1}{f}] \to B[\frac{1}{f}]$ soit surjectif définit le plus grand ouvert au-dessus duquel $\operatorname{Spec} B \to \operatorname{Spec} A$ est une immersion fermée.

Nous tirons le lemme suivant de [Petersen 2010] :

**14.5. Lemme.** *Soit $\varphi : A \to B$ un morphisme de $k$-algèbres (où, ici, $k$ est un anneau quelconque) : le morphisme* $\operatorname{Spec} B \to \operatorname{Spec} A$ *est une immersion ouverte si et seulement si l'idéal engendré dans $B$ par l'idéal $P$ des $f \in A$ tels que $A[\frac{1}{f}] \to B[\frac{1}{f}]$ soit un isomorphisme, est l'idéal unité de $B$. De plus, si c'est le cas, l'image de cette immersion est l'ouvert complémentaire du fermé de* $\operatorname{Spec} A$ *défini par l'idéal $P$.*

*Démonstration.* Si $\operatorname{Spec} B \to \operatorname{Spec} A$ est une immersion ouverte, son image dans $\operatorname{Spec} A$ est la réunion d'ouverts principaux $D(f_i)$ pour certains $f_i \in A$ (qui *a priori* pourraient ne pas être en nombre fini), pour chacun d'entre eux $A[\frac{1}{f_i}] \to B[\frac{1}{f_i}]$ est un isomorphisme puisque le morphisme de schémas correspondant en est un, et comme les $D(f_i)$ recouvrent $\operatorname{Spec} B$, les images des $f_i$ dans $B$ engendrent l'idéal unité de $B$.

Réciproquement, supposons donnés un certain nombre (que cette fois on peut d'emblée supposer fini) de $f_i$ dans $A$ qui engendrent l'idéal unité de $B$ et tels que les $A[\frac{1}{f_i}] \to B[\frac{1}{f_i}]$ soient des isomorphismes : alors le morphisme $\operatorname{Spec} B \to \operatorname{Spec} A$ est un isomorphisme au-dessus de chacun des ouverts principaux $D(f_i)$ de $\operatorname{Spec} A$, et leurs images réciproques recouvrent $\operatorname{Spec} B$ : il s'agit donc d'un isomorphisme de $\operatorname{Spec} B$ sur l'ouvert réunion des $D(f_i)$ dans $A$. Cet ouvert est bien le complémentaire du fermé défini par les $f_i$. □

En particulier, si $A$ et $B$ sont des $k$-algèbres de type fini décrites par une présentation, alors on peut algorithmiquement tester si $\varphi$ définit une immersion ouverte $\operatorname{Spec} B \to \operatorname{Spec} A$.

**14.6. Morphismes finis.** Soit $A = k[X_1, \dots, X_m]/I$ et $B = k[X_1, \dots, X_m, Y_1, \dots, Y_n]/\hat{J}$ avec $I \subseteq \hat{J}$ (comme en 14.2) la présentation d'un morphisme $A \to B$ entre $k$-algèbres de type fini. Alors on peut tester algorithmiquement si ce morphisme est *fini* (c'est-à-dire si $B$ est un $A$-module de type fini). En effet, c'est le cas si et seulement si les images modulo $\hat{J}$ de chacun des $Y_i$ sont entières sur $A$, ce qui d'une part permet de se ramener au cas où ($n = 1$, c'est-à-dire) $B = k[X_1, \dots, X_m, Y]/\hat{J}$, et d'autre part c'est le cas si et seulement si $\hat{J}$ contient un polynôme dont le monôme initial, pour un ordre monomial pour lequel $Y$ est supérieur à tout monôme en les $X_1, \dots, X_m$, est une puissance de $Y$. Autrement dit, sous ces hypothèses, c'est le cas lorsque l'idéal initial de $\hat{J}$ (pour un tel ordre ; c'est-à-dire l'idéal engendré par les monômes initiaux des éléments de $\hat{J}$) contient une puissance de $Y$. Or les monômes initiaux de la base de Gröbner de $\hat{J}$ (pour l'ordre considéré) engendrent son idéal initial : il s'agit donc simplement de tester si la base contient un élément dont le monôme initial est une puissance de $Y$.

Notons qu'on a alors une présentation explicite de $B$ comme un $A$-module de type fini, dont les générateurs sont les monômes sur les $Y_i$. Dans le cas où $n = 1$, il suffit d'aller jusqu'à la puissance donnée par le degré de l'équation entière satisfaite par $Y$, moins 1).

## 15. Algèbre commutative effective

Dans cette section (en fait, à partir de 15.2) et dans la suite de cette partie, même si elles ne sont pas partout indispensables, nous ferons implicitement les hypothèses 12.7 sur le corps $k$ auquel on a affaire.

**15.1. Fonction de Hilbert et dimension.** Soit $S = k[Z_1, \dots, Z_d]$ où $k$ est un corps, que nous voyons comme une $k$-algèbre *graduée* (par le degré total). Si $M = \mathrm{Coker}(G \to F)$ est un $S$-module gradué explicitement présenté (c'est-à-dire que $F$ et $G$ sont des modules libres de type fini gradués, i.e. des sommes directes finies de $S[n_i]$ où $S[n_i]$ désigne l'algèbre $S$ où le degré est décalé de $n_i$ ; et où la flèche $G \to F$ qui décrit $M$ est homogène de degré 0), alors on sait algorithmiquement calculer la fonction de Hilbert de $M$ (qui à $n$ associe la dimension sur $k$ de l'espace vectoriel des éléments de $M$ homogènes de degré $n$) : c'est-à-dire que non seulement on sait calculer sa valeur en chaque degré $n$ donné, mais on sait aussi calculer un rang à partir duquel cette fonction est un polynôme, et quel est ce polynôme. (Pour ce fait, nous renvoyons à [Eisenbud 1995, théorème 15.26] et [Cox, Little et O'Shea 2007, chapitre 9, § 2–3].)

Dans ce même contexte, on sait calculer une résolution libre graduée *finie* de $M$ (c'est-à-dire de la forme $0 \to F_r \to \cdots \to F_2 \to F_1 \to F_0 \to 0$ avec $F_i$ libre de type fini gradué et les flèches homogènes de degré 0) : cf. [Eisenbud 1995, corollaire 15.11]. (Notons que ce calcul ne fait pas appel à des recherches non bornées : notamment, la longueur $r$ de la résolution est majorée par le nombre $d$ de variables.)

Si $S = R/I$ où $I$ est un idéal homogène de $R$, la fonction de Hilbert d'un $S$-module gradué $M$ est celle du $R$-module gradué sous-jacent à $M$. (En revanche, il n'existe pas en général de résolution libre finie de $M$ comme $S$-module.)

**15.2. Décomposition primaire.** Dans cette section, où on se penche sur la calculabilité d'une décomposition primaire d'un idéal $I$ (dans un anneau de polynômes), on entend par ce terme le calcul de couples $(M_i, \mathfrak{p}_i)$ d'idéaux (du même anneau) tels que $I = \bigcap_i M_i$, l'intersection étant irredondante, avec $\mathfrak{p}_i$ des idéaux premiers deux à deux distincts et $M_i$ (pour chaque $i$) un idéal $\mathfrak{p}_i$-primaire. Nous cherchons donc à la fois à calculer les $M_i$ (qui ne sont pas uniques) et les $\mathfrak{p}_i$ (qui le sont). Remarquons en particulier que la décomposition primaire recouvre le calcul du radical d'un idéal (intersection des $\mathfrak{p}_i$), et permet de tester si un idéal est radical, ou s'il est premier.

**15.2.1.** *Dimension zéro.* Soit $k$ un corps (dont, conformément à la convention 12.7, on suppose qu'il est calculable et dispose d'un algorithme de factorisation et une $p$-base finie explicite). Alors, si $I$ est un idéal de $k[Z_1, \dots, Z_d]$ de dimension 0 (voir la section suivante pour le cas plus général), on peut algorithmiquement calculer le radical, et plus généralement une décomposition primaire, de $I$ : pour cela, on renvoie soit à [Gianni, Trager et Zacharias 1988, § 6], ou bien, pour une présentation peut-être plus simple, à [Becker et Weispfenning 1993, théorème 8.22] combiné à [Steel 2005] pour lever les difficultés liées à l'inséparabilité (cette dernière référence se place dans un cas plus restreint, mais il est aisé de voir que ces hypothèses additionnelles ne servent que pour obtenir celle que nous avons faite de factorisation dans les extensions finies).

**15.2.2.** *Dimension arbitraire.* Soit $k$ un corps (dont, comme dans le paragraphe précédent, et conformément à la convention 12.7, on suppose qu'il est calculable et dispose d'un algorithme de factorisation et une $p$-base finie explicite). Alors, si $I$ est un idéal de $k[Z_1, \dots, Z_d]$, on peut algorithmiquement calculer une décomposition primaire de $I$ en ramenant ce problème à celui de la dimension 0 : on renvoie pour cela à [Gianni, Trager et Zacharias 1988, § 8 et § 9] et [Becker et Weispfenning 1993, théorème 8.101] (cf. aussi [Steel 2005, § 5.3]).

(On pourra remarquer au passage que l'algorithme IDEALDIV2 décrit en [Becker et Weispfenning 1993, p. 268], est primitivement récursif.)

**15.3. Algèbre de Rees et gradué associé.** Soit $R$ une $k$-algèbre de type fini, disons $R = S/I$ où $S = k[Z_1, \dots, Z_d]$ est un anneau de polynômes, et soit $J$ un idéal de $R$ engendré par des éléments $f_1, \dots, f_r \in R$.

L'**algèbre de Rees** associée à cette situation est la sous-algèbre $R[Jt] = R \oplus Jt \oplus J^2t^2 \oplus \cdots$ de l'algèbre $R[t]$ des polynômes en une indéterminée $t$ sur $R$ formée des polynômes dont le coefficient de degré $i$ appartient à $J^i$. On peut calculer une présentation de $R[Jt]$ comme $R$-algèbre, donc aussi comme $k$-algèbre, de la manière suivante. Soit $L$ l'idéal de $S[T_1, \dots, T_r]$ défini comme l'intersection de ce dernier avec l'idéal de $S[T_1, \dots, T_r, t]$ engendré par $I$ et par les $T_i - f_it$ : alors on sait calculer $L$ par un algorithme d'élimination (cf. [Eisenbud 1995, § 15.10.4]). Or $R[Jt]$ est isomorphe au quotient $S[T_1, \dots, T_r]/L$, l'isomorphisme envoyant $a \in R$ sur la classe modulo $L$ de n'importe quel $\hat{a} \in S$ qui le représente (remarquer que $L$ contient $I$), et $f_{u_1} \cdots f_{u_i} t^i \in J^i t^i$ sur la classe de $T_{u_1} \cdots T_{u_i}$ modulo $L$. (Cf. [Vasconcelos 2005, proposition 1.5].)

Le quotient de l'algèbre de Rees $R[Jt]$ par l'idéal $J$ de $R$ définit l'**algèbre graduée associée** à $J$ dans $R$, soit $\mathrm{gr}_J(R) = (R/J) \oplus (J/J^2)t \oplus (J^2/J^3)t^2 \oplus \cdots$ (les $t^i$, qui servent simplement à étiqueter les degrés, sont souvent omis de cette description). D'après ce qui précède, on peut aussi en calculer une présentation comme $(R/J)$-algèbre ou comme $k$-algèbre, à savoir $S[T_1, \dots, T_r]/(L + (f_1, \dots, f_r))$.

Expliquons comment ceci s'adapte au cas d'un $R$-module $M$ explicitement présenté (cf. 13.3) pour obtenir une présentation explicite de $M[Jt] = M \oplus JMt \oplus J^2Mt^2 \oplus \cdots$ comme module sur l'algèbre de Rees $R[Jt]$, ainsi donc que du gradué $\mathrm{gr}_J(M) = (M/JM) \oplus (JM/J^2M)t \oplus \cdots$ comme module sur $\mathrm{gr}_J(R)$. À partir d'une présentation $M = S^n/Q$ de $M$ comme $S$-module (où $Q \subseteq S^n$ contient $I^n$), on définit encore le $S[T_1, \dots, T_r]$-module $L$ — tel que $M[Jt]$ soit isomorphe à $S[T_1, \dots, T_r]^n/L$ — comme l'intersection de $S[T_1, \dots, T_r]^n$ avec le sous-$S[T_1, \dots, T_r, t]$-module de $S[T_1, \dots, T_r, t]^n$ engendré par $Q$ et les produits de $T_i - f_it$ par les éléments de la base canonique de $S[T_1, \dots, T_r, t]^n$. On peut calculer une présentation de $L$ car la théorie de l'élimination fonctionne encore pour les sous-modules des modules libres de type fini sur les anneaux de polynômes ([Eisenbud 1995, remarque suivant la proposition 15.29 et exercice 15.37]).

**15.4. Calculs de longueurs et de multiplicités.** Soit comme dans la section précédente (15.3) $R$ une $k$-algèbre de type fini, et soit maintenant $\mathfrak{p}$ un idéal *premier* de $R$. (On rappelle que grâce à 15.2 on sait tester si un idéal de $R$ est premier.) Remarquons que le corps $\kappa_\mathfrak{p} := R_\mathfrak{p}/\mathfrak{p}R_\mathfrak{p} = \mathrm{Frac}(R/\mathfrak{p})$ est calculable, et on a même les propriétés (12.7) d'admettre un algorithme de factorisation et une $p$-base finie explicite en vertu de 12.5 et des remarques 12.6. Si $M$ est un $R$-module explicitement présenté, on a obtenu ci-dessus une présentation explicite de $\mathrm{gr}_\mathfrak{p} R$ comme algèbre de type fini sur $R/\mathfrak{p}$ et de $\mathrm{gr}_\mathfrak{p} M$ comme module sur $\mathrm{gr}_\mathfrak{p} R$ : ceci donne donc également une présentation explicite de $\mathrm{gr}_\mathfrak{p} R_\mathfrak{p}$ comme algèbre de type fini sur $\kappa_\mathfrak{p}$ et de $\mathrm{gr}_\mathfrak{p} M_\mathfrak{p}$ comme module sur cette algèbre de type fini. D'après 15.1, on peut calculer la fonction de Hilbert de $\mathrm{gr}_\mathfrak{p} M_\mathfrak{p}$, c'est-à-dire la fonction $\dim_{\kappa_\mathfrak{p}} \mathfrak{p}^i M_\mathfrak{p}/\mathfrak{p}^{i+1} M_\mathfrak{p}$ (« calculer » au sens où on peut à la fois calculer le polynôme avec lequel cette fonction coïncide pour $i$ assez grand, et expliciter une borne à partir de laquelle elle coïncide avec lui).

En particulier, on sait calculer la *multiplicité* $\mathrm{mult}_{R_\mathfrak{p}}(\mathfrak{p}, M_\mathfrak{p})$ de $\mathfrak{p}$ dans $M_\mathfrak{p}$ (cf. [Serre 1965, V(A)2] et [Eisenbud 1995, remarque suivant le corollaire 12.5]).

Expliquons comment, dans le cas où $\mathfrak{a}$ est seulement supposé être un idéal $\mathfrak{p}$-primaire de $R$ (avec $\mathfrak{p}$ un idéal premier), on peut *majorer* la multiplicité $\mathrm{mult}_{R_\mathfrak{p}}(\mathfrak{a}, M_\mathfrak{p})$ (nous n'aurons pas besoin d'un calcul exact). En calculant la fonction de Hilbert de $\mathrm{gr}_\mathfrak{p} R/\mathfrak{a}$, on calcule (sans effectuer de recherche non bornée) un $r$ assez grand pour que $\mathfrak{p}^r \subseteq \mathfrak{a}$, ce qui entraîne $\mathfrak{p}^{ri} M \subseteq \mathfrak{a}^i M$ pour tout $i$ : la longueur sur $R_\mathfrak{p}$ de $M_\mathfrak{p}/\mathfrak{a}^i M_\mathfrak{p}$ est donc majorée par celle de $M_\mathfrak{p}/\mathfrak{p}^{ri} M_\mathfrak{p}$, or les paragraphes précédents montrent qu'on sait calculer la longueur sur $R_\mathfrak{p}$ de $\mathfrak{p}^j M_\mathfrak{p}/\mathfrak{p}^{j+1} M_\mathfrak{p}$ (qui est la dimension de ce $\kappa_\mathfrak{p}$-espace vectoriel), donc de $M_\mathfrak{p}/\mathfrak{p}^j M_\mathfrak{p}$.

**15.5. Normalisation.**

**15.5.1.** Soient $A$ un anneau nœthérien réduit, $K$ son anneau total des fractions et $\widetilde{A}$ le normalisé de $A$ dans $K$. La détermination constructive de $\widetilde{A}$ est un problème classique, qui fait l'objet d'une littérature abondante : citons notamment les articles [Stolzenberg 1968], [T. de Jong 1998], [Singh et Swanson 2009] (en caractéristique positive uniquement) et les livres [Vasconcelos 2005, chap. 6] et [Huneke et Swanson 2006, chap. 15]. Pour la commodité du lecteur, nous rappelons ici brièvement un argument (tiré de [T. de Jong 1998]). Soit $I$ un idéal de $A$ contenant un élément $i$ non diviseur de zéro et posons $A' = \mathsf{End}_A(I)$. Le morphisme évident $A \to A'$ est entier ; il est injectif, de même que le morphisme $A' \to K$, $\varphi \mapsto \varphi(i)i^{-1}$, qui est indépendant du choix de $i$. Notons que la structure d'anneau de $A'$ est facilement explicitable, d'abord comme $A$-module puis comme anneau ; cf. **13.4** et **13.6.2** ou bien [ibid., § 3]. Il résulte de ce qui précède que si $A$ est normal, le morphisme $A \to A'$ est un isomorphisme. On a la réciproque suivante, moyennant des hypothèses supplémentaires.

15.5.2. **Proposition** (Grauert-Remmert). *Soient $k$ un corps* parfait *et $A = k[x_1, \dots, x_n]/(f_1, \dots, f_r)$ un anneau intègre de normalisé $\widetilde{A}$ dans son corps des fractions. Notons $I$ le radical de l'idéal de Fitting $J = \mathrm{Fitt}_d(\Omega^1_{A/k})$ engendré par l'image dans $A$ des déterminants des sous-matrices jacobiennes de taille $n-d$. Alors :*

- ***(i)*** *l'idéal $J$ est non nul et inclus dans l'idéal conducteur $\mathfrak{c} \coloneqq \mathrm{Ann}_A(\widetilde{A}/A)$ ;*
- ***(ii)*** *le morphisme $A \to \mathsf{End}_A(I)$ est un isomorphisme si et seulement si $A$ est normal.*

Pour (i), voir par exemple [Singh et Swanson 2009, rem. 1.5], où il est implicitement fait usage du théorème de normalisation de Noether « génériquement étale » ([Eisenbud 1995, cor. 16.18], [Zariski et Samuel 1975, V, § 3, th. 8]). (Voir aussi [Huneke et Swanson 2006, exercice 12.12].) Pour (ii), voir par exemple [Grauert et Remmert 1984, VI, § 5],

L'anneau $A$ étant japonais, la suite croissante $A \subseteq A' \subseteq A'' \subseteq \cdots$ obtenue en itérant la construction $A \mapsto \mathsf{End}_I(A)$, pour $I = \sqrt{\mathrm{Fitt}(\Omega^1_{A/k})}$, est stationnaire. D'après (ii), sa limite est le normalisé $\widetilde{A}$ que l'on souhaite calculer.

**15.5.3**. La méthode précédente ne dit rien sur le nombre d'opérations à faire : on s'arrête simplement lorsque le morphisme d'inclusion d'un terme dans le suivant est un isomorphisme, condition que l'on sait tester (**14.2**). L'existence de bornes *a priori*, mais non calculables, est connue ([Dries et K. Schmidt 1984, § 3]) mais inutile ici. On se propose ici de montrer que l'on peut contrôler cette terminaison par un calcul de multiplicité et ainsi calculer la normalisation sans faire de recherche non bornée.

Rappelons ([Serre 1965, IV, th. 11]) qu'un anneau est normal si et seulement si il est $R_1$, c'est-à-dire régulier en codimension 1, et satisfait la condition $S_2$ de Serre ([ÉGA IV$_2$, 5.7.2]).

La condition $S_2$ est facile à satisfaire : si $P$ est une sous-$k$-algèbre de polynômes de $A$ telle que $P \to A$ soit fini, et que l'on note $D(-) = \mathsf{Hom}_P(-, P)$, le bidual $B = D(D(A))$ est $S_2$ et est la $S_2$-ification de $A$, c'est-à-dire que $A \to B$ est la plus petite extension finie contenue dans l'anneau total des fractions de $A$ qui soit $S_2$. Voir [Vasconcelos 2005, prop. 6.21], [Huneke et Swanson 2006, démonstration du théorème 15.3.3 et exercices 15.9 et 15.12]. (La multiplication sur $B$ utilisée ici provient de son plongement dans le corps des fractions de $A$ ; mais elle peut aussi se définir de façon intrinsèque : si $\xi, \eta \in D(D(A))$ sont vues comme des formes $P$-linéaires sur $D(A)$, leur produit est la forme $P$-linéaire sur $D(A)$ qui à $\varphi \in D(A)$ associe $\eta(x \mapsto \xi(y \mapsto \varphi(xy)))$ — cette multiplication du bidual, appelée « multiplication d'Arens » par les analystes dans le contexte des algèbres de Banach, cf. [Palmer 1974], n'est pas commutative en général, mais il facile de se convaincre qu'elle l'est, et qu'elle coïncide bien avec la restriction de la multiplication sur $\mathrm{Frac}(A)$, dans le cas où on s'est placé, cf. [Madore 2014].) Voir aussi [ÉGA IV$_2$, § 5.10] et [Achar et Sage 2009] pour une présentation intrinsèque de $S_2$-ification. Or $P$ est calculable en utilisant une démonstration explicite du lemme de normalisation de Noether ([Eisenbud 1995, théorème 13.3 en utilisant le lemme 13.2(a)] ou [Serre 1965, III(D)2, théorème 2]), on sait décrire explicitement $A$ comme un $P$-module de type fini d'après la remarque faite en **14.6**, on peut en déduire la structure de $B$ comme $P$-module d'après **13.4**, et bien sûr comme algèbre (**13.6.1**, puisque la multiplication sur $D(D(A))$ est calculable).

D'autre part, si $A$ est $R_1$, il en est de même de $B$ (en fait, $A \to DD(A)$ est un isomorphisme en codimension 1 ; cf. [Huneke et Swanson 2006, exercice 15.11]).

Ceci nous ramène donc à normaliser (=régulariser) en codimension 1. D'après la proposition précédente (i), on peut trouver un élément non nul $f$ dans le conducteur $\mathfrak{c}$. Notons $\mathfrak{p} \in \mathrm{Spec}(A)$ un point maximal (de codimension 1) de $V(f)$, calculé par décomposition primaire de $f$. Soit $B$ le localisé de $A$ en $\mathfrak{p}$ et notons $\widetilde{B}$ le normalisé de $B$. La longueur de toute chaîne strictement croissante $B \subsetneq B' \subsetneq \dots \subsetneq \widetilde{B}$ est majorée par l'entier $\mathrm{long}_B(\widetilde{B}/B)$ lui-même inférieur ou égal à $\mathrm{long}_B(\widetilde{B}/\mathfrak{c})$. Rappelons que ce dernier est égal à la multiplicité $\mathrm{mult}(\mathfrak{c}, B)$ de $\mathfrak{c}$ dans $B$, à son tour inférieur ou égal à $\mathrm{mult}((f), B)$, que l'on sait majorer (cf. **15.4**). Ceci montre que l'on a peut calculer une borne sur le nombre d'étapes pour rendre $A$ régulier en codimension 1 [①].

Justifions brièvement l'égalité $\mathrm{long}_B(\widetilde{B}/\mathfrak{c}) = \mathrm{mult}(\mathfrak{c}, B)$. Soit $n \geq 1$ un entier ; on a $\mathrm{long}_B(\widetilde{B}/\mathfrak{c}^n) = \sum_{\mathfrak{q}} [\kappa(\mathfrak{q}) : \kappa(\mathfrak{p})] \cdot \mathrm{long}_{\widetilde{B}_{\mathfrak{q}}}(\widetilde{B}_{\mathfrak{q}}/\mathfrak{c}^n)$, où $\mathfrak{q}$ parcourt le spectre maximal (fini) de l'anneau de Dedekind semi-local $\widetilde{B}$. Par régularité des anneaux $\widetilde{B}_{\mathfrak{q}}$ et nouvelle application de cette formule (pour $n = 1$), on obtient l'égalité $\mathrm{long}_B(\widetilde{B}/\mathfrak{c}^n) = n \cdot \mathrm{long}_B(\widetilde{B}/\mathfrak{c})$ d'où $\mathrm{mult}(\mathfrak{c}, \widetilde{B}) = \mathrm{long}_B(\widetilde{B}/\mathfrak{c})$. Comme d'autre part $\mathrm{long}_B(\widetilde{B}/B)$ est finie, on a $\mathrm{mult}(\mathfrak{c}, \widetilde{B}) = \mathrm{mult}(\mathfrak{c}, B)$.

**15.6. Lissité et étalitude.** Soient $A$ un anneau, $P = A[X_1, \dots, X_n]$ un anneau de polynômes sur $A$ et $B$ le quotient de $P$ par un idéal de type fini $I$. Considérons le $B$-module $\mathrm{Ext}^1_B(\mathbb{L}_{B/A}, I/I^2)$, où $\mathbb{L}_{B/A}$ est le complexe cotangent. D'après [Illusie 1971-1972, III.1.2.9.1], ce module est naturellement isomorphe au quotient $\mathsf{End}_B(I/I^2)/d^\star\mathsf{Hom}_P(\Omega^1_{P/A}, I/I^2)$, où $d$ désigne la dérivation $I/I^2 \to \Omega^1_{P/A} \otimes_P B$. On peut donc en déterminer la structure (cf. **13.4**). Notons $H_A(P, I)$ son idéal annulateur en tant que $P$-module ; d'après une variante du critère jacobien ([Gabber et Ramero 2003, 5.4.2]), le lieu lisse (à la source) du morphisme $\mathrm{Spec}(B) \to \mathrm{Spec}(A)$ est l'ouvert complémentaire du fermé $\Sigma$ défini par l'idéal image de $H_A(P, I)$ dans $B$. (Un avantage de cette description est que, contrairement à la description plus classique du lieu singulier par un idéal de Fitting du $B$-module $\Omega^1_{B/A}$ (cf. **15.5**), elle ne fait pas d'hypothèse d'équidimensionalité ou de platitude.) On trouve en [ibid., § 5.4.5] et [Elkik 1973, § 0.2] une variante moins intrinsèque mais plus explicite qui, donnés des générateurs $\{f_1, \dots, f_r\}$ de $I$, produit — par dérivation de ces générateurs et opérations élémentaires sur les idéaux de $P$ (cf. **13.2** et **13.3**) — un idéal contenu dans $H_A(P, I)$ et définissant également le fermé $\Sigma$.

Il résulte de ce qui précède que l'on peut tester si le morphisme de présentation finie $\mathrm{Spec}(B) \to \mathrm{Spec}(A)$ est *lisse*, par exemple lorsque $A$ est une algèbre de type fini sur un corps. Pour vérifier s'il est *étale*, il suffit de vérifier si les fibres au-dessus des points maximaux de $\mathrm{Spec}(A)$ (c'est-à-dire les points génériques des composantes irréductibles, calculables d'après **15.2**) sont vides ou de dimension nulle.

**15.7. Présentation d'un pincement.**

**15.7.1. Définition.** Soit $A$ un anneau et $I$ un idéal de $A$. On définit une $A$-algèbre $B$ en munissant le $A$-module $A \oplus I$ de la multiplication donnée par $(a, u) \cdot (b, v) = (ab, av + bu + uv)$. On dit que $B$ est la $A$-algèbre obtenue par **pincement de $A$ le long de $I$**.

Remarquons que $B$ peut aussi se voir comme la $A$-algèbre $A \times_{A/I} A$ — en identifiant le couple $(a, u)$ de $A \oplus I$ à $(a, a + u)$ de $A \times_{A/I} A$ — ou comme la $A$-algèbre obtenue en ajoutant une unité à l'idéal $I$ vu comme une $A$-algèbre-non-unitaire.

**15.7.2. Proposition.** *Soient $k$ un corps, $A = k[t_1, \dots, t_d]/(u_1, \dots, u_q)$ une $k$-algèbre de type fini explicitement présentée et $I$ l'idéal de $A$ engendré par des éléments $x_1, \dots, x_r$ de $A$, images d'éléments $\hat{x}_1, \dots, \hat{x}_r$ de $k[t_1, \dots, t_d]$). Alors on peut trouver algorithmiquement une présentation explicite de l'algèbre $B = A \oplus I$ obtenue par pincement de $A$ le long de $I$ (cf. **15.7.1**).*

*Démonstration.* Le $A$-module $B$ est engendré par $1, x_1, \dots, x_r$. D'après **13.3**, on sait trouver algorithmiquement des générateurs $y_1, \dots, y_s$ des relations linéaires entre les $x_1, \dots, x_r$ (donc entre les $1, x_1, \dots, x_r$). Quitte à écrire chaque $x_i x_j$ comme combinaison de $x_1, \dots, x_r$, on trouve également des relations quadratiques $q_{i,j}$ telles que définies en **13.6.1**, qui assure alors que les relations $y_i$ et les $q_{i,j}$ définissent l'algèbre $B$ comme un quotient de $A[t'_1, \dots, t'_r]$, donc, quitte à prendre leur relèvement et y ajouter les $u_1, \dots, u_q$, comme un quotient de $k[t_1, \dots, t_d, t'_1, \dots, t'_r]$. □

## 16. Schémas de type fini sur un corps : description algorithmique

On rappelle que les hypothèses **12.7** sont implicitement faites sur le corps $k$. Par ailleurs, nous ne parlerons ici que de schémas de type fini sur $k$, omettant fréquemment les mots « de type fini ».

[①] Cette technique nous a été suggérée par O. Gabber.

**16.1. Schémas affines et quasi-affines, morphismes entre iceux.** On représentera un *schéma affine* (sous-entendu : de type fini) $X$ sur un corps $k$ par une algèbre de type fini $R$ dont il est le spectre, cette algèbre $R$ étant elle-même représentée comme un quotient $k[Z_1, \dots, Z_d]/I$ d'une algèbre de polynômes, c'est-à-dire par la donnée d'un ensemble (fini !) de générateurs de $I$, équations de $X$ dans l'espace affine $\mathbb{A}^d_k$. Remarquons qu'on peut tester algorithmiquement si un tel schéma est vide (il s'agit exactement de tester si $I$ est l'idéal unité).

On représentera un morphisme $X \to Y$ de $k$-schémas affines, où $X = \mathrm{Spec}(R)$ et $Y = \mathrm{Spec}(S)$ sont les spectres de deux $k$-algèbres de type fini $R, S$, comme un morphisme $S \to R$ de $k$-algèbres (cf. 14.2).

On représentera un schéma *quasi-affine*, c'est-à-dire un ouvert $U$ d'un schéma affine $X$, au moyen d'un fermé dont il est le complémentaire (décrit par des équations $f_i$ de ce fermé : ceci revient à écrire $U$ comme la réunion des ouverts principaux $D(f_i)$). Remarquons que, à l'intérieur d'un schéma affine $X$ fixé, on sait tester l'inclusion ou l'égalité entre des ouverts (cela revient à tester l'inclusion entre les radicaux des idéaux définissant les fermés : cf. 15.2 pour le calcul du radical).

Un morphisme d'un schéma affine $X$ vers un schéma quasi-affine $V$ ouvert complémentaire de $Z$ dans un schéma affine $Y$ sera décrit comme un morphisme $X \to Y$ qui se factorise par $V$ : ce fait est testable algorithmiquement en testant si l'image réciproque de $Z$ par $X \to Y$ est vide (c'est-à-dire si l'algèbre produit tensoriel de celles de $X$ et $Z$ au-dessus de celle de $Y$ est nulle).

Un morphisme d'un schéma quasi-affine $U$ ouvert d'un schéma affine $X$, vers un schéma quasi-affine $V$, sera décrit comme une collection de morphismes $U_i \to V$ qui coïncident sur $U_i \cap U_j$, où les $U_i$ sont les ouverts affines principaux recouvrant $U$ (c'est-à-dire les $D(f_i)$ avec $f_i$ parcourant des équations d'un fermé dont $U$ est le complémentaire dans $X$ : on sait écrire des équations de $D(f_i)$ comme schéma affine en le considérant comme une hypersurface d'équation $T f_i - 1$ au-dessus de $X$, où $T$ est une nouvelle indéterminée ; et on peut représenter l'ouvert affine $U_i \cap U_j$ comme $D(f_i f_j)$). De nouveau, on peut calculer la composée de tels morphismes entre schémas quasi-affines, et tester l'égalité de deux d'entre eux (même si le quasi-affine $U$ de départ n'est pas représenté par le même recouvrement par des $D(f_i)$ : il suffit de prendre un raffinement commun entre deux recouvrements, ce qui est facile).

**16.2. Description des schémas et de leurs morphismes.** On représentera un *schéma* (sous-entendu : de type fini) $X$ sur un corps $k$ par la donnée d'un nombre fini de schémas affines $U_i$ et, pour chacun, d'un recouvrement $V_{ii'}$ par des ouverts quasi-affines (les variables $i$ et $i'$ parcourent ici le même ensemble fini) et, pour chaque paire $i, i'$, d'un morphisme $\varphi_{ii'} : V_{ii'} \to V_{i'i}$ vérifiant la condition de compatibilité que ($V_{ii} = U_i$ et $\varphi_{ii} = \mathrm{Id}_{U_i}$ et que) $\varphi_{i_2 i_3} \circ \varphi_{i_1 i_2}$ et $\varphi_{i_1 i_2}$ coïncident là où tous deux sont définis (en particulier, les $\varphi_{ii'}$ sont des isomorphismes). Toutes ces conditions sont bien testables algorithmiquement, et le schéma $X$ défini est alors le recollement des $U_i$ en identifiant l'ouvert $V_{ii'}$ de $U_i$ avec l'ouvert $V_{i'i}$ de $U_{i'}$ au moyen de $\varphi_{ii'}$. On dit aussi que les $U_i$ (avec les autres données les accompagnant) constituent un **atlas** affine de $X$. Un **raffinement** d'un tel atlas est un atlas obtenu en remplaçant chaque $U_i$ par un recouvrement de celui-ci par des ouverts affines principaux $U_{ij}$ (ici le $j$ parcourt un ensemble fini qui peut dépendre de $i$), avec les données évidemment déduites de ce recouvrement.

On représentera un morphisme de schémas $X \to Y$, décrits par des atlas $U_i$ pour $X$ et $V_j$ pour $Y$, en se donnant un raffinement $U_{ij}$ de l'atlas initial de $X$ et des morphismes $U_{ij} \to V_j$ de schémas affines, qui se recollent aux intersections décrites par l'atlas.

Notons qu'on peut algorithmiquement calculer la composée de morphismes de schémas ainsi décrits, et par ailleurs que, donnés deux morphismes $X \to Y$ entre les deux mêmes schémas décrits par les *mêmes* atlas, on peut tester leur égalité (ceci se fait en prenant un raffinement commun aux deux atlas de $X$ qui décrivent les morphismes à comparer, ce qu'on peut faire puisqu'il s'agit de raffiner des recouvrements de mêmes schémas affines $U_i$).

Remarquons aussi que si dans la définition d'un schéma on ne suppose plus les $U_i$ affines mais que ce sont des schémas plus généraux (autrement dit, si on cherche à recoller un atlas formé de schémas non nécessairement affines), on peut encore se ramener algorithmiquement à la situation où ils sont affines (quitte à remplacer chaque schéma $U_i$ par un atlas affine qui le décrit).

Nous ne prétendons pas qu'il soit possible de tester l'égalité (l'isomorphisme) de deux schémas décrits par des atlas (ceci sera néanmoins possible dans le cas étale : cf. 17.3 ci-dessous). Remarquons à ce sujet que dans la suite si nous écrivons par exemple « si $Z$ est affine » il faut comprendre « si on s'est donné une description de $Z$ comme schéma affine » et pas « si on s'est donné de $Z$ une description comme un schéma général, et qu'il s'avère que $Z$ est affine [chose qu'on ne sait pas tester] » : ceci ne devrait pas prêter à confusion.

En fait, à ce stade de la description, nous ne savons même pas encore tester si un morphisme *donné* entre schémas est un isomorphisme (ceci, en revanche, sera bien décidable : cf. **16.5**).

**16.3. Produits fibrés de schémas.** Donnés $X, Y, Z$ trois schémas (i.e., $k$-schémas de type fini) décrits comme précédemment, et donnés $X \to Z$ et $Y \to Z$ deux morphismes, on peut algorithmiquement calculer le produit fibré $X \times_Z Y$. En effet, si $X, Y, Z$ sont affines (disons $X = \operatorname{Spec} R$, $Y = \operatorname{Spec} S$ et $Z = \operatorname{Spec} A$), il s'agit de calculer un produit tensoriel de $k$-algèbres de type fini, or on a vu en **14.2** qu'on pouvait calculer des présentations finies de $R$ et $S$ comme $A$-algèbres, auquel cas leur produit tensoriel se calcule simplement en réunissant les générateurs et les relations. La démonstration dans le cas général suit celle de [ÉGA I, 3.2.6] : si $Z$ est toujours supposé affine mais que $X, Y$ ne le sont plus, on obtient un atlas de $X \times_Z Y$ comme en [ÉGA I, 3.2.6.3], en prenant les $U_i \times_Z V_j$ pour $U_i$ parcourant un atlas de $X$ et $V_j$ de $Y$ ; si $Z$ n'est plus supposé affine, étant donné un atlas $W_i$ de $Z$, la donnée même des morphismes $X \to Z$ et $Y \to Z$ fournit des atlas de $X$ et $Y$ appropriés à les représenter, c'est-à-dire par des $U_i \to W_i$ où $U_i$ est un ouvert de $X$ (non nécessairement affine, mais réunion d'ouverts affines) et $V_i \to W_i$ de même, on peut donc calculer les produits fibrés $U_i \times_{W_i} V_i$ comme on vient de le dire, et les recoller à leur tour.

Dans le cas où $X = Y$, en recollant les morphismes $R \otimes_A R \to R$ (de multiplication, qui ne posent pas de difficulté à décrire algorithmiquement), on obtient une description de la diagonale $X \to X \times_Z X$ d'un morphisme $X \to Z$ quelconque.

**16.4. Pincement (cas non nécessairement affine).** On renvoie à [FERRAND 2003] pour la question générale de l'existence des pincements dans la catégorie des schémas, dont on tire notamment (théorème 7.1(B)) le fait que la somme amalgamée (= « pincement ») $X \amalg_F X$ est représentable dans la catégorie des schémas lorsque $X$ est un schéma et $F$ un sous-schéma fermé de $X$, et plus exactement, représentée par la même somme amalgamée dans la catégorie des espaces annelés.

On a vu en **15.7.2** comment calculer algorithmiquement ce pincement si $X$ est un schéma affine de type fini sur un corps $k$ (si $X = \operatorname{Spec} A$ et $F = \operatorname{Spec}(A/I)$ avec $I$ un idéal de $A$, alors $X \amalg_F X = \operatorname{Spec}(A \oplus I)$ où $A \oplus I$ est muni de la structure d'algèbre décrite en **15.7.1**). Remarquons par ailleurs que la construction est fonctorielle : si $\varphi\colon A \to A'$ est un morphisme d'algèbres et $I' = \varphi(I)\,A'$ l'idéal engendré dans $A'$ par l'image de l'idéal $I$ de $A$, alors $(A \oplus I) \to (A' \oplus I')$ défini par $(a, u) \mapsto (\varphi(a), \varphi(u))$ est bien un morphisme d'algèbres.

Si $X$ est un schéma de type fini sur $k$ et $U$ un ouvert affine de $X$, alors d'après [ibid., lemme 4.4] $U \amalg_V U$ est un ouvert de $X \amalg_F X$, où $V$ désigne l'ouvert $F \cap U$ de $F$ : ces ouverts sont calculables algorithmiquement d'après ce qu'on vient de dire (il est clair que $F \cap U$ est calculable, au besoin d'après **16.3**), ils recouvrent $X \amalg_F X$, et si $U, U'$ sont deux ouverts affines de $X$, l'ouvert quasi-affine $(U \amalg_V U) \cap (U' \amalg_{V'} U') = (U \cap U') \amalg_{V \cap V'} (U \cap U')$ se décrit comme réunion des $U'' \amalg_{V''} U''$ pour $U'' \subseteq U \cap U'$ et s'envoie vers $U \amalg_V U$ et $U' \amalg_{V'} U'$ par les morphismes de fonctorialité (cf. ci-dessus). Ceci fournit donc une description de $X \amalg_F X$ au sens de **16.2**.

**16.5. Immersions ouvertes et isomorphismes de schémas.** Donné un schéma $X$ décrit par un atlas d'ouverts affines $U_i = \operatorname{Spec} R_i$ (dont on notera $V_{ij} = U_i \cap U_j$ les intersections deux à deux que l'atlas identifie) et un schéma $Y$ que nous supposons dans un premier temps affine $Y = \operatorname{Spec} S$, on peut tester algorithmiquement si un morphisme $X \to Y$ est une immersion ouverte : en effet, pour ceci, il faut et il suffit que chacun des $U_i \to Y$ soient des immersions ouvertes, et que l'intersection (dans $Y$) des images de $U_i$ et $U_j$ coïncide avec $V_{ij}$ en tant qu'ouvert de $U_i$ ; or d'après **14.5** on sait tester si $U_i \to Y$ est une immersion ouverte, on peut calculer $U_i \cap U_j$ dans $Y$ comme $\operatorname{Spec}(R_i \otimes_S R_j)$, et tester si l'ouvert en question de $U_i = \operatorname{Spec} R_i$ coïncide bien avec $V_{ij}$.

Le cas où le schéma $Y$ cible n'est plus supposé affine ne pose pas de difficulté particulière : le morphisme est décrit dans des atlas adaptés, et il est une immersion ouverte si et seulement si sa restriction à chaque carte de la cible est une immersion ouverte.

On peut également tester algorithmiquement si un morphisme $X \to Y$ de schémas est un isomorphisme et, le cas échéant, calculer sa réciproque. En effet, en considérant d'abord le cas où $Y$ est affine, ceci se fait en testant d'abord s'il s'agit d'une immersion ouverte comme décrit par le pénultième paragraphe dont nous reprenons les notations, puis, si c'est bien le cas, en vérifiant que les $U_i$ recouvrent $Y$, ce qui se fait en calculant (toujours par **14.5**) un fermé complémentaire de $U_i$ dans $Y$, ce qui permet de tester si l'intersection de ces fermés est vide (les équations qui les décrivent engendrent l'idéal unité de $S$) ; le cas échéant, les ouverts principaux dont les $U_i$ sont écrits comme réunion forment un atlas de $Y$ pour lequel l'écriture de l'isomorphisme réciproque est claire. Le cas où $Y$ n'est pas affine ne présente pas de difficulté nouvelle (il s'agit, de nouveau, d'être un isomorphisme en restriction à chaque carte de la cible).

**16.6. Réduit, composantes irréductibles, composantes connexes.** On peut tester algorithmiquement si un schéma est réduit, ou bien calculer son réduit : en effet, $X$ est réduit si et seulement si chacun des ouverts affines constituant un atlas de $X$ l'est, ce qui ramène à tester si une $k$-algèbre de type fini est réduite, donc si un idéal d'un anneau de polynômes est radical (cf. **15.2**) ; et réduire $X$ se fait en réduisant chaque carte d'un atlas de $X$ (les morphismes de recollement passent au réduit grâce à la fonctorialité de celui-ci).

De même, on peut tester algorithmiquement si un schéma est irréductible, ou bien calculer ses composantes irréductibles.

Enfin, en testant parmi les composantes irréductibles celles dont l'intersection est vide (on a déjà observé ci-dessus qu'on pouvait tester si un schéma est vide), on peut calculer algorithmiquement les composantes connexes d'un schéma.

**16.7. Image fermée, immersions fermées.** Si $X \to Y$ est un morphisme de schémas, on peut en calculer algorithmiquement l'image fermée schématique (c'est-à-dire, [ÉGA I, 9.5], le plus petit sous-schéma fermé $Y'$ de $Y$ par lequel le morphisme $X \to Y$ se factorise, et bien sûr on cherche aussi à calculer la factorisation en question) : en effet, si on suppose d'abord que $Y = \operatorname{Spec} S$ est affine, et si $X$ est décrit par un atlas d'ouverts affines $U_i = \operatorname{Spec} R_i$, alors $Y' = \operatorname{Spec}(S/N)$ où $N$ est l'intersection des noyaux de tous les $S \to R_i$ décrivant le morphisme (cf. [ÉGA I, 9.5.2]). Si $Y$ n'est plus supposé affine, de nouveau, il suffit d'effectuer cette construction localement.

Grâce à cette construction, on peut tester si un morphisme de schémas $X \to Y$ est une immersion fermée : en effet, cela revient à tester si la factorisation $X \to Y'$, où $Y'$ est l'image fermée schématique décrite ci-dessus, est un isomorphisme, et on a vu comment tester ce fait.

On peut également tester si un morphisme de schémas est dominant (c'est-à-dire d'image dense), puisque cela signifie précisément que son réduit est schématiquement dominant (c'est-à-dire que son image fermée schématique, définie ci-dessus, est toute la cible).

**16.8. Image d'un morphisme, surjectivité.** Si $X \to Y$ est un morphisme de schémas, on peut tester algorithmiquement s'il est surjectif, ou même calculer son image (une partie constructible de $Y$, décrite comme combinaison booléenne de fermés de $Y$). En effet, cette image est caractérisée par l'ensemble de ses points à valeurs dans la clôture algébrique de $k$ (on rappelle que nous ne considérons ici que des schémas de type fini sur $k$ ; cf. [ÉGA $\mathrm{I}_{\mathrm{Spr}}$, proposition 7.1.8 et sa démonstration]), or pour $k$ algébriquement clos il est algorithmique de calculer l'image d'une partie constructible de $k^n$ par la projection sur un sous-ensemble de ses coordonnées ([Fried et Jarden 2008, théorème 9.3.1], cf. aussi **23.1**), ce qui suffit à calculer l'image de $X(k)$ par $f$ comme la projection du graphe de $f$.

**16.9. Morphismes et schémas séparés, radiciels.** On peut tester si un morphisme de schémas est séparé, resp. radiciel, en testant si sa diagonale (qu'on sait calculer d'après **16.3**) est une immersion fermée ([ÉGA I, 5.4.1]), resp. une surjection ([ÉGA $\mathrm{IV}_1$, 1.8.7.1]), ce qu'on sait tester d'après **16.7**, resp. **16.8**.

**16.10. Détection de points.** Expliquons brièvement pourquoi, donné un schéma $X$ non vide sur $k$, on peut algorithmiquement en expliciter un point géométrique, c'est-à-dire un point sur « la » clôture algébrique de $k$. Comme il suffit de trouver un point d'un ouvert affine de $X$, on peut évidemment supposer que $X$ est affine. Le cas où $X$ est un fermé de la droite affine $\mathbb{A}^1_k$ est trivial (s'il est décrit comme l'ensemble $\{f = 0\}$ des zéros d'un polynôme $f$, on considère une racine de $f$ dans la clôture algébrique de $k$), et le cas où $X$ est un ouvert de $\mathbb{A}^1_k$ ne l'est pas moins (si l'ouvert en question est d'écrit comme $\{f \neq 0\}$, on peut considérer le fermé $\{f = 1\}$ qui y est contenu, se ramenant ainsi au cas précédent). Dans le cas général, on procède par récurrence sur la dimension de l'espace affine dans lequel $X$ est inclus : si $\pi$ est la projection sur une coordonnée, alors $\pi(X)$ est calculable d'après 16.8, on peut trouver un point géométrique de $\pi(X) \subseteq \mathbb{A}^1_k$ d'après ce qui vient d'être dit, et on est ramené au même problème dans la fibre de $\pi$ au-dessus de ce point.

**16.11. Passage à la limite.** Soient $k_0$ un corps, $X_0$ un $k_0$-schéma séparé de type fini et $k/k_0$ une extension algébrique (non nécessairement finie). On note $X$ le produit fibré $X_0 \times_{k_0} k$ et on suppose donné un morphisme de type fini $f : T \to X$. Alors, on peut construire une extension *finie* $k_1/k_0$ et un morphisme $f_1 : T_1 \to X_1 = X_0 \times_{k_0} k_1$ induisant le morphisme $f$ par changement de base : il suffit de considérer le sous-corps de $k$ engendré des coefficients définissant $f$. Si $f$ est étale, on peut supposer qu'il en est de même de $k_1/k_0$ (par invariance topologique du site étale).

## 17. Géométrie algébrique effective

On rappelle que les hypothèses 12.7 sont implicitement faites sur le corps $k$, et que les schémas sur $k$ sont supposés de type fini.

**17.1. Dimension.** Grâce au calcul de la fonction de Hilbert (cf. 15.1), il est possible de calculer la dimension d'un schéma (cf. [Cox, Little et O'Shea 2007, § 9.3]), ou la dimension de ses composantes connexes et irréductibles.

**17.2. Lissité et étalitude.** Les propriétés d'être lisse ou étale étant locales, le fait de pouvoir tester cette propriété sur un morphisme d'algèbres (15.6) permet de tester si un morphisme de schémas est lisse, resp. étale.

**17.3. Sections et isomorphismes de morphismes étales.** On peut décider algorithmiquement si un morphisme étale séparé $X \to S$ admet une section. En effet, d'après [ÉGA IV$_4$, 17.4.9] (ou bien [Milne 1980, corollaire I.3.12]), c'est le cas si et seulement si sa restriction à une réunion de composantes connexes de $X$ est un isomorphisme.

Par conséquent, on peut aussi décider si deux morphismes étales finis $X \to S$ et $Y \to S$ sont $S$-isomorphes. Pour s'en convaincre, constatons d'abord qu'on peut calculer le morphisme $\underline{\mathrm{Isom}}_S(X, Y) \to S$ où $\underline{\mathrm{Isom}}$ est le schéma paramétrant les $S$-isomorphismes $X \to Y$ : ses équations s'écrivent explicitement en fonction de celles de $X$ et $Y$. Plus exactement, on peut décrire $\underline{\mathrm{Isom}}_S(X, Y)$ comme le fermé des $(u, v) \in \underline{\mathrm{Hom}}_S(X, Y) \times_S \underline{\mathrm{Hom}}_S(Y, X)$ défini par les équations $v \circ u = \mathrm{Id}_X$ et $u \circ v = \mathrm{Id}_Y$, où $\underline{\mathrm{Hom}}_S(X, Y)$ (aussi noté $\mathfrak{R}_{X/S}(Y \times_S X)$, où $\mathfrak{R}$ désigne la restriction à la Weil) est le schéma paramétrant les morphismes $X \to Y$ : voir [Bosch, Lütkebohmert et Raynaud 1990, § 7.6, notamment théorème 4], dont la démonstration fournit une description explicite de ce schéma (cf. aussi [Debarre 2001, chapitre 2]). Comme $\underline{\mathrm{Isom}}_S(X, Y) \to S$ est lui-même étale fini quand $X$ et $Y$ le sont, et que savoir si $X$ et $Y$ sont isomorphes sur $S$ revient à savoir s'il a une section, ce qui nous ramène à la question précédente.

**17.4. Proj.** Soient $r$ un entier, $k$ un corps et $S$ l'algèbre graduée de polynômes $k[x_0, \dots, x_r]$, où les variables $x_i$ sont de degré 1. On va montrer que pour chaque un $S$-module gradué de type fini $M$, on peut calculer le $k$-module de type fini $\mathrm{H}^0(\mathbb{P}^r_k, \mathscr{M})$, où $\mathscr{M}$ désigne le faisceau quasi-cohérent naturellement associé à $M$ ([ÉGA II, 2.5]). Rappelons ([Serre 1955, ch. III, § 4, ¶ 69, cor. 2]), que l'on a $\mathrm{H}^0(\mathbb{P}^r_k, \mathscr{M}) = \mathrm{colim}_\nu \mathsf{Hom}_S((x_0^\nu, \dots, x_r^\nu), M)_0$. (L'indice 0 indique que l'on ne considère que les morphismes de degré nul.) Étant donné une résolution libre (de type fini) $L_\bullet$ de $M$ — que l'on peut déduire d'une présentation de $M$ comme conoyau d'un morphisme $L_1 \to L_0$ entre $S$-modules gradués libres de type fini —, il résulte de [ibid., ch. III, § 3, ¶ 63, prop. 3 (a)] que l'on a égalité $\mathrm{H}^0(\mathbb{P}^r_k, \mathscr{M}) = \mathsf{Hom}_S((x_0^\nu, \dots, x_r^\nu), M)_0$ dès que $\nu + r$ est supérieur à chaque entier $n$ tel que $S(-n)$ apparaisse comme facteur direct des $L_i$. Rappelons qu'un $S$-module gradué libre est somme directe de modules $S(n)$, où $n \in \mathbb{Z}$ et $S(n)$ est le $S$-module $S$ muni de la graduation $S(n)_i = S_{n+i}$. (Un tel $\nu$ est aussi lié à la « régularité », au sens de Mumford-Castelnuevo, de $M$, cf. [Bayer et Mumford 1993, définition 3.2].)

Pour le calcul d'une résolution libre, cf. 15.1.

**17.5. Application.** Soient $k$ un corps et $X = V(I) \subseteq \mathbb{P}^r_k$ un schéma projectif. Il résulte de ce qui précède (et des résultats de la section 13) que l'on peut calculer $\mathrm{H}^0(X, \mathscr{O}_X) = \mathrm{H}^0(\mathbb{P}^r_k, \mathscr{O}_{\mathbb{P}^r_k}/I)$. En particulier, on peut vérifier si la flèche naturelle $k \to \mathrm{H}^0(X, \mathscr{O}_X)$ est un isomorphisme : ceci fournit une approche alternative à **16.6** pour le calcul des composantes connexes géométriques d'une variété projective.

La même approche, en remplaçant $k$ par une $k$-algèbre $A$ de type fini (qu'on peut supposer être une algèbre de polynômes) permettrait de calculer la factorisation de Stein d'un morphisme projectif $f : Y \to \operatorname{Spec} A$ en calculant la $A$-algèbre $f_*\mathscr{O}_Y$ : on pourrait en déduire la même chose pour $f : Y \to X$ sur une base non nécessairement affine.

## III. Modèle de calcul « universel »

### Motivation

Le but de cette partie est de présenter un formalisme permettant de définir le concept de fonction calculable sur un corps en évitant de passer par les corps calculables (**12.1**), et qui conduit à une notion plus restrictive (voir le corollaire **21.4** et la remarque **21.5**) qui permet d'obtenir automatiquement des énoncés tels que **0.7**.

Ce formalisme se passe d'hypothèse de type « calculabilité » sur le corps, et revient à travailler sur tous les corps algébriquement clos (ou même simplement parfaits) de caractéristique fixée : autrement dit, nous allons définir une notion de fonction calculable $\mathbb{N}^r \times \mathbb{A}^n_{\mathbb{F}_p} \to \mathbb{N}$ ou $\mathbb{N}^r \times \mathbb{A}^n_{\mathbb{F}_p} \to \mathbb{A}^1_{\mathbb{F}_p}$. On peut procéder de deux manières, dont nous prouverons l'équivalence.

L'une considère des machines semblables aux machines de Turing qui peuvent manipuler, outre les entiers naturels, les éléments d'un corps (parfait de caractéristique $p \geq 0$) dont elles ne savent rien, sous forme de « boîtes noires » sur lesquelles elles peuvent effectuer les opérations algébriques (addition, soustraction, multiplication, inverse, Frobenius inverse) et le test d'égalité : c'est ce que nous verrons en 21. C'est cette approche qui permet de se convaincre que les fonctions qu'on a naturellement envie de qualifier de « calculables », et notamment celles qualifiées de telles dans la partie II, sont bien calculables au sens du modèle décrit ici.

L'autre approche (que nous utiliserons essentiellement comme définition en **20.1** et **20.4**, et dont on peut trouver un énoncé simple pour les fonctions totales en **22.4**) consiste à définir les fonctions calculables par une stratification calculable en parties constructibles de $\mathbb{A}^n_{\mathbb{F}_p}$ sur chacune desquelles parties la fonction prend une valeur constante (si elle doit renvoyer une valeur dans $\mathbb{N}$) ou définie par une fonction rationnelle composée avec $\mathrm{Frob}_p^{-i}$ (si elle doit renvoyer une valeur dans $\mathbb{A}^1_{\mathbb{F}_p}$). C'est ainsi que, par exemple, savoir calculer si une variété algébrique affine est, disons, géométriquement normale, revient exactement, pour le modèle de calcul ici décrit, à savoir calculer quelles fibres d'un morphisme entre variétés algébriques affines le sont.

Ceci fournit donc un point de vue effectif sur la constructibilité des propriétés géométriques tel qu'exposé en [ÉGA IV$_3$, § 9] (voir notamment [ÉGA IV$_3$, corollaire 9.9.5]) : *le fait de montrer qu'une propriété géométrique est calculable (dans le sens présenté ici) montre automatiquement qu'elle est constructible* (au moins pour des schémas de type fini).

Les résultats de cette partie ne sont pas difficiles, et peuvent s'approcher soit par la géométrie algébrique classique soit par la théorie des modèles (le point central pouvant reposer sur la compacité de la topologie constructible sur $\mathbb{A}^n_{\mathbb{F}_p}$, cf. **18.3**, ou sur le théorème de compacité de la logique du premier ordre ; de même, les résultats de la section 23 peuvent s'expliquer par le théorème de Chevalley ou l'élimination des quantificateurs dans la théorie des corps algébriquement clos). Nous avons cherché à évoquer les deux points de vue, qui se complètent. Chacun des deux serait d'ailleurs sans doute susceptible de généralisations (pour ne pas considérer uniquement des corps mais des schémas arbitraires ou des modèles d'une théorie quelconque qui réalise l'élimination des quantificateurs[①]) : nous n'avons pas cherché à poursuivre cette piste mais montrons dans la dernière section (26) comment le modèle peut s'adapter au cas où on ne veut pas fixer la caractéristique $p$.

[①] Un exemple évident serait la théorie des corps réel-clos, correspondant aux fonctions calculables par des machines pouvant en manipuler les éléments dans des boîtes noires permettant les opérations algébriques et la comparaison pour l'ordre : le rôle de $\mathbb{A}^n_{\mathbb{F}_p}$ serait alors joué par l'espace des types de cette théorie.

Signalons par ailleurs que le formalisme de cette partie peut être considéré soit comme un cas particulier de la récursion sur les structures générales tel que défini par exemple dans [Moschovakis 1974] ou [Fitting 1981] (mais nous utiliserons ici les spécificités du cas des corps, à commencer par l'élimination des quantificateurs), soit comme une réinterprétation « classique » de l'algèbre constructive telle qu'exposée par exemple dans [Mines, Richman et Ruitenburg 1988] ou [Lombardi et Quitté 2011] (mais l'algèbre constructive est formulée dans le cadre de la logique intuitionniste : nous avons préféré travailler dans le cadre plus familier des mathématiques classiques ; au sujet de ce lien, voir aussi la parenthèse suivant le corollaire 22.4).

On rappelle que l'*exposant caractéristique* d'un corps désigne sa caractéristique lorsque celle-ci est différente de 0, et 1 lorsque la caractéristique est 0.

Nous attirons par ailleurs l'attention du lecteur sur le fait que, nous serons dans cet annexe fréquemment amenés à considérer des fonctions partielles[1] : une fonction partielle $f : E \to F$ est simplement une fonction [totale] $E_0 \to F$ pour une partie $E_0 \subseteq E$, et on dit que $f(x)$ est *définie* (resp. non définie) si $x \in E_0$ (resp. $x \notin E_0$). Nous avons pris garde à utiliser l'adjectif « partiel » ou « total » pour qualifier les fonctions dès qu'il pouvait y avoir ambiguïté.

## 18. Types et parties constructibles

**18.1. Définition.** Soit $p$ un nombre premier ou bien 1, et soit $\mathbb{F}_p$ le corps premier de cet exposant caractéristique (c'est-à-dire le corps fini à $p$ éléments si $p > 1$ et $\mathbb{Q}$ lorsque $p$ vaut 1). On appelle **formule booléenne** en les indéterminées $X_1, \dots, X_n$ sur $\mathbb{F}_p$ une combinaison booléenne, c'est-à-dire une combinaison par les connecteurs $\wedge$ (et logique), $\vee$ (ou logique) et $\neg$ (négation logique) d'expressions de la forme $f(X_1, \dots, X_n) = 0$ où $f \in \mathbb{F}_p[X_1, \dots, X_n]$. On dira indifféremment qu'une telle formule est « vraie », « valide », « vérifiée » ou « satisfaite » en des éléments $(x_1, \dots, x_n)$ d'un corps $k$ d'exposant caractéristique $p$ lorsque sa valeur de vérité évaluée sur les éléments en question est vraie.

À une telle formule $\phi$ on associe une **partie constructible** $E_\phi$ de l'ensemble $\mathbb{A}^n_{\mathbb{F}_p}$ des idéaux premiers de $\mathbb{F}_p[X_1, \dots, X_n]$ de la manière évidente : c'est-à-dire que la partie associée à $f(X_1, \dots, X_n) = 0$ est le fermé de Zariski ayant cette équation, et les parties associées à des combinaisons par les connecteurs $\wedge, \vee, \neg$ sont obtenues respectivement par intersection, union et complémentaire des parties associées aux formules connectées. De même, si $k$ est un corps d'exposant caractéristique $p$, on associe à $\phi$ une partie (dite constructible) $E_\phi(k)$ de $k^n$ définie de façon analogue : c'est l'ensemble des $n$-uplets d'éléments de $k$ qui vérifient la formule $\phi$.

Par analogie avec le cas de $x \in k^n$ pour $k$ un corps d'exposant caractéristique $p$ (où il y a par définition équivalence entre $x \in E_\phi(k)$ et $\phi(x)$), si $\mathfrak{x} \in \mathbb{A}^n_{\mathbb{F}_p}$ et si $\phi$ est une formule booléenne, on écrira parfois encore abusivement que $\mathfrak{x}$ vérifie $\phi$, ou que $\phi(\mathfrak{x})$ est vraie, pour dire $\mathfrak{x} \in E_\phi$.

On aura parfois tendance à identifier abusivement la partie constructible $E_\phi$ de $\mathbb{A}^n_{\mathbb{F}_p}$ et le foncteur de même nom qui à un corps $k$ d'exposant caractéristique $p$ (éventuellement supposé algébriquement clos) associe la partie $E_\phi(k)$ de $k^n$ ; cette confusion est bénigne en raison de la proposition suivante :

**18.2. Proposition.** *Soient $\varphi, \psi$ deux formules booléennes en les indéterminées $X_1, \dots, X_n$ sur $\mathbb{F}_p$, soient $E$ et $F$ les parties constructibles de $\mathbb{A}^n_{\mathbb{F}_p}$ qu'elles définissent, et pour tout corps $k$ d'exposant caractéristique $p$ soient $E(k)$ et $F(k)$ les parties constructibles de $k^n$ qu'elles définissent. Il y a équivalence entre :*

- (i) *$E \subseteq F$ en tant que parties de $\mathbb{A}^n_{\mathbb{F}_p}$ (i.e., ensembles d'idéaux premiers de $\mathbb{F}_p[X_1, \dots, X_n]$),*
- (ii) *$E(k) \subseteq F(k)$ pour tout corps $k$ d'exposant caractéristique $p$,*
- (ii′) *$E(k) \subseteq F(k)$ pour tout corps $k$ algébriquement clos d'exposant caractéristique $p$,*
- (iii) *l'énoncé $\forall x_1 \cdots \forall x_n(\varphi \Rightarrow \psi)$ est un théorème de la théorie des corps d'exposant caractéristique $p$,*
- (iii′) *l'énoncé $\forall x_1 \cdots \forall x_n(\varphi \Rightarrow \psi)$ est un théorème de la théorie des corps algébriquement clos d'exposant caractéristique $p$.*

[1] Même si on s'intéresse finalement surtout aux fonctions calculables totales, ici définies en § 22, toute définition raisonnable de calculabilité doit nécessairement considérer des fonctions calculables partielles, à cause de l'indécidabilité du problème de l'arrêt : cf. [Odifreddi 1989, prop. II.2.1].

*On a, bien sûr, les mêmes équivalences en remplaçant partout «$\subseteq$» par «$=$» et «$\Rightarrow$» par «$\Leftrightarrow$».*

Pour le sens de l'expression « être un théorème de la théorie $\mathscr{T}$ », et la définition des théories considérées, cf. par exemple [Fried et Jarden 2008, § 8.1, exemple 7.3.1 et chapitre 9]. L'équivalence entre (ii) et (iii) (resp. entre (ii′) et (iii′)) découle du théorème de complétude de Gödel ([Marker 2002, théorème 2.1.2] ou [Fried et Jarden 2008, corollaire 8.2.6]) ; les autres équivalences sont de démonstration immédiate.

Il sera souvent commode d'identifier une formule booléenne en $n$ variables en la formule équivalente en $n' > n$ variables (ne faisant pas intervenir les nouvelles variables).

**18.3**. Les parties constructibles de $\mathbb{A}^n_{\mathbb{F}_p}$ sont la base d'ouverts d'une topologie sur $\mathbb{A}^n_{\mathbb{F}_p}$ plus fine que la topologie de Zariski (puisque tous les fermés de Zariski sont à la fois ouverts et fermés pour cette nouvelle topologie). Celle-ci s'appelle la **topologie constructible** (ou **topologie de Stone**) sur $\mathbb{A}^n_{\mathbb{F}_p}$, qui est alors dans la terminologie des logiciens l'**espace des types** (en $n$ variables, et à coefficients dans $\mathbb{F}_p$) de la théorie des corps algébriquement clos d'exposant caractéristique $p$, et cette topologie est *compacte* et *totalement discontinue* : cf. [ÉGA IV$_1$, 1.9.15] et [Marker 2002, lemme 4.1.8 et exemple 4.1.14]. Plus généralement, on peut associer à tout schéma affine $X = \mathrm{Spec}(A)$, un schéma $X^{\mathrm{cons}}$, naturellement isomorphe au spectre de l'anneau absolument plat universel $T(A)$ associé, dont l'espace topologique sous-jacent est l'ensemble $X$ muni de la topologie constructible (cf. [Olivier 1966-1968, prop. 5]). Pour chaque $x \in X$, l'anneau local $\mathscr{O}_{X^{\mathrm{cons}},x}$ est naturellement isomorphe au *corps* résiduel de $x$.

**18.4**. **Remarques.** Il importe aussi de souligner que la notion de multiplicité n'existe pas ici : les formules booléennes $f = 0$ et $f^2 = 0$ sont toujours équivalentes en ce qu'elles définissent la même partie constructible. Un fermé de Zariski — sans notion de « multiplicité » ni structure de schéma — équivaut à une partie constructible définie par une formule booléenne *positive* (c'est-à-dire sans utiliser la négation logique) : en effet, quitte à utiliser la distributivité de $\vee$ sur $\wedge$, une telle formule peut s'écrire comme une conjonction de disjonctions ($\bigwedge_i \bigvee_j f_{ij} = 0$, dite aussi « forme normale conjonctive »), or la disjonction $f_{i,1} = 0 \vee \cdots \vee f_{i,s} = 0$ équivaut à $f_i \coloneqq f_{i,1} \cdots f_{i,s} = 0$ de sorte que la conjonction de ces clauses définit un fermé de Zariski $f_1 = \cdots = f_r = 0$. (La réciproque est claire : le fermé $V(f_1, \ldots, f_r)$ est défini par $\bigvee_i f_i = 0$.)

On peut reformuler le fait qu'une partie constructible de $\mathbb{A}^n_{\mathbb{F}_p}$ est réunion finie de parties *localement fermées* (pour la topologie de Zariski) ainsi : toute formule booléenne comme disjonction de formules *purement conjonctives*, c'est-à-dire de la forme $f_1 = 0 \wedge \cdots \wedge f_r = 0 \wedge g_1 \neq 0 \wedge \cdots \wedge g_s \neq 0$ ; on peut d'ailleurs supposer $s = 1$ quitte à remplacer la fin de la formule par $g_1 \cdots g_s \neq 0$. On pourra cette fois parler de « forme normale disjonctive », et on peut réécrire une formule booléenne sous cette forme en utilisant la distributivité de $\wedge$ sur $\vee$ (ainsi que les lois de De Morgan pour traiter la négation).

**18.5**. **Observation.** Les formules booléennes sur $\mathbb{F}_p$ peuvent être manipulées par une machine de Turing.

Il est algorithmique de décider si la partie constructible associée à une formule booléenne est vide : en effet, comme expliqué en **18.4**, on peut se ramener à une formule purement conjonctive, c'est-à-dire de la forme $f_1 = 0 \wedge \cdots \wedge f_r = 0 \wedge g_1 \neq 0 \wedge \cdots \wedge g_s \neq 0$, or tester si une telle formule est universellement fausse revient, d'après le Nullstellensatz, à tester si l'un des $g_j$ (ou 1 dans le cas $s = 0$) appartient au radical de l'idéal engendré par les $f_i$. Ou encore que l'idéal engendré par les $f_i$ et les $1 - yg_j$ contienne 1, où $y$ est une nouvelle variable. Ce fait est décidable, par exemple en utilisant les algorithmes de bases de Gröbner rappelés en **13.1**.

Par conséquent, il est décidable de tester l'inclusion ou l'égalité entre deux telles parties constructibles et on peut considérer que les machines de Turing peuvent manipuler directement les parties constructibles (à coefficients dans $\mathbb{F}_p$).

## 19. Fonctions F-rationnelles

**Motivation.** Même si le modèle de calcul élaboré ici s'intéresse surtout aux fonctions calculables à valeurs entières, qu'on va définir en **20.1**, il est utile (et même indispensable pour pouvoir formuler l'approche « boîtes noires » en 21) de disposer d'une notion de fonctions calculables à valeurs dans $\mathbb{A}^1_{\mathbb{F}_p}$ (c'est-à-dire dans « le corps non spécifié sur lequel on travaille ») : pour cette notion (qui sera définie en 20.4), il est très commode (même si ce n'est pas absolument indispensable) d'admettre l'inverse du Frobenius comme opération calculable — ceci est justifié notamment par le souhait d'obtenir le lien simple entre fonctions à valeurs entières et fonctions à valeurs dans $\mathbb{A}^1_{\mathbb{F}_p}$ qui sera explicité en 24.3. C'est ce qui explique l'apparition de $\mathrm{Frob}_p^{-1}$ dans cette section et ailleurs.

**19.1. Définition.** Soit $p$ un nombre premier ou bien 1. Si $k$ est un corps d'exposant caractéristique $p$, on note $\mathrm{Frob}_p : k \to k$ le morphisme $x \mapsto x^p$. Si $k$ est parfait (et notamment si $k$ est algébriquement clos), on note $\mathrm{Frob}_p^{-1}$ le morphisme réciproque de $\mathrm{Frob}_p$.

Remarquons notamment que $\mathrm{Frob}_p : \mathbb{F}_p(X_1, \dots, X_n) \to \mathbb{F}_p(X_1, \dots, X_n)$ (où $\mathbb{F}_p(X_1, \dots, X_n)$ désigne le corps des fractions rationnelles en $n$ variables sur le corps premier $\mathbb{F}_p$) peut se voir soit comme l'élévation à la puissance $p$ soit comme la composition par la substitution $(X_1, \dots, X_n) \mapsto (X_1^p, \dots, X_n^p)$.

On appellera **fonction F-rationnelle** (sous-entendu : à coefficients dans $\mathbb{F}_p$) un élément de la clôture parfaite de $\mathbb{F}_p(X_1, \dots, X_n)$, qui est aussi la colimite du système inductif indexé par $\mathbb{N}$ dont le $i$-ième terme $K_i$ est $\mathbb{F}_p(X_1, \dots, X_n)$, et les flèches sont les Frobenius itérés $(\mathrm{Frob}_p)^{j-i} : K_i \to K_j$ si $i \leq j$.

Autrement dit, une fonction F-rationnelle est simplement une expression de la forme $\mathrm{Frob}_p^{-i} \circ f$ où $f \in \mathbb{F}_p(X_1, \dots, X_n)$, avec les identifications évidentes, ou si on préfère une fraction rationnelle en les $X_1^{1/p^i}, \dots, X_n^{1/p^i}$ pour un certain $i$. La *valeur* d'une telle fonction F-rationnelle en $(x_1, \dots, x_n) \in k^n$, si $k$ est un corps parfait d'exposant caractéristique $p$, est définie comme $\mathrm{Frob}_p^{-i}(f(x_1, \dots, x_n))$ (c'est-à-dire qu'elle est définie si et seulement si $f(x_1, \dots, x_n)$ l'est, et lorsque c'est le cas elle a pour valeur $\mathrm{Frob}_p^{-i}(f(x_1, \dots, x_n))$).

(Dans [Marker 2002, définition 3.2.13], ces fonctions sont appelées « quasirationnelles ».)

Comme pour les formules booléennes, il sera souvent commode d'identifier une fonction F-rationnelle en $n$ variables $X_1, \dots, X_n$ en la fonction F-rationnelle de même expression en $n' > n$ variables (ne faisant pas intervenir les nouvelles variables).

**19.2. Remarque.** Le *graphe* $\Gamma_u$ d'une fonction F-rationnelle $u = \mathrm{Frob}_p^{-i} \circ f$ est la partie constructible (en fait, localement fermée) de $\mathbb{A}^n_{\mathbb{F}_p} \times \mathbb{A}^1_{\mathbb{F}_p}$ définie par l'équation $y^{p^i} = f(x_1, \dots, x_n)$ (c'est-à-dire $vy^{p^i} - u = 0$ et $v \neq 0$ où $f = u/v$ est la forme irréductible de la fonction rationnelle $f$ ; ici, $x_1, \dots, x_n$ sont les coordonnées sur $\mathbb{A}^n_{\mathbb{F}_p}$ et $y$ celle sur $\mathbb{A}^1_{\mathbb{F}_p}$).

Il est clair qu'une fonction F-rationnelle est complètement déterminée par son graphe (quitte à appliquer une puissance de $\mathrm{Frob}_p$ on se ramène au cas des fractions rationnelles).

Remarquons d'autre part (cf. [ÉGA IV$_2$, 2.4.5]) que la projection $\Gamma_u \to D_u$ où $D_u \subseteq \mathbb{A}^n_{\mathbb{F}_p}$ est l'ouvert de définition de $u$ (ou de façon équivalente, de $f$) est un homéomorphisme universel de schémas (pour la topologie de Zariski) : ceci permet de voir une fonction F-rationnelle $u$ comme une application ensembliste $D_u \to \mathbb{A}^1_{\mathbb{F}_p}$, donc une fonction *partielle* $\mathbb{A}^n_{\mathbb{F}_p} \to \mathbb{A}^1_{\mathbb{F}_p}$.

Si $E$ est une partie constructible de $\mathbb{A}^n_{\mathbb{F}_p}$, on appellera « fonction F-rationnelle sur $E$ » la restriction à $E$ d'une fonction F-rationnelle (ceci n'impose pas que la fonction soit définie sur tout $E$).

La proposition suivante a pour objet d'aider l'intuition du lecteur. Nous ne faisons qu'en esquisser la démonstration puisque nous ne l'utiliserons pas et que nous prouverons un résultat plus général (plus effectif) en 24.2.

**19.3. Proposition.** *Soit $Z$ une partie localement fermée de $\mathbb{A}^{n+1}_{\mathbb{F}_p}$ (munie de la structure de schéma réduit induite) telle que la projection $\pi : Z \to \mathbb{A}^n_{\mathbb{F}_p}$ soit radicielle (c'est-à-dire que pour tout corps $k$ d'exposant caractéristique $p$ et tout $x \in \mathbb{A}^n_{\mathbb{F}_p}(k)$ il existe au plus un $z \in Z(k)$ se projetant en $x$). Alors on peut écrire $\pi(Z)$ comme réunion finie disjointe de parties localement fermées $X_1, \dots, X_n$ telles que sur chacune d'entre elles la restriction $Z_i = \pi^{-1}(X_i)$ soit le graphe d'une fonction F-rationnelle sur $X_i$.*

*Esquisse de démonstration.* Le morphisme $\pi : Z \to \mathbb{A}^n_{\mathbb{F}_p}$ est quasi-fini : quitte à écrire son image comme réunion finie disjointe de localement fermés $X_i$, qu'on peut même supposer affine, on peut s'arranger pour que $\pi$ soit fini (entier et de type fini) sur chaque $X_i$. On peut supposer qu'il n'y a qu'un seul $X_i$, notons-le $X$. Mettons que $X = \mathrm{Spec}\, A$ et $Z = \mathrm{Spec}\, B$, de sorte que $B$ est une $A$-algèbre de type fini entière ; l'immersion fermée $Z \hookrightarrow X \times \mathbb{A}^1_{\mathbb{F}_p}$ définit une surjection de $A$-algèbres $A[t] \twoheadrightarrow B$, donc $B$ est engendré par un unique élément $y$ (image de $t$) qui vérifie une équation radicielle, $y^{p^i} = x$ pour $x \in A$ qu'on peut représenter par une fonction rationnelle. □

(On pourra comparer avec [Marker 2002, proposition 3.2.14].)

**19.4. Observation.** Il est clair que les fonctions F-rationnelles peuvent être manipulées par une machine de Turing (on peut par exemple le déduire de 12.5(i,v)) : il est algorithmique de calculer des sommes, différences, produits et quotients dans ce corps, ou de tester la nullité d'un élément, ou encore de calculer des composées (et notamment de tester si une composée est bien définie). On se convainc de même facilement qu'on peut par exemple calculer algorithmiquement le graphe d'une fonction F-rationnelle (comme partie constructible) ou tester l'égalité de deux fonctions F-rationnelles sur une partie constructible.

## 20. Fonctions partielles calculables

Nous cherchons maintenant à définir la notion de fonction partielle $\mathbb{N}^r \times \mathbb{A}^n_{\mathbb{F}_p} \to \mathbb{N}$ ou $\mathbb{N}^r \times \mathbb{A}^n_{\mathbb{F}_p} \to \mathbb{A}^1_{\mathbb{F}_p}$ *calculable* au sens de Church-Turing. L'intuition de la définition qui suit est que, pour évaluer une telle fonction $f$ en des arguments $(i_1, \ldots, i_r, x_1, \ldots, x_n) \in \mathbb{N}^r \times k^n$ (où $k$ est un corps quelconque d'exposant caractéristique $p$), une machine de Turing va produire une suite de fonctions booléennes $\varphi_j$ (indexées par les entiers naturels $j$, et pouvant dépendre des paramètres entiers $i_1, \ldots, i_r$) et de valeurs $v_j$, et la valeur de la fonction en $x_1, \ldots, x_n$ est le $v_j$ accompagnant la première formule $\varphi_j$ qui soit vraie en $x_1, \ldots, x_n$. Formellement :

**20.1. Définition.** Soit $p$ un nombre premier ou bien 1. On appelle **fonction partielle** $f : \mathbb{N}^r \times \mathbb{A}^n_{\mathbb{F}_p} \to \mathbb{N}$ **calculable au sens de Church-Turing** la donnée d'une fonction partielle calculable au sens de Church-Turing (au sens usuel[①]) qui à un $(r+1)$-uplet d'entiers naturels $(\underline{i}, j) = (i_1, \ldots, i_r, j)$ associe une formule booléenne $\varphi_{\underline{i},j}$ sur $\mathbb{F}_p$ en des indéterminées $X_1, \ldots, X_n$ ainsi qu'un entier naturel $v_{\underline{i},j}$.

Si $(x_1, \ldots, x_n)$ sont des éléments d'un corps $k$ quelconque d'exposant caractéristique $p$, ou bien une notation abusive pour un point de $\mathbb{A}^n_{\mathbb{F}_p}$, et si $(i_1, \ldots, i_r)$ sont des entiers naturels, alors la *valeur* de $f$ en $(i_1, \ldots, i_r, x_1, \ldots, x_n)$ est définie comme valant $v_{\underline{i},j}$, où $j$ est le plus petit entier naturel tel que

- la formule $\varphi_{\underline{i},j}$ est définie, et $\varphi_{\underline{i},j}(x_1, \ldots, x_n)$ est vraie, et
- pour tout $j' < j$, la formule $\varphi_{\underline{i},j'}$ est définie, et $\varphi_{\underline{i},j'}(x_1, \ldots, x_n)$ n'est pas vraie.

Si un tel $j$ n'existe pas, alors la valeur de $f$ n'est pas définie en ces points.

(Pour alléger la notation, on omettra souvent les indices $i_1, \ldots, i_r$ sur $\varphi$ et $v$, notant ainsi simplement $\varphi_j$ et $v_j$.)

On définit aussi la notion de fonction partielle $f : \mathbb{N}^r \times (\coprod_{n=0}^{+\infty} \mathbb{A}^n_{\mathbb{F}_p}) \to \mathbb{N}$ calculable au sens de Church-Turing : il s'agit d'une fonction partielle calculable qui à un $(r+2)$-uplet $(i_1, \ldots, i_r, n, j)$ d'entiers naturels associe un couple $(\varphi_j, v_j)$, la valeur de $f$ en $(i_1, \ldots, i_r, (x_1, \ldots, x_n))$ étant définie de la même manière.

Il faut comprendre cette définition de la manière suivante : pour définir une fonction partielle calculable $\mathbb{N}^r \times \mathbb{A}^n_{\mathbb{F}_p} \to \mathbb{N}$, on donne une machine de Turing qui produit, pour chaque valeur des paramètres entiers $(i_1, \ldots, i_r)$, une suite récursive de couples $(\varphi_j, v_j)$ qui signifient « si les paramètres $(x_1, \ldots, x_n)$ vérifient $\varphi_j$ (et aucune formule précédente), alors retourner la valeur $v_j$ ». (Si pour un certain $j$ la machine de Turing ne termine pas, ou bien sûr si elle renvoie autre chose qu'un couple formé d'une formule booléenne et d'un entier naturel, la valeur de la fonction est considérée comme indéfinie, sauf pour les valeurs qui ont déjà été attrapées par les $\varphi_{j'}$ précédents. De même, si aucune formule n'est vérifiée, la fonction n'est pas définie au point considéré.) Pour une fonction $f : \mathbb{N}^r \times (\coprod_{n=0}^{+\infty} \mathbb{A}^n_{\mathbb{F}_p}) \to \mathbb{N}$, la dimension $n$ doit également être passée en paramètre à la machine de Turing.

À titre d'exemple trivial (avec $r = 0$ et $n = 1$), la machine de Turing qui à $j = 0$ associe le couple (“$x = 0$”, 0) et à $j = 1$ le couple (“$\neg(x = 0)$”, 1) définit une fonction $\mathbb{A}^1_{\mathbb{F}_p} \to \mathbb{N}$ qui, sur n'importe

[①] C'est-à-dire définie, par exemple, par une machine de Turing, avec la convention que lorsque l'exécution de la machine de Turing sur l'entrée fournie ne termine pas, la fonction n'est pas définie : [Odifreddi 1989, définition II.1.1, théorème II.1.8 et discussion p. 132–133].

quel corps, vaut (l'entier naturel) 0 en (l'élément du corps) 0 et 1 ailleurs. On verra en 22.1 que cette fonction est « totale » : dans le cas des fonctions totales, la définition des fonctions calculables se simplifie considérablement, comme on le verra en 22.4.

**20.2. Remarque.** Pour être plus rigoureux, la définition faite en 20.1 est celle de la *description* d'une fonction partielle calculable $\mathbb{N}^r \times \mathbb{A}^n_{\mathbb{F}_p} \to \mathbb{N}$ : on considérera que deux fonctions calculables $\mathbb{N}^r \times \mathbb{A}^n_{\mathbb{F}_p} \to \mathbb{N}$ sont égales lorsque leurs valeurs sont les mêmes en tout $(i_1, \ldots, i_r, x_1, \ldots, x_n)$ pour n'importe quel corps $k$ d'exposant caractéristique $p$ — d'après 18.2 il revient au même de faire cette hypothèse pour $k$ algébriquement clos, ou encore pour tout $(i_1, \ldots, i_r, x) \in \mathbb{N}^r \times \mathbb{A}^n_{\mathbb{F}_p}$. On peut donc bien considérer les fonctions partielles calculables $\mathbb{N}^r \times \mathbb{A}^n_{\mathbb{F}_p} \to \mathbb{N}$ comme des fonctions partielles $\mathbb{N}^r \times \mathbb{A}^n_{\mathbb{F}_p} \to \mathbb{N}$ particulières.

**20.3. Remarque.** Dans la définition 20.1, on peut, si on le souhaite, supposer que les $\varphi_j$ sont deux à deux exclusives (c'est-à-dire que les parties constructibles qu'elles définissent sont deux à deux disjointes) : en effet, il suffit pour cela de remplacer $\varphi_j$ par sa conjonction avec la négation de tous les $\varphi_{j'}$ antérieurs. (Cette transformation est évidemment algorithmique.)

De même, si on le souhaite, et de nouveau par une transformation algorithmique, on peut supposer que les $\varphi_j$ sont des formules purement conjonctives (donc définissent des parties localement fermées, cf. 18.4). On se convainc même qu'on peut faire ces deux hypothèses simultanément (car toute formule booléenne peut se réécrire — algorithmiquement — comme disjonction *exclusive* de formules purement conjonctives, quitte à distinguer tous les cas possibles d'égalité ou de non-égalité à 0 des polynômes qui y interviennent).

La notion de fonction partielle calculable à valeurs dans $\mathbb{A}^1_{\mathbb{F}_p}$ est définie de façon analogue à celle de fonction partielle calculable à valeurs dans $\mathbb{N}$, mais cette fois les valeurs $v_j$ retournées par la machine de Turing sont des fonctions F-rationnelles qui sont évaluées en $(x_1, \ldots, x_n)$. Précisément :

**20.4. Définition.** On appelle **fonction partielle** $f : \mathbb{N}^r \times \mathbb{A}^n_{\mathbb{F}_p} \to \mathbb{A}^1_{\mathbb{F}_p}$ **calculable au sens de Church-Turing** une fonction comme en 20.1 mais où cette fois les $v_j$ sont censés être des fonctions F-rationnelles (cf. 19.1) en $n$ variables $X_1, \ldots, X_n$ ; la valeur de la fonction $f$ en $(i_1, \ldots, i_r, x_1, \ldots, x_n)$, où $x_1, \ldots, x_n$ appartiennent à un corps parfait $k$ d'exposant caractéristique $p$, est alors définie comme la valeur en $(x_1, \ldots, x_n)$ de la fonction F-rationnelle $v_j$ (et indéfinie si $v_j$ n'est pas définie en $(x_1, \ldots, x_n)$).

On définit de même les fonctions partielles calculables $f : \mathbb{N}^r \times (\coprod_{n=0}^{+\infty} \mathbb{A}^n_{\mathbb{F}_p}) \to \mathbb{A}^1_{\mathbb{F}_p}$. Et on définit sans difficulté les fonctions partielles calculables à valeurs dans $\mathbb{N}^s \times \mathbb{A}^m_{\mathbb{F}_p}$ (ou même vers $\coprod_{m=0}^{+\infty} \mathbb{A}^m_{\mathbb{F}_p}$, en convenant que pour spécifier un élément de celui-ci on spécifie un $m$ et une valeur dans $\mathbb{A}^m_{\mathbb{F}_p}$, c'est-à-dire un $m$-uplet de fonctions F-rationnelles).

En particulier, les fonctions F-rationnelles elles-mêmes sont des fonctions calculables.

Comme on l'a convenu en 20.2 pour les fonctions partielles calculables $\mathbb{N}^r \times \mathbb{A}^n_{\mathbb{F}_p} \to \mathbb{N}$, les fonctions partielles calculables $\mathbb{N}^r \times \mathbb{A}^n_{\mathbb{F}_p} \to \mathbb{A}^1_{\mathbb{F}_p}$ sont définies par leurs *valeurs* (on remarquera que c'est le cas pour les fonctions F-rationnelles : cf. 19.2).

**20.5. Remarque.** Dans le contexte de la définition 20.4, il pourra parfois être commode de supposer que la fonction F-rationnelle $v_j$ est nécessairement définie quand on cherche à l'évaluer (autrement dit, que la non-définition de la fonction $f$ ne peut être due qu'aux formules $\varphi_{j'}$) : ceci est très simple à assurer, il suffit de remplacer $\varphi_j$ par sa conjonction avec la formule qui assure que le dénominateur de $v_j$ ne s'annule pas.

## 21. Équivalence avec le modèle « boîte noire »

Le modèle de calcul que nous venons de définir a l'avantage qu'il définit automatiquement des fonctions valables à la fois sur tous les corps d'exposant caractéristique $p$ (remarquons au passage que si la fonction dispose de $n = 0$ variables d'entrées dans le corps, alors la seule sortie qu'elle puisse produire est un élément de $\mathbb{F}_p$).

Il n'est cependant peut-être pas évident que cette notion de calcul soit naturelle et décrive bien toutes les manipulations qu'on souhaite pouvoir effectuer algorithmiquement sur les éléments d'un corps. Pour répondre à cette crainte, nous allons à présent montrer que les fonctions partielles calculables dans le modèle qu'on vient de définir coïncident précisément avec celles qui sont calculables par une machine de Turing qui manipulerait les éléments d'un corps non spécifié comme des « boîtes noires » sur lesquelles elle peut effectuer uniquement les opérations du corps (addition, opposé, multiplication, inverse, comparaison à zéro).

Montrons d'abord un lemme qui énonce le fait que dans les définitions **20.1** et **20.4**, il revient au même de supposer que la machine de Turing qui produit la valeur finale $v_j$ de la fonction peut poser des questions sur la satisfaction de formules booléennes des variables, y compris en fonction des réponses aux questions précédemment posées (c'est-à-dire, sur la base d'un arbre de décision) :

**21.1. Lemme.** *On obtient des définitions équivalentes à **20.1** et **20.4** si la machine de Turing considérée prend en entrée un multi-indice $\underline{i} = (i_1, \ldots, i_r)$ et un mot fini $\gamma$ sur l'alphabet $\{\top, \bot\}$ et renvoie* soit *une formule booléenne $\psi_\gamma$* soit *une valeur finale $w_\gamma$ (qui est un entier naturel ou bien une fonction F-rationnelle selon le type censé être retourné par $f$), la valeur de $f$ en $(i_1, \ldots, i_r, x_1, \ldots, x_n)$ étant alors définie comme $w_\gamma$ où $\gamma$ est l'unique mot tel que*

- *la machine renvoie pour $\gamma$ une valeur finale (qui est $w_\gamma$), et*
- *pour tout préfixe strict $\gamma'$ de $\gamma$, la machine renvoie une formule $\psi_{\gamma'}$, et le symbole suivant immédiatement $\gamma'$ dans $\gamma$ est $\top$ ou $\bot$ selon que $\psi_{\gamma'}(x_1, \ldots, x_n)$ est vraie ou fausse.*

(La donnée de la machine qui à $\gamma$ associe une formule ou une valeur finale comme ci-dessus s'appellera un **arbre de décision** pour la fonction qu'elle calcule. Plus exactement, le mot vide s'appelle la **racine** de l'arbre de décision, à chaque fois que la machine renvoie une valeur $w_\gamma$, le mot $\gamma$ est appelé une **feuille** de l'arbre de décision, tandis que si elle renvoie une formule $\psi_\gamma$, on dira que $\gamma$ est un **nœud** ayant pour fils les deux mots obtenus en concaténant respectivement $\top$ et $\bot$ à la fin de $\gamma$.)

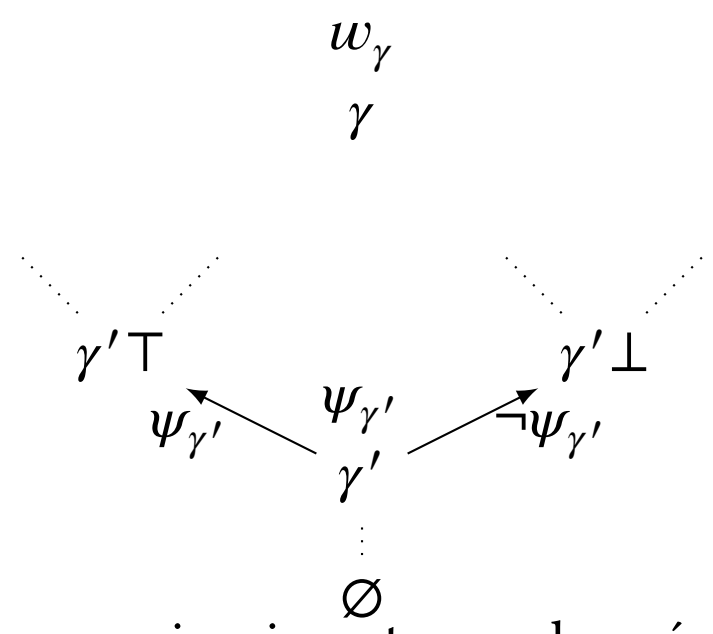

*Démonstration.* Il est clair que cette nouvelle définition est au moins aussi puissante que la précédente : donnée une suite $(\varphi_j, v_j)$ calculée comme en **20.1** ou **20.4**, on définit une machine $M$ qui calcule la même fonction pour la définition du présent lemme : si $\gamma$ est un mot formé de $j$ signes $\bot$ suivi d'un signe $\top$ alors $M$ renvoie pour $\gamma$ une valeur $w_\gamma = v_j$, et si $\gamma$ est formé de $j$ signes $\bot$ alors $M$ renvoie pour $\gamma$ une formule $\psi_\gamma = \varphi_j$ (et pour tout mot $\gamma$ qui n'est pas formé d'une suite de $\bot$ puis d'un éventuel signe $\top$, ce que renvoie $M$ est sans importance).

Réciproquement, donnons-nous une machine $M$ qui pour tout mot $\gamma$ renvoie soit une formule $\psi_\gamma$ soit une valeur $w_\gamma$ comme dans l'énoncé du présent lemme. Il s'agit de définir une nouvelle machine de Turing $M'$ qui énumère une suite $(\varphi_j, v_j)$ en fonction de ce que renvoie la machine $M$ calculant les $\psi_\gamma$ ou $w_\gamma$. Pour cela, on va effectuer un parcours en parallèle de l'arbre de tous les calculs possibles par $M$. Précisément, la machine $M'$ commence par exécuter $M$ sur le mot vide, puis, si celle-ci termine et renvoie une formule booléenne, exécute $M$ *en parallèle* sur les deux mots $\top$ et $\bot$, et ainsi de suite récursivement : à chaque fois qu'une exécution de $M$ sur un mot $\gamma$ termine en renvoyant une formule booléenne, la machine $M'$ exécute en parallèle $M$ sur les deux mots obtenus en ajoutant à $\gamma$ un symbole $\top$ ou $\bot$ à la fin ; en revanche, si la machine $M$ termine en produisant une valeur finale $w_\gamma$, alors la machine $M$ énumère pour prochaine valeur $j$ (c'est-à-dire, le premier $j$ qui n'ait pas encore été utilisé) le couple $(\varphi_j, v_j)$ où $v_j = w_\gamma$ et $\varphi_j$ est la conjonction de toutes les formules $\psi_{w'}$ qui ont été calculées par $M$ sur les préfixes stricts $w'$ de $w$. □

**21.2. Proposition.** *Les fonctions partielles calculables $\mathbb{N}^r \times \mathbb{A}^n_{\mathbb{F}_p} \to \mathbb{N}$ et $\mathbb{N}^r \times \mathbb{A}^n_{\mathbb{F}_p} \to \mathbb{A}^1_{\mathbb{F}_p}$ au sens des définitions **20.1** et **20.4** :*

- *contiennent les fonctions constantes (à valeurs dans les entiers ou les éléments de $\mathbb{F}_p$);*
- *contiennent la fonction successeur $i \mapsto i + 1$ sur les entiers naturels, la fonction de test à zéro (qui vaut l'entier naturel $0$ en $0 \in \mathbb{A}^1_{\mathbb{F}_p}$ et $1$ ailleurs), ainsi que les opérations d'addition ($(x, x') \mapsto x + x'$), opposée ($x \mapsto -x$), multiplication ($(x, x') \mapsto xx'$) et inverse ($x \mapsto x^{-1}$ non définie en $0$), et la fonction $\mathrm{Frob}_p^{-1}$ ;*
- *contiennent les projections vers n'importe quelle variable d'entrée;*

- *sont stables par composition : si $h$ et $g_1, \dots, g_N$ sont calculables et que la composition a un sens, alors $h(g_1, \dots, g_N)$ est calculable;*
- *sont stables par récursion primitive : si $g$ et $h$ sont calculables, le premier argument de $h$ étant un entier naturel, alors la fonction $f$, dont le premier argument est un entier naturel, définie par $f(0, \dots) = g(\dots)$ et $f(u+1, \dots) = h(u, f(u, \dots), \dots)$, est calculable — ici, les points de suspension désignent n'importe quelle combinaison de sortes de variables;*
- *sont stables par recherche non bornée (ou « opérateur $\mu$ de Kleene ») : si $f$ est calculable, son premier argument étant un entier naturel, alors la fonction qui à des paramètres d'entrée $(\dots)$ associe le plus petit $u$ tel que $f(u, \dots) = 0$ et que $f(u', \dots) \neq 0$ pour tout $u' < u$ (ce qui sous-entend que toutes ces valeurs $f(u', \dots)$ sont définies), est calculable.*

(Comparer avec [Odifreddi 1989, définition II.1.1]. Voir la proposition suivante pour la réciproque.)

*Démonstration.* Les trois premiers points sont évidents.

Les suivants découlent du lemme précédent, car il est facile de fabriquer les arbres de décision des fonctions décrites d'après ceux des fonctions qui les définissent.

Par exemple, pour la composée $h \circ g$, on fabrique l'arbre de décision de $h \circ g$ en attachant à la chaque feuille de l'arbre de décision de $g$ (c'est-à-dire, dans la notation du lemme, à chaque fois qu'au mot $\gamma$ il affecte une valeur finale $w_\gamma$ et non une formule) une copie de l'arbre de décision de $h$, dans lequel les formules et éventuelles fonctions booléennes ont été composées par $w_\gamma$ (ce qui signifie, si $g$ renvoie des entiers, qu'on a calculé l'arbre de décision de $h$ pour cet entier $w_\gamma$, et si $g = (g_1, \dots, g_N)$ prend ses valeurs dans $\mathbb{A}^N_{\mathbb{F}_p}$, qu'on a composé toutes les formules booléennes, et éventuellement valeurs de retour si ce sont des fonctions F-rationnelles, dans cet arbre, par la fonction $w_\gamma$). □

**21.3. Proposition.** *Les fonctions partielles calculables $\mathbb{N}^r \times \mathbb{A}^n_{\mathbb{F}_p} \to \mathbb{N}$ et $\mathbb{N}^r \times \mathbb{A}^n_{\mathbb{F}_p} \to \mathbb{A}^1_{\mathbb{F}_p}$ au sens des définitions **20.1** et **20.4** sont le plus petit ensemble de telles fonctions qui possèdent les stabilités énumérées à la proposition précédente.*

Ces fonctions pourraient donc également être appelées **fonctions calculables à partir des opérations de corps et $\mathrm{Frob}_p^{-1}$ comme boîtes noires**.

*Démonstration.* On vient de voir une implication. Il s'agit maintenant de montrer que les fonctions calculables au sens des définitions **20.1** et **20.4** sont effectivement calculables à partir des opérations de corps et $\mathrm{Frob}_p^{-1}$ comme boîtes noires.

Mais les fonctions calculables en ce dernier sens contiennent au moins les fonctions calculables au sens de Church-Turing, et par ailleurs contiennent toutes les fonctions F-rationnelles et les formules booléennes (considérées comme fonctions à valeurs dans $\{0, 1\}$). Il reste donc simplement à montrer qu'à partir d'une suite récursive $(\varphi_j, v_j)$ on peut calculer la valeur finale de $f$ : mais c'est clair en utilisant la primitive récursion et la recherche (du plus petit $j$ tel que $\varphi_j$ soit vrai et tous les précédents soient faux). □

On en déduit que les fonctions partielles calculables dans notre modèle sont calculables quand on les applique à n'importe quel « corps calculable » :

**21.4. Corollaire.** *Si $\mathfrak{K}$ est un corps calculable (cf. **12.1**) d'exposant caractéristique $p$, alors pour toute fonction partielle calculable $\mathbb{N}^r \times \mathbb{A}^n_{\mathbb{F}_p} \to \mathbb{N}$ (resp. $\mathbb{N}^r \times \mathbb{A}^n_{\mathbb{F}_p} \to \mathbb{A}^1_{\mathbb{F}_p}$ en supposant $\mathfrak{K}$ parfait) au sens des définitions **20.1** ou **20.4**, la fonction partielle $\mathbb{N}^r \times \mathfrak{K}^n \to \mathbb{N}$ qu'elle définit est calculable (au sens usuel) (resp. il existe une fonction calculable $\mathbb{N}^r \times \mathfrak{K}^n \to \mathfrak{K}$ qui la représente).*

*Démonstration.* Il s'agit simplement de voir que les fonctions énumérées en **21.2** sont calculables sur $\mathfrak{K}$. Le seul cas qui n'est pas totalement évident est celui de $\mathrm{Frob}_p^{-1}$ : mais comme pour le cas de fonctions $\mathbb{N}^r \times \mathbb{A}^n_{\mathbb{F}_p} \to \mathbb{A}^1_{\mathbb{F}_p}$ on a supposé $\mathfrak{K}$ parfait, et que dans le cas de fonctions $\mathbb{N}^r \times \mathbb{A}^n_{\mathbb{F}_p} \to \mathbb{N}$ on peut construire et travailler dans la clôture parfaite de $\mathfrak{K}$ d'après le début de la démonstration de **12.5**(v)), on peut calculer $\mathrm{Frob}_p^{-1}$ (au pire[①] d'après **12.8**). □

[①] Si on voulait remplacer la notion de calculabilité par celle de fonction *primitivement* récursive, la bonne notion de corps « primitivement récursivement parfait » supposerait, bien sûr, que $\mathrm{Frob}_p^{-1}$ fait partie des opérations de corps.

**21.5. Remarque.** La réciproque *n'est pas vraie*. Par exemple, si on considère la clôture algébrique $\mathbb{F}_p(T)^{\text{alg}}$ de $\mathbb{F}_p(T)$ construite de la façon habituelle (cf. 12.5), la fonction qui à $x$ associe 0 si $x$ est algébrique sur $\mathbb{F}_p$ et 1 si $x$ est transcendant, sera calculable (ceci n'est pas surprenant : il est en effet facile, en considérant l'écriture d'une fonction algébrique sur $\mathbb{F}_p$, de voir si l'indéterminée apparaît réellement) ; en revanche, la fonction $\mathbb{A}^1_{\mathbb{F}_p} \to \mathbb{N}$ qui répond à la même définition *n'est pas* calculable au sens de la définition 20.1, comme il résultera par exemple du corollaire 22.4 (l'ensemble des éléments transcendants sur $\mathbb{F}_p$ n'étant, de toute évidence, pas une partie constructible de $\mathbb{A}^1_{\mathbb{F}_p}$). Plus généralement, la différence entre la notion de fonction calculable sur un corps calculable $\mathfrak{K}$ et celle de fonction calculable au sens « universel » défini dans la présente partie est que les premières peuvent utiliser les détails du codage des éléments de $\mathfrak{K}$ alors que les « boîtes noires » ne permettent de réaliser que des opérations algébriques.

## 22. Fonctions totales

**22.1. Définition.** On dit qu'une fonction partielle calculable $\mathbb{N}^r \times \mathbb{A}^n_{\mathbb{F}_p} \to \mathbb{N}$ ou $\mathbb{N}^r \times \mathbb{A}^n_{\mathbb{F}_p} \to \mathbb{A}^1_{\mathbb{F}_p}$ (au sens des définitions 20.1 et 20.4) est **totale** lorsque sa valeur est définie pour *toutes* valeurs $(i_1, \ldots, i_r, x_1, \ldots, x_n)$, les $x_1, \ldots, x_n$ appartenant à n'importe quel corps parfait $k$ d'exposant caractéristique $p$.

En vertu de 18.2, il revient au même de demander que la valeur soit définie pour toutes valeurs $(i_1, \ldots, i_r, x_1, \ldots, x_n)$, les $x_1, \ldots, x_n$ appartenant à n'importe quel corps *algébriquement clos* $k$ d'exposant caractéristique $p$, ou encore, en tout point de $\mathbb{N}^r \times \mathbb{A}^n_{\mathbb{F}_p}$ où $\mathbb{A}^n_{\mathbb{F}_p}$ est vu comme l'ensemble des idéaux premiers de $\mathbb{F}_p[X_1, \ldots, X_n]$.

**22.2. Contre-exemple :** La fonction qui à $x \in \mathbb{A}^1_{\mathbb{F}_p}$ associe le degré de $x$ sur $\mathbb{F}_p$, c'est-à-dire le plus petit degré d'un polynôme $f \in \mathbb{F}_p[X]$ tel que $f(x) = 0$, est partielle calculable (il suffit de produire une machine de Turing qui énumère les polynômes irréductibles $f_j$ sur $\mathbb{F}_p$, et renvoie le couple formé de la formule booléenne $f_j = 0$ et l'entier $\deg f_j$), mais même si elle est définie sur tout élément de $\mathbb{F}_p^{\text{alg}}$, elle n'est pourtant pas totale, car appliquée à l'élément transcendant $T$ de $\mathbb{F}_p(T)$ elle n'est pas définie (ou, si on préfère voir $\mathbb{A}^1_{\mathbb{F}_p}$ comme l'ensemble des idéaux premiers de $\mathbb{F}_p[T]$, la fonction n'est pas définie au point générique). Il résulte de la proposition qui suit (ou de son corollaire) qu'il n'existe aucune manière de prolonger cette fonction en une fonction calculable totale.

**22.3. Proposition.** *(A) Pour qu'une fonction partielle calculable $\mathbb{N}^r \times \mathbb{A}^n_{\mathbb{F}_p} \to \mathbb{N}$ ou $\mathbb{N}^r \times \mathbb{A}^n_{\mathbb{F}_p} \to \mathbb{A}^1_{\mathbb{F}_p}$ soit totale, il (faut et il) suffit qu'elle soit définie pour toutes les valeurs $(i_1, \ldots, i_r, x_1, \ldots, x_n)$, où les $x_1, \ldots, x_n$ appartiennent à une clôture algébrique fixée de $\mathbb{F}_p(T_1, \ldots, T_n)$.*

*(B) De plus, pour une fonction totale, dans la définition 20.1 ou 20.4, il existe un $N$ (dépendant de $i_1, \ldots, i_r$ et, évidemment, de la fonction considérée) telle qu'on puisse tronquer la suite $(\varphi_j, v_j)$ à partir de $N$ sans changer la fonction ainsi définie.*

*Démonstration.* Le point (A) est conséquence du fait que les corps résiduels des points de $\mathbb{A}^n_{\mathbb{F}_p}$ sont de degré de transcendance au plus $n$ sur $\mathbb{F}_p$. Détaillons dans le cas d'une fonction $\mathbb{N}^r \times \mathbb{A}^n_{\mathbb{F}_p} \to \mathbb{N}$ ; le cas d'une fonction $\mathbb{N}^r \times \mathbb{A}^n_{\mathbb{F}_p} \to \mathbb{A}^1_{\mathbb{F}_p}$ ne présentant pas de difficulté supplémentaire grâce à la remarque 20.5. Si la fonction n'est pas totale, il existe $(i_1, \ldots, i_r) \in \mathbb{N}^r$ et un $n$-uplet $(x_1, \ldots, x_n)$ d'éléments d'un certain corps $k$ d'exposant caractéristique $p$, qu'on peut certainement supposer algébriquement clos, tels que (1) pour un certain $j_0$, les $\varphi_0, \ldots, \varphi_{j_0-1}$ soient définis et faux (=non vérifiés) pour $(x_1, \ldots, x_n)$ mais que $\varphi_{j_0}$ ne soit pas défini (la machine censée le produire ne termine pas), ou bien (2) les $\varphi_j$ soient définis pour tout $j$ et tous faux pour $(x_1, \ldots, x_n)$. (Dans le cas inintéressant, c'est-à-dire (1), définissons comme $0 = 1$ tous les $\varphi_j$ pour $j \geq j_0$.) Par hypothèse et construction, l'intersection des parties constructibles $E_j := E_{\neg\varphi_j}$ est non vide. Le corps résiduel d'un point quelconque de cette intersection se plonge dans le corps de l'énoncé.

Montrons (B) : le fait que la fonction soit totale implique qu'il n'existe pas de $n$-uplet $(x_1, \ldots, x_n)$ d'éléments d'un corps algébriquement clos d'exposant caractéristique $p$ qui vérifieraient $\neg\varphi_j$ pour tout $j$. Autrement dit, l'intersection des parties constructibles $E_j = E_{\neg\varphi_j}$ est *vide*. L'espace topologique $\mathbb{A}^{n\,\text{cons}}_{\mathbb{F}_p}$ étant compact et les $E_j$ fermés, il existe un $N$ tel que $\bigcap_{j<N} E_j = \varnothing$. On peut alors omettre tous les $\varphi_j$ pour $j \geq N$ dans la définition de la fonction car l'une des conditions $\varphi_j$, $j < N$ est toujours satisfaite. □

Notons, bien que nous n'en ferons pas usage que si dans (A) l'entier $n$ n'est pas fixé, on peut tester si une fonction est totale en considérant des variables à valeurs dans une clôture algébrique de $\mathbb{F}_p(T_i, i \in \mathbb{N})$[①].

**22.4. Corollaire.** *Une fonction totale calculable $\mathbb{A}^n_{\mathbb{F}_p} \to \mathbb{N}$ est exactement la donnée d'une fonction* continue $\mathbb{A}^{n\,\mathrm{cons}}_{\mathbb{F}_p} \to \mathbb{N}$*, c'est-à-dire d'une partition de $\mathbb{A}^n_{\mathbb{F}_p}$ en un nombre fini de parties constructibles, sur chacune desquelles la fonction prend une valeur constante entière (en particulier, la fonction est bornée).*

*Une fonction totale calculable $\mathbb{A}^n_{\mathbb{F}_p} \to \mathbb{A}^m_{\mathbb{F}_p}$ est exactement la donnée d'une partition de $\mathbb{A}^n_{\mathbb{F}_p}$ en un nombre fini de parties constructibles, sur chacune desquelles la fonction coïncide avec une certaine fonction F-rationnelle (définie partout sur la partie en question).*

*Une fonction totale calculable $\mathbb{N}^r \times \mathbb{A}^n_{\mathbb{F}_p} \to \mathbb{N}$ ou $\mathbb{N}^r \times \mathbb{A}^n_{\mathbb{F}_p} \to \mathbb{A}^m_{\mathbb{F}_p}$ est exactement la donnée d'une fonction calculable au sens usuel qui à toute valeur $(i_1, \ldots, i_r)$ des $r$ premiers paramètres associe une donnée comme décrite dans les deux paragraphes précédents.*

*Démonstration.* Ceci découle de manière simple du (B) de la proposition précédente. On notera que le $N$ de cet énoncé peut être calculé en fonction de $(i_1, \ldots, i_r)$ en vertu de l'observation **18.5** (car on sait tester si $E_{\varphi_1} \cup \cdots \cup E_{\varphi_N}$ recouvrent $\mathbb{A}^n_{\mathbb{F}_p}$). □

(On pourra comparer avec [Marker 2002, proposition 3.2.14]. Par ailleurs, l'uniformité affirmée dans le premier paragraphe du corollaire ci-dessus est la même que celle énoncée dans [Galligo, Gonzalez-Vega et Lombardi 2003, introduction, p. 80] — qui ne considère cependant que des entrées dans $\mathbb{A}^n$.)

**22.5. Remarque.** La proposition et son corollaire valent également (*mutatis mutandis*) pour des fonctions calculables totales définies sur $\mathbb{N}^r \times E$ où $E$ est une partie constructible de $\mathbb{A}^n_{\mathbb{F}_p}$, ou même sur une partie $E$ de $\mathbb{N}^r \times \mathbb{A}^n_{\mathbb{F}_p}$ dont la fonction indicatrice est calculable (on pourrait dire une partie récursive-constructible de $\mathbb{N}^r \times \mathbb{A}^n_{\mathbb{F}_p}$). En effet, il suffit d'appliquer ce qu'on a dit à la fonction qui coïncide avec la fonction considérée sur la partie $E$ et qui, ailleurs que sur $E$, prenne une valeur spéciale.

On pourra donc considérer que la fonction $x \mapsto x^{-1}$ est « totale » si besoin est, puisqu'elle est totale sur une partie constructible (le complémentaire de $\{0\}$) : les fonctions « *essentiellement* partielles » (celles qu'on ne peut pas trivialement considérer comme totales) sont celles pour lesquelles on est réellement confronté à un calcul qui pourrait ne pas terminer en temps fini, ne serait-ce que le contre-exemple présenté en **22.2**.

## 23. Élimination des quantificateurs

Commençons par rappeler un fait bien connu, dont on trouvera par exemple une démonstration dans [Fried et Jarden 2008, théorème 9.3.1], et pour lequel notre modèle de calcul fournira une preuve différente :

**23.1. Proposition.** *Si $Z$ est une partie constructible de $\mathbb{A}^{n+1}_{\mathbb{F}_p}$ alors sa projection sur $\mathbb{A}^n_{\mathbb{F}_p}$ est encore constructible, et, de plus, une formule booléenne la définissant se calcule algorithmiquement en fonction d'une formule booléenne définissant $Z$.*

*Démonstration.* Comme remarqué en **18.4**, on peut supposer que $Z$ est une partie localement fermée, définie par une formule purement conjonctive, c'est-à-dire la conjonction de $f_1 = \cdots = f_r = 0$ et de $g \neq 0$ (où $f_1, \ldots, f_r$ et $g$ sont des polynômes en $x_1, \ldots, x_n, y$).

[①] Une autre approche à (A), utilisant la théorie des modèles, repose sur le théorème de Löwenheim-Skolem ([Marker 2002, théorème 2.3.7]) : avec les notations de la démonstration, la suite $\neg\varphi_0, \neg\varphi_1, \neg\varphi_2, \ldots$, réunie à la théorie des corps algébriquement clos d'exposant caractéristique $p$, admet un modèle de degré de transcendance $\aleph_0$ sur $\mathbb{F}_p$, qui permet de tester si une fonction est totale. Pour (B), on pourrait substituer au théorème de compacité de $\mathbb{A}^{n\,\mathrm{cons}}_{\mathbb{F}_p}$ celui de la logique du premier ordre ([ibid., théorème 2.1.4]).

Observons maintenant le fait suivant : si cette fois $f_1, \dots, f_r$ et $g$ sont des polynômes en *une seule* variable $Y$ sur un *corps* $k$ algébriquement clos d'exposant caractéristique $p$, il existe un algorithme décidant si les équations $f_1 = \dots = f_r = 0$ ont une solution commune telle que $g \neq 0$. En effet, il suffit pour cela de calculer le PGCD $u$ (unitaire) de $f_1, \dots, f_r$ selon l'algorithme d'Euclide, et de vérifier qu'il ne divise pas $g^{\deg(u)}$ (l'exposant $\deg u$ garantit que les racines de $u$ dans une clôture algébrique sont un sous-ensemble de celles de $g$ si et seulement si $u$ divise $g^{\deg u}$).

Cet algorithme n'utilise manifestement que des opérations des corps et des tests d'égalité à 0, donc en vertu des résultats de la section 21, il est, en fait, calculable dans notre sens « universel », et comme la fonction calculée est totale, d'après **22.4**, ceci signifie qu'il existe une fonction calculable totale qui à des entiers naturels $r$, $D$ associe la formule d'une partie constructible $\Pi$ de l'espace affine $\mathbb{A}^{(r+1)(D+1)}_{\mathbb{F}_p}$ telle que des polynômes $f_i = c_{i,0} + c_{i,1}Y + \dots + c_{i,D}Y^D$ et $g = c'_0 + c'_1 Y + \dots + c'_D Y^D$ dans un corps $k$ algébriquement clos d'exposant caractéristique $p$ admettent une solution à $f_1 = \dots = f_r = 0$ et $g \neq 0$ si et seulement si $(c_{0,0}, \dots, c'_D) \in \Pi$.

On obtient alors une formule pour la projection de $Z$ en remplaçant les coefficients en $Y$ de nos $f_1, \dots, f_r, g$ de départ dans la formule booléenne définissant $\Pi$ (pour un $D$ assez grand). □

Autrement dit, ceci signifie que donnée une formule booléenne $\varphi$ on peut trouver, algorithmiquement, une formule booléenne (donc sans quantificateur !) qui équivaut sur tout corps algébriquement clos d'exposant caractéristique $p$ à la formule avec quantificateur $\exists y(\varphi(y, \dots))$.

(Rappelons que l'*existence* d'une telle formule sans quantificateur est conséquence du théorème de Chevalley : l'image d'un morphisme $\mathrm{Spec}(B) \to \mathrm{Spec}(A)$ de présentation finie est constructible ; cf. par exemple [ÉGA IV$_1$, § 1.8] ou [Olivier 1966-1968, § 3].)

Naturellement, ceci permet aussi de dire que la partie représentée par $\forall y(\varphi(y, \dots))$, qui est le complémentaire de la projection du complémentaire de la partie définie par $\varphi$, se calcule aussi effectivement.

La proposition suivante assure que, si une fonction $f$ ($\mathbb{N}^r \times \mathbb{A}^{n+1}_{\mathbb{F}_p} \to \mathbb{N}$ ou $\mathbb{N}^r \times \mathbb{A}^{n+1}_{\mathbb{F}_p} \to \mathbb{A}^1_{\mathbb{F}_p}$) est calculable, alors il est possible de tester *toutes* les valeurs d'une de ses variables $y$ (toutes, c'est-à-dire dans n'importe quel corps algébriquement clos d'exposant caractéristique $p$ contenant les autres variables) pour détecter si l'une ou toutes les valeurs de la fonction $f$ vaut 0 ; néanmoins, si une des valeurs n'est pas définie, alors la fonction résultante ne sera pas non plus définie. Précisément :

**23.2. Proposition.** *Si $f$ est une fonction partielle calculable $\mathbb{N}^r \times \mathbb{A}^{n+1}_{\mathbb{F}_p} \to \mathbb{N}$ ou $\mathbb{N}^r \times \mathbb{A}^{n+1}_{\mathbb{F}_p} \to \mathbb{A}^1_{\mathbb{F}_p}$ et $y$ une des coordonnées de $\mathbb{A}^{n+1}_{\mathbb{F}_p}$ intervenant dans $f$, alors il existe une fonction partielle calculable $f^{\forall y}$ (resp. $f^{\exists y}$) $\mathbb{N}^r \times \mathbb{A}^n_{\mathbb{F}_p} \to \mathbb{N}$ (manifestement unique) telle que, si $k$ est un corps algébriquement clos d'exposant caractéristique $p$, et si $\underline{z} \in \mathbb{N}^r \times \mathbb{A}^n_{\mathbb{F}_p}(k)$, alors $f^{\forall y}(\underline{z})$ (resp. $f^{\exists y}(\underline{z})$) vaut :*

- 0 *si, pour tout $y$ dans n'importe quel corps algébriquement clos $k'$ contenant $k$, la valeur $f(\underline{z}, y)$ est définie et vaut toujours* 0 *(resp. et prend au moins une fois la valeur* 0*),*
- 1 *si, pour tout $y$ dans n'importe quel corps algébriquement clos $k'$ contenant $k$, la valeur $f(\underline{z}, y)$ est définie et prend au moins une fois une valeur* $\neq 0$ *(resp. et ne prend que des valeurs* $\neq 0$*),*
- *non définie si pour au moins un $y$ dans un corps algébriquement clos $k'$ contenant $k$, la valeur $f(\underline{z}, y)$ n'est pas définie.*

(Ces fonctions $f^{\forall y}$ et $f^{\exists y}$ sont évidemment totales quand $f$ l'est. Par ailleurs, on peut bien sûr remplacer les mots « algébriquement clos » par « parfait » partout dans tout cet énoncé, et pour une fonction $\mathbb{N}^r \times \mathbb{A}^{n+1}_{\mathbb{F}_p} \to \mathbb{N}$ on peut les omettre.)

*Démonstration.* On peut évidemment supposer que $f$ est à valeurs dans $\mathbb{N}$ (quitte à la composer par la fonction de test à zéro, qu'on sait être calculable). Mettons qu'on cherche à représenter $f^{\forall y}$ (le cas de $f^{\exists y}$ étant complètement analogue). Soient $(\varphi_j, v_j)$ définissant $f$ avec les notations habituelles (**20.1**). Comme on l'a signalé en **20.3**, quitte à remplacer $\varphi_j$ par sa conjonction avec la négation de tous les $\varphi_{j'}$ antérieurs, on peut supposer que les parties constructibles qu'elles définissent sont deux à deux disjointes.

On parcourt la liste des $(\varphi_j, v_j)$ en calculant, au fur et à mesure, deux suites $(\psi_j)$ et $(\psi_j^{(0)})$ de formules booléennes en $n$ variables. (L'idée est que la partie constructible de $\mathbb{A}^n_{\mathbb{F}_p}$ définie par $\psi_j$ sera celle sur laquelle on sait que $f$ est définie pour tout $y$, et celle définie par $\psi_j^{(0)}$ sera celle sur laquelle on sait que $f$ est définie et vaut 0 pour tout $y$.) Précisément, pour chaque $j$, on calcule (d'après ce qu'on a rappelé en 23.1) une formule booléenne $\psi_j$ équivalente à $\forall y\,(\varphi_1(\ldots,y) \vee \cdots \vee \varphi_j(\ldots,y))$ et une autre formule booléenne $\psi_j^{(0)}$ équivalente à $\forall y\,(\varphi_{j_1}(\ldots,y) \vee \cdots \vee \varphi_{j_s}(\ldots,y))$ où les $j_u$ sont les indices entre 1 et $j$ tels qu'on ait $v_{j_u} = 0$. Enfin, pour représenter $f^{\forall y}$, on énumère alternativement $(\psi_j^{(0)}, 0)$ et $(\psi_j \wedge \neg\psi_j^{(0)}, 1)$.

Si pour un certain $\underline{z}$, la fonction $f(\underline{z}, y)$ est définie pour tous les $y$, c'est-à-dire que $\mathbb{A}^1_k$ est la réunion des parties constructibles définies par les $\varphi_j$ (où $\underline{z}$ a été substitué par sa valeur), par compacité de la topologie constructible il existe un $j$ tel que les parties constructibles définies par $\varphi_1, \ldots, \varphi_j$ suffisent à recouvrir $\mathbb{A}^1_k$, autrement dit, $\underline{z}$ vérifie $\psi_j$, et la fonction $f^{\forall y}$ est bien définie en $\underline{z}$ et y a manifestement la bonne valeur. □

Exactement comme en 22.3(A), on peut éviter d'avoir à considérer des éléments dans un corps $k'$ plus gros que $k$ en ne considérant qu'un seul corps $K$, fixé une fois pour toutes, algébriquement clos d'exposant caractéristique $p$ et de degré de transcendance infini (cf. [Marker 2002, proposition 4.3.4 et exemple 4.3.10]). (Dans le cas de fonctions calculables totales, on peut se contenter d'un corps algébriquement clos d'exposant caractéristique $p$ quelconque.)

**23.3. Corollaire.** *Si $f$ est une fonction partielle calculable $\mathbb{N}^r \times \mathbb{A}^{n+1}_{\mathbb{F}_p} \to \mathbb{N}$ et $y$ une des coordonnées de $\mathbb{A}^{n+1}_{\mathbb{F}_p}$ intervenant dans $f$, et $\max_y f$ (resp. $\min_y f$) désigne la fonction qui à $\underline{z} \in \mathbb{N}^r \times \mathbb{A}^n_{\mathbb{F}_p}$ associe le maximum (resp. minimum) des $f(\underline{z}, y)$ si ceux-ci sont définis pour tout $y$, alors $\max_y f$ et $\min_y f$ sont calculables. (Elles sont évidemment totales quand $f$ l'est.)*

*Démonstration.* La fonction $\max_y f$ peut se définir comme le plus petit $u$ tel que $\forall y(f(\underline{z}, y) \leq u)$, or la fonction qui à $u$ associe la valeur de vérité de $\forall y(f(\underline{z}, y) \leq u)$ est calculable d'après la proposition, et d'après 21.2 (stabilité par recherche), la fonction qui calcule le plus petit $u$ pour lequel cette valeur de vérité est vraie est elle-même calculable. □

**23.4. Remarque.** Dans la lignée de la remarque 22.5, on peut calculer le maximum d'une fonction calculable pour $y$ parcourant les éléments d'une partie constructible de $\mathbb{A}^n_{\mathbb{F}_p}$. (En effet, on peut définir modifier la fonction pour valoir $-\infty$ en-dehors de cette partie constructible, ou adapter facilement la ou les démonstrations ci-dessus.)

Le même procédé (soit en imitant la preuve de la proposition 23.2, soit en l'appliquant comme dans le corollaire qui précède), permet de montrer la calculabilité de toutes sortes de fonctions où on parcourt l'ensemble de tous les $y$ pour vérifier une propriété, par exemple on pourrait facilement construire une fonction calculable $f^{\exists! y}$, comme dans l'énoncé de la proposition 23.2, qui teste l'existence d'un unique $y$ pour lequel $f(\underline{z}, y) = 0$. Dans la section qui suit, on va voir qu'on peut, en fait, retourner la valeur de ce $y$.

## 24. Sélection d'un élément unique

Le résultat suivant a pour objet d'éclaircir le lien entre fonctions calculables au sens de 20.1 et 20.4 en ramenant les secondes aux premières. Trivialement, si $f : \mathbb{N}^r \times \mathbb{A}^n_{\mathbb{F}_p} \to \mathbb{A}^1_{\mathbb{F}_p}$ (partielle ou totale) est calculable, alors la fonction $h : \mathbb{N}^r \times \mathbb{A}^{n+1}_{\mathbb{F}_p} \to \mathbb{N}$ (partielle ou totale selon ce que $f$ est) définie par $h(i_1, \ldots, i_r, x_1, \ldots, x_n, y) = 1$ lorsque $y = f(i_1, \ldots, i_r, x_1, \ldots, x_n, y)$ et $= 0$ sinon (et non définie si $f$ ne l'est pas), est calculable. La proposition suivante, à rapprocher de 19.3, et son corollaire, montrent que la réciproque est vraie en un certain sens ; on commence par montrer un lemme facile qui est un cas particulier de ce qu'on veut prouver :

**24.1. Lemme.** *La fonction (totale) qui à des polynômes univariés $f_1, \ldots, f_r, g \in k[Y]$ à coefficients dans un corps algébriquement clos d'exposant caractéristique $p$ associe $(1, \xi)$ s'il existe un unique $\xi \in k$ tel que $f_1(\xi) = \cdots = f_r(\xi) = 0$ et $g(\xi) \neq 0$, et $(0, 0)$ sinon, est une fonction calculable totale, au sens du modèle de calcul universel (c'est-à-dire de la définition 20.4).*

(Il faut bien sûr comprendre qu'on prend en entrée un couple $(r, D) \in \mathbb{N}^2$, où $r$ est le nombre de polynômes et $D$ une borne sur leurs degrés, ainsi que leurs coefficients vus comme un élément de $\mathbb{A}_{\mathbb{F}_p}^{(r+1)(D+1)}$, et on renvoie le couple annoncé — ou, si on veut, qu'on dispose de deux fonctions, l'une qui teste s'il existe un $\xi$ unique comme indiqué et l'autre qui renvoie sa valeur.)

*Démonstration.* Le raisonnement est analogue à celui dans la démonstration de 23.1 : on effectue le calcul suivant *dans notre modèle de calcul universel.*

Donnés des polynômes univariés $f_1, \ldots, f_r, g \in k[Y]$ à coefficients dans un corps $k$ d'exposant caractéristique $p$, on calcule le PGCD unitaire $u$ de $f_1, \ldots, f_s$ (au moyen de l'algorithme d'Euclide), ce qui fournit un polynôme dont les racines sont les racines communes à $f_1, \ldots, f_s$. Puis on calcule le quotient $v$ de $u$ par le PGCD unitaire de $u$ avec $g^{\deg(u)}$, ce qui fournit un polynôme dont les racines sont celles de $u$ qui ne sont pas des racines de $g$ (comme dans la démonstration de 23.1, l'exposant $\deg u$ sert à garantir qu'il est « assez grand »).

Puis on répète (au plus $\deg(v)$ fois) la boucle suivante : si le degré de $v$ vaut 0, on renvoie $(0, 0)$ ; si $v$ est de degré 1 et vaut $aY + b$ on renvoie $(1, -b/a)$. Si la dérivée $v'$ de $v$ vaut 0, on remplace $v$ par sa racine $p$-ième (qu'on peut calculer car notre modèle de calcul permet d'effectuer $\mathrm{Frob}_p^{-1}$ ; remarquons que seuls les coefficients de degré multiple de $p$ sont non-nuls ici) ; sinon, on remplace $v$ par le quotient de celui-ci par le PGCD de $v$ et $v'$ ; si ces différentes opérations ne changent pas le polynôme (c'est-à-dire qu'il est à racines simples et de degré $> 1$), on renvoie $(0, 0)$. Il est clair que ces opérations renvoient bien l'unique racine $\xi$ de $v$ lorsque celle-ci existe et est unique. □

**24.2. Proposition.** *Si $h$ est une fonction partielle calculable $\mathbb{N}^r \times \mathbb{A}_{\mathbb{F}_p}^{n+1} \to \mathbb{N}$, alors il existe une fonction partielle calculable $f : \mathbb{N}^r \times \mathbb{A}_{\mathbb{F}_p}^{n} \to \mathbb{A}_{\mathbb{F}_p}^{1}$ telle que, si $k$ est un corps algébriquement clos d'exposant caractéristique $p$, et si pour $(\underline{i}, \underline{x}) \in \mathbb{N}^r \times \mathbb{A}_{\mathbb{F}_p}^{n}(k)$ il existe un* unique *$y$ tel que $h(\underline{i}, \underline{x}, y) = 1$, alors $f(\underline{i}, \underline{x})$ vaut $y$ (toute autre valeur de $f$ est non spécifiée).*

*Démonstration.* Soient $(\varphi_j, v_j)$ définissant $h$ avec les notations habituelles. Comme on l'a signalé en 20.3, on peut supposer que les $\varphi_j$ sont mutuellement exclusives et sont des formules purement conjonctives, et on le fera.

On énumère successivement les $(\varphi_j, v_j)$ en ignorant celles pour lesquelles $v_j \neq 1$. Considérons maintenant une formule $\varphi_j$ pour laquelle $v_j = 1$. Comme remarqué en 18.4, on peut supposer qu'elle s'écrit comme la conjonction de $f_1 = 0, \ldots, f_r = 0$ et de $g \neq 0$ (où $f_1, \ldots, f_r$ et $g$ sont des polynômes en $x_1, \ldots, x_n, y$). D'après le lemme (et le corollaire 22.4), on peut calculer un nombre fini (disons $s$) de parties constructibles disjointes de l'espace affine (de dimension $(r+1)(D+1)$ où $D$ est une borne sur les degrés des polynômes) et de fonctions F-rationnelles sur chacune d'elles, qui, sur n'importe quel corps $k$ algébriquement clos d'exposant caractéristique $p$, si des polynômes univariés ont un unique zéro commun qui ne soit pas zéro d'un autre, *calculent* le zéro en question. On substitue dans les équations de ces parties et fonctions les coefficients de nos $f_1, \ldots, f_r$ en la variable $y$ (qui sont des polynômes en $x_1, \ldots, x_n$) : ceci définit les $(\psi_{j,1}, w_{j,1}), \ldots, (\psi_{j,s}, w_{j,s})$ qui participent à la définition de la fonction $f$. □

**24.3. Corollaire.** *Une fonction (partielle ou totale) $f : \mathbb{N}^r \times \mathbb{A}_{\mathbb{F}_p}^{n} \to \mathbb{A}_{\mathbb{F}_p}^{1}$ est calculable, si et seulement si la fonction $h : \mathbb{N}^r \times \mathbb{A}_{\mathbb{F}_p}^{n+1} \to \mathbb{N}$ (partielle ou totale selon ce que $f$ est) définie par $h(i_1, \ldots, i_r, x_1, \ldots, x_n, y) = 1$ lorsque $y = f(i_1, \ldots, i_r, x_1, \ldots, x_n)$ et $= 0$ sinon (et non définie si $f$ ne l'est pas), est calculable*

*Démonstration.* On a déjà remarqué que l'implication « seulement si » est triviale. Pour ce qui est du « si », il découle facilement de la proposition (sauf éventuellement que $f$ pourrait être définie en trop de points, mais il suffit d'évaluer $h(\underline{i}, \underline{x}, 0)$ avant de faire tout autre calcul pour éviter ce phénomène). □

## 25. Factorisation dans le modèle de calcul

**25.1. Remarque.** Dans le cadre de la proposition 23.2, le modèle de calcul que nous avons adopté permet de tester l'existence d'un $x$ tel que $f(x, \underline{z}) = 0$, mais ne permet pas, *stricto sensu*, d'en *calculer* un. Par exemple, si $f$ (fonction d'une seule variable $x$, et à valeurs entières) est définie comme valant 0 si $x^2 + 1 = 0$ et 1 si $x^2 + 1 \neq 0$, alors $\exists x(f(x) = 0)$ est toujours vérifié, c'est-à-dire que la fonction $f^{\exists x}$ est égale à la constante 0, et pourtant, en caractéristique 0 ou bien $\equiv 3 \pmod 4$, aucun élément du corps premier, donc aucune fonction calculable constante (cf. 22.4), ne peut témoigner de ce fait.

Ceci est raisonnable si on pense que notre modèle de calcul ne « sait » rien du corps sur lequel il travaille (autrement que sa caractéristique), même pas si ce corps est algébriquement clos.

Les résultats qui vont suivre ont pour objet de montrer que, même si le modèle de calcul ne permet pas de factoriser un polynôme (comme l'illustre l'exemple de la remarque ci-dessus), il permet cependant de calculer tout ce qu'on pourrait calculer à partir d'une factorisation.

Nous commençons par considérer l'action de $\mathfrak{S}_n$ (groupe symétrique sur $n$ objets) sur $\mathbb{A}^{m+n}_{\mathbb{F}_p}$ (dont les coordonnées seront notées $Z_1, \ldots, Z_m, X_1, \ldots, X_n$) par permutation des $n$ dernières coordonnées. Cette action correspond à une action sur $\mathbb{F}_p[Z_1, \ldots, Z_m, X_1, \ldots, X_n]$ par permutation des variables $X_i$ : on dira comme d'habitude qu'un polynôme est « totalement symétrique » en les $n$ dernières variables lorsqu'il est invariant par cette action : on sait bien alors qu'il peut s'écrire, de façon unique, comme un polynôme en $Z_1, \ldots, Z_m$ et en les fonctions symétriques élémentaires de $X_1, \ldots, X_n$. Nous montrons deux lemmes cruciaux quoique intuitivement évidents.

**25.2. Lemme.** *Soit $F$ un fermé de Zariski de $\mathbb{A}^{m+n}_{\mathbb{F}_p}$ tel que $\sigma(F) = F$ pour tout $\sigma \in \mathfrak{S}_n$, où le groupe $\mathfrak{S}_n$ symétrique sur $n$ objets opère en permutant les $n$ dernières coordonnées. Alors $F$ peut être défini par la formule $f_1 = \cdots = f_r = 0$ où les $f_i$ sont eux-mêmes totalement symétriques en les $n$ dernières variables. De plus, on peut calculer les $f_i$ algorithmiquement en fonction d'équations définissant $F$.*

Rappelons qu'un « fermé de Zariski » n'est pas un *sous-schéma* fermé mais bien une *partie* constructible définie par une formule booléenne positive, cf. **18.3**.

*Démonstration.* Soient $g_1 = \cdots = g_r = 0$ les équations initiales de $F$. Quitte à ajouter à ces équations toutes celles qui s'obtiennent en permutant de manière arbitraire les $n$ dernières variables d'un des $g_i$, on peut supposer que l'ensemble des $g_i$ est stable par $\mathfrak{S}_n$ (agissant par permutation des $n$ dernières variables). Soient $f_1, \ldots, f_r$ les $r$ fonctions symétriques élémentaires appliquées à $g_1, \ldots, g_r$ (i.e., $f_1 = g_1 + \cdots + g_r$, $f_2 = g_1 g_2 + g_1 g_3 + \cdots + g_{r-1} g_r$, ..., $f_r = g_1 \cdots g_r$). Alors les $f_i$ sont manifestement totalement symétriques en les $r$ dernières variables, et les équations $f_1 = \cdots = f_r = 0$ sont équivalentes (ensemblistement !) à $g_1 = \cdots = g_r = 0$. ☐

(On aurait bien sûr pu formuler ce lemme dans un contexte beaucoup plus général. Remarquons par ailleurs qu'il n'est plus vrai avec des multiplicités : le fermé $X_1 = \cdots = X_n = 0$ de $\mathbb{A}^n_{\mathbb{F}_p}$ ne peut pas être défini par des équations totalement symétriques en les $X_i$.)

Montrons maintenant le résultat pour des combinaisons booléennes non nécessairement positives :

**25.3. Lemme.** *Soit $\varphi$ une formule booléenne (cf. **18.1**) en les indéterminées $Z_1, \ldots, Z_m, X_1, \ldots, X_n$ qui soit totalement symétrique en les $n$ dernières variables, c'est-à-dire que $\varphi$ soit équivalente à $\varphi(Z_1, \ldots, Z_m, X_{\sigma(1)}, \ldots, X_{\sigma(n)})$ pour toute permutation $\sigma \in \mathfrak{S}_n$. Alors on peut écrire $\varphi$ comme combinaison booléenne de formules de la forme $f = 0$ où chaque $f$ est lui-même symétrique en les $n$ dernières variables. De plus, cette transformation peut être faite de manière algorithmique.*

*Démonstration.* Soit $P$ l'ensemble de tous les polynômes qui apparaissent (positivement ou négativement) dans la formule booléenne $\varphi$, ainsi que tous les polynômes obtenus en appliquant un $\sigma \in \mathfrak{S}_n$ quelconque à l'un de ceux-ci (c'est-à-dire en permutant ses $n$ dernières variables par $\sigma$). Ainsi, $\mathfrak{S}_n$ opère sur $P$, et donc aussi sur l'ensemble des parties de cet ensemble (en préservant, évidemment, le cardinal des parties).

Pour chaque partie $I$ de $P$, on définit deux formules booléennes : $\pi_I$ sera la conjonction de tous les $f = 0$ pour $f \in I$, et $\delta_I$ sera la conjonction de tous les $f = 0$ pour $f \in I$ et des $f \neq 0$ pour $f \notin I$. Il est ainsi évident que les $\delta_I$ sont deux à deux exclusifs (la conjonction de deux d'entre eux est fausse) et recouvrent tous les cas possibles (leur disjonction est vraie) ; intuitivement, les $\delta_I$ représentent tous les cas possibles de nullité ou non-nullité des $f$ (donc tous les cas de figure auxquels on sera éventuellement amenés à avoir affaire), mais il se peut, bien sûr, que certains des $\delta_I$ soient universellement faux.

Remarquons que $\pi_I$ peut s'écrire comme $\delta_I \vee \bigvee_J \delta_J$ où la disjonction $\bigvee_J \delta_J$ est prise sur toutes les $J$ contenant strictement $I$.

Il est clair que $\varphi$ peut s'écrire comme disjonction des $\delta_I$ pour l'ensemble $\mathscr{C}$ de parties de $I$ telles que $\delta_I \Rightarrow \varphi$ au sens de 18.2 : il est également clair qu'on peut calculer $\mathscr{C}$ algorithmiquement. Naturellement, $\mathfrak{S}_n$ opère sur $\mathscr{C}$, c'est-à-dire que ce dernier est une réunion d'orbites.

Si $\mathscr{X}$ est une réunion d'orbites de $\mathfrak{S}_n$ agissant sur les parties de $P$, on notera encore $\pi_{\mathscr{X}}$ (resp. $\delta_{\mathscr{X}}$) la disjonction des $\pi_I$ (resp. des $\delta_I$) pour $I \in \mathscr{X}$. D'après le lemme précédent, les $\pi_{\mathscr{X}}$ (puisqu'ils définissent des fermés de Zariski invariants par $\mathfrak{S}_n$) peuvent s'écrire comme conjonction de formules de la forme $f = 0$ avec $f$ symétrique en les $X_i$, et on peut calculer ces $f$ algorithmiquement.

Montrons maintenant que si $\mathscr{X}$ est une réunion d'orbites de $\mathfrak{S}_n$ agissant sur les parties de $P$, alors $\delta_{\mathscr{X}}$ est une combinaison booléenne de formules de la forme $f = 0$ avec $f$ symétrique en les $X_i$ (en appliquant à $\mathscr{X} = \mathscr{C}$, ceci donnera le résultat voulu, vu que $\delta_{\mathscr{C}} = \varphi$). Pour cela, d'après ce qu'on vient de dire, il suffit de montrer que $\delta_{\mathscr{X}}$ est une combinaison booléenne des $\pi_{\mathscr{Y}}$. On peut évidemment supposer que $\mathscr{X}$ est une unique orbite. On procède alors par récurrence décroissante sur le cardinal commun $s$ des éléments de $\mathscr{X}$.

En écrivant $\pi_{\mathscr{X}} = \bigvee_{I\in\mathscr{X}} \pi_I$ et chaque $\pi_I$ comme $\delta_I \vee \bigvee_J \delta_J$ (comme expliqué ci-dessus), on fait apparaître $\pi_{\mathscr{X}}$ comme la disjonction (exclusive !) de $\delta_{\mathscr{X}}$ et de certains $\delta_J$ pour $J$, tous de cardinal $> s$, parcourant un $\mathscr{Z}$ qui est une réunion d'orbites. Par l'hypothèse de récurrence, $\delta_{\mathscr{Z}}$ s'écrit comme combinaison booléenne des $\pi_{\mathscr{Y}}$, donc il en va de même de $\delta_{\mathscr{X}} = \pi_{\mathscr{X}} \wedge \neg\delta_{\mathscr{Z}}$. □

**25.4. Proposition.** *Soit $h\colon \mathbb{N}^r \times \mathbb{A}^{m+n}_{\mathbb{F}_p} \to \mathbb{N}$ une fonction partielle calculable au sens de la définition 20.1, et telle que pour toute permutation $\sigma \in \mathfrak{S}_n$ on ait $h(i_1,\dots,i_r,z_1,\dots,z_m,x_1,\dots,x_n) = h(i_1,\dots,i_r,z_1,\dots,z_m,x_{\sigma(1)},\dots,x_{\sigma(n)})$.*

*Pour $k$ un corps algébriquement clos (variable) d'exposant caractéristique $p$, considérons l'application partielle $f_k$ qui à $(i_1,\dots,i_r,z_1,\dots,z_m,c_1,\dots,c_n) \in \mathbb{N}^r \times k^{m+n}$ associe $h(i_1,\dots,i_r,z_1,\dots,z_m,\xi_1,\dots,\xi_n)$, où $\xi_1,\dots,\xi_n$ sont les racines, avec multiplicités, du polynôme $X^n + c_1X^{n-1} + \dots + c_n \in k[X]$. Alors les $f_k$ définissent une fonction partielle calculable $f\colon \mathbb{N}^r \times \mathbb{A}^{m+n}_{\mathbb{F}_p} \to \mathbb{N}$ (évidemment unique : cf. 20.2).*

*Démonstration.* Suivant la définition 20.1, la fonction $h$ est définie par une suite (énumérée par une machine de Turing) de parties constructibles $H_j = \{\phi_j\}$ de $\mathbb{A}^{m+n}_{\mathbb{F}_p}$ et de valeurs $v_j \in \mathbb{N}$ associées. Quitte à remplacer $\phi_j$ par la disjonction des $\phi_j(z_1,\dots,z_m,x_{\sigma(1)},\dots,x_{\sigma(n)})$ pour $\sigma \in \mathfrak{S}_n$ (ce qui ne change pas la valeur de $h$) d'après l'hypothèse faite sur celle-ci), on peut supposer que $H_j$ vérifie $\sigma(H_j) = H_j$ pour tout $\sigma \in \mathfrak{S}_n$ (où $\sigma$ opère sur $\mathbb{A}^{m+n}_{\mathbb{F}_p}$ en permutant les $n$ dernières coordonnées).

Il s'agit donc de voir que si $H$ est une partie constructible de $\mathbb{A}^{m+n}_{\mathbb{F}_p}$ telle que $\sigma(H) = H$ pour tout $\sigma \in \mathfrak{S}_n$, alors le foncteur $F$ qui à un corps $k$ algébriquement clos d'exposant caractéristique $p$ associe l'ensemble des $(z_1,\dots,z_m,c_1,\dots,c_n) \in k^n$ tels que $(z_1,\dots,z_m,\xi_1,\dots,\xi_n) \in H(k)$ où $\xi_1,\dots,\xi_n$ sont les racines, avec multiplicité, du polynôme $X^n + c_1X^{n-1} + \dots + c_n \in k[X]$, définit bien une partie constructible de $\mathbb{A}^{m+n}_{\mathbb{F}_p}$ et que, de plus, cette partie constructible $F$ se calcule effectivement en fonction de $H$ (c'est-à-dire qu'on peut décrire une machine de Turing qui, en fonction d'une formule définissant $H$, produit une formule définissant $F$).

Or, en vertu des lemmes qui précèdent, on peut se ramener au cas où $H$ est un fermé de Zariski défini par une unique équation $g = 0$ avec $g$ un polynôme totalement symétrique en les $n$ dernières variables. On sait alors algorithmiquement écrire $g$ comme un polynôme (à coefficients dans $\mathbb{F}_p$) en les $Z_1,\dots,Z_m$ et les fonctions symétriques élémentaires en $X_1,\dots,X_n$, ce qui donne les équations du fermé $F$ recherché. □

**25.5. Proposition.** *Le même énoncé que 25.4 vaut pour une fonction $h\colon \mathbb{N}^r \times \mathbb{A}^{m+n}_{\mathbb{F}_p} \to \mathbb{A}^1_{\mathbb{F}_p}$.*

*Démonstration.* Il suffit d'appliquer 24.3 pour se ramener à 25.4. □

**25.6. Remarque.** Dans ces deux propositions, la fonction $h$ est supposée vérifier $h(i_1,\dots,i_r,z_1,\dots,z_m,x_1,\dots,x_n) = h(i_1,\dots,i_r,z_1,\dots,z_m,x_{\sigma(1)},\dots,x_{\sigma(n)})$ pour toute permutation $\sigma$ des $n$ dernières variables. On peut évidemment se passer de cette hypothèse si on n'exige d'avoir $f(i_1,\dots,i_r,z_1,\dots,z_m,c_1,\dots,c_n) = h(i_1,\dots,i_r,z_1,\dots,z_m,\xi_1,\dots,\xi_n)$ que lorsque le membre de droite ne dépend pas de la permutation des $\xi_i$ — en effet, il suffit de remplacer, le cas échéant, la fonction $h$ par une fonction $h'$ qui calcule $h(i_1,\dots,i_r,z_1,\dots,z_m,x_{\sigma(1)},\dots,x_{\sigma(n)})$ pour tout $\sigma$, renvoie la valeur commune si toutes ces valeurs sont toutes égales, et renvoie 42 si elles ne le sont pas.

On peut donc résumer ces résultats ainsi : dans le contexte de notre modèle de calcul universel, *on peut supposer qu'une fonction calculable peut calculer les racines d'un polynôme $X^n+c_1X^{n-1}+\cdots+c_n \in k[X]$ lorsque les calculs qu'elle fait avec ces racines sont invariants par toute permutation de celles-ci.*

**25.7. Factorisation par simulation.** Expliquons une autre manière de comprendre ces résultats, qui en fournit une démonstration différente, dans l'optique de la présentation « boîte noire » du modèle de calcul. Il s'agit donc d'expliquer pourquoi le fait d'ajouter à une machine de Turing, qui manipule dans des boîtes noires les éléments d'un corps $k$ (parfait d'exposant caractéristique $p \geq 1$) non spécifié, la faculté de factoriser un polynôme dans sa clôture algébrique, n'augmente pas les fonctions qu'elle est capable de calculer.

(La description faite ci-dessous est à rapprocher de la méthode décrite en [Lombardi et Quitté 2011, chapitre VII] pour se passer de clôtures algébriques en algèbre constructive.)

Imaginons donc une machine de Turing qui effectue des calculs en manipulant les éléments de $k$ dans des boîtes noires et qui dispose, en plus de la possibilité de calculer des additions, soustractions, multiplications, divisions, et Frobenius inverse (comme en 21), de la faculté de factoriser un polynôme. On souhaite expliquer comment effectuer les mêmes calculs sans cette faculté : il s'agit donc de simuler la possibilité de factoriser les polynômes. Dans la discussion qui suit, on appellera « programme simulé » la machine de Turing qui possède la faculté de factoriser les polynômes, et « simulateur » celui que nous décrivons, qui doit effectuer la même tâche sans factoriser de polynômes.

Le simulateur maintient une $k$-algèbre artinienne (c'est-à-dire de dimension finie comme $k$-espace vectoriel) et étale (c'est-à-dire géométriquement réduite, ou, ce qui revient au même puisque $k$ est un corps parfait, réduite), que nous appellerons « algèbre courante » : celle-ci est vue comme quotient d'une algèbre de polynômes sur $k$ par un idéal de celle-ci, et décrit par exemple par une base de Gröbner (cf. **13.1**) ; au cours de l'exécution du simulateur, on pourra ajouter des indéterminées à cette « algèbre courante », ou des relations (c'est-à-dire faire grossir l'idéal), mais il y aura toujours un morphisme naturel de l'ancienne algèbre courante vers la nouvelle. Les éléments manipulés par le programme simulé (sous forme de boîtes noires) sont représentés, dans le simulateur, par des éléments de cette algèbre courante (et si l'algèbre courante est modifiée par l'ajout de nouvelles indéterminées ou de nouvelles relations, les éléments en question sont transportés par le morphisme naturel, c'est-à-dire qu'ils sont représentés par les mêmes polynômes dans les anciennes indéterminées).

Initialement, l'algèbre courante ne comporte pas d'indéterminées et pas de relations (elle est donc égale à $k$). Lorsque le programme simulé effectue une somme, une différence, ou un produit, le simulateur effectue la même opération dans l'algèbre courante.

Lorsque le programme simulé demande à calculer l'inverse d'un élément, représenté dans le simulateur par un élément $x$ de l'algèbre courante, le simulateur ajoute une nouvelle indéterminée $Y$ à l'algèbre courante, et la relation $xY-1=0$ (c'est-à-dire ajoute l'élément $xY-1$ à ceux qui engendrent l'idéal des relations de l'algèbre courante), et il renvoie l'élément représenté par l'indéterminée $Y$. (Remarquons que l'algèbre courante peut devenir nulle, si $x$ valait 0 : on s'apercevra ci-dessous du problème.)

Lorsque le programme simulé demande à calculer $\mathrm{Frob}_p^{-1}$ d'un élément représenté dans le simulateur par un élément $x$ de l'algèbre courante, le simulateur ajoute une nouvelle indéterminée $Y$ et la relation $Y^p - x = 0$ ; puis il calcule le réduit de cette algèbre (cf. par exemple [Becker et Weispfenning 1993, théorème 8.22], qui utilise le « lemme 92 » de Seidenberg, [ibid., lemme 8.13] ; on utilise bien sûr ici le fait que $k$ est parfait et que le simulateur peut calculer $\mathrm{Frob}_p^{-1}$ dans $k$), et renvoie la classe de $Y$ dans l'algèbre courante.

Lorsque le programme simulé demande à tester la nullité d'un élément représenté dans le simulateur par un élément $x$ de l'algèbre courante, le simulateur ajoute une nouvelle indéterminée $Y$ et la relation $Y - x = 0$ à l'algèbre courante, puis il calcule le polynôme $h \in k[Y]$ obtenu par l'élimination de toutes les autres variables de l'algèbre courante, c'est-à-dire le polynôme qui engendre l'intersection de $k[Y]$ avec l'idéal de $k[Y, \ldots]$ définissant l'algèbre courante (où les points de suspension représentent d'autres indéterminées) : cf. [Eisenbud 1995, § 15.10.4] pour un tel algorithme. Si le degré de $h$ vaut 0, cela signifie que l'algèbre courante est nulle, ce qui ne peut se produire que si le programme simulé a fait calculer l'inverse de 0 : dans ce cas, le simulateur exécute une boucle infinie (de façon à ne pas terminer). Sinon, et si $h$ est de la forme $Y^r$ (i.e., $Y$ est nilpotent dans l'algèbre courante), le simulateur ajoute la relation $Y = 0$ et répond que $x$ vaut 0 ; sinon, le simulateur ajoute la relation $h(Y)/Y^r = 0$ où $r$ est l'ordre (=la valuation en $Y$) de $h$ et répond que $x$ ne vaut pas 0. De même, si le programme simulé retourne $x$ comme valeur finale, le simulateur effectue la même opération d'élimination : si le polynôme $h$ est de degré 1, le simulateur retourne son unique zéro comme valeur finale ; s'il est de degré 0 (on a effectué une division par 0) ou $> 1$ (le programme simulé a violé la promesse de calculer une valeur qui ne dépende pas de la numérotation des racines des polynômes factorisés), le simulateur ne termine pas.

Enfin, si le programme simulé demande à factoriser un polynôme $P = X^n + c_1 X^{n-1} + \cdots + c_n$, on ajoute à l'algèbre courante les indéterminées et relations de l'algèbre « de décomposition universelle » de $P$ (cf. par exemple [Lombardi et Quitté 2011, III, § 4, et VII, § 4]), c'est-à-dire l'algèbre $A$ (finie sur $k$) quotient de $k[Z_1, \ldots, Z_n]$ par les relations du type $e_i(Z_1, \ldots, Z_n) - (-1)^i c_i$ où les $e_i$ sont les fonctions symétriques élémentaires. (Rappelons que si $P = \prod_{i=1}^{n}(X - \xi_i)$, avec les $\xi_i$ deux à deux distincts, est un polynôme séparable sur un corps, alors $A$ est, en fait, isomorphe à $k^{\mathfrak{S}_n}$ par l'isomorphisme qui envoie $Z_i$ sur l'ensemble des $\sigma(\xi_i)$ pour tous les $\sigma \in \mathfrak{S}_n$.) Puis on réduit l'algèbre courante, comme expliqué pour le calcul de $\mathrm{Frob}_p^{-1}$. Le simulateur renvoie comme racines de $P$ les indéterminées $Z_1, \ldots, Z_n$ qui viennent d'être ajoutées.

Pour montrer la correction de cet algorithme, il suffit de considérer les points de l'algèbre courante sur une clôture algébrique de $k$ : il s'agit d'un nombre fini de tuples de valeurs de cette clôture algébrique, à savoir les valeurs que peuvent prendre les indéterminées pour une certaine numérotation des racines des polynômes factorisés. Dès lors que le programme simulé calcule une valeur ne dépend pas de cette numérotation, le simulateur effectue bien le calcul voulu.

**25.8. Extensions de corps dans le modèle de calcul.** Si $k$ est un corps algébriquement clos d'exposant caractéristique $p \geq 1$, dont les éléments sont manipulés sous forme de « boîtes noires » (c'est-à-dire dans le cadre du modèle de calcul développé ici), on peut encore manipuler les éléments de tous les corps obtenus à partir de $k$ par itération des extensions énumérées en 12.5, à savoir

(i) l'ajout d'un transcendant : $L = K(T)$ où $T$ est une indéterminée,
(ii) l'ajout d'un élément algébrique : $L = K[X]/(f)$ où $f \in K[X]$ est irréductible (non supposé séparable), donné,
(iii) le passage à « la » clôture algébrique $L = K^{\mathrm{alg}}$ de $K$,
(iv) le passage à « la » clôture séparable $L = K^{\mathrm{sép}}$ de $K$,
(v) (dans le cas où $K$ est de caractéristique $p > 0$) le passage à la clôture parfaite $L = K^{1/p^\infty}$ de $K$,

en utilisant essentiellement les mêmes algorithmes qu'expliqués en *loc. cit.* L'algorithme de Rabin pour construire la clôture algébrique n'est pas applicable tel quel, mais on peut le remplacer par celui, essentiellement équivalent, décrit en 25.7 qui permet de manipuler les éléments de « la » clôture algébrique d'un corps $K$ (ou plutôt, de toutes ses clôtures algébriques simultanément) si on sait manipuler ceux de $K$ (la construction suppose $K$ parfait, mais ceci n'est pas une difficulté puisqu'on peut de toute façon passer d'abord à la clôture parfaite de $K$).

## 26. Indépendance de la caractéristique (esquisse)

Toute la description faite jusqu'ici du modèle de calcul « universel » supposait fixée la caractéristique du corps manipulé. Esquissons maintenant la manière dont il faudrait modifier les choses pour obtenir l'uniformité en $p$. L'intuition sous-jacente, dans la présentation « boîte noire » est que la machine qui manipule des boîtes noires ne peut jamais être sûre d'être en caractéristique 0, elle peut seulement tester si $p \neq 0$ pour différentes valeurs de $p$, et donc « minorer » la caractéristique ; par ailleurs, elle ne peut appliquer $\mathrm{Frob}_p^{-1}$ qu'à condition que la caractéristique vaille effectivement $p$ (sous peine que la fonction ne soit pas définie).

**26.1. Types et parties constructibles.** On appelle formule booléenne en $X_1, \ldots, X_n$ sur $\mathbb{Z}$ une combinaison booléenne d'expressions de la forme $f(X_1, \ldots, X_n) = 0$ où $f \in \mathbb{Z}[X_1, \ldots, X_n]$ ; on lui associe une partie constructible de $\mathbb{A}^n_{\mathbb{Z}}$, ou de $k^n$ pour $k$ un corps quelconque, de la façon évidente. L'analogue de **18.2** vaut encore, en remplaçant « $\mathbb{F}_p$ » par « $\mathbb{Z}$ » et en supprimant les mots « d'exposant caractéristique $p$ ». On remarquera, par exemple, que la formule $p = 0$ définit la partie pleine sur les corps de caractéristique $p$ et la partie vide sur les corps de caractéristique $\neq p$.

**26.2. Fonctions F-rationnelles.** On définit la notion de « fonction F-rationnelle sur $\mathbb{Z}$ » de façon stupide : une telle fonction est *soit* la donnée d'une fonction rationnelle sur $\mathbb{Q}$ (c'est-à-dire le quotient de deux polynômes à coefficients dans $\mathbb{Z}$, qu'on peut supposer sans facteur premier non trivial), auquel cas elle est définie aux points où son dénominateur réduit ne s'annule pas, *soit* la donnée d'une fonction F-rationnelle sur $\mathbb{F}_p$ pour un $p > 0$ (c'est-à-dire pouvant faire intervenir $\mathrm{Frob}_p^{-1}$), auquel cas elle n'est définie que sur les corps de caractéristique $p$. Pour déplaisante que soit cette définition, la notion de fonction calculable qui en résulte est raisonnable.

**26.3. Fonctions partielles calculables.** On définit la notion de fonction partielle calculable $f : \mathbb{N}^r \times \mathbb{A}^n_{\mathbb{Z}} \to \mathbb{N}$ ou $f : \mathbb{N}^r \times \mathbb{A}^n_{\mathbb{Z}} \to \mathbb{A}^1_{\mathbb{Z}}$ comme dans les définitions **20.1** et **20.4**, en remplaçant « $\mathbb{F}_p$ » par « $\mathbb{Z}$ » et en supprimant les mots « d'exposant caractéristique $p$ ».

À titre d'exemple, la fonction (*partielle !*) qui renvoie la caractéristique du corps $k$ si celle-ci est non nulle, et non définie si $k$ est de caractéristique nulle, est calculable : il suffit d'énumérer les formules $p = 0$ (pour $p$ parcourant les nombres premiers), chacune étant associée à la valeur $p$ (autrement dit, la machine teste successivement les nombres premiers et retourne la valeur $p$ si elle constate que $p = 0$).

L'équivalence avec le modèle « boîte noire » se montre alors sans difficulté nouvelle : on souligne que la machine peut appliquer une fonction $\mathrm{Frob}_p^{-1}$ pour n'importe quel $p$ (mais dès qu'elle le fait, la fonction cesse d'être définie si la caractéristique ne valait pas $p$).

**26.4. Fonctions totales.** On dira qu'une fonction partielle calculable $f : \mathbb{N}^r \times \mathbb{A}^n_{\mathbb{Z}} \to \mathbb{N}$ ou $f : \mathbb{N}^r \times \mathbb{A}^n_{\mathbb{Z}} \to \mathbb{A}^1_{\mathbb{Z}}$ est **totale** lorsque sa valeur est définie pour toutes valeurs $(i_1, \dots, i_r, x_1, \dots, x_n)$, les $x_1, \dots, x_n$ appartenant à n'importe quel corps parfait $k$ (quelle que soit sa caractéristique).

L'analogue de **22.3**(A) est qu'une fonction est totale lorsque, pour chaque $p \geq 1$, elle est définie en tous les $(i_1, \dots, i_r, x_1, \dots, x_n)$ où $x_1, \dots, x_n$ appartenant à $K_p$, corps algébriquement clos d'exposant caractéristique $p$ et de degré de transcendance infini sur $\mathbb{F}_p$. L'énoncé **22.3**(B) reste valable sans modification. Le corollaire **22.4** devient l'énoncé suivant (exactement analogue, mais suffisamment notable pour être répété intégralement) :

**26.5. Corollaire.** *Une fonction totale calculable $\mathbb{A}^n_{\mathbb{Z}} \to \mathbb{N}$ est exactement la donnée d'une partition de $\mathbb{A}^n_{\mathbb{Z}}$ en un nombre* fini *de parties constructibles, sur chacune desquelles la fonction prend une valeur constante entière (en particulier, la fonction est bornée).*

*Une fonction totale calculable $\mathbb{A}^n_{\mathbb{Z}} \to \mathbb{A}^m_{\mathbb{Z}}$ est exactement la donnée d'une partition de $\mathbb{A}^n_{\mathbb{Z}}$ en un nombre* fini *de parties constructibles, sur chacune desquelles la fonction coïncide avec une certaine fonction F-rationnelle (définie partout sur la partie en question : donc soit une fonction rationnelle sur $\mathbb{Q}$, soit une fonction F-rationnelle sur $\mathbb{F}_p$ uniquement sur une partie constructible du fermé $\mathbb{A}^n_{\mathbb{F}_p}$ de $\mathbb{A}^n_{\mathbb{Z}}$).*

*Une fonction totale calculable $\mathbb{N}^r \times \mathbb{A}^n_{\mathbb{Z}} \to \mathbb{N}$ ou $\mathbb{N}^r \times \mathbb{A}^n_{\mathbb{Z}} \to \mathbb{A}^m_{\mathbb{Z}}$ est exactement la donnée d'une fonction calculable au sens usuel qui à toute valeur $(i_1, \dots, i_r)$ des $r$ premiers paramètres associe une donnée comme décrite dans les deux paragraphes précédents.*

On peut faire la même remarque que **22.5** : nous avons ici considéré des fonctions calculables totales sur $\mathbb{N}^r \times \mathbb{A}^n_{\mathbb{Z}}$, mais on pourrait énoncer les résultats analogues pour des fonctions calculables totales sur $\mathbb{N}^r \times E$ où $E$ est une partie constructible de $\mathbb{A}^n_{\mathbb{Z}}$ (en particulier, il est clair que tout ce que nous disons s'applique à l'ouvert $\mathbb{A}^n_{\mathbb{Z}[\frac{1}{N}]}$ de $\mathbb{A}^n_{\mathbb{Z}}$).

À titre d'exemple, la fonction qui calcule la caractéristique de $k$ *n'est pas* calculable (puisqu'elle est totale et ne répond pas aux conditions de ce corollaire) ; en revanche, la fonction qui calcule la caractéristique de $k$ si celle-ci est $\leq p_0$ (pour un certain $p_0$), et 0 si la caractéristique est $> p_0$, est bien calculable.

Notons la conséquence immédiate suivante :

**26.6. Corollaire.** *Soit $f$ une fonction totale calculable $\mathbb{N}^r \times \mathbb{A}^n_{\mathbb{Z}} \to \mathbb{N}$, et soient $i_1, \dots, i_r \in \mathbb{N}$ et $x_1, \dots, x_n \in \mathbb{Z}$. Alors il existe $p_0$ tel que la valeur de $f$ en $(i_1, \dots, i_r, x_1, \dots, x_n)$ ne dépende pas de $p$ si*

*$p \geq p_0$ ou $p = 1$, où on considère $x_1, \ldots, x_n$ comme des éléments d'un corps d'exposant caractéristique $p$ (cette valeur ne dépend pas non plus du corps choisi). De plus, ce $p_0$ peut être calculé algorithmiquement.*

*On a un résultat analogue pour une fonction totale calculable $\mathbb{N}^r \times \mathbb{A}^n_{\mathbb{Z}} \to \mathbb{A}^m_{\mathbb{Z}}$ : pour tous $i_1, \ldots, i_r \in \mathbb{N}$ et $x_1, \ldots, x_n \in \mathbb{Z}$, il existe $p_0$, calculable algorithmiquement, et des rationnels $y_1, \ldots, y_m$ dont le dénominateur n'est divisible par aucun nombre premier $\geq p_0$, tel que la valeur de $f$ en $(i_1, \ldots, i_r, x_1, \ldots, x_n)$, où on considère $x_1, \ldots, x_n$ comme des éléments d'un corps d'exposant caractéristique $p$, soit égale à $\underline{y}$ modulo $p$ si $p \geq p_0$, et à $\underline{y}$ si $p = 1$.*

*Démonstration.* Considérons le cas d'une fonction $\mathbb{N}^r \times \mathbb{A}^n_{\mathbb{Z}} \to \mathbb{N}$. L'évaluation de $f$ en $(i_1, \ldots, i_r, x_1, \ldots, x_n)$ définit une fonction calculable totale de 0 variables, c'est-à-dire $\mathbb{A}^0_{\mathbb{Z}} \to \mathbb{N}$ (ici, évidemment, $\mathbb{A}^0_{\mathbb{Z}}$ est $\operatorname{Spec} \mathbb{Z}$) elle est donc décrite, comme ci-dessus, par une partition de $\operatorname{Spec} \mathbb{Z}$ en un nombre fini de parties constructibles, calculables algorithmiquement : il suffit de renvoyer un $p_0$ plus grand que tous les $p$ intervenant dans les formules booléennes (qui ne peuvent être que $p = 0$) décrivant ces parties constructibles.

Le second cas est complètement analogue : la donnée d'une fonction calculable $\mathbb{A}^0_{\mathbb{Z}} \to \mathbb{A}^m_{\mathbb{Z}}$ revient à se donner des rationnels sur un nombre fini de parties constructibles de $\operatorname{Spec} \mathbb{Z}$. □

## Bibliographie

### Sigles

## Références